\documentclass[11pt, reqno
, final
]{amsart}

\usepackage[latin1]{inputenc}
\usepackage{amssymb, amsfonts, amsmath, amscd, amsthm}
\usepackage[scr]{rsfso}
\usepackage{hyperref}
\hypersetup{
    colorlinks=true,
    linkcolor=blue
}

\usepackage{mathtools}
\usepackage{a4wide, enumitem}
\usepackage{tikz-cd}
\usepackage{graphicx}
\usepackage{caption}
\usepackage{subcaption} 
\usepackage{float}
\usepackage{dsfont}
\usetikzlibrary{arrows, bending}
\usepackage[new]{old-arrows}
\usetikzlibrary{shapes.geometric, calc}
\usetikzlibrary{decorations.markings}
\usetikzlibrary{decorations.pathmorphing}
\usetikzlibrary{shapes.misc}
\usetikzlibrary{bbox} % to get bounding box
\usepackage{outlines}
\usepackage{stmaryrd}
\usepackage{xcolor} %\pagecolor{white}
\usepackage{comment}
\usepackage{tikz-3dplot}
\usetikzlibrary{decorations.pathreplacing} %% package for the braces above/below interval in picture 
\usepackage{pagecolor} %to get white page-color
\usepackage{textcomp} %To get interrobang = !+?

\usepackage%
{standalone}% requires -shell-escape, to get standalone tikzpictures

\usepackage{subfiles}% to be able to compile subfiles separately

\usepackage{import}
\usepackage{pdfpages}
\usepackage{transparent}
\usepackage{overpic} %% package to layer latex-text over pictures

\usepackage{cleveref}
\usepackage{autonum}%only number equations that are actually referenced

\usepackage{macros/utilities}%collect some general LaTeX-functions
\usepackage{macros/macros}
\usepackage{macros/tikz-macros}%Eilind's macros for drawing pictures
\usepackage{macros/environments}%our custom enviroments
\usepackage{macros/diagram-macros}%Tashi's macros for squares,cubes,etc.

\usepackage[
  a4paper
, centering
, margin=3cm
]{geometry}

\numberwithin{equation}{subsection}%number equations within each section
\pagecolor{white} %gives eilind's compiler actual white pages and not grey :)

\makeatletter \def\l@subsection{\@tocline{2}{0pt}{2pc}{6pc}{}} \makeatother

\begin{document}

\title{Assembly of constructible factorization algebras}

\author{Eilind Karlsson}
\address[Eilind Karlsson]{TU M{\"{u}}nchen, Boltzmannstra{\ss}e 3, 85748 Garching, Germany}
\email{eilind.karlsson@ma.tum.de}
\author{Claudia I.~Scheimbauer}
\address[Claudia I.~Scheimbauer]{TU M{\"{u}}nchen, Boltzmannstra{\ss}e 3, 85748 Garching, Germany}
\email{scheimbauer@tum.de}
\author{Tashi Walde}
\address[Tashi Walde]{Universit{\"{a}}t Regensburg, Fakult{\"{a}}t f{\"{u}}r Mathematik,  93040 Regensburg, Germany}
\email{tashi.walde@ur.de}

\maketitle
\thispagestyle{empty}%to remove page number from first page

\begin{abstract}
  We provide a toolbox
  of extension, gluing, and assembly techniques for factorization algebras.

  Using these tools, we fill various gaps in the literature
  on factorization algebras on stratified manifolds,
  the main one being that constructible factorization algebras form
  a sheaf of symmetric monoidal $\infty$-categories.
  Additionally,
  we explain how to assemble constructible factorization algebras
  from the data on the individual strata
  together with module structures associated to the relative links;
  thus answering a question by Ayala.

  Along the way, we give detailed proofs of the following facts
  which are also of independent interest:
  constructibility is a local condition;
  the \infy-category of disks is a localization of any sufficiently fine poset of disks;
  constructibility implies the Weiss condition on disks;
  constructible factorization algebras are algebras for the \infy-operad of embedded disks.
  For each of these, variants or special cases already existed,
  but they were either incomplete or not general enough.
\end{abstract}

\tableofcontents
\thispagestyle{empty}%to remove page number from second page

% \listoftodos{}

\newpage

\section{Introduction}

The notion of a factorization algebra was first developed by Beilinson and Drinfeld in the algebro-geometric context. 
In the topological context, Lurie defined a variant in \cite{LurHA} under the name factorizable cosheaves, and Costello and Gwilliam \cite{CG} defined 
a variant suitable for (perturbative quantization of) observables of a field theory. Constructible factorization algebras live on stratified spaces  
and enjoy the property that their values on certain inclusions of disks (which ``look the same with respect to the stratification'') are equivalences.  
A  special case is that of an empty stratification, which reflects a field theory being topological.
A closely related concept is factorization homology as developed by Lurie, Ayala, Francis, and Tanaka: indeed, a key source 
of constructible factorization algebras is provided by factorization homology. We discuss this below in more detail.

Informally, a factorization algebra on a space \(X\)
valued in a symmetric monoidal \infy-category \(\targetcat\)
is an assignment
which sends an open set $U$ to an object $\Fa(U)\in \targetcat$
and to any inclusion of disjoint opens
$U_1 \disjun \cdots \disjun U_n \subseteq V$
it assigns a morphism in $\targetcat$
\begin{equation}
  f_{U_1,\cdots, U_n; V}\colon \Fa(U_1)\otimes\cdots \otimes \Fa(U_n) \longrightarrow \Fa(V)
\end{equation}
respecting composition.
Moreover, this assignment is required to be multiplicative,
which says that the maps
\begin{equation}
  f_{;\emptyset}\colon \monunit \xrightarrow{\simeq} \Fa(\emptyset)
  \quad
  \text{and}
  \quad
  f_{U, V; U\disjun V}\colon \Fa(U)\otimes\Fa(V) \xrightarrow{\simeq} \Fa(U\disjun V)
\end{equation}
are equivalences for disjoint opens $U$ and $V$,
and is required to satisfy descent
for so-called Weiss covers\footnote{%
  The descent condition is a standard descent condition in topos theory.
  However, there are variations appearing in the literature,
  so we encourage even the reader already familiar with factorization algebras
  to look at \Cref{rem:Fact_alg_variants} and \Cref{rem:terminology}.}.
To encode this formally with all the correct coherences,
one uses the operadic language:
A factorization algebra on a space $X$
is an algebra in \(\targetcat\) for an operad $\Open{X}$ of disjoint open sets in $X$,
which satisfies a multiplicativity condition
and descent for Weiss covers. 

Basic examples include associative algebras,
viewed as factorization algebras on the real line,
bimodules (or left/right modules) with a chosen basepoint,
viewed as constructible factorization algebras on the reals
(or non-positive/non-negative reals)
with the point 0 as the stratification,
and \Ealgebra{n}s, viewed as factorization algebras on \(\reals^n\).

The main goal of this article is to provide a toolbox of techniques to prove various statements about factorization algebras.
While the three main applications deal with constructible factorization algebras, many of the intermediate steps work without this assumption.
Before delving deeper into the toolbox we first explain these applications.

The first application of this toolbox is the following:
Although widely expected (and sometimes claimed),
there is no proof in the literature that factorization algebras can be glued;
for example, that given an open cover $X=U\cup V$
and factorization algebras $\Fa$ on $U$ and $\Ga$ on $V$
together with an isomorphism \(\Fa\restrict{U\cap V}\simeq \Ga\restrict{U\cap V}\),
we can produce a unique factorization algebra on $X$
that restricts to \(\Fa\) and \(\Ga\).
In other words, factorization algebras form a sheaf of categories.
We fill this gap with this article in the constructible setting, that is, for smooth conical manifolds\footnote{A smooth conical manifold is what Ayala--Francis--Tanaka call a ``conically smooth stratified space''. Furthermore, in all main Theorems we require the existence of ``enough good disks'', see \Cref{def:enough_disks} and \Cref{rem:existence_good_disks}.}.

Our first main theorem states that constructible factorization algebras can be glued from constructible factorization algebras on the opens in any cover.
\begin{thmx}[\Cref{thm:main-gluing}] \label{thmx:main-gluing}
  Let $\Ua$ be an open cover of a smooth conical manifold $X$. Then the map induced by restricting  factorization algebras
  \begin{equation}
    \FactCstr[\targetcat]{X} \longrightarrow \lim_{\{U_1,\ldots, U_n\}\subseteq \Ua}  \FactCstr[\targetcat]{U_1\cap\cdots\cap U_n}
  \end{equation}
  is an equivalence of \infy-categories. In other words, $\Fact[\targetcat]{(-)}$ is a sheaf of \infy-categories.
\end{thmx}
In the unstratified situation, i.e., for locally constant factorization algebras,
this was already shown in \cite{Matsuoka}.

We also record an easy consequence of our gluing results.
For our very basic examples above, note that both algebras and (bi)modules admit a symmetric monoidal structure, given by endowing the tensor product of the underlying vector spaces (or object in the underlying category) with the desired structure. The analogous statement for constructible factorization algebras holds. This result was also long believed to be true and has been stated several times in the literature, albeit without proof. 
\begin{thmx}[\Cref{thm:symmetric_monoidal} and \Cref{prop:symm_mon_gluing}]
  For every smooth conical manifold $X$,
  the functor
  \begin{equation}
    \FactCstrO[\targetcat]{X}\colon \Fin\to \Catinfty,
    \quad\quad
    I_+\mapsto \FactCstr[\targetcat]{X^{\amalg I}}
  \end{equation}
  exhibits a monoid object in \((\Catinfty,\times)\),
  i.e., a symmetric monoidal structure
  on the $\infty$-category \(\FactCstr[\targetcat]{X}\) of constructible factorization algebras on $X$.
  
  In fact, \Cref{thmx:main-gluing} can be upgraded to the symmetric monoidal version: $\Fact[\targetcat]{(-)}$ is a sheaf of symmetric monoidal \infy-categories.
\end{thmx}

\begin{rem}
These two results are essential ingredients for building higher Morita categories using constructible factorization algebras as was exploited in work by the second author \cite{Scheimbauer, JFS}. As these results had been widely claimed in the literature, the authors relied on them, and a purpose of this article is to fill this gap.
\end{rem}

The last application of our toolbox is to unravel the data required for a constructible factorization algebra on a smooth conical manifold. The simplest version is for cones, and answers a question of David Ayala.
The same result was also obtained by Brav--Rozenblyum~\cite{brav-rozenblyum}
using very different techniques, see also \Cref{rem:brav-rozenblyum}.

\begin{thmx}[\Cref{thm:cones_as_modules}] \label{thmx:cones_as_modules}
  Let $Z$ be a compact smooth conical manifold.
  We have a pullback square of \(\infty\)-categories
  \begin{equation}
    \begin{tikzcd}
      \FactCstr[\targetcat]{ \topcone{Z} }  \arrow{d}{\pf{p}} \arrow{r}
      \isCartesian
      &
      \FactCstr[\targetcat]{Z\times (0,\infty) }   \arrow{d}{\pf{p}}
      &
      \Aa
      \ar[d,mapsto]
      \\
      \FactCstr[\targetcat]{[0,\infty) } \arrow{r}
      &
      \FactCstr[\targetcat]{ (0,\infty) }
      &
      \int_Z\Aa
    \end{tikzcd}
  \end{equation}
  where the vertical functors are pushforward
  along the quotient map \(p\colon \topcone{Z}\to [0,\infty)\)
  and the horizontal functors are restriction to open subspaces.
\end{thmx}

The lower right hand side is equivalent to algebra objects in $\targetcat$,
and the left bottom side encodes a pointed module for such an algebra;
see \Cref{ex:ass-and-rmod}.
Moreover, there should be a version of Dunn's additivity for constructible factorization algebras, and combining this with the above-mentioned result gives that the top right corner encodes an algebra object in constructible factorization algebras on~$Z$.

We elaborate on how to globalize this statement for a general smooth conical manifold in \Cref{constr:step-unlinking} and \Cref{remark:globalizing_unlinking}.

\subsection{Opening the toolbox}

To prove the main theorems, we develop a toolbox of techniques for (constructible) factorization algebras. 
A key step in the toolbox is to formulate the statements at a level general enough to be useful in other settings as well.

Let us briefly outline the strategy of the proof of gluing of constructible factorization algebras in the special case of an open cover of a topological space of \(X\) with two elements \(\Ua=\{U, V\}\) and use this to explain our toolbox.

{\em Step 1:} The first step, which is done in \Cref{subsect:FirstStepGluing},  is to compute the pushout of $\infty$-operads along the natural inclusions
\begin{equation}
  \begin{tikzcd}
    \Open{U\cap V} \arrow{r}\arrow{d} &\Open{U} \arrow[dashed]{d}\\
    \Open{V} \arrow[dashed]{r} &\operad{\cdown{\Ua}}\,.
  \end{tikzcd}
\end{equation}
We explain why this can be computed as the (strict) colimit of dendroidal sets. %In our example of having just two open sets $U$ and $V$ in the cover, we have that
Here, \(\operad{\cdown{\Ua}}\) is simply the \infy-operad corresponding to the collection of opens which are either in $U$ or in $V$, and \(\operad{\cdown{\Ua}}\)-algebras and hence also \(\operad{\cdown{\Ua}}\)-factorization algebras can be evaluated on precisely those opens. From this pushout we immediately see that (constructible) \(\operad{\cdown{\Ua}}\)-factorization algebras are glued from 
(constructible) \(\Open{U}\)- and \(\Open{V}\)-factorization algebras, see \Cref{prop:GlueFactOnCover} and \Cref{cor:GlueFactCstrOnCover}. 
In fact, we explain how to compute the colimit of certain $\infty$-operads of open sets in more generality, which is also used in 
proving \Cref{thmx:cones_as_modules}.

From these partially defined, constructible factorization algebras on $X$ we extend to factorization algebras defined on all opens in $X$ by left Kan extending along two inclusions.
In each step, there are two substeps:
the first (a) uses our toolbox, namely performing certain Kan extensions.
The second (b) uses ``constructibility tools'', which are generalizations of statements already (mostly) existing in the literature.

{\em Step 2a:} First, in \Cref{subsect:ExtendingToDisjUniCompletion} we enlarge \(\Bf=\operad{\cdown{\Ua}}\) by adding disjoint unions of opens which are already contained in \(\Bf\) to obtain \(\disjunComp{\Bf}\). In our  example of only two open sets $U$ and $V$ in the cover, this means adding open sets of the form $W_1 \disjun W_2$, where $W_1\subset U$ and $W_2 \subset V$. We show that there is an equivalence of categories between multiplicative prefactorization algebras for \(\Bf\) and for \(\disjunComp{\Bf}\), given by left Kan extending along the inclusion, see \Cref{prop:MultAlgToDisjCompletion}.

{\em Step 2b:} To guarantee Weiss descent we restrict to constructible factorization algebras,
since in this case the Weiss condition is automatically satisfied
as long as we consider algebras Kan extended from a suitable poset of disks.
This is proven in \Cref{prop:final_step_in_gluing} and relies on a generalization of a proof by Ayala-Francis, see \Cref{cor:constructible=>Weiss}.

{\em Step 3a:} This is the final step, which is done in \Cref{subsect:ExtendingUFactAlgs}, to reach factorization algebras defined on all opens. 
In fact, \(\disjunComp{\Bf}\) is a so-called factorizing basis, so it is enough to show that a factorization algebra can be recovered from its values on a factorizing basis. This was already shown in \cite{CG} for $\targetcat$ being chain complexes by an explicit construction.
For general $\targetcat$,  we again perform an operadic left Kan extension, this time along the inclusion of the factorizing basis into all opens, see \Cref{lemma:FactBasisExtends}.  
This last step already gives an equivalence at the level of prefactorization algebras satisfying only Weiss descent (\Cref{prop:AlgWeissEquivFactBasis}).

{\em Step 3b:}
The above equivalence restricts to constructible factorization algebras (\Cref{prop:extension-from-sieve-cover}) since constructibility is local, which we prove in \Cref{thm:constructle-local} generalizing a proof of Ginot.

Summarizing, we have the following chain of equivalences.
\begin{equation}
  \begin{tikzcd}
    \FactCstr{{X}} \arrow[r, "\sim", "\ref{prop:extension-from-sieve-cover}"']
    & \FactCstr{\cdisj{\cdown{\Ua}}} \arrow[r, "\sim", "\ref{prop:final_step_in_gluing}"']
    & \FactCstr{\cdown{\Ua}} \arrow[r, "\sim", "\ref{prop:GlueFactOnCover}"']
    & \displaystyle \lim_{\{U_1,\ldots, U_n\}\subset \Ua} \FactCstr{U_1\cap\cdots\cap U_n}\,.	
  \end{tikzcd}
\end{equation}

\subsection{Constructibility tools}

In \Cref{sec:external-tools} we collect, complete, and generalize some folklore results (partially) existing in the literature needed in parts (b) of the steps above. These are interesting in their own right and may be considered a constructibility toolbox.

First we generalize a result of Ginot which proves that local constancy is in fact a local condition.
\begin{thmx} [\Cref{thm:constructle-local}]
Let \(M\) be a conical manifold
  and \(\Aa\) a factorization algebra on \(M\).
  Then \(\Aa\) is constructible if and only if it is locally constructible.
\end{thmx}
The proof is an adaptation of Ginot's proof to the constructible setting and allows checking constructibility locally, for instance using a cover.

The second constructibility tool deals with localizing at isotopy equivalences. 
This relates the localization of the discrete category of (certain) disks with the topological/\infy-category \(\ddisk\) of disks and is a main result of \cite{AFT-fh-stratified}. 
Although versions of the missing argument had appeared before in various places
(see \Cref{rem:localizing_gap}),
here we provide the details;
in the unstratified setting the original outline
was also completed independently by Arakawa~\cite{Arakawa}
at the same time this paper was being written.
Moreover, we need a more general version of this result, which allows localizing a smaller poset \(\Bf\subseteq\diskX\) of disks:
\begin{thmx}[\Cref{thm:Disk-localized}]
  Let \(X\) be a smooth conical manifold and
  \(\Bf\subseteq\diskX\) be a decomposable\footnote{
    decomposable = closed under disjoint decompositions;
    see \Cref{defn:DecomposablePresieve}
  }
  and multiplicative\footnote{
    multiplicative = closed under disjoint unions;
    see \Cref{defn:mult_fact_basis}
  } basis of \(X\)
  Then the tautological functor
  \begin{equation}
    \Bf	\to
    \ddiskofX
  \end{equation}
  is an \(\infty\)-categorical localization.
\end{thmx}
We apply this result in our main \Cref{thmx:cones_as_modules}, and given its flexibility we expect this to be applicable to other situations as well.

For a long time it has been folklore that constructibility automatically implies the Weiss condition.
We prove that this is true on disks. 
\begin{thmx}
  [\Cref{cor:constructible=>Weiss}, \Cref{cor:extending-from-disk-basis}, and \Cref{cor:extending-from-disks}]
  \label{thmx:extending-from-disk-basis}
  Let \(X\) be a smooth conical manifold with enough good disks
  and \(\Bf\subseteq\diskX\) a decomposable disk-basis 
  which is factorizing\footnote{
    factorizing = closed under disjoint unions and intersections;
    see \Cref{defn:mult_fact_basis}
  }.
  \begin{enumerate}[label=(\arabic*), ref=(\arabic*)]
  \item
    Every multiplicative prefactorization algebra on \(\Bf\)
    is automatically a Weiss cosheaf, hence a factorization algebra.
  \item
    \label{it:inthmx-extending-from-disk-basis}
    Restriction yields equivalences
    \begin{equation}
      \FactCstr[\targetcat]{X}
      \xrightarrow{\simeq}
      \AlgMCstr[\targetcat]{\diskX}
      \xrightarrow{\simeq}
      \FactCstr[\targetcat]{\Bf}=\AlgMCstr[\targetcat]{\Bf}
    \end{equation}
    of \(\infty\)-categories.
  \end{enumerate}
\end{thmx}

We also prove an operadic extension of the aforementioned localization result
and use it to compare the two main different ways that
``constructible factorization algebras''
are usually defined in the literature.

\begin{thmx}
  [\Cref{prop:Disk-operad-localization} and \Cref{cor:FactX-AlgDiskX}]
  Let \(X\) be a smooth conical manifold.
  \begin{enumerate}[label=(\arabic*), ref=(\arabic*)]
  \item
    For every decomposable multiplicative disk-basis \(\Bf\) of \(X\),
    the tautological map
    \begin{equation}
      \operad{\Bf}\to\DDiskof{X}
    \end{equation}
    is a localization of \infy-operads.
  \item
    If \(X\) has enough good disks,
    we construct an equivalence
    \begin{equation}
      \FactCstr[\targetcat]{X}\simeq \AlgM[\targetcat]{\ddiskofX}
    \end{equation}
    between constructible factorization algebras on \(X\)
    and multiplicative algebras for the \infy-operad
    of embeddings of disks in \(X\).
  \end{enumerate}
\end{thmx}

Note that in part \ref{it:inthmx-extending-from-disk-basis} of
\Cref{thmx:extending-from-disk-basis},
we prove the first equivalence by first proving that the composite is an equivalence.
For this, it is crucial that the collection \(\Bf\) of disks
is closed under intersections.
In particular, the collection \(\Bf\coloneqq \disk{X}\) of \emph{all} disks
does not satisfy this, 
but for example taking \emph{convex} disks of an unstratified smooth manifold
with respect to a chosen metric does.
In the unstratified setting,
an analogous result was also obtained
in Carmona's thesis~\cite{Carmona-thesis},
who also needs the assumption on good disks.

We contrast this with \cite{AFprimer}
where in the unstratified setting the analogous statement is made for \(\Bf=\disk{X}\);
see also \Cref{rem:AF_analog_extension}.

A corollary of these results is that given a fixed (framed) $n$-manifold and an $E_n$-algebra $\Aa$,
assigning to any open $U$ the factorization homology $\int_U \Aa$ is a factorization algebra.
As stated, this was already proven in \cite{GTZ} under the assumption of the gluing \Cref{thm:main-gluing}
and has long been folklore to be a corollary of the localization result.
Here we record this as a corollary of the previous result in the stratified setting.
\begin{thmx}[\Cref{thm:fact_hom_is_fact_alg}]
  Let \(\cB\) be a type of tangential structure.
  For $M$ a $\cB$-structured manifold and \(\Aa\) a $\ddiskB$-algebra,
  the prefactorization algebra
  \begin{equation}
    \Fa_\Aa \colon \Open{M}   \to \MfldB \to
    \MMfldB \xrightarrow{\int_{-}\Aa} \targetcatOT\, 
  \end{equation}
  given by precomposing the factorization homology functor
  is a factorization algebra.
\end{thmx}

\subsection{Conventions and notation}
\label{sec:guide}

\begin{itemize}
\item
  We freely use the language of \infy-categories
  ((co)limits, Kan extensions, localizations, (co)cartesian fibrations etc.)
  as developed in \cite{LurHTT}; another good reference is \cite{Cisinski}.
  None of our arguments make use of any of the specific features
  of the Joyal--Lurie implementation in quasicategories,
  so that we allow ourselves to use a model-agnostic language.
\item
  We view posets and (ordinary) categories as \infy-categories
  whose mapping spaces are \((-1)\)-truncated and \(0\)-truncated, respectively.
  Thus we do not need to explicitly apply the nerve construction
  to treat a poset or category as an \infy-category.
\item
  We adopt the usual conventions regarding set-theoretic size issues:
  we fix Grothendieck universes of ``small'' and ``large'' sets.
  Sets and spaces are always small by default;
  \infy-categories are large but have small mapping spaces.
\item
  We also make free use of the theory of \infy-operads and algebras thereover,
  which we view as implemented within the theory of \infy-categories as in \cite{LurHA}.
  We make heavy use of operadic colimits and left operadic Kan extensions;
  for the convenience of the reader we summarize its key features
  in \Cref{app:OLKE} specializing and summarizing Lurie's more general theory.
\item
  We write \(\Gpdinfty\), \(\Catinfty\), \(\SMCatinfty\), and \(\Opdinfty\)
  for the \infy-categories of (small)
  spaces/\infy-groupoids, \infy-categories, symmetric-monoidal \infy-categories,
  and \infy-operads, respectively.
\item
  The corresponding full subcategories of (small)
  sets, categories, symmetric monoidal categories, and operads
  are denoted by
  \(\Set\), \(\Cat\), \(\SMCat\), and \(\Opd\), respectively.
  Note that this makes \(\Cat\) and \(\Opd\) into \((2,1)\)-categories
  whose 2-cells are the natural isomorphisms.
\item
  In contrast, we denote by \(\strictCat\) and \(\strictOpd\)
  the \((1,1)\)-categories of categories and operads, respectively.
  This means that we allow no non-trivial natural isomorphisms
  and that invertible arrows are \emph{isomorphisms} of categories/operads---%
  a notion that is usually too strict to be useful.
\item
  We use the symbols \(X,Y,Z\) for topological spaces,
  which we always assume to be Hausdorff.
  Often they will be (topological or smooth) manifolds
  or stratified manifolds; we recall the relevant definitions in
  \Cref{sec:conical-manifold}.
\item
  When talking about \infy-operads, we write \(\Oa^\otimes\) for the operad itself
  and \(\Oa\) for its underlying \infy-category of 1-ary operations.
  Conversely, some \infy-categories \(\Ua\) can be naturally
  enhanced to an \infy-operad
  (for example if they carry a partial symmetric monoidal structure)
  which we then denote by \(\operad{\Ua}\).
\item
  All our algebras take values in \((\targetcat,\otimes)\),
  which is an \(\otimes\)-presentable symmetric monoidal \infy-category
  that we fix throughout the paper.
  We write \(\Alg[\targetcat]{\Oa}\) for the \infy-category of \(\Oa^\otimes\)-algebras
  and often omit ``\(\targetcat\)'', since it is always fixed.
 \item
   We write \(\openX\) for the poset of open subsets of \(X\)
   and by \(\cdisk{X}\) its full subposet consisting of embedded disks.
   In contrast, \(\diskX\) denotes the full subposet of embedded multidisks,
   i.e., finite disjoint unions of embedded disks;
   see \Cref{not:disks}.
   We sometimes call the elements of  \(\diskX\) just ``disks''
   and say ``contractible disks'' to emphasize when an object lies in \(\cdisk{X}\).
\item
  We write \(\Alg{X}\) and \(\Fact{X}\) for the \infy-categories
  of prefactorization algebras and factorization algebras on \(X\), respectively;
  see \Cref{subsect:FactAlgebras} for the definitions.
  Our prefactorization algebras are lax by default,
  while our factorization algebras are non-lax (=multiplicative).
\item
  In operads coming from posets of open subsets,
  we usually talk about ``\pcocartesian{}'' (with a \pcocartdecoration{})---%
  rather than ``cocartesian''---morphisms;
  see \Cref{def:pcocartesian} and \Cref{rem:pcocartesian}.
  This distinction only exists for technical reasons and can safely be ignored by most readers;
  for most cases of practical relevance the two notions agree.
\item
  We use the superscripts ``$\lc$'' and ``$\cstr$'' to denote \infy-categories
  of locally constant or constructible (pre-)factorization algebras, respectively;
  see \Cref{sec:ConstructibleFAs} for the definitions.
\item
  We write \(\Snglr\) and \(\SSnglr\) for the ordinary category
  and the \infy-category of
  conically smooth stratified spaces of \cite{AFT-local-structures},
  respectively,
  which we call (smooth) conical manifolds;
  \(\ddisk\) refers to the full sub-\infy-category of multidisks.
  We recall the definitions in \Cref{sec:smooth-conical-manifold}.
\item
  We also use the versions
  \(\mfldB\) and \(\mmfldB\) where all smooth conical manifolds
  are equipped with a tangential structure of type \(\cB\).
\item
  A copresheaf is simply a functor to \(\targetcat\),
  usually defined on some poset of opens of \(X\),
  see \Cref{def:copresheaf}
\item
  We use a simplified version of the general theory of descent and cosheaves,
  which is self-contained and well-adapted for our setting.
  We provide a detailed account in the first background section,
  \Cref{sec:background}.
\item
  In particular, we use certain terms such as
  ``sieve'', ``cover'', ``refinement'', ``hypercover''
  in a slightly non-standard way.
  We explain how these relate to the more standard notions in 
  \Cref{rem:terminology}.
\item
  We write \(U\cup V\) for the union of two subsets of \(X\),
  and \(U\disjun V\) if this union is disjoint.
  We distinguish this from \(U\amalg V\) which refers to the abstract
  external coproduct, unrelated to the inclusions into \(X\).
\item
  If \(\Ua\) is a poset of open subsets of \(X\),
  we denote by \(\ccaps{\Ua}\), \(\cdown{\Ua}\) and \(\cdisj{\Ua}\)
  its closure under binary intersections, open subsets and finite disjoint unions,
  respectively.
  \item
  We write \(\topcone Z\) for the cone of $Z$ and \( \topcone[t]X \)
  for the cone of \(X\) of radius \(t\), see \Cref{cstr:stratified-cone}.
\end{itemize}

\subsection{Acknowledgements}
We thank David Ayala for posing the question that led to \Cref{thm:cones_as_modules}.
Pelle Steffens pointed out the magic of the Seifert--Van Kampen Theorem and suggested applying it for \Cref{thm:Disk-localized}.
CS thanks her coauthors for joining in finally filling these gaps in the literature.

All authors were supported by the SFB 1085: Higher Invariants from the Deutsche Forschungsgemeinschaft
(DFG).
CS was supported by the Simons Collaboration on Global Categorical Symmetries (1013836).

\newpage
\section{Background on Weiss cosheaves}
\label{sec:background}

In this section we will set up the foundations and background needed to be precise in later arguments. In particular, we will first explain how we think of posets of opens in different ways. Then we will discuss the notion of locality, before defining Weiss covers and cosheaves with respect to Weiss covers. We also comment on the terminology used here before lastly defining the notion of a factorizing basis.

We fix
\begin{itemize}
\item
  a topological space \(X\) which we always assume to be Hausdorff;
\item
  a \(\otimes\)-presentable symmetric monoidal \(\infty\)-category
  \((\targetcat,\otimes)\).
\end{itemize}

\subsection{Sieves, (hyper)covers and refinements}
In this paper, we will be juggling many different sets of open subsets of \(X\).
Sometimes, we will think of such a set as a poset
(inheriting the inclusion-ordering of \(\openX\))
while other times we just treat them as subsets.

To minimize potential confusion, we summarize our conventions here,
before defining all the terms;
see also \Cref{rem:hypercover-nonstandard} below for more explanation
justifying this terminology.
\begin{itemize}
\item
  ``(Pre)covers'' refer to subsets of \(\openX\);
  they are denoted by mathcal-letters such as \(\Ua,\Va,\Wa\).
\item
  ``Hyper(pre)covers'' and ``(pre)sieves'' refer to (full) subposets of \(\openX\);
  they are denoted by mathfrak-letters such as \(\Bf,\Uf,\Wf\).
\end{itemize}

\begin{defn}
  \begin{itemize}
  \item[]
  \item
    A \emph{precover} of \(X\) is
    a set \(\Ua\subseteq \openX\) of open subsets of \(X\).
  \item
    Given two such precovers \(\Va,\Ua\) of \(X\) we say that
    \(\Va\) refines \(\Ua\)
    if for every \(V\in \Va\) there is an \(U\in \Ua\) with \(V\subseteq U\).
    We say that we have a \emph{refinement} \(\Va\to \Ua\).
  \item
    Clearly refinements compose.
    We write \(\PCov{X}\) for the category of precovers and refinements.
  \item
    We call a refinement of the form \(\Va\to\{U\}\) (with singleton codomain)
    a precover. We say that \(\Va\) is a precover of the open set \(U\)
    and just write \(\Va\to U\).
  \item
    Such a precover is called a \emph{cover} if \(\collunion \Va=U\).
  \end{itemize}
  When we just say that \(\Ua\) is a precover,
  we implicitly mean that it is precover of \(X\) unless stated otherwise.
\end{defn}

We distinguish different kinds of subposets of \(\openX\).

\begin{defn}
  When we say that \(\Uf\) is a \emph{\posetofopens{}} (of \(X\)),
  we mean that \(\Uf\subseteq\openX\) is a full subposet.
  It is called
  \begin{itemize}
  \item
    a \emph{presieve}, if it is closed under binary intersections,
    i.e.\ if \(U,U'\in\Uf\) implies \(U\cap U'\in \Uf\);
  \item
    a \emph{sieve} if it is downward closed,
    i.e.\ if \(U\in \Uf\) and \( U'\subseteq U\) imply \(U'\in\Uf\).
  \end{itemize}
\end{defn}

Precovers and \posetsofopens{} are related by the following constructions:

\begin{constr}
  \begin{itemize}
  \item[]
  \item
    Every \posetofopens{} \(\Uf\) has an underlying precover
    \(\fgtposet{\Uf}\)
    by just forgetting the poset structure.
  \item   
    For every precover \(\Ua\)
    we have the associated presieve
    \begin{equation}
      \ccaps{\Ua}\coloneqq \{U_1\cap\dots\cap U_n\mid n\geq 1, U_i\in \Ua\}
    \end{equation}
    obtained by closing under binary intersections
    and the associated sieve
    \begin{equation}
      \cdown{\Ua}\coloneqq \{U'\in\openX\mid \exists U\in \Ua: U'\subseteq U\}
    \end{equation}
    obtained by closing it downward.
  \end{itemize}
\end{constr}

Throughout this paper,
we will be working with \(\targetcat\)-valued copresheaves on \(X\),
i.e.\ functors \(\Aa\colon\openX\to \targetcat\).
Moreover, it will be convenient to consider copresheaves
defined not on all opens of \(X\) but only on a smaller selection.

\begin{defn}
  \label{def:copresheaf}
  Let \(\Uf\subseteq\openX\) be a \posetofopens{}.
  A (\(\targetcat\)-valued) \emph{copresheaf on \(\Uf\)}
  is a functor \(\Aa\colon \Uf\to\targetcat\).
\end{defn}

\begin{notation}
  We say that \(\Aa\) is a
  \emph{copresheaf defined on \(\Uf\)}
  if \(\Aa\colon \Uf'\to \targetcat\) is a copresheaf
  with \(\Uf'\supseteq \Uf\) some \posetofopens{}
  which we wish to leave unspecified.
  Whenever we use this turn of phrase, we are making the implicit claim that
  any statement we then make about \(\Aa\) does not depend on the choice of \(\Uf'\).
\end{notation}

\begin{constr}
  \label{cstr:copresheaf-on-presieve}
  Let \(\Uf\) be a \posetofopens{} and
  \(\Aa\colon \Uf\to \targetcat\) be a copresheaf.
  Then we can evaluate \(\Aa\) not just on elements of \(\Uf\),
  but also on any precover \(\Va\in\PCov{X}\)
  (and refinements between them)
  by left Kan extending along the inclusion
  \begin{equation}
    \Uf\subseteq \openX\hookrightarrow \PCov{X},
    \quad\quad U\mapsto \{U\}.
  \end{equation}
  Explicitly, the value of \(\Aa\) at \(\Va\) is given by the colimit
  \begin{equation}
    \label{eq:def-A-of-V-with-U}
    \AonPCov[\Uf]{\Va} \coloneqq\colim \Aa\restrict{\Uf\cap\cdown{\Va}}.
  \end{equation}
\end{constr}

We will also need the following more general notion.

\begin{defn}
  A \emph{\hyperprecover{} of \(X\)} is a full subposet \(\Uf\subseteq \openX\).
  Given two \hyperprecover{}s \(\Uf\) and \(\Wf\),
  we say that \emph{\(\Uf\) hyperrefines \(\Wf\)}, written \(\Uf\hyprefines\Wf\)
  if for each \(U\in \Uf\)
  the poset
  \(\undercat{\Wf}{U}\coloneqq\{W\in \Wf\mid U\subseteq W\}\) is weakly contractible.
\end{defn}

\begin{rem} \label{rem:hyprefines-gives-refinement-onfgtposet}
  If \(\Uf\) hyperrefines \(\Wf\)
  then for each \(U\in\Uf\) the poset \(\undercat{\Wf}{U}\) is in particular non-empty,
  so that \(\fgtposet{\Uf}\) refines \(\fgtposet{\Wf}\).
\end{rem}

\begin{rem}
  Although the phrases
  ``\posetofopens{}''
  and
  ``\hyperprecover{}''
  are synonyms,
  we use them differently:
  A \posetofopens{} is used as the \emph{domain of definition} of a co(pre)sheaf
  \(\Aa\colon \Uf\to \targetcat\),
  while a \hyperprecover{} is something a co(pre)sheaf
  can be \emph{evaluated on}
  (see \Cref{defn:A-on-hyperprecover}
  and \Cref{defn:A-on-hyprefinement} below).
  Of course every \posetofopens{} can be interpreted as a \hyperprecover{}:
  the value of \(\Aa\colon \Uf\to\targetcat\) on \(\Uf\)
  corresponds to the global sections of \(\Aa\).
\end{rem}

\begin{notation}
  When we have a hyperrefinement \(\Wf\hyprefines\{U\}\) of a singleton,
  we say that \(\Wf\) is a \emph{\hyperprecover{} of} \(U\)
  and just write \(\Wf\hyprefines U\).
  Note that this happens precisely if we have a precover
  \(\fgtposet{\Wf}\to U\).

  A \hyperprecover{} \(\Wf\hyprefines U\) is called a \emph{hypercover} if,
  for each finite list \(W_1,\dots,W_n \in \Wf\),
  the precover
  \begin{equation}
    \{W\in\fgtposet{\Wf}\mid W\subseteq W_1\cap\dots\cap W_n\}\to W_1\cap\dots\cap W_n
  \end{equation}
  is a cover.
  Here the intersection is taken in \(U\)
  so that the case \(n=0\) says precisely that \(\fgtposet{\Wf}\to U\) is a cover.
\end{notation}

\begin{rem} \label{rem:Inclusions-are-hyprefinements}
  Every inclusion \(\Uf\hookrightarrow \Wf\) of \hyperprecover{}s
  is also a hyperrefinement
  because for each \(U\in\Uf\)
  the poset \(\undercat{\Wf}{U}\) has an initial object,
  namely \(U\) itself.
  In particular, the hyperrefinement-relation is reflexive.
\end{rem}

\begin{war}
  In general hyperrefinements do \emph{not} compose; i.e.\ given \(\Uf \hyprefines \Vf \hyprefines \Wf\) it does not follow that \(\Uf \hyprefines \Wf\). 
  Consider for example \(\Uf = \{U\}, \Vf = \{V\} \) and \(\Wf = \{W, W'\}\) such that we have
  \begin{equation}
    \begin{tikzcd}
      U \arrow[rrd, hook] \arrow[r, hook]& V \arrow[r, hook] & W
      \\ && W'
    \end{tikzcd}
  \end{equation}
  and no other inclusions between them.
  By construction we have \(\undercat{\Wf}{U} = \{W, W'\}\) with no non-trivial morphisms,
  i.e.\ it is not weakly contractible and hence \(\Uf\) does not hyperrefine \(\Wf\).
\end{war}

\begin{rem} \label{rem:SomeHyprefsCompose}
  Although hyperrefinements do not compose in general, there are special cases where they do. The cases relevant in this article are:
  \begin{enumerate}
  \item Inclusions of \hyperprecover{}s compose: Given \(\Uf \hookrightarrow \Vf \hookrightarrow \Wf\) it is clear that \(\Uf \hookrightarrow \Wf\). 
  \item Hyperrefinements into a singleton compose: Given \(\Uf \hyprefines \Vf \hyprefines \{W\}\), it follows that \(\Uf \hyprefines \{W\}\). By \Cref{rem:hyprefines-gives-refinement-onfgtposet} the given hyperrefinements gives refinements \(\fgtposet{\Uf} \ra \fgtposet{\Vf}\) and \(\fgtposet{\Vf} \ra \{W\}\). Since refinements do compose it follows that \(\fgtposet{\Uf} \ra \{W\}\) is a refinement, which is equivalent to \(\Uf \hyprefines \{W\}\) exactly because the target is a singleton. 
  \end{enumerate} 
\end{rem}

We want to evaluate any copresheaf \(\Aa\colon \Uf\to \targetcat\) on \hyperprecover{}s \(\Wf\), as well as hyperrefinements between them, in a similar way to \Cref{cstr:copresheaf-on-presieve}. However, since \hyperprecover{}s and hyperrefinements do not form a category we, a bit unsatisfyingly, need to pointwise define what the evaluation is.

\begin{defn}
  \label{defn:A-on-hyperprecover}
  % \label{cstr:evaluate-on-hyperprecover}	
  Let \(\Aa\colon \Uf\to \targetcat\) be a copresheaf. For every \hyperprecover{} \(\Wf \subseteq \Uf\) we define the evaluation of \(\Aa\) on \(\Wf\) to be 
  \begin{equation}  \label{eq:AOnHPCovInU}
    	\AonHPCov[\Uf]{\Wf} \coloneqq \displaystyle \colim_{W\in \Wf} \Aa(W).
  \end{equation}
  Note that this definition does not depend on the ambient \posetofopens{} \(\Uf\), 
  so we often omit it from the notation and simply write \(\AonHPCov{\Wf}\).
\end{defn}

We can now also define how to evaluate a copresheaf on certain hyperrefinements: 
\begin{defn} \label{defn:A-on-hyprefinement}
  Let \(\Vf, \Wf\) be \hyperprecover{}s with a hyperrefinement \(\Vf \hyprefines \Wf\).
  Let \(\Aa\) be a copresheaf defined on \(\Vf\) and \(\Wf\).
  The evaluation of \(\Aa\) on this hyperrefinement is defined to be the composite
  \begin{equation}
    \label{eq:A-on-hyperrefinements}
    \colim_{V\in \Vf} \Aa(V)
    \xleftarrow{\simeq}
    \colim_{V\subseteq W}\Aa(V)
    \simeq
    \colim_{W\in \Wf}\colim_{V\in \Vf\cap \cdown{W}}\Aa(V)
    \to
    \colim_{W\in \Wf} \Aa(W).
  \end{equation}
  In the middle we compute the colimit over
  \(\{V\subseteq W\}\coloneqq\{( V,W )\in \Vf\times \Wf\mid V\subseteq W\}\)
  or, equivalently, a double colimit over \(W\) and \(V\) separately.
  Note that since \(\Vf\) hyperrefines \(\Wf\),
  the cartesian fibration \(\{V\subseteq W\} \to \Vf\) has weakly contractible fibers,
  hence is colimit cofinal;
  this is why the first map of \eqref{eq:A-on-hyperrefinements} is an equivalence. 
\end{defn}

\begin{lemma}
  \label{lem:hypref-functorial}
  Let \(\Uf, \Vf, \Wf\) be \hyperprecover{}s with
  hyperrefinements \(\Uf \hyprefines \Vf \hyprefines \Wf\)
  and let \(\Aa\) be a copresheaf defined on all of them.
  Assume that we also have the hyperrefinement \(\Uf \hyprefines \Wf\).
  Then the induced diagram
  \begin{equation}
    \begin{tikzcd}[cells={font=\everymath\expandafter{\the\everymath\displaystyle}}] % puts every cell in displaystyle instead of manually writing it everywhere
      \AonHPCov{\Uf} \simeq \colim_{\Uf} \Aa  \arrow[r] \arrow[dr] &  \colim_{\Vf} \Aa \simeq \AonHPCov{\Vf} \arrow[d]
      \\ 		&  \colim_{\Wf} \Aa \simeq \AonHPCov{\Wf}
    \end{tikzcd}
  \end{equation}
  commutes in \(\targetcat\).
\end{lemma}
\begin{proof}
  Using Equation \eqref{eq:A-on-hyperrefinements} we get the outer square of the diagram
  \begin{equation} \label{eq:"functoriality"-on-hyprefinements}
    \begin{tikzcd}[cells={font=\everymath\expandafter{\the\everymath\displaystyle}}, column sep=8ex] 
      \colim_{U\in \Uf} \Aa(U)  & \arrow[l, "\simeq" '] \colim_{U\subseteq V} \Aa(U) \arrow[r] & \colim_{V \in \Vf} \Aa(V)
      \\  & \colim_{U\subseteq V \subseteq W} \Aa(U) \arrow[u, "\circ" '] \arrow[dl, "\bullet"] \arrow[r] & \colim_{V\subseteq W} \Aa(V) \arrow[u, "\simeq" '] \arrow[d]
      \\ \arrow[uu, "\simeq"] \arrow[rr]\colim_{U\subseteq W} \Aa(U) &  & \colim_{W\in \Wf} \Aa(W)
    \end{tikzcd}.
  \end{equation}
  Firstly, note that the inner 3 squares commute. 
  Secondly, the arrow marked with a circle (``\(\circ\)'') is an equivalence: we have a pullback (of posets)
  \begin{equation}
    \begin{tikzcd}
      \{U\subseteq V \subseteq W\} \arrow[r] \arrow[d] \isCartesian & \{V\subseteq W\} \arrow[d]
      \\ \{U\subseteq V\} \arrow[r] 	& \{V \in \Vf\},
    \end{tikzcd}
  \end{equation}
  where \(U, V, W\) ranges over \(\Uf, \Vf, \Wf\) respectively. Here, the right vertical map is a cartesian fibration, and its fibers are weakly contractible exactly because \(\Vf \hyprefines \Wf\). 
  It follows that the left vertical map is also a cartesian fibration with weakly contractible fibers, and hence colimit cofinal. 
    
  From commutativity of the left upper square in \eqref{eq:"functoriality"-on-hyprefinements} it follows that the arrow marked with a filled circle (``\(\bullet\)'') is also an equivalence. 
  Thus, the diagram witnesses the evaluation of \(\Aa\) on the hyperrefinement \(\Uf \hyprefines \Wf\) as the composition of \(\Aa\) evaluated on \(\Uf \hyprefines \Vf\) and \(\Vf  \hyprefines \Wf\), and we are done. 
\end{proof}

\begin{lemma}
  \label{lem:formulas-AU-HPCov}
  Let \(\Aa\) be a copresheaf defined on \(\Uf\).
  We have the following explicit formulas:
  \begin{enumerate}
  \item
    Evaluating on a \hyperprecover{} of the form \(\Wf\hyprefines V \)
    yields the canonical map
    \begin{equation}
      \colim_{W\in \Wf}\Aa(W)\to \Aa(V),
    \end{equation}
    induced by the inclusions \(W\hookrightarrow V\).
  \item \label{it:A(presieve)=A(precover)}
    For any precover \(\Wa\) with \(\ccaps{\Wa}\subseteq \Uf\)
    the canonical map
    \begin{equation}
      \AonHPCov{\ccaps{\Wa}}\coloneqq \colim\Aa\restrict{\ccaps{\Wa}}
      \xrightarrow{\simeq}
      \colim\Aa\restrict{\Uf\cap\cdown{\Wa}}
      \eqqcolon
      \AonPCov[\Uf]{\Wa}
    \end{equation}
    is an equivalence.
  \item 
    For any precover \(\Wa\) with \(\ccaps{\Wa}\subseteq \Uf\),
    we have
    \(\AonPCov[\Uf]{\Wa}\simeq\AonPCov[\ccaps{\Wa}]{\Wa}\).
  \end{enumerate}
\end{lemma}

\begin{proof}
  \begin{enumerate}
  \item[]
  \item
    Direct simplification of the formula \eqref{eq:A-on-hyperrefinements}
    in the case where \(\Vf=\{V\}\) is a singleton.
  \item
    We claim that the inclusion
    \(\ccaps{\Wa}\hookrightarrow \Uf\cap\cdown{\Wa}\)
    is colimit cofinal.
    Indeed, for every \(U\in \Uf\cap\cdown{\Wa}\)
    the slice \(\undercat{(\ccaps{\Wa})}{U}\) is nonempty
    and closed under finite intersections,
    hence downward directed, hence weakly contractible.
  \item
    Follows directly from part \ref{it:A(presieve)=A(precover)}
    since both sides are just equivalent to \(\AonHPCov{\ccaps{\Wa}}\).
    \qedhere
  \end{enumerate}
\end{proof}

\begin{notation}
  \label{not:indep_of_poo}
  \Cref{lem:formulas-AU-HPCov} tells us that the value
  on precovers \(\Wa\)
  is independent of the ambient \posetofopens{} \(\Uf\),
  as long as \(\ccaps{\Wa}\subseteq \Uf\).
  Thus we may abbreviate
  \begin{equation}
    \AonPCov{\Wa}\coloneqq \AonPCov[\Uf]{\Wa}
  \end{equation}
  without ambiguity
  as long as \(\Aa\) is defined at least on \(\ccaps{\Wa}\).
\end{notation}

\begin{war}
  \label{war:cover-vs-hypercover}
  For a \hyperprecover{} \(\Wf\subseteq\Uf\)
  with underlying precover \(\Wa\coloneqq\fgtposet{\Wf}\),
  much headache can arise by conflating
  \begin{equation}
    \AonHPCov{\Wf}
    \quad
    \text{and}
    \quad
    \AonPCov[\Uf]{\Wa}.
  \end{equation}
  We emphasize that these two do \emph{not} agree in general.
  In other words, when evaluating a copresheaf on a collection of opens,
  we have to carefully distinguish
  whether we view the collection as a poset (=\hyperprecover{})
  or just a set (=precover).

  Luckily, \Cref{lem:formulas-AU-HPCov} states that if \(\Wf\) is a presieve,
  then the two ways of evaluating agree;
  hence at least in that case we do not have to worry about this distinction.
\end{war}

It is sometimes convenient to compute the value on a precover
as a colimit of its \v{C}ech nerve,
which in our context most naturally takes the form of a punctured cube.

\begin{lemma}
  \label{lem:A(U)-cubical}
  Let \(\Uf\) be a \posetofopens{} 
  and \(\Va\) a precover with \(\ccaps{\Va}\subseteq \Uf\).
  Let \(\Aa \colon \Uf\to\targetcat\) be a copresheaf.
  Let \(\beta\) be an ordinal and
  \(V_\bullet\colon \beta\twoheadrightarrow \Va\)
  a surjection, so that we can write \(\Va=\{V_i\}_{i<\beta}\).
  Then the punctured cube
  \begin{equation}
    \label{eq:cube-for-presieve-value}
    \Pfinop{\beta}\setminus\{\emptyset\}\to \Uf\cap\cdown{\Va},
    \quad\quad I\mapsto V_I\coloneqq \bigcap_{i\in I} V_i
  \end{equation}
  (which is well defined because \(\Uf\) contains \(\ccaps{\Va}\))
  is colimit cofinal.
  In particular, it induces an equivalence
  \begin{equation}
    \colim_{\beta\supseteq I\supsetneq \emptyset}\Aa(V_I)
    \xrightarrow{\simeq}
    \colim_{U\in \Uf\cap\cdown{\Va}} \Aa(U)
    =
    \AonPCov[\Uf]{\Va}
    .
  \end{equation}
\end{lemma}

\begin{proof}
  Fix an open \(U\in \Uf\cap\cdown{\Va}\) and consider the slice poset
  \begin{equation}
    K\coloneqq 
    \{I\in \Pfinop{\beta}\setminus\{\emptyset\}
    \mid
    U\subseteq V_I\}
  \end{equation}
  under \(U\).
  This poset is nonempty because by definition there exists a \(V\in \Va\)
  with \(U\subseteq V\), so that \(\{i\}\in K\) for every \(i<\beta\)
  with \(V_i=V\); such \(i\) exists because \(V_\bullet\) is surjective.
  The poset \(K\) is also downward directed
  since for every \(I,I'\in K\),
  we also have \(I\cup I'\in K\).
  Every directed poset is weakly contractible,
  so the result follows.
\end{proof}

\subsection{Locality and descent}
In this subsection we first define what we mean by locality for refinements and hyperrefinements. Then we collect some general results regarding situations where locality follows, which will be useful later on.

\begin{defn} \label{defn:LocalRef}
  Let \(\Uf\) be a \posetofopens{} and \(\Aa\colon \Uf\to \targetcat\) a copresheaf.
  \begin{itemize}
  \item
    Let \(\Wa'\to \Wa\) be a refinement with \(\Wa',\Wa\subseteq \Uf\).
    We say that the refinement \(\Wa'\to\Wa\) is \emph{\(\Aa\)-local},
    if the induced map \(\AonPCov[\Uf]{\Wa'}\to\AonPCov[\Uf]{\Wa}\)
    is an equivalence in \(\targetcat\).
  \item
    A hyperrefinement \(\Wf' \hyprefines \Wf\)
    of \hyperprecover{}s \(\Wf',\Wf\subseteq \Uf\)
    is called \emph{\(\Aa\)-local},
    if the induced map \(\AonHPCov{\Wf'} \to \AonHPCov{\Wf}\) is an equivalence. 
  \end{itemize}
\end{defn}

\begin{rem}
  Commonly one says that \(\Aa\) \emph{satisfies descent}
  with respect to a (pre-)cover \(\Ua\to U\) (or a hyper(pre-)cover \(\Uf\hyprefines U\))
  if that (pre-)cover (or hyper(pre-)cover) is \(\Aa\)-local.
\end{rem}

\begin{rem}
  By definition, a refinement \(\Wa'\to\Wa\) is \(\Aa\)-local
  if and only if
  the inclusion
  \(\Uf\cap\cdown{\Wa'}\hookrightarrow\Uf\cap\cdown{\Wa}\)
  of \hyperprecover{}s
  is an \(\Aa\)-local hyperrefinement.
  Note that this depends on the ambient \posetofopens{} \(\Uf\)
  unless \(\Wa'\) and \(\Wa\) are presieves.
\end{rem} 

The easiest way to produce local hyperrefinements is to exhibit an inverse:

\begin{lemma}
  Let \(\Vf,\Wf\) be two hyperprecovers with hyperrefinements
  \(\Vf\hyprefines \Wf\) and \(\Wf\hyprefines \Vf\).
  Then both those hyperrefinements are \(\Aa\)-local for every copresheaf
  defined on \(\Vf\) and \(\Wf\).
\end{lemma}

\begin{proof}
  Noting that the identity hyperrefinements \(\Wf\hyprefines \Wf\)
  and \(\Vf\hyprefines \Vf\) yield the identity map on colimits,
  we can apply \Cref{lem:hypref-functorial}
  to the hyperrefinements
  \begin{equation}
    \begin{tikzcd}
      \Vf\ar[dr,"="]\ar[r,"\mathrm{h}" description]&\Wf\ar[d,"\mathrm{h}" description] \ar[dr,"="]
      \\
      &\Vf\ar[r,"\mathrm{h}" description]& \Wf
    \end{tikzcd}
  \end{equation}
  to see that the induced map \(\AonHPCov{\Wf}\to \AonHPCov{\Vf}\)
  has both a left and a right inverse, hence is an equivalence.
  And analogously for \(\Vf\hyprefines\Wf\).
\end{proof}

\begin{ex}
  \label{rem:iso-of-refinements}
  Some useful examples of invertible hyperrefinements are: 
  \begin{itemize}
  \item
    For every \(V\in\Uf\),
    the inclusion \(\{V\}\hookrightarrow \Uf\cap\cdown{V}\)
    has an inverse
    \(\Uf\cap \cdown{V} \hyprefines V\).
  \item
    For every presieve \(\Vf\),
    the inclusion \(\Vf\hookrightarrow \cdown{\Vf}\)
    has an inverse
    \(\cdown{\Vf}\hyprefines{\Vf}\)
    since for every \(U\in \cdown{\Vf}\),
    the poset \(\{V\in \Vf\mid U\subseteq V\}\)
    is non-empty by definition and downward directed by \(\cap\).
  \end{itemize}
\end{ex}

\begin{notation}
  For two precovers \(\Ua\) and \(\Va\), we write
  \begin{equation}
    \Ua\wedge \Va\coloneqq \{U\cap V\mid U\in \Ua, V\in \Va\}
  \end{equation}
  for (the specific choice of) a coarsest common refinement.
  In terms of the associated sieves we have
  \begin{equation}
    \cdown{(\Ua\wedge\Va)}=\cdown{\Ua}\cap\cdown{\Va}.
  \end{equation}
  We use the same notation \(\Uf\wedge\Wf\)
  also for \posetsofopens{}.
\end{notation}

A useful lemma for detecting local refinements is the following:

\begin{lemma}
  \label{lem:coll-restrict-local}
  Let \(\Ua\) be a precover and \(V\subseteq X\) an open subset.
  Let \(\Aa\) be a copresheaf defined both on \(\ccaps{\Ua}\)
  and on \(\ccaps{\Ua}\wedge \{V\}\).
  Assume that for each \(U\in\ccaps{\Ua}\)
  the inclusion \(U\cap V\to U\) is \(\Aa\)-local.
  Then the refinement
  \begin{equation}
    \Ua\wedge\{V\}\coloneqq \{U\cap V\mid U\in\Ua\} \to \Ua
  \end{equation}
  is also \(\Aa\)-local.
\end{lemma}

\begin{proof}
  Choose a surjection
  \(U_\bullet\colon\beta\twoheadrightarrow\Ua\),
  which then yields the surjection
  \(U_\bullet\cap V\colon \beta\twoheadrightarrow\Ua\wedge \{V\}\).
  With the notation \(U_I\coloneqq \bigcap_{i\in I}U_i\) for all finite \(I\subset \beta\),
  we have a commutative square
  \begin{equation}
    \cdsquareNA
    {\colim\limits_{\beta\supseteq I\supsetneq\emptyset}\Aa(U_I\cap V)}
    {\colim\limits_{\beta\supseteq I\supsetneq\emptyset}\Aa(U_I)}
    {\AonPCov{\Ua\wedge \{V\}}}{\AonPCov{\Ua}},
  \end{equation}
  using \(\bigcap_{i\in I}(U_i\cap V)=U_I\cap V\) in the top left corner.
  The vertical maps are equivalences by \Cref{lem:A(U)-cubical}
  and the top horizontal map is an equivalence by assumption
  (because each \(U_I\) lies in \(\ccaps{\Ua}\)).
  Hence the lower horizontal map is also an equivalence, as desired.
\end{proof}

\begin{defn}
  \label{def:strongly-local}
  A hyperrefinement \(\Vf\hyprefines\Uf\) is called \emph{strongly \(\Aa\)-local}
  if for every \(U\in \Uf\)
  the induced \hyperprecover{} \(\Vf\cap\cdown{U}\hyprefines U\) is \(\Aa\)-local.
  This definition makes sense whenever \(\Aa\) is a copresheaf defined
  on \(\Vf\) and \(\Uf\).
\end{defn}

\begin{lemma}
  \label{lem:strongly-local-local}
  Let \(\Aa\) be a copresheaf defined on all opens that appear.
  \begin{enumerate}
  \item
    An inclusion \(\Vf\hookrightarrow \Uf\) of \posetsofopens{} is
    a strongly \(\Aa\)-local hyperrefinement
    if and only if
    the copresheaf \(\Aa\restrict{\Uf}\)
    is a left Kan extension of its restriction \(\Aa\restrict{\Vf}\).
  \item
    Every strongly \(\Aa\)-local hyperrefinement is \(\Aa\)-local.
  \end{enumerate}
\end{lemma}

\begin{proof}
  \begin{enumerate}
  \item[]
  \item
    By the pointwise formula,
    \(\Aa\restrict{\Uf}\) is a left Kan extension of \(\Aa\restrict{\Vf}\)
    if and only if each comparison map
    \(\colim \Aa\restrict{\Vf\cap \cdown{U}}\to \Aa(U) \)
    is an equivalence
    if and only if
    each \hyperprecover{}
    \(\Vf\cap\cdown{U}\hyprefines U\) is \(\Aa\)-local;
    that is the definition of strong \(\Aa\)-locality.
  \item
    Let \(\Wf\hyprefines \Vf\) be a strongly \(\Aa\)-local refinement.
    Then by definition, for each \(V\in \Vf\) the
    \hyperprecover{} \(\Wf\cap\cdown{V}\) is \(\Aa\)-local,
    which says precisely that the right map
    in the formula \eqref{eq:A-on-hyperrefinements}
    is a colimit of equivalences, hence an equivalence.
    Therefore the map \(\Aa(\Wf)\to\Aa(\Vf)\) is an equivalence, as desired.
    \qedhere
  \end{enumerate}
\end{proof}

\begin{lemma}
  \label{lem:strongly-local-to-sieve}
  If \(\Vf\hookrightarrow \Uf\) is a strongly \(\Aa\)-local inclusion of
  \posetsofopens{} and \(\Wa\) is a precover,
  then the inclusion \(\Vf\cap \cdown{\Wa}\hookrightarrow \Uf\cap\cdown{\Wa}\)
  is also strongly \(\Aa\)-local.
\end{lemma}

\begin{proof}
  Assume that \(\Vf\hookrightarrow\Uf\) is strongly \(\Aa\)-local.
  For every \(U\in\Uf\cap\cdown{\Wa}\)
  we have \(\cdown{U}\subseteq \cdown{\Wa}\),
  hence the \hyperprecover{}
  \begin{equation}
    (\Vf\cap\cdown{\Wa})\cap \cdown{U}=\Vf\cap\cdown{U}\hyprefines U
  \end{equation}
  is \(\Aa\)-local by assumption.
\end{proof}

\subsection{Weiss cosheaves and hypercosheaves}

In the standard Grothendieck topology on \(\openX\) a cover
corresponds to an ordinary open cover.
However, for factorization algebras we do
not want a cosheaf condition for ordinary open covers but rather so-called
\emph{Weiss covers}.
In \cite{CG},
factorization algebras are defined to be Weiss cosheaves
without an extra multiplicativity condition we will impose,
so one could alternatively call them
\emph{non-multiplicative factorization algebras}.

\begin{defn} \label{defn:JinftyCover}
  A cover \(\Wa\to U\) is called a \emph{Weiss cover}
  if for each finite subset \(S\subset U\)
  there exists a \(W \in \Wa\) with \(S\subseteq W\).
\end{defn}

\begin{ex}
  A cover \(\Wa\to U\) is called degenerate if \(U\in \Wa\).
  Every degenerate cover is trivially Weiss.
\end{ex}

If a Weiss cover is not degenerate, then it is necessarily infinite.
In this sense, Weiss covers are typically very big.
There are two tautological ways of constructing Weiss covers for every space \(X\),
which in some sense embody two different extremes:

\begin{ex}
  Choose an infinite set \(\Ka\) of pairwise disjoint closed subsets of \(X\).
  Then \(\{X\setminus K \mid K\in \Ka\}\) is a Weiss cover of \(X\).
  For example, \(\Ka\coloneqq \{\{x\}\mid x \in X\}\)
  could simply be the set of all singletons of \(X\).
\end{ex}

\begin{ex}\label{ex:WeissHausdorff}
  Let \(\Bf\) be a basis for the
  topology of \(X\). Fix a subset \(U\subseteq X\).
  For each finite set \(S\subset U\) we may choose
  pairwise disjoint opens \(B_S^s \in \Bf\) (for \(s\in S\))
  with \(s \in B_S^s \subset U\)
  (since \(X\) is Hausdorff and \(\Bf\) is a basis);
  then we set \(B_S\coloneqq \bigdisjun_{s\in S} B_S^s\).
  By construction, the set \(\Wa\coloneqq \{B_S\mid S\subset U \text{ finite}\}\)
  is a Weiss cover of \(U\).
\end{ex}

\begin{rem}
  If \(\Wa\to U\) is a Weiss cover,
  then \(\Wa\wedge\{V\}\to V\) is again a Weiss cover for every open subset \(V\subseteq U\).
\end{rem}

Given a \posetofopens{} \(\Uf\) we say that a precover \(\Wa\to U\)
or \hyperprecover{} \(\Wf\hyprefines U\) are \emph{in \(\Uf\)}
if all the involved open sets are contained in \(\Uf\),
i.e., if \(\Wa\cup \{U\}\subseteq \fgtposet{\Uf}\)
or \(\Wf\cup\{U\}\subseteq \Uf\), respectively.

\begin{defn}
  \label{defn:JinftyCondition} %= only used to reference Weiss cosheaves, not weiss hypercosheaves!
  Let \(\Uf\) be a \posetofopens{} and \(\Aa\) a copresheaf defined on \(\Uf\).
  We say that \(\Aa\) is \emph{a Weiss cosheaf on \(\Uf\)}
  if every Weiss cover \(\Wa\to U\) in \(\Uf\) is \(\Aa\)-local,
  i.e.\ if we have an equivalence
  \begin{equation} \label{eq:WeissCosheafConditionExplicit}
    \AonPCov[\Uf]{\Wa} \coloneqq
    \colim \Aa\restrict{\Uf\cap\cdown{\Wa}} \xrightarrow{\simeq} \Aa(U). 
  \end{equation}    
  We denote the full \(\infty\)-category of
  Weiss cosheaves on \(\Uf\) by 
  \begin{equation}
    \cShvWeiss{\Uf} \subset \Fun(\Uf, \targetcat).
  \end{equation}
\end{defn}

It will be convenient to also work with hyper-analogues of Weiss covers and cosheaves. 

\begin{defn}\label{defn:weiss-hypercover-and-strongly-weiss}
\begin{itemize}
	\item[]
	\item A \hyperprecover{} \(\Wf\hyprefines U\) is called a \emph{Weiss hypercover} if 
		the precover \(\fgtposet{\Wf} \to U\) is a Weiss cover 
		and for each \(W \in \ccaps{\Wf}\) we have a Weiss cover 
		\(\fgtposet{(\Wf\cap\cdown{W})}\to W\). 
	\item A hyperrefinement \(\Vf\hyprefines \Uf\)
    		is called \emph{strongly Weiss} if for each \(U\in \Uf\)
    		we have a Weiss hypercover \(\Vf\cap\cdown{U}\hyprefines U\).
\end{itemize}
\end{defn}

\begin{rem}
  Note that for each \(W\in \ccaps{\Wf}\) and
  \(W'\in \ccaps{(\Wf\cap\cdown{W})}=\ccaps{\Wf}\cap\cdown{W}\)
  we have \((\Wf\cap\cdown{W})\cap \cdown{W'}=\Wf\cap\cdown{W'}\).
  Hence for a Weiss hypercover \(\Wf\hyprefines U\), each
  \(\Wf\cap\cdown{W}\hyprefines W\) does not just have an underlying Weiss cover,
  but is even itself a Weiss hypercover again.

  Thus one could define the class of Weiss hypercovers coinductively:
  A \hyperprecover{} \(\Wf \hyprefines U\) is a Weiss hypercover if
	the precover \(\fgtposet{\Wf} \to U\) is a Weiss cover 
	and for each \(W\in \ccaps{\Wf}\) we have a Weiss hypercover
  \(\Wf\cap\cdown{W}\hyprefines W\). 
\end{rem}

\begin{rem}
  \label{rem:presieve-hyper-Weiss}
  Let \(\Wf\hyprefines U\) be a \hyperprecover{}.
  If it is a Weiss hypercover, then by definition
  the underlying precover \(\fgtposet{\Wf}\to U\) is a Weiss cover.

  The converse does \emph{not} hold in general;
  but it \emph{does} hold when \(\Wf\) is a presieve
  because then for each \(W\in \ccaps{\Wf}\) we have a degenerate cover
  \(\Wf\cap\cdown{W}\to W\).
\end{rem}

Before we define Weiss hypercosheaves we give a short lemma showing how Weiss covers and Weiss hypercovers interact. 
\begin{lemma}
  \label{lem:Weiss-to-hyper}
  Let \(\Wa\to U\) be a Weiss cover.
  \begin{enumerate}
  \item
    We have a Weiss hypercover \(\ccaps{\Wa}\hyprefines U\).
  \item
    The inclusion \(\cdown{\Wa}\hookrightarrow \cdown{U}\) is strongly Weiss.
  \end{enumerate}
\end{lemma}

\begin{proof}
  \begin{enumerate}
  \item[]
  \item
    Follows from \Cref{rem:presieve-hyper-Weiss}.
  \item
    For every \(V\subseteq U\) we have a restricted Weiss cover
    \begin{equation}
      \label{eq:sieve-cap=sieve-wedge}
      \cdown{\Wa}\cap\cdown{V}=\cdown{\Wa}\wedge\{V\}\to V,
    \end{equation}
    which by \Cref{rem:presieve-hyper-Weiss}
    yields the desired Weiss hypercover
    \(\cdown{\Wa}\cap\cdown{V}\hyprefines V\)
    since it is a presieve (even a sieve).
    \qedhere
  \end{enumerate}
\end{proof}

\begin{rem}
  Note that we could not replace \(\cdown{W}\) by \(\ccaps{W}\)
  in the second statement since the equation in \eqref{eq:sieve-cap=sieve-wedge}
  would not hold in that case.
\end{rem}

\begin{defn}
  Let \(\Uf\) be a \posetofopens{} and \(\Aa\) a copresheaf defined on \(\Uf\).
  We say that \(\Aa\) is 
  \emph{a Weiss hypercosheaf on \(\Uf\)}
  if every Weiss hypercover \(\Wf\hyprefines U\) in \(\Uf\) is \(\Aa\)-local, i.e. if 
  \begin{equation}
    \AonHPCov[\Uf]{\Wf} \coloneqq \displaystyle \colim_{W\in \Wf} \Aa(W) \xrightarrow{\simeq} \Aa(U).
  \end{equation}
  We denote the full \(\infty\)-category of
  Weiss hypercosheaves on \(\Uf\) by
  \begin{equation}
    \cHShvWeiss{\Uf} \subset \Fun(\Uf, \targetcat).
  \end{equation}
\end{defn}

\begin{rem}
  \label{rem:Weiss-cosheaf-Weiss-presieves}
  Let \(\Uf\) be a presieve.
  Then by \Cref{lem:formulas-AU-HPCov}
  a Weiss cover \(\Wa\to U\) in \(\Uf\) evaluates the same as
  the Weiss hypercover \(\ccaps{\Wa}\hyprefines U\).
  So while a Weiss hypercosheaf is characterized by \emph{all}
  Weiss hypercovers \(\Wf\hyprefines U\) being \(\Aa\)-local,
  to be a Weiss cosheaf only requires it for Weiss \emph{presieves},
  i.e., Weiss hypercovers which are closed under intersection.
  In particular,
  we have a (fully faithful) inclusion
  \begin{equation}
    \cHShvWeiss{\Uf}\subseteq\cShvWeiss{\Uf}.
  \end{equation}
\end{rem}

\begin{war}
  If \(\Uf\) is not a presieve then
  it might happen that not every Weiss hypercosheaf on \(\Uf\)
  is a Weiss cosheaf.
  For example, let \(X\coloneqq \reals\)
  and \(\Uf\coloneqq\{X\}\cup\{X\setminus\{x\}\mid x\in X\}\).
  Then there are no nondegenerate hypercovers in \(\Uf\),
  so that the hypercosheaf condition is vacuous.
  However, the tautological Weiss cover
  \(\Uf\setminus\{X\}\to X \)
  imposes the non-trivial constraint
  \(\Aa(X)\simeq\coprod_{x\in X}\Aa(X\setminus\{x\})\)
  on Weiss cosheaves \(\Aa\). 
\end{war}

\subsection{Multiplicative and factorizing bases}
In this subsection we will introduce a few notions of bases that prove convenient for example when extending factorization algebras to bigger \posetsofopens{} (or all of \(\openX\)). 
The notion of a factorizing basis is due to Costello--Gwilliam,
see \cite[\S 7 Definition~2.1.1]{CG}.
We need the following variants:

\begin{defn}
  \label{defn:mult_fact_basis}
  Let \(\Uf\) be a \posetofopens{}.
  A \emph{multiplicative basis} of \(\Uf\) is a subposet \(\Bf\subseteq \Uf\)
  which is
  \begin{itemize}
  \item
    a basis for the topology of \(\collunion \Uf\),
    i.e., for each \(x\in V\in \cdown{\Uf}\) there is a \(B \in \Bf\)
    with \(x \in B \subseteq V\);
  \item
    closed under disjoint unions subordinate to \(\Uf\);
    this means
    \(\emptyset\in \Bf\) and
    for each \(B,B'\in \Bf\) with \(B\disjun B' \in \cdown{\Uf}\)
    we have \(B\disjun B' \in \Bf\).
  \end{itemize}
  A multiplicative basis \(\Bf\) of \(\Uf\) is called a \emph{factorizing basis}
  if it is additionally
  \begin{itemize}
  \item
    closed under intersections,
    i.e., \(B\cap B'\in \Bf\) for all \(B,B'\in \Bf\).
  \end{itemize}
\end{defn}

\begin{notation}
  When we refer to a multiplicative/factorizing basis of an open \(U\subseteq X\),
  we mean a multiplicative/factorizing basis of \(\Uf\coloneqq\cdown{U}\).
\end{notation}

\begin{constr}
  For any \posetofopens{} \(\Uf\), we denote by
  \begin{equation}
    \cdisj{\Uf}\coloneqq
    \{U_1\disjun\cdots\disjun U_n
    \mid
    n\in \naturals; \, U_1,\dots,U_n\in \Uf \text{ pairwise disjoint}\}
  \end{equation}
  the completion under finite disjoint unions.
  Note that \(\cdisj{\Uf}\) always contains the empty set
  \(\emptyset\) which is the disjoint union of the empty list.
\end{constr}

\begin{ex}
  Let \(X\) be a smooth manifold and denote by
  \(\Bf\coloneqq\diskX\subset \open{X}\)
  the \posetofopens{}
  consisting of all disjoint unions of embedded disks;
  see also \Cref{not:disks} below.
  Then \(\Bf\) is a multiplicative basis of \(X\)
  but \emph{not} a factorizing basis because the intersection
  of two embedded disks need not be a disjoint union of embedded disks again;
  for example the intersection of the embedded disks
  \(\reals^3\setminus \reals_{\geq 0}\)
  and
  \(\reals^3\setminus \reals_{\leq 0}\)
  in \(\reals^3\)
  is \(\reals^3\setminus \reals= \reals\times (\reals^2\setminus \{0\})\)
  which has non-trivial fundamental group hence isn't a disjoint union of disks.
\end{ex}

It is sometimes convenient to work with \posetsofopens{}
which are closed under disjoint decompositions.
\begin{defn}
  \label{defn:DecomposablePresieve}
  We say that a \posetofopens{} \(\Uf\) is \emph{decomposable}
  if for any open \(U \disjun U' \in \Uf\)
  with \(U\neq \emptyset\neq U'\)
  it follows that \(U\in \Uf\) and \(U'\in \Uf\). 
\end{defn}

\begin{ex}
	Both \(\cdiskX\) and \(\diskX\) are examples of decomposable \posetsofopens{}. 
	Recall that their opens are contractible embedded disks,
  and finite disjoint unions of embedded disks, respectively.
\end{ex}

\begin{rem}
  If \(\Uf\) is a decomposable \posetsofopens{}
  and \(\Bf\subset\Uf\) is a multiplicative/factorizing basis,
  then we can always add all non-empty disjoint union components of elements of \(\Bf\)
  to obtain a multiplicative/factorizing basis which is decomposable.

  Thus we will sometimes assume without loss of generality that
  our multiplicative/factorizing bases are decomposable.
\end{rem}

The main selling point of multiplicative bases
is that they provide a large supply of Weiss (hyper)covers.

\begin{lemma}
  \label{lem:basis-to-cover}
  Let \(\Uf\) be a \posetofopens{} and 
  let \(\Bf\subseteq \Uf\) be a multiplicative basis for \(\Uf\).
  Then for each \(U\in\cdown{\Uf}\)
  we have a Weiss hypercover \(\Bf\cap\cdown{U}\hyprefines U\).
\end{lemma}

\begin{proof}
  It suffices to check that for each
  \(V \in \cdown{\Uf}\)
  we have a Weiss cover \(\fgtposet{(\Bf\cap\cdown{V})}\to V\);
  we can then apply this to \(V\coloneqq U\)
  and each \(V \in \ccaps{(\Bf\cap\cdown{U})}\).
  
  For this, let \(S\subset V\) be a finite set.
  Choose pairwise disjoint open subsets \(B_s\in \Bf\cap \cdown{V}\)
  with \(s\in B_s\)
  (indexed by \(s\in S\)),
  which is possible because \(\Bf\) is a basis of \(\Uf\) and \(X\) is Hausdorff.
  Then the disjoint union \(B\coloneqq \bigdisjun_{s\in S} B_s\)
  is contained in \(V\),
  hence lies in \(\Bf\cap\cdown{V}\) by multiplicativity.
  By construction, it contains \(S\), which concludes the proof.
\end{proof}

More generally, we can input a Weiss cover to get more Weiss (hyper)covers.

\begin{lemma}
  \label{lem:factorizing-is-Weiss}
  Let \(\Bf\) be a multiplicative basis for \(\Uf\).
  Let \(\Wa\to U\) be a Weiss cover in \(\cdown{\Uf}\).
  \begin{enumerate}
  \item
    \label{it:fact-basis-Weiss-cover}
    We have a Weiss hypercover \(\Bf\cap \cdown{\Wa}\hyprefines U\).
  \item
    \label{it:fact-basis-Jinfty-inclusion-cover}
    The inclusion
    \(\Bf\cap \cdown{\Wa}\hookrightarrow \cdown{U}\)
    is strongly Weiss.
  \end{enumerate}
\end{lemma}

\begin{proof}
  \begin{enumerate}
  \item[]
  \item %[\ref{it:fact-basis-Weiss-cover}]
    Let \(S\subset U\) be a finite set.
    Since \(\Wa\to U\) is a Weiss cover,
    we find a \(W\in \Wa\subset \cdown{U}\) with \(S\subset W\).
    Since \(\fgtposet{\left(\Bf\cap\cdown{W}\right)}\to W\) is a Weiss cover of \(W\)
    by \Cref{lem:basis-to-cover},
    we find a \(B\in \Bf\cap\cdown{W}\subseteq \Bf\cap\cdown{\Wa}\)
    with \(S\subset B\).
    This shows that we have a Weiss cover \(\fgtposet{(\Bf\cap\cdown{\Wa})}\to U\).
    Moreover, for each
    \(W\in \ccaps{(\Bf\cap\cdown{\Wa})}=\ccaps{\Bf}\cap\cdown{\Wa}\)
    we have the Weiss cover
    \begin{equation}
      \fgtposet{((\Bf\cap\cdown{\Wa})\cap \cdown{W})}
      =\fgtposet{(\Bf\cap\cdown{W})}\to W
    \end{equation}
    (again by \Cref{lem:basis-to-cover}),
    hence we have the desired Weiss hypercover.
  \item %[\ref{it:fact-basis-Jinfty-inclusion-cover}]
    Let \(V\in \cdown{U}\).
    We can restrict the Weiss cover \(\Wa\to U\)
    to a Weiss cover
    \(\Wa\wedge \{V\}\to V\).
    Hence we have the desired Weiss hypercover
    \begin{equation}
      (\Bf\cap\cdown{\Wa})\cap\cdown{V}
      =
      \Bf\cap\cdown{(\Wa\wedge\{V\})}
      \hyprefines
      V
    \end{equation}
    by part~\ref{it:fact-basis-Weiss-cover}.
    \qedhere
  \end{enumerate}
\end{proof}

Finally, we come to the main result of this subsection which essentially outlines how factorizing bases interact nicely with the Weiss condition when left Kan extending. This will be important in \Cref{subsect:ExtendingUFactAlgs} where we prove that factorization algebras defined on a factorizing basis extends. 

\begin{prop}[Unique extension from a basis]
  \label{prop:UniqueExtFromBasis}
  Let \(\Aa\colon \Uf\to \targetcat\) be a copresheaf
  on the \posetofopens{} \(\Uf\).
  Let \(\Bf\subseteq\Uf\) be a factorizing basis of \(\Uf\).
  \begin{enumerate}[label=(\arabic*), ref=(\arabic*)]
  \item
    \label{it:Weiss-is-LKE-from-Weiss}
    Every Weiss cosheaf on \(\Uf\)
    is a Weiss cosheaf on \(\Bf\)
    and a left Kan extension along \(\Bf\hookrightarrow \Uf\).
  \item
    \label{it:LKE-from-Weiss-is-Weiss}
    Assume that \(\Uf\) is a presieve.
    Left Kan extension along \(\Bf\hookrightarrow\Uf\)
    from a Weiss cosheaf on \(\Bf\)
    yields a Weiss cosheaf on \(\Uf\).
  \end{enumerate}
\end{prop}

\begin{proof}
  \begin{enumerate}
  \item[]
  \item
    Since \(\Bf\) is a presieve,
    the evaluation of \(\Aa\) on a Weiss cover
    \(\Wa\to U\) computed in \(\Bf\)
    agrees with the evaluation in \(\Uf\)
    (see \Cref{not:indep_of_poo});
    hence restriction of a Weiss cosheaf is a Weiss cosheaf again.

    For every \(U\in \Uf\),
    the pointwise left Kan extension condition amounts to
    the hyperrefinement
    \(\Bf\cap\cdown{U}\hyprefines{U}\) being \(\Aa\)-local.
    This is a Weiss hypercover by \Cref{lem:basis-to-cover}
    and also a presieve because \(\Bf\) is a presieve.
    Therefore it is local for every Weiss cosheaf
    (see \Cref{rem:Weiss-cosheaf-Weiss-presieves}).
  \item
    Assume that \(\Aa\restrict{\Uf}\) is a left Kan extension
    of a Weiss cosheaf \(\Aa\restrict{\Bf}\).
    Since \(\Uf\) is a presieve,
    we only have to show that every Weiss presieve
    \(\Wf\hyprefines U\) is \(\Aa\)-local
    (see \Cref{rem:Weiss-cosheaf-Weiss-presieves}).
    We have the following square of hyperrefinements
    \begin{equation}
      \cdsquareOpt
      {\Bf\cap\cdown{\Wf}}
      {\Bf\cap\cdown{U}}
      {\Wf}
      {U}
      {hookrightarrow}
      {"\mathrm{h}" description}
      {"\mathrm{h}" description}
      {"\mathrm{h}" description}
    \end{equation}
    where the left vertical map is the composite
    \(\Bf\cap\cdown{\Wf}\hookrightarrow \cdown{\Wf}\hyprefines{\Wf}\)
    which exists since \(\Wf\) is a presieve. 
    Note that by \Cref{rem:SomeHyprefsCompose} the composite is again a hyperrefinement 
    making the diagram commute, and \Cref{lem:hypref-functorial} 
    ensures that applying \(\Aa\) yields a commuting diagram in \(\targetcat\). 
       \begin{itemize}
    \item
      The right vertical map is
      \(\Aa\)-local because \(\Aa\) is
      a pointwise left Kan extension along \(\Bf\hookrightarrow \Uf\)
    \item
      We claim that the upper horizontal inclusion is strongly \(\Aa\)-local.
      By \Cref{lem:factorizing-is-Weiss},
      we have the pointwise Weiss hypercovers
      \(\Bf\cap\cdown{(\Wf\wedge \{V\})}\hyprefines{V}\)
      in \(\Bf\)
      for all \(V\in \cdown{U}\), hence a fortiori for all \(V\in\Bf\cap\cdown{U}\).
      Since they are sieves relative to \(\Bf\),
      they take the same value
      as their underlying Weiss cover
      (computed in \(\Bf\)).
      Hence they are \(\Aa\)-local because \(\Aa\restrict{\Bf}\) is a Weiss cosheaf.
    \item
      The left vertical hyperrefinement is strongly \(\Aa\)-local,
      because each
      \(\Bf\cap\cdown{W}\hyprefines{W}\)
      is again an instance of the pointwise left Kan extension.
      \qedhere
    \end{itemize}
  \end{enumerate}
\end{proof}

\begin{rem}
  In both conclusions of part~\ref{it:Weiss-is-LKE-from-Weiss}
  of \Cref{prop:UniqueExtFromBasis}
  we make crucial use of the fact that \(\Bf\) is intersection-closed.
  For example, the value on a Weiss cover \(\Wa\to U\)
  might be different when computed in \(\Bf\) or in \(\Uf\)
  so that the restriction of a Weiss cosheaf might not be a Weiss cosheaf again.
  However, the statement~\ref{it:Weiss-is-LKE-from-Weiss}
  of \Cref{prop:UniqueExtFromBasis}
  still holds with exactly the same proof for an arbitrary
  (not necessarily intersection-closed)
  multiplicative basis \(\Bf\subset \Uf\)
  if we replace ``Weiss cosheaf'' by ``Weiss hypercosheaf''
  both in the hypothesis and in the conclusion.

  One might also expect a version of part~\ref{it:LKE-from-Weiss-is-Weiss}
  for Weiss hypercosheaves.
  We do not know of such a statement
  and note that our proof can not be adapted to yield one:
  indeed if \(\Wf\hyprefines U\) is an arbitrary Weiss hypercover
  which is not a presieve,
  there does not usually exist the hyperrefinement
  \(\Bf\cap\cdown{\Wf}\hyprefines \Wf\) on which our proof crucially relies.
\end{rem}

\subsection{Remarks about terminology}
\label{rem:terminology}

A few remarks are in order about the terminology we chose to adopt,
and how it relates to and differs from other sources.

\begin{rem}
  \label{rem:hypercover-nonstandard}
  In topos theory one can express sheaf or hypersheaf conditions
  on a site \(H\)
  via three different types of objects:
  sets/families of morphisms, simplicial diagrams, and sieves.
  \begin{itemize}
  \item
    The covering families are certain sets \(\{u_i\to x\}\) of morphisms.
  \item
    A hypercover is defined to be a simplicial diagram
    in a site, such that all matching maps are covers;
    a \v{C}ech diagram is a special type of hypercover
    in which all but the \(0\)th matching map are equivalences.
  \item
    A sieve on \(x\) is a subfunctor of the hom-functor \(H(-,x)\),
    which just translates to a collection of maps \(u\to x\)
    closed under precomposition.
  \end{itemize}
  In this paper we implicitly work with topologies on \(\openX\)
  (mainly the Weiss topology),
  but express the various notions in a slightly more ad-hoc way
  using sets and posets of opens;
  completely avoiding simplicial diagrams.
  With this adjustment,
  (what we call) presieves are the analog of \v{C}ech nerves
  and (what we call) hypercovers
  are the analog of the usual hypercovers.
  We refer to Section~4.4 of \cite{DuggerIsaksen}, 
  where this relation is explained in more detail;
  note that what we call hypercovers corresponds to their complete covers. 
\end{rem}

\begin{rem} \label{rem:Equivalent-colimits-Weiss-condition}
  In the literature on (co)sheaves and factorization algebras,
  the reader might encounter formulas that look similar,
  but a bit different from ours.
  To preempt this confusion,
  we quickly sketch how some of these possible different formulas relate to each other.
  Let \(\Uf\) be a category and \(\Aa\colon \Uf\to\targetcat\)
  a copresheaf on \(\Uf\).
  Of course for us \(\Uf\subseteq \open{X}\) is always a \posetofopens{} of \(X\),
  but the discussion below works in greater generality.
 
  \begin{enumerate}
    \newcommand{\sieve}{\mathrm{sieve}}
    \newcommand{\Cechraw}[1]{\check{C}(#1)}
  \item \label{it:covering-sieve-condition}
    Let \(\Va=\{V_i\to U\}_{i<\beta}\) be a precover of some \(U\in\Uf\).
    From the perspective of topos theory,
    the ``official'' way of evaluating \(\Aa\) on \(\Va\to U\) is as a colimit
    over the associated sieve,
    which is defined as the full subcategory
    \(\sieve(\Va)\subseteq \overcat{\Uf}{U}\)
    spanned by those maps that factor through some \(V_i\to U\).
    Of course, in our setting \(\Uf\) is just a poset,
    so that this full subcategory is precisely
    \(\sieve(\Va)=\Uf\cap \cdown{\Va}\),
    thus justifying our original formula \eqref{eq:def-A-of-V-with-U}.
  \item \label{it:unordered-cech-nerve-condition}
    In the presence of sufficient fiber products,
    one can replace\footnote{
      One way to see this is to observe that in the presheaf-topos
      \(\Fun(\Uf^\op,\Gpdinfty)\)
      the \((-1)\)-truncated approximation of the map
      \(\coprod_{i<\beta}\Uf(-, V_i)\to \Uf(- , U)\)
      can either be computed explicitly pointwise to yield the subobject
      \(\sieve(\Va)\subseteq \Uf(-,U) \)
      of those \(f\colon O\to U\) which factor through some \(V_i\to U\),
      or as the colimit of its \v{C}ech nerve $\Cechraw{\Va}$
      (see \cite[Prop.\ 6.2.3.4]{LurHTT}).
      Applying the colimit preserving extension
      \(\Aa\colon \Fun(\Uf^\op,\Gpdinfty)\to \targetcat\)
      to the two ways of presenting this approximation,
      one obtains the desired comparison
      \(\colim \Cech{\Aa}{\Va}\simeq \colim \Aa\restrict{\sieve(\Va)}\).
    } the colimit over the sieve
    with the colimit over the (unordered) \v{C}ech nerve of the cover,
    defined as the simplicial object
    \begin{equation}
      \label{eq:unordered-Cech-nerve}
      \Cech{\Aa}{\Va}\colon \Delta^\op \ni [n] \mapsto
      \coprod_{i_0,\dots,i_n\in \beta}\Aa(V_{i_0}\times_{U}\dots\times_{U} V_{i_n}).
    \end{equation}
    In our setting, these fiber products are just intersections
    and so this \v{C}ech nerve is well defined as soon as \(\ccaps{\Va}\subseteq \Uf\).
  \item
    One can normalize the unorderd \v{C}ech nerve by removing all
    the degenerate cofactors,
    yielding the semi-simplicial object
    \begin{equation}
      \label{eq:normalized-ordered-Cech-nerve}
      \Delta^\op_{\inj}\ni [n]\mapsto
      \coprod_{\substack{i_0,\dots,i_n\in \beta\\\forall_{j}:\, i_{j-1}\neq i_j}}
      \Aa(V_{i_0}\times_U\dots \times_U V_{i_n}),
    \end{equation}
    where we only allow those indexing sequences with no adjacent duplicate indices.
    To recover the original \v{C}ech nerve one simply left Kan extends
    along the inclusion \(\Delta^\op_{\inj}\to \Delta^\op\),
    so they have the same colimit.
  \item
    Sometimes it is convenient to work with the punctured cubical diagram
    \begin{equation}
      \label{eq:cubical-Cech}
      \Pfinop{\beta}\setminus \{\emptyset\} \ni \{i_0,\dots,i_n\}
      \mapsto \Aa(V_{i_0}\times_U\dots\times_U V_{i_n})
    \end{equation}
    rather than the simplicial diagram.
    In our setting, this is justified by \Cref{lem:A(U)-cubical}
    which makes use of the fact that the sieve
    \(\mathrm{sieve}(\Va)=\Uf\cap\cdown{\Va}\) is a poset.
    More generally,
    a similar proof shows that the map
    \begin{equation}
      \Pfinop{\beta}\setminus\{\emptyset\} \to \mathrm{sieve}(\Va)
    \end{equation}
    is colimit cofinal as long as every map \(V_i\to U\) of the precover
    is a monomorphism (i.e.\ \((-1)\)-truncated).
  \item
    Another way of encoding the same colimit is via the
    \emph{strictly ordered}
    \v{C}ech nerve, defined as the semisimplicial object
    \begin{equation}
      \soCech{\Aa}{\Va}\colon
      \Delta_{\mathrm{inj}}^\op\to \targetcat,
      \quad\quad [n] \mapsto
      \coprod_{i_0< \dots< i_n <\beta}
      \Aa(V_{i_0}\times_U\dots\times_U V_{i_n}).
    \end{equation}
    or the \emph{weakly ordered}
    \v{C}ech nerve, defined as the simplicial object
    \begin{equation}
      \oCech{\Aa}{\Va}\colon
      \Delta^\op\to \targetcat,
      \quad\quad [n] \mapsto
      \coprod_{i_0\leq \dots\leq i_n <\beta}
      \Aa(V_{i_0}\times_U\dots\times_U V_{i_n}).
    \end{equation}
    These arise from the punctured cubical diagram \eqref{eq:cubical-Cech}
    by left Kan extension along
    the tautological map
    \begin{equation}
      \Pfinop{\beta} \setminus\{\emptyset\} \to\Delta^\op_{\mathrm{inj}}\to \Delta^\op
    \end{equation}
    that sends \(\emptyset\neq I\subset \beta\)
    to the nonempty set \(I\) with the linear order induced from \(\beta\).

    Since this left Kan extension does not change the colimit,
    the respective colimits over the ordered \v{C}ech nerves again
    compute the value of \(\Aa\) on the precover \(\Va\).
    Note that since the comparison between ordered and unorderd \v{C}ech
    nerves passes through the punctured cubical diagram,
    it only holds if all maps \(V_i\to U\) of the cover are monomorphisms.
  \item \label{it:stable-totalizations}
    \newcommand{\NDK}{\mathrm{N}_{\mathrm{DK}}}
    \newcommand{\tot}{\mathrm{tot}}
    In the case where the target \(\infty\)-category \(\targetcat\)
    is additive
    (for example, the derived category of an abelian category)
    one can invoke the Dold--Kan correspondence
    (see \cite{WaldeDK})
    \begin{equation}
      \NDK\colon
      \Fun(\Delta^\op,\targetcat)
      \xrightarrow{\simeq}
      \Chconn{\targetcat},
    \end{equation}
    which maps a simplicial object
    to its connective chain complex of normalized chains.
    Computing the colimit of an (ordered or unordered) \v{C}ech nerve
    then just amounts to totalizing its associated chain complex.%
    \footnote{
      One way to see this is to note that
      the Dold--Kan equivalence identifies
      the full subcategories of constant simplicial objects
      and 
      \(0\)-truncated chain complexes
      (see \cite[Corollary~7.1.7]{WaldeDK}),
      hence intertwines the left adjoints
      \(\colim_\Delta\colon \Fun(\Delta^\op,\targetcat)\to \targetcat\)
      and \(\tot\colon \Chconn{\targetcat}\to \targetcat\)
      to the respective inclusions.
    }
  \end{enumerate}
\end{rem}

\begin{war} \label{war:CGDefinitionNotQuiteRight}
  It might seem more natural to normalize the unordered \v{C}ech nerve
  as the semi-simplicial object
  \begin{equation}
    \Delta^\op_{\mathrm{inj}}\ni [n]\mapsto
    \coprod_{\substack{i_0,\dots,i_n\in \beta\\\text{pairwise distinct}}}
    \Aa(V_{i_0}\times_U\dots \times_U V_{i_n}),
  \end{equation}
  i.e., by throwing out index sequences with \emph{any} duplicated indices,
  not just adjacent ones.
  However, this would yield the \emph{wrong} colimit:
  for example, if \(\Va=\{V_1,V_2\}\) is the degenerate 2-element cover of \(U\)
  with \(V_1=U=V_2\),
  this formula would yield the colimit over the diagram
  \begin{equation}
    \begin{tikzcd}
      \Aa(V_1\cap V_2)
      \ar[d]
      \ar[dr]
      &
      \Aa(V_2\cap V_1)
      \ar[d]
      \ar[dl]
      \\
      \Aa(V_1)
      &
      \Aa(V_2)
    \end{tikzcd}
  \end{equation}
  which is \(\sphere{1}\otimes \Aa(U)\)
  instead of the correct value \(\Aa(\Va)=\Aa(U)\).
\end{war}

\begin{rem}
  \label{rem:compare-Jinfty}
  \newcommand{\Jinfty}{J_{\infty}}
  In this paper we chose to work with Weiss cosheaves
  (rather than with \emph{hyper}cosheaves)
  as is more standard in the literature surrounding factorization algebras. 
  We want to warn the reader that this convention is not universal,
  so that there is a high potential for confusion,
  exacerbated by the relative subtle (and often not explicitly addressed) distinction
  between covers and hypercovers when both just amount to a collection of opens;
  see also \Cref{war:cover-vs-hypercover}.

  For example, we want to compare our setting with that of Eric Berry,
  whose PhD thesis~\cite{Eric} on the additivity of locally constant factorization algebras
  uses many techniques similar to ours
  (and in fact we learned some of them from there!).
  Berry's factorization algebras \(\Aa\) are ``\(\Jinfty\)-cosheaves''
  (see \cite[Definition~2.0.10]{Eric})
  which amounts to the condition
  \begin{equation}
    \label{eq:Berry-Jinfty-cosheaf-condition}
    \colim_{U\in\Uf}\Aa(U)\xrightarrow{\simeq} \Aa(O).
  \end{equation}
  for every ``\(\Jinfty\)-cover'' \(\Uf\) of an open \(O\)
  (see~\cite[Definition~2.0.7]{Eric}).
  Since Berry's \(\Jinfty\)-covers
  are precisely our Weiss hypercovers \(\Uf\hyprefines O\),
  his ``\(\Jinfty\)-cosheaf'' condition translates 
  to what we call the ``Weiss hypercosheaf'' condition.
\end{rem}

\newpage

\section{Background on (constructible) factorization algebras}
\label{sec:(Constructible)FactAlgs}

In this section we define the \infy-category \(\Fact{\PS}\) of factorization algebras on a \posetofopens{} \(\PS\). In addition we define what it means for a
 factorization algebra to be ``constructible'' or ``locally constant with respect to a stratification''.

 \subsection{Basics on factorization algebras}
 \label{subsect:FactAlgebras}
 We give a brief introduction to factorization algebras in this section.
 Most of the definitions below can be found in \cite{CG} in the special case of chain complex valued, non-multiplicative factorization algebras.
 See also \cite{Ginot} for a nice introduction.

 Both \cite{CG} and \cite{Ginot} work in a 1-categorical context
 using techniques from homological algebras to ensure the correct homotopical notions.
 Instead, we work directly in the intrinsic \infy-categorical setting of \infy-operads.
 This is similar to Lurie's theory of factorizable cosheaves \cite{LurHA}
 and is also the approach of \cite{AFprimer} and \cite{Matsuoka}.
 There are other related approaches, such as the model categorical approach taken in \cite{Carmona-thesis}.
 A \(G\)-equivariant version is developed in \cite{Murray-thesis}.

 Let \(\PS\) be a \posetofopens{}, i.e.\ a full subposet \(\PS \subseteq \open{X}\).
 We turn this into a colored operad (in \(\Set\)) as follows: the colors are all the opens of \(\PS\).
 Let \((b_i)_{i\in n}\), where \(n\in \naturals\), be an \(n\)-tuple of colors.
 The set of multimorphisms out of this tuple into some color \(c\) is the singleton set if the \(b_i\)'s are pairwise disjoint and each \(b_i\) is a subset of \(c\);
 otherwise it is the empty set. 

 There is a standard way to turn a colored operad into an \infy-operad (\cite[Construction 2.1.1.7]{LurHA}), and we spell out this construction for \(\PS\) below. Recall that \(\Fin\) is the category of finite pointed sets and point-preserving maps. 

 \begin{constr} 
	The \infy-operad \(\operad{\PS}\) (corresponding to the \posetofopens{} \(\PS\)) has as objects pairs \((I_+, (b_i))\), where \(I_+ \in \Fin\) and \((b_i)\) is an \(I\)-indexed collection of opens \(b_i \in \PS\). A morphism \(f: (I_+, (b_i)_{i\in I}) \rightarrow (J_+, (b'_j)_{j\in J})\) is a morphism \(f: I_+ \rightarrow J_+\) in \(\Fin\) such that for each \(j\in J\), the tuple \((b_i\, |\, i\in f^{-1}(j))\) is a collection of pairwise disjoint open subsets of \(b'_j\). 
\end{constr} 

\begin{figure}[H] \centering 
\scalebox{0.8}{
    \begin{tikzpicture}
%	\draw[step=1cm,gray,very thin] (0,0) grid (8,4); %grid for help if you want to change something
	\draw[color=red] plot [smooth cycle] coordinates {(3,2.5) (4,3.5) (5,3.5) (4.5,3) (5,1.5)}; \node at (4,2.8){\(U_2\)}; 
	\draw[color=blue] plot [smooth cycle] coordinates {(4.5,2) (3.5,1) (4.5,0) (5.5,1)}; \node at (4.5,1){\(U_3\)};
	\draw[color=blue] plot [smooth cycle] coordinates {(7.5,1) (7.5,2) (9,0) (7,-1) }; \node at (8,0.5){\(U_4\)};
	\draw[color=red] plot [smooth cycle] coordinates {(0,0) (-1,2) (2,2) (2,0.5) (1,1)}; \node at (0,1.5){\(U_1\)};
	\draw[color=black] plot [smooth cycle] coordinates {(11,4) (12,3) (11,2)}; \node at (11.5,3){\(U_5\)};
	%bigger open around the first two smaller opens
	\draw[color=red] plot [smooth cycle] coordinates {(-1,-1) (4,0) (5.5,1.5) (4.8, 4) (2,3) (-1, 3)}; \node at (2, 3.5){\(V_1\)};
	%second bigger open
	\draw[color=blue] plot [smooth cycle] coordinates {(10,0) (9,3) (3,2) (3, -1) (5.5,-0.5) (7, -2)}; \node at (8,3.4){\(V_2\)};
    \end{tikzpicture}	}
  \caption{An example of a morphism \(f\colon (5_+, (U_i)) \ra (2_+, (V_j))\) in \(\Open{\reals^2}\).
    The map \(f\colon 5_+ \ra 2_+\) is the map in \(\Fin\) defined by \(1, 2 \mapsto 1\);
    \(3, 4 \mapsto 2\) and \(5 \mapsto +\).
    This requires that the opens \(U_1\) and \(U_2\) are disjoint and contained in \(V_1\)
    and that the opens \(U_3\) and \(U_4\) are disjoint and contained in \(V_2\).
    There are no other conditions, so for example \(U_2\) and \(U_3\) may overlap
    and the open \(U_5\) need not be included in any \(V_j\).
  } \label{fig:MorpInOpenX}
\end{figure}

\begin{defn}
  \label{def:pcocartesian}
  A morphism $ f\colon (I_+, (b_i)) {\longrightarrow} ( J_+, (b'_j)) $
  in \(\operad{\PS}\) is called \emph{\pcocartesian{}} (with a \pcocartdecoration{}) if and only if
  for each \(j\in J\) we have \(\bigdisjun_{i\in f^{-1}(j)}b_i=b'_j\).
\end{defn}

\begin{rem}
  \label{rem:pcocartesian}
  A morphism
  $ f\colon (I_+, (b_i)) {\longrightarrow} ( J_+, (b'_j)) $
  in \(\operad{\PS}\) is cocartesian if and only if
  for each \(j\in J\) we have \(\sup(b_i\mid i \in f^{-1}(j))=b'_j\),
  with the supremum taken in the poset \(\PS\) (if it exists).
  Of course if the disjoint union \(\bigdisjun_{i} b_i\) lies in \(\PS\),
  then it is automatically the supremum;
  this case---%
  which is the only one we really care about in this article---%
  is exactly the one where the morphism \(f\) is \pcocartesian{}.

  We chose to use such a subtle typographical distinction
  because the two notions coincide in the cases that are most often considered,
  namely \(\PS=\open{X}\)
  or more generally when \(\PS\) is closed under subordinate disjoint unions,
  e.g., any multiplicative or factorizing basis;
  in these cases it is standard in the factorization algebra literature
  to just talk about cocartesian morphisms
  and we did not want to use a new word for the same concept.
  If there is ever a distinction between the two,
  we will never care about the actual cocartesian morphisms (without the \pcocartdecoration{}).
  But for the most part, all but the most detail-interested readers
  are probably best served by just ignoring the distinction.

  For an explicit example where the two notions differ,
  let \(\PS\subset\open{\reals}\)
  be the poset of open intervals.
  This poset has all finite disjoint joins/suprema
  (in fact, even the non-disjoint ones),
  but it is given by the convex hull and not by disjoint union.
  For example, the operation \((0,1),(1,2)\to (0,2)\)
  corresponds to a cocartesian morphism \(\operad{\PS}\),
  which is not \pcocartesian{}.
\end{rem}

We remind the reader that we work with an \(\otimes\)-presentable symmetric monoidal \infy-category \((\targetcat,\otimes)\).  
In particular, this means that \(\targetcat\) has all small colimits
and that the tensor product \(\otimes\) preserves them in each variable.

\begin{rem} 
  Here and in the rest of the paper,
  we view the symmetric monoidal \(\infty\)-category \((\targetcat,\otimes)\)
  as encoded in a cocartesian fibration \(\targetcatOT \rightarrow \Fin\)
  as in \cite[Chapter~2]{LurHA}.
  Its tensor products are defined uniquely (up to contractible choice)
  by cocartesian transfer,
  so that the cocartesian morphisms in \(\targetcatOT\) are
  (tautologically) the ones of the form 
	\begin{equation}
		(I_+, (V_i)) \stackrel{f}{\longrightarrow} \Bigg( J_+, \Big(  \bigotimes_{i \in f^{-1}(j)} V_i \Big) \Bigg).
	\end{equation}
\end{rem}

Let us now give the first definition on the way towards factorization algebras.
\begin{defn}
  Let \(\PS\) be a \posetofopens{}.
  A \emph{prefactorization algebra on \(\PS\)} (with values in \(\targetcat\))
  is a morphism of \infy-operads 
	\begin{equation}
	\Aa \colon \operad{\PS} \ra \targetcatOT.
	\end{equation}
  We equivalently say that \(\Aa\) is a \(\PS\)-prefactorization algebra in this
  situation.
  The \infy-category of \(\PS\)-prefactorization algebras is
  \(\Alg{\PS} \coloneqq\Alg[\targetcat]{\PS}\), where we usually leave the
  target category implicit unless necessary.
  For any open \(U\), we abbreviate \(\Alg{U}\coloneqq\Alg{\cdown{U}}\).
\end{defn}

\begin{ex}
  \label{ex:algebra}
  Let $X=\reals$ and $\Bf$
  the poset of nonempty open intervals.
  Let $\algebra$ be an associative unital algebra in $\targetcat = \Vect_\field$ of $\field$-vector spaces.
  We can define a $\Bf$-prefactorization algebra $\Aa$ as follows.
  Any interval $b$ is sent to $\algebra$, and,
  given an inclusion of intervals \(\bigdisjun_{i=1}^n b_i \hookrightarrow c\),
  we assign:
  \begin{itemize}
  \item if $n=1$, the identity $A\xrightarrow{=} A$;
  \item if $n\geq 2$, the composite
    \begin{equation}
      \Aa(b_1)\otimes \cdots \otimes \Aa(b_n)
      \xrightarrow[\rho]{\cong}
      \Aa(b_{\rho(1)})\otimes \cdots \otimes \Aa(b_{\rho(n)})
      =
      A\otimes \cdots \otimes A
      \xrightarrow{\mu^n} \algebra ;
    \end{equation}
    where \(\mu^n\) is the \(n\)-uple multiplication map
    and \(\rho\) is the unique permutation of \(\{1,\dots,n\}\)
    such that the disjoint intervals \(b_{\rho(1)}< \cdots < b_{\rho(n)}\)
    are in ascending order with respect to the ordering induced from \(\reals\).
  \item if $n=0$, the unit viewed as a map $\field \to \algebra$.
  \end{itemize}
  This assignment being a prefactorization algebra translates precisely to the algebra $\algebra$ being associative and unital. 
\end{ex}

\begin{ex}\label{ex:algebra_on_circle}
Now let $X=\sphere{1}$ and choose an orientation.
Let $\Bf$ be the poset of nonempty open intervals in $X$.
Again, let $\algebra$ be an associative unital algebra in $\targetcat = \Vect_\field$.
We define a $\Bf$-prefactorization algebra $\Aa_{\sphere{1}}$ exactly as in the previous example, but now using the unique permutation \(\rho\) of \(\{1,\dots,n\}\) such that when identifying $c$ with $\reals$ using the orientation, 
    the disjoint intervals \(b_{\rho(1)}< \cdots < b_{\rho(n)}\)
    are in ascending order with respect to the ordering induced from \(\reals\).
\end{ex}

\begin{ex}
  \label{ex:bimod}
  Let $X=\reals$ and $\Bf$ the poset of nonempty open intervals.
  Let $\algebra_1, \algebra_2$ be  associative unital algebras in $\Vect_\field$, $\bimod$ an $(\algebra_1, \algebra_2)$-bimodule, and $m\in \bimod$.
  We define a $\Bf$-prefactorization algebra $\Ma$ as follows:
  \begin{equation}
    \Ma(b) =
    \begin{cases}
      \algebra_1, & \text{if }b\subset (-\infty, 0),\\
      \bimod, & \text{if }0\in b,\\
      \algebra_2, & \text{if }b\subset (0, \infty).\\
    \end{cases}
  \end{equation}
  Let  \(\bigdisjun_{i=1}^n b_i \hookrightarrow c\) be an inclusion of intervals 
  and let \(\rho\) be the unique permutation of \(\{1,\dots,n\}\)
  such that the disjoint intervals \(b_{\rho(1)}<\dots<b_{\rho(n)}\)
  are in ascending order with respect to the ordering induced from \(\reals\).
  Let \(1\leq k\leq n\) be the largest number such that
  \(b_{\rho(k)}\not\subseteq(0,\infty)\).

  Then to this inclusion we assign the composite
  \begin{align}
    &\bigotimes_{i\in n}\Ma(b_i)
    \xrightarrow[\rho]{\cong}
    \Ma(b_{\rho(1)})
    \otimes
    \cdots
    \otimes
    \Ma(b_{\rho(n)})
    \\
    &
      \begin{cases}
        \xrightarrow{=}A_1^{\otimes{k-1}} \otimes M \otimes A_2^{\otimes n-k}
        \xrightarrow{\mu^{k-1}_1\otimes \id_M\otimes\mu^{n-k}_2}
        A_1\otimes M\otimes A_2
        \xrightarrow{\beta}
        M, & \text{if } 0\in b_{\rho(k)}
        \\
        \xrightarrow{=}A_1^{\otimes{k}} \otimes A_2^{\otimes n-k}
        \xrightarrow{\mu^k_1\otimes\mu^{n-k}_2}
        A_1\otimes A_2
        \xrightarrow{\beta(-,m,-)}
        M,
           & \text{if } b_{\rho(k)}\subset(-\infty,0)
      \end{cases}
  \end{align}
  where
  \(\mu_i^l\colon A_i^{\otimes l}\to A_i\) are the \(l\)-uple multiplication maps
  (which is the identity for \(l=1\) and the unit \(\field\to A_i\) for \(l=0\)),
  and where \(\beta\colon A_1\otimes M\otimes A_2\to M\)
  is the bimodule action map.
  Note that if \(n=0\) then this is just the map \(\field\xrightarrow{\cdot m} M\).
\end{ex}

\begin{ex}\label{ex:algebra_on_disjoints}
We continue our previous examples, but now extending from \(\Bf\) to \(\disjunComp{\Bf}\), 
i.e.\ the poset of finite disjoint unions of non-empty open intervals.
We can extend the $\Bf$-prefactorization algebras $\Fa=\Aa$, $\Fa=\Aa_{\sphere{1}}$, and $\Fa=\Ma$ from the previous examples to $\disjunComp{\Bf}$-prefactorization algebras by setting
\begin{equation}
   \Fa(b_1\disjun\cdots \disjun b_n) \coloneqq  \Fa(b_1)\otimes \cdots \otimes \Fa(b_n)\,.
\end{equation}
This includes the empty set (corresponding to the empty disjoint union),
to which is assigned the empty tensor product, i.e., \(\Fa(\emptyset)=\field\).
\end{ex}

We impose two conditions on a prefactorization algebra to get a factorization algebra.
First, a factorization algebra satisfies a local-to-global condition. This is expressed by the Weiss cosheaf condition recalled in \Cref{defn:JinftyCondition}.

\begin{defn} \label{defn:WeissAlgebrasCat}
  A \emph{Weiss algebra on \(\PS\)} is a \(\PS\)-prefactorization algebra
  whose underlying copresheaf \(\Bf\to \targetcat\)
  is a Weiss cosheaf.
	Define the full \infy-subcategory of Weiss algebras on \(\PS\) to be the pullback 
 	\begin{equation}
 		\cdsquareOpt[pb]
    {\AlgWeiss[\targetcat]{\PS}}
    {\Alg[\targetcat]{\PS}}
    {\cShvWeiss[\targetcat]{\PS} }
    {\Fun(\PS, \targetcat)}
    {hookrightarrow}
    {}{}
    {hookrightarrow}.
 	\end{equation}
\end{defn}

The second condition is that disjoint unions of opens should be sent to tensor products. 
\begin{defn}
  \label{defn:multiplicative-prefact-alg}
	A \(\PS\)-prefactorization algebra \(\Aa\) is said to be \emph{multiplicative} if it sends \pcocartesian{} morphisms in \(\operad{\PS}\) to cocartesian morphisms in \(\targetcatOT\). Denote by
	 \begin{equation}
	 	\AlgM[]{\PS} \hookrightarrow \Alg[]{\PS}
	\end{equation}	
	the full \infy-subcategory whose objects are the multiplicative \(\PS\)-prefactorization algebras. 
\end{defn}

\begin{rem}
  The \pcocartesian{} morphisms in \(\operad{\PS}\)
  are generated under concatenation of tuples by the ones of the form
  \((b_1,\dots,b_n)\to b\),
  where \(b=b_1\disjun\cdots \disjun b_n\).
  Thus a prefactorization algebra on \(\PS\) is multiplicative
  if and only if for each such decomposition,
  the induced map
  \begin{equation}\label{eqn:multiplicative}
    \Aa(b_1)\otimes \cdots \otimes \Aa(b_n)\to \Aa(b_1\disjun\cdots \disjun b_n)
  \end{equation}
  in \(\targetcat\) is an equivalence.
  Of course by induction it suffices to check this for \(n=0\) and \(n=2\).
\end{rem}

Combining the two conditions we arrive at our main definition.
\begin{defn}\label{defn:UFactAlg}
  Let \(\PS\) be a \posetofopens{} (of \(X\)).
  The \(\infty\)-category of \(\PS\)-factorization algebras
  (with values in \(\targetcat\)) is the intersection
  \begin{equation}
    {\Fact[\targetcat]{\PS}}
    \coloneqq
    {\AlgM[\targetcat]{\PS}}
    \cap
    {\AlgWeiss[\targetcat]{\PS} }
  \end{equation}
  of full subcategories of
  \(\Alg[\targetcat]{\PS}\).
	Explicitly, a \emph{\(\PS\)-factorization algebra} is a morphism of \infy-operads
  \(\Aa \colon \operad{\PS} \ra \targetcatOT \)
	that restricts to a Weiss cosheaf
  and sends \pcocartesian{} morphisms in \(\operad{\PS}\) to cocartesian morphisms in \(\targetcatOT\). 
\end{defn}

\begin{rem}
  In \cite{CG}, ``factorization algebras'' are defined to be what we call ``Weiss algebras'',
  without imposing multiplicativity.
  They write ``multiplicative factorization algebra'' for what we just call ``factorization algebra''.
  Most other references use our convention.
\end{rem}

\begin{ex}\label{ex:algebra_is_fact}
  Consider the prefactorization algebra $\Aa$ from \Cref{ex:algebra}.
  It is tautologically multiplicative, since the only \pcocartesian{} morphisms in \(\operad{\PS}\) lie over $1_+ \in \Fin$.
  To see that it is a Weiss cosheaf,
  first observe that every Weiss cover \(\Wa\to U\) in \(\PS\)
  must contain an ascending chain of intervals $(w_i)_{i\in\naturals}$ covering $U$.
  We explain this for $U=(a, \infty)$; the other cases are similar. 
  For each \(i\in\naturals\), let $S_i\coloneqq\{a+\frac{1}{i},i\}$.
  Since \(\Wa\) is a Weiss cover, there is a $w_i \in \Wa$ containing $S_i$;
  thus we have \(U=\bigcup_{i\in\naturals}w_i\).
  After potentially passing to a subsequence, we may assume without loss of generality
  that the \(w_i\) are ascending, thus yielding a colimit cofinal map
  \(w\colon \naturals \to \Bf\cap\cdown{\Wa}\).

  Finally, we conclude that we have an equivalence
  \begin{equation}
    \colim_{\naturals} \algebra \simeq \colim_{\Bf\cap\cdown{\Wa}} \algebra
    = \colim \Aa\restrict{\Bf\cap\cdown{\Wa}}  \xrightarrow{\simeq} \Aa(U) = \algebra,
  \end{equation}
  because \(\naturals\) is weakly contractible and all maps in the colimit are equivalences.
\end{ex}

\begin{ex}\label{ex:bimod_is_fact}
  Consider the prefactorization algebra $\Ma$ from \Cref{ex:bimod}.
  Again, there is nothing to show for multiplicativity.
  To see that it is a Weiss algebra, we use the same argument as in the previous example.
  Given a Weiss cover \(\Wa\to U\) in \(\Uf\), we find \((w_i)_{i\in\naturals} \subset \Bf\cap\cdown{\Wa}$ which is cofinal. 
  So
  \begin{equation}
    \colim_{i\in\naturals} \Ma(w_i)  \simeq \colim \Ma\restrict{\Bf\cap\cdown{\Wa}}  \xrightarrow{\simeq} \Ma(U),
  \end{equation}
  where we use that either \(\Ma(w_i)=A_1=\Ma(U)\)
  or \(\Ma(w_i)=A_2=\Ma(U)\) (for all \(i\in \naturals\)) if \(0\notin U\),
  or \(\Ma(w_i)=M=\Ma(U)\) (for all \(i\gg 0\)) if \(0\in U\).
  \end{ex}

  Since there are several slightly different definitions of (multiplicative) factorization algebras in the literature
  we now remark on how some of them compare. 

\begin{rem}\label{rem:Fact_alg_variants}
	Consider factorization algebras on \(\PS\), and let \(\Wa\) be a Weiss cover such that \(\ccaps{\Wa} \subseteq \PS\).
	We then have several equivalent ways of giving the Weiss cosheaf condition as thoroughly explained in \Cref{rem:Equivalent-colimits-Weiss-condition}.
	In particular, the definition used in this article (see Equation \eqref{eq:WeissCosheafConditionExplicit}) 
	exactly corresponds to taking the colimit over the Weiss covering sieves as in \Cref{rem:Equivalent-colimits-Weiss-condition} (\ref{it:covering-sieve-condition}). 
	This is also the condition appearing in \cite[Definition~2.20]{AFprimer} (for \(\PS = \openX\)). 

	To compare our definition to that of \cite[Section~4.1]{Ginot} we first note that his ``factorizing bases" only become 
	Weiss covers after adding all disjoint unions, and one needs to impose multiplicativity for the comparison.
  Modulo these 
	differences his definition corresponds to taking a colimit over the unordered \v{C}ech nerve as in Equation \eqref{eq:unordered-Cech-nerve}. 

  Regarding \cite[\S6 Definition~1.4.1]{CG},
  it is clear from the context (see e.g.\ \cite[App.~A, Definition~4.3.1]{CG}),
  that Costello-Gwilliam mean the \v{C}ech 
	complex associated to any of the equivalent \v{C}ech nerves described in
  \Cref{rem:Equivalent-colimits-Weiss-condition}.
  However, we warn the reader that as written
	their formulas can easily be misinterpreted to yield the wrong result,
  as explained in \Cref{war:CGDefinitionNotQuiteRight}. 

  Related but slightly different is \cite[Definition~2.0.12]{Eric},
  who works with Weiss hypercosheaves rather than Weiss cosheaves;
  see also \Cref{rem:compare-Jinfty}.
  Since evaluation on a Weiss hypercover is just its colimit,
  no version of the \v{C}ech nerve needs to be considered in his context.
\end{rem}

We also introduce the following notation to match the literature on factorization algebras.
\begin{notation}
	Whenever we are in the situation where the \posetofopens{} \(\PS\)
  consists of all opens,
  i.e., \(\PS=\cdown{X} = \open{X}\),
  we simply omit referencing \(\PS\) and often also leave \(X\) implicit.
  Explicitly, \(\Fact{X}\) is called the \(\infty\)-category of
  factorization algebras on \(X\).
\end{notation}

\subsection{Conical manifolds}
\label{sec:conical-manifold}
\label{sec:smooth-conical-manifold}

We recall the key notions of conical manifolds and smooth conical manifolds,
introduced by Ayala--Francis--Tanaka~\cite{AFT-local-structures}
under the names ``\(C^0\) stratified spaces'' and ``conically smooth stratified spaces''.
These are the analogs of topological and smooth manifolds in the conically stratified setting.

\begin{defn}
  A \emph{stratified space} \(X=(X\to P)\) consists of a topological space \(X\),
  together with a surjective continuous map \(X\to P\),
  where \(P\) is a poset equipped with the upwards-closed topology
  \begin{equation}
    \{
    U\subseteq P
      \mid
    \forall u\in U, \, \forall p\in P\colon u<p \implies p\in U\
    \}.
  \end{equation}
  A map \((X\to P)\to (Y\to Q)\) of stratified spaces
  is just a commutative square.
\end{defn}

\begin{constr}
  \label{cstr:stratified-cone}
  For every topological space \(X\) and every \(t\in (0,\infty]\),
  we write
  \begin{equation}
    \topcone[t]X \coloneqq \{0\}\amalg_{\{0\}\times X}[0,t)\times X
  \end{equation}
  for \emph{the cone of \(X\) of radius \(t\)}.
  We abbreviate \(\topcone{X}\coloneqq \topcone[\infty]{X}\).
  If \(X\) is equipped with a stratification \(\sigma\colon X\to P\),
  then we equip \(\topcone[t]{X}\) with the stratification
  \begin{equation}
    \topcone[t]{X}\to \{-\infty\}\star P,\quad (s,x)\mapsto
    \begin{cases}
      -\infty\text{, for } s=0\\
      \sigma(x)\text{, for } s>0\\
    \end{cases}
  \end{equation}
  where \(\{-\infty\}\star P\)
  arises from the poset \(P\) by adding a new minimal element.
\end{constr}

\begin{ex}
	Let \(X=\sphere{1} \ra P= \{0 < 1\} \) be the circle stratified with 3 points as in \Cref{fig:MarkedCircle}. 
	The corresponding cone \(\topcone{X} \ra \{-\infty\} \star P\), where \(\{-\infty\} \star P = \{-\infty < 0 <1\}\), is illustrated in \Cref{fig:Cone}. 
	Each of the marked points trace out a line in the cone \(\topcone{X}\), and the added minimal element, i.e. \(-\infty\), corresponds to the marked point at the apex of the cone. 
	
\begin{figure}[H]
\begin{subfigure}{0.45\linewidth}
\begin{minipage}{0.9\linewidth}
  % \includestandalone[width=\textwidth]{figures/GluingMarkedCircle} 
  \includegraphics*[width=\textwidth]{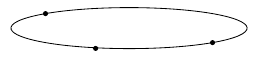}
  \caption{A circle stratified with 3 points}
  \label{fig:MarkedCircle}
\end{minipage} 
\end{subfigure} \hfill
\begin{subfigure}{0.45\linewidth}
\begin{minipage}{0.9\linewidth}
  % \includestandalone[width=\textwidth]{figures/ConePicture} 
  \includegraphics*[width=\textwidth]{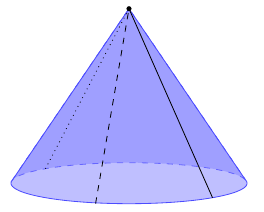} 
  \caption{The corresponding cone}
  \label{fig:Cone}
\end{minipage}
\end{subfigure}
\label{fig:StratCircleCone}
\caption{An example of a stratified circle and the corresponding cone.}
\end{figure}
\end{ex}

\begin{rem}
  \label{rem:cone-sphere}
  As a topological space, we have a canonical homeomorphism
  \(\topcone{\sphere{n-1}}\cong\reals^{n}\)
  which restricts to the corresponding homeomorphism
  \begin{equation} \label{eq:homeo-reals-cone}
  \topcone[t]{\sphere{n-1}}\cong \ball{t}{0} 
  \end{equation}
  for each \(t>0\).
  In this article we always consider \(\reals^{n}\) with its trivial stratification unless otherwise specified. 
  Hence, the homeomorphism \eqref{eq:homeo-reals-cone} is \emph{not} an isomorphism of stratified spaces,
  since \(\reals^n\) has no \(0\)-dimensional stratum \(\{0\}\)
  corresponding to the cone point of \(\topcone{\sphere{n-1}}\).
\end{rem}

Now we introduce conical manifolds, the stratified analog of a topological manifold.
These were introduced by Ayala--Francis--Tanaka under the name
``\(C^0\) stratified spaces'',
see \cite[Definition~2.1.15 and Lemma~2.2.2]{AFT-local-structures}.

\begin{defn}
  \label{defn:conical-manifold}
  Conical manifolds and conical disks are defined by mutual induction:
  \begin{enumerate}
  \item
    A stratified space \(X\)
    is called a \emph{conical manifold}
    if it is second countable Hausdorff and
    the collection of embedded conical disks
    forms a basis for the topology of \(X\).

    If \(X\) is nonempty
    then we say that \(X\) is \emph{of dimension} \(N\) or \(\leq N\)
    if all these embedded conical disks are of dimension
    \(N\) or \(\leq N\), respectively.
    We declare the empty conical manifold%
    \footnote{
      The empty space \(\emptyset\) is a conical manifold with empty stratification.
      Its set of embedded conical disks (which is empty) is vacuously a basis.
    }
    to be of dimension \(-1\).
  \item
    A \emph{conical disk} is a stratified space %this corresponds to what AFT calls a \(C^0\) basic, c.f. Def 2.2.1
    isomorphic to
    \begin{equation}
      \stratdisk{X}{n}\coloneqq  \topcone{X}\times \reals^n,
    \end{equation}
    where \(n\in\naturals\)
    and \(X\) is a compact conical manifold.
    If \(X\) is of dimension \(k\) or \(\leq k\)
    then we declare \(\stratdisk{X}{n}\) to be of dimension
    \(k+n+1\) or \(\leq k+n+1\), respectively.
  \end{enumerate}
\end{defn}

\begin{rem}
  Unlike the case of ordinary manifold,
  even a connected conical manifold might not have a well-defined dimension.
  For example \(\topcone{(\sphere{1}\amalg \sphere{0})}\) is homeomorphic
  to \(\reals^2\vee\reals\);
  it is a conical manifold of dimension \(\leq 2\)
  but it includes points of local dimension \(1\) and \(2\).

  A conical manifold might even have unbounded dimension,
  in which case it is not compact.
  Since we define conical disks only to arise from cones of \emph{compact}
  conical manifolds, they are always of bounded dimension.
  This in turn means that conical manifolds at least
  have \emph{locally} bounded dimension.
\end{rem}

\begin{ex}
  A conical manifold with trivial stratification
  is the same as an ordinary (topological) manifold,
  because in this case the only local disks that can appear are
  of the form \(\stratdisk{\emptyset}{n}=\reals^n\).

  Since we declare the empty conical manifold \(\emptyset\) to be of dimension \(-1\),
  the disks \(\reals^n=\stratdisk{\emptyset}{n}\) are of dimension \(n\),
  as expected.
  This means that in the unstratified case we recover the usual dimension of a manifold.
\end{ex}

\begin{ex}\label{ex:line_with_pts}
The simplest example of a non-trivial stratification is a line with finitely many points.
Take $X=\reals \to P=\{0<1\}$,  with the preimage of 0 a finite set of points $S\subset X$.
This space is covered by conical disks isomorphic to $\topcone{\sphere{0}}$ and $\reals$.
\end{ex}

\begin{rem}
Every conical manifold \(X\) admits a Weiss cover by (finite disjoint unions of) conical disks. This is analogous to the argument in \Cref{ex:WeissHausdorff}, using that the conical disks form a basis for the topology of \(X\) and that \(X\) is Hausdorff. 
\end{rem}

While annoyingly inductive, \Cref{defn:conical-manifold} is relatively
straightforward, since (just as for topological manifolds)
one does not need to put any extra conditions
on the transition maps appearing between the charts by conical disks.
The situation for \emph{smooth} conical manifolds is significantly more delicate
and requires an elaborate inductive definition.
In this paper we never work directly with the definition but only with
follow-up facts proved by Ayala--Francis--Tanaka;
hence we allow ourselves
to only sketch the definition and refer to \cite{AFT-local-structures} for the details.

Note that we call ``smooth conical manifolds''
what Ayala--Francis--Tanaka call ``conically smooth stratified spaces''.

\begin{defn}[Sketch of {\cite[Section~3.2]{AFT-local-structures}}]
  The following notions are defined by mutual recursion.
  \begin{itemize}
  \item
    A \emph{smooth conical manifold} \(X\)
    is a paracompact and Hausdorff stratified space equipped with
    a maximal collection\footnote{
      Unlike the set of opens in euclidean space
      (which are used to define ordinary smooth manifolds),
      the set of smooth conical disks is not necessarily small.
      Thus, to avoid set-theoretic issues,
      one has to interpret maximality up to isomorphism as follows:
      For every \(I\subseteq I'\) and any larger collection
      $\{\psi_{i'}\}_{i'\in I'}$ of the same kind,
      for every \(i'\in I'\) there exists an \(i\in I\)
      and a \emph{conically smooth isomorphism} \(\stratdiskNA_i\cong\stratdiskNA_{i'}\)
      intertwining \(\psi_i\) and \(\psi_{i'}\).
      Alternatively, one could inductively choose a small set of representatives
      for all types of smooth conical disks, only allow embeddings of those,
      and then use maximality with respect to plain inclusion.
    }
    \(\{\psi_i\colon \stratdiskNA_i\hookrightarrow X\}_{i\in I}\)
    (called a conically smooth atlas)
    of open embeddings of \emph{smooth conical disks} \(\stratdiskNA_i\),
    such that for each \(i,j\in I\) and each \(x\in \psi_i(\stratdiskNA_i)\cap\psi_j(\stratdiskNA_j)\)
    there is a \(k\in I\) with \(x\in \stratdiskNA_k\) and a commutative square
    \begin{equation}
      \begin{tikzcd}
        \stratdiskNA_k\ar[r,hookrightarrow]\ar[dr,hookrightarrow,"\psi_k"]
        \ar[d,hookrightarrow]
        &
        \stratdiskNA_i\ar[d,"\psi_i",hookrightarrow]
        \\
        \stratdiskNA_j\ar[r,hookrightarrow,"\psi_j"]
        &X
      \end{tikzcd}
    \end{equation}
    in which the left vertical map and upper horizontal map
    are \emph{conically smooth embeddings}.
  \item
    A \emph{conically smooth map} \(f\colon X'\to X\)
    between smooth conical manifolds
    with smooth conical atlas \(\{\psi'_i\}\) and \(\{\psi_j\}\), respectively,
    is a map of stratified spaces
    which induces conically smooth maps \(\psi^{-1}_j\circ f\circ \psi'_i\)
    on all charts.
  \item
    A \emph{smooth conical disk} is a stratified space of the form
    \begin{equation}
      \stratdisk{X}{n}\coloneqq \topcone{X}\times \reals^n,
    \end{equation}
    where \(X\) is a compact \emph{smooth conical manifold}.
  \item
    \emph{Conical smoothness} for a map \(\stratdisk{X}{n}\to\stratdisk{Y}{m}\)
    of \emph{smooth conical disks} is defined in terms
    of \emph{conical smoothness} of maps \(X\to Y\)
    (plus the ordinary notion of smoothness between euclidean spaces).
    We omit the precise definition here.
  \end{itemize}
\end{defn}

\begin{rem}
  Every smooth manifold can be viewed as a smooth conical manifold 
  with the trivial stratification.
  Between smooth manifolds, the conically smooth maps are precisely the ordinary smooth maps.
\end{rem}

\begin{notation}
  \label{not:link}
  For every point \(x\in X\) in a smooth conical manifold
  there is a unique compact smooth conical manifold \(Z\) (up to isomorphism)
  and a unique \(n\in \naturals\)
  such that embedded disks of the type \(\stratdisk{Z}{n}\)
  form a neighborhood basis around \(x\);
  see \cite[Theorem~4.3.1 and Corollary~4.3.2]{AFT-local-structures}.
  We call \(Z\) the \emph{link} at \(x\) in \(X\),
  and write \(\link{X}{x}\coloneqq Z\).

  For each stratum \(S\) we have the link bundle \(\Link{X}{S}\to S\)
  whose fiber over \(s\in S\) is the link at \(s\) in \(X\);
  see \cite[Definition~7.3.3]{AFT-local-structures}.
  If \(S\) is connected then the fiber
  of the link bundle is well defined up to isomorphism.
  In this case we write \(\link{X}{S}\) for this fiber
  and call it the \emph{link} at \(S\) in \(X\).
\end{notation}

\begin{notation}
  \label{not:snglr-cat}
  Smooth conical manifolds and open embeddings assemble into a category \(\Snglr\)
  and, more subtly,
  into an \(\infty\)-category \(\SSnglr\) whose homotopies are
  given by conically smooth isotopies.
  We refer to \cite[Section~4.1]{AFT-local-structures} for the precise definition.

  Note that every arrow in the $1$-category \(\Snglr\) is a monomorphism,
  so that all its slices are posets.
  Thus we are justified in writing
  \(\overcat{\Snglr}{M}=\open{M}\)
  and identifying an open embedding with its image.
  (Of course, the analogous statement is false for \(\SSnglr\).)
\end{notation}

\begin{rem}
  Nocera--Volpe~\cite{NoceraVolpe} proved that any Whitney stratified space
  admits a canonical smooth conical structure. Hence, all classical examples of
  Whitney stratified spaces from differential topology like analytic varieties are
  also examples of smooth conical manifolds.
\end{rem}

\subsection{Constructible factorization algebras}
\label{sec:ConstructibleFAs}

Having a good understanding of the stratifications relevant for us,
we now want to define constructibility.
This will be a condition checked on conical disks,
hence we first define when it makes sense to ask for such a condition to begin with. 

\begin{notation}
  \label{not:disks}
  We write \(\diskX\subset \openX\) for the full poset of
  finite disjoint unions of embedded disks,
  where ``embedded disks'' needs to be interpreted according
  to the type of object we view \(X\)
  as. More specifically, \(\diskX\) refers to
  \begin{itemize}
  \item
    the full subposet of finite disjoint unions of embedded conical disks
    if \(X\) is (viewed as) a conical manifold;
  \item
    the full subposet of finite disjoint unions of
    conically smoothly embedded smooth conical disks
    if \(X\) is (viewed as) a smooth conical manifold.
  \end{itemize}
  As a special case we get the usual notion of continuously or smoothly
  embedded disks in a topological or smooth (unstratified) manifold.
\end{notation}

\begin{defn}
	Let \(\Uf\) be a \posetofopens{}. We say that a subposet \(\Bf\subseteq \Uf \) is a \emph{disk-basis of \(\Uf\)} if \(\Bf\subset \diskX\) and \(\Bf\) is a basis for the topology of \(\bigcup \Uf\).
\end{defn}

\begin{defn}\label{def:enough_disks}\
  A \posetofopens{} \(\Uf\subset \openX\) is said to \emph{have enough good disks}
  if there exists a factorizing disk-basis \(\Bf\) of \(\Uf\) (cf.~\Cref{defn:mult_fact_basis}).
  
    We say that \emph{\(X\) has enough good disks} if the \posetofopens{} \(\openX\)
  has enough good disks witnessed by a factorizing disk-basis \(\Bf\subset\diskX\).
\end{defn}

\begin{rem}
  If \(X\) has enough good disks, every sieve \(\cdown{\Ua}\subset\openX\) also has
  enough good disks,
  because \(\Bf\cap \cdown{\Ua}\) is still a factorizing disk-basis of \(\cdown{\Ua}\).
\end{rem}

\begin{rem}[Existence of enough good disks] \label{rem:existence_good_disks}
  In general it is not true that \(\Bf\coloneqq\diskX\) is a factorizing basis of \(X\),
  because the intersection of two embedded disks need not be a disk.
  Thus a priori, there might not exist any \posetofopens{} of \(X\)
  which has enough good disks.

  For (unstratified) smooth manifolds there is a standard trick to get around this:
  one can choose a Riemannian metric on \(X\)
  and then let \(\Bf\subset \diskX\) consist only of disjoint unions of
  \emph{geodesically convex} disks.
  Since the intersections of convex disks is a convex disk again,
  this \(\Bf\) does indeed yield a factorizing basis of \(X\).

  For this paper, we choose not to address this question of existence further.
  Instead we just ``define the issue away'':
  when needed, we simply assume that our (smooth) conical manifolds
  have enough good disks.
\end{rem}

\begin{rem}
  \label{rem:enough-disks-prefact}
  Let \(\Uf\) be a \posetofopens{} which is closed under subordinate disjoint unions,
  i.e.\ \(U\disjun U' \in \Uf\) whenever
  \(U,U'\in \Uf\) are disjoint opens and \(U\disjun U' \in \cdown{\Uf}\).
  Then to establish enough good disks it suffices to find
  a \emph{prefactorizing disk-basis} \(\Bf\) of \(\Uf\),
  i.e., a subposet \(\Bf\subset\Uf\) which
  \begin{itemize}
  \item
    is a basis of the topology of \(\collunion \Uf\),
  \item
    is closed under intersection,
  \item
    consists of (disjoint unions of) disks
  \end{itemize}
  (but is not necessarily closed under disjoint unions).
  Indeed, one can then obtain a factorizing disk-basis of \(\Uf\)
  by simply closing \(\Bf\) under disjoint unions subordinate to \(\Uf\).
\end{rem}

\begin{defn}\label{defn:ConstructibleUFA}
	Let \(X\) be a conical manifold
  and let \(\Uf\subseteq \open{X}\) be a \posetofopens{} 
  with enough good disks.
  Let \(\Aa\colon \operad{\Uf} \ra \targetcatOT\)
  be a \(\Uf\)-factorization algebra.
  We say that \(\Aa\) is \emph{constructible}
  or \emph{locally constant with respect to the stratification} if every inclusion
  \(U\hookrightarrow V\) in \(\Uf\) of (abstractly) isomorphic conical disks
  is \(\Aa\)-local.
  We write
  \begin{equation}
    \FactCstr{\Uf} \subset \Fact{\Uf},
  \end{equation}
  for the full \infy-subcategory of constructible factorization algebras on \(\Uf\).
\end{defn}

Note that without requiring the \posetofopens{} \(\Uf\) above to have enough disks one might get an entirely vacuous condition. In the special case of ordinary manifolds we
abbreviate the name:

\begin{defn}
  Let \(X\) be a manifold and \(\Uf\subseteq \openX\) a \posetofopens{}
  with enough good disks.
  A \(\Uf\)-factorization algebra is called
  \emph{locally constant}
  if it is constructible with respect to the trivial stratification.
\end{defn}

\begin{ex}
  We return to \Cref{ex:algebra} and \Cref{ex:algebra_is_fact}.
  In this case, $\Aa$ is locally constant as a $\Bf$-factorization algebra, since every inclusion of intervals is even sent to the identity.
  Of course it remains constructible if
  we refine the stratification, such as by adding finitely many points as in
  \Cref{ex:line_with_pts}.
\end{ex}

\begin{ex}
  We return to \Cref{ex:bimod} and \Cref{ex:bimod_is_fact}. In this case, $\Ma$ is not locally constant, since the inclusion of an interval not containing the point $0$ into one containing the point $0$ is not an equivalence.
  However, it is a constructible factorization algebra
  if we view \(\reals\) as the smooth conical manifold
  \(\reals=\topcone{\sphere{0}}\),
  i.e., with the stratification \(\{0\}\subset \reals\).
\end{ex}

\subsection{Restriction and pushforward}
We recall the two operations for transferring factorization algebras
from one space to another.

\begin{constr}[Restriction and pushforward]
  \begin{itemize}
  \item[]
  \item
    Given an open subspace \(U\subseteq X\),
    the inclusion of operads
    \(\Open{U}\hookrightarrow \Open{X}\)
    induces the \emph{restriction} functor
    \begin{equation}
      -\restrict{U}\colon \Alg{X} \to \Alg{U}
    \end{equation}
    of algebras by precomposition.
  \item
    For each continuous map \(f\colon X\to Y\), 
    we have an induced functor
    \begin{equation}
      \preim{f}\colon \Open{Y}\to \Open{X}
    \end{equation}
    of operads,
    which induces the \emph{pushforward} functor
    \begin{equation}
      \pf{f}\colon \Alg{X}\to\Alg{Y}
    \end{equation}
    of algebras by precomposition.
  \end{itemize}
\end{constr}

\begin{lemma}
  \label{lem:restriction-preserves}
  Restriction along \(U\subseteq X\) preserves (separately)
  multiplicativity, Weiss cosheaves, and constructibility.
  In particular, it induces restriction functors
  \begin{equation}
    {-}\restrict{U}\colon\Fact{X}\to\Fact{U}
    \quad
    \text{and}
    \quad
    {-}\restrict{U}\colon\FactCstr{X}\to\FactCstr{U}
  \end{equation}
  on (constructible) factorization algebras.
\end{lemma}

\begin{proof}
  This is straightforward,
  since every instance in \(U\) of a
  multiplicativity, Weiss cosheaf, or constructibility
  condition is also an instance in \(X\).
  Note that here we use that \(\cdown{U}\subseteq\open{X}\)
  is a presieve (in fact, a sieve),
  so that Weiss covers in \(U\) are evaluated the same as Weiss covers in \(X\).
\end{proof}

\begin{lemma}
  \label{lem:pushforward-preserves}
  Pushforward along a continuous map \(f\colon X\to Y\) preserves (separately)
  multiplicative algebras and Weiss cosheaves.
  In particular, it induces a pushforward functor
  \begin{equation}
    \pf{f}\colon\Fact{X}\to\Fact{Y}
  \end{equation}
  on factorization algebras.
\end{lemma}

\begin{proof}
  Every instance
  \(V=V_1\disjun \cdots\disjun V_n\)
  of the multiplicativity condition on \(Y\)
  translates to an instance
  \(\preim[V]{f}=\preim[V_1]{f}\disjun\cdots\disjun\preim[V_n]{f}\)
  of the multiplicativity condition on \(X\);
  thus multiplicative algebras are preserved by \(\pf{f}\)

  Similarly, to prove that Weiss cosheaves are preserved by \(\pf{f}\),
  it suffices to observe that \(\preim{f}\) preserves intersections
  (so that \(\pf{f}\) is compatible with the evaluation of (pre-)covers),
  and to show that the preimage of a Weiss cover \(\Wa\to V\) in \(Y\)
  is again a Weiss cover;
  indeed, for every finite set \(S\subset \preim[V]{f}\)
  we can find a \(W\in\Wa\) covering the finite set \(f(S)\)
  so that \(\preim[W]{f}\in\preim[\Wa]{f}\) covers \(S\) as desired.
\end{proof}

\begin{ex} 
  If \(X\) and \(Y\) are smooth (conical) manifolds,
  even pushforward along (conically) smooth maps need not preserve
  constructible factorization algebras.
  For example consider an associative algebra \(A\)
  viewed as a locally constant factorization algebra \(\Aa\) on \(\reals\) with \(\Aa(\reals)=A\).
  On finite disjoint unions of intervals, this is simply the $\disjunComp{\Bf}$-factorization algebra\footnote{
  One can extend this to all open sets using the techniques from this paper.
   However, for our purposes in this example, the $\disjunComp{\Bf}$-factorization algebra is enough.
}
 from \Cref{ex:algebra_on_disjoints}.
    Consider the smooth map \(f\colon \reals\to\reals\) given by \(f(x)=x^2\).
  Then the inclusion \(\reals_{>0}\subset \reals\)
  of abstractly isomorphic disks
  need not be local for \(\pf{f}(\Aa)\),
  since it evaluates as the multiplication map
  \begin{equation}
    \label{eq:counterexample-x^2}
    A\otimes A = \Aa(\reals\setminus\{0\})\to\Aa(\reals)=A,
  \end{equation}
  which of course need not be an equivalence.

  This defect can be corrected by refining the stratification:
  Indeed, pushforward along the map \(f\colon x\mapsto x^2\) above
  does preserve constructible factorization algebras if the codomain
  is equipped not with the trivial stratification,
  but instead with the stratification \(\{0\}\subset \reals\).
  In this case the two conical disks
  \(\reals_{>0}\) and \(\reals\cong \topcone{\sphere{0}}\)
  are no longer abstractly isomorphic,
  which excludes the problematic non-equivalence
  \eqref{eq:counterexample-x^2}
  from the picture.
\end{ex}

\begin{rem}
  There are many maps between smooth conical manifolds along
  which the pushforward \emph{does} preserve constructible factorization algebras,
  for example locally trivial fibrations of (unstratified) manifolds
  or so-called ``adequately stratified'' maps of conical manifolds as defined by Ginot;
  see \cite[Sections~5.1~\&~6.1]{Ginot}.
\end{rem}

In this paper we only need the following very degenerate example,
for which the proof is straightforward.

\begin{lemma}\label{lemma:fold_is_constr}
  For each \(n\geq 0\),
  pushforward along the \(n\)-fold fold map
  $\nabla^n\colon X^{\amalg n}\to X$
  (which is the empty map for \(n=0\))
  preserves constructible factorization algebras.
\end{lemma}

\begin{proof}
  Let $\Fa$ be a constructible factorization algebra on $X^{\amalg n}$;
  we need to check that $\pf{\nabla^n}(\Fa)$ is constructible.
  To see this, start with an inclusion of abstractly isomorphic conical disks $U\subseteq V\subset X$.
  Then we have the commutative square
  \begin{equation}
    \begin{tikzcd}
      \pf{\nabla^n}\Fa(U)  \arrow{r}\arrow[equal]{d} & \pf{\nabla^n}\Fa(V)\arrow[equal]{d}\\
      \Fa(U^{\amalg n}) \arrow{r} & \Fa(V^{\amalg n}) \\
      \Fa(U)^{\otimes n}  \arrow{u}{\simeq}  \arrow{r}{\simeq}  & \Fa(V)^{\otimes n} \arrow{u}{\simeq}
    \end{tikzcd}
  \end{equation}
  where the lower vertical and horizontal maps are equivalences
  by multiplicativity and constructibility of $\Fa$, respectively;
  hence the top horizontal map is also an equivalence as desired.
\end{proof}

\begin{lemma}
  \label{lemma:union_is_constr}
  Let $f_i \colon X_i\to Y_i$ for $i\in I$.
  If pushforward along each \(f_i\) preserves constructible factorization algebras,
  then so does pushforward along their disjoint union \(\coprod_if_i\).
\end{lemma}

\begin{proof}
  This follows directly from the fact that
  any inclusion of abstractly isomorphic conical disks $U\subseteq V$
  in the disjoint union $\coprod_i Y_i$ already lies in one of the $Y_i$.
\end{proof}

\newpage

\section{The toolbox}

%local macros
\newcommand{\PSTwo}{\mathfrak{R}}

The goal of this section is to build a toolbox of techniques and results related to gluing of (constructible) factorization algebras.
First of all we motivate the different tools in \Cref{sec:MotivationForToolbox} by giving a detailed outline of how gluing of constructible factorization algebras will be proven.
In \Cref{subsect:FirstStepGluing} we use the theory of dendroidal sets to extend factorization algebras on opens \(\operad{U_I}\) of some cover \(\Ua\) to factorization algebras on \(\operad{\cdown{\Ua}}\).
Then, in \Cref{subsect:ExtendingUFactAlgs} we prove that factorization algebras defined on some factorizing basis can be extended.
Finally, in \Cref{subsect:ExtendingToDisjUniCompletion} we prove that multiplicative algebras on some decomposable \posetofopens{} \(\Bf\) can be extended to multiplicative algebras on the disjoint union completion \(\disjunComp{\Bf}\).

\subsection{Overview}
\label{sec:MotivationForToolbox}
Let \(\Ua\) be an open cover of a space \(X\). Let \(\beta\) be an ordinal and \(V_\bullet \colon \beta \twoheadrightarrow \Ua\) be a surjection, 
so we can write \(\Ua = \{U_i\}_{i < \beta}\). For any finite subset \(I \subset \beta\), denote \(U_I \coloneqq\bigcap_{i\in I} U_i\). 

Our motivating application (\Cref{thm:main-gluing}) was to show that constructible factorization algebras glue; that is, that the natural map 
\begin{equation}
\FactCstr{X} \xrightarrow{\simeq} \displaystyle \lim_{\beta\supset I \supsetneq \emptyset} \FactCstr{U_I},
\end{equation}
which can be found in the first horizontal row in the diagram \eqref{fig:diagram-outlining-proof-of-gluing}, is an equivalence.
The proof of this relies on three main steps corresponding to the columns in the diagram \eqref{fig:diagram-outlining-proof-of-gluing}. 

The leftmost column of equivalences are all proven in \Cref{subsect:FirstStepGluing}. 
The key insight is to understand the \infy-operad \(\operad{\cdownR{\Ua}}\) as a colimit of the \infy-operads \(\operad{\cdownR{U_I}}\), 
where \(\cdownR{\Ua}\) respectively \(\cdownR{U_I}\) denotes a restricted version of the sieve \(\cdown{\Ua}\) respectively \(\cdown{U_I}\)
which excludes the empty set (see \Cref{defn:restricted-sieve}). 
To be able to work explicitly with this colimit we work in the model of dendroidal sets. 
This key insight gives rise to an equivalence of the corresponding \infy-categories of algebras (\Cref{cor:GluingAlgebras}), 
which immediately restricts to the subcategories of multiplicative algebras (\Cref{lemma:alg-open-cover-eq-descends-to-mult}). 
Combining this with \Cref{prop:MultAlgToDisjCompletion} we get an equivalence also on the multiplicative 
algebras of the corresponding \infy-operads from the honest sieves, i.e.\ with the empty set added in (\Cref{cor:equi-multiplicative-limit-with-emptyset}). 
Moreover, this equivalence immediately restricts to the subcategories of factorization algebras and constructible factorization algebras (\Cref{prop:GlueFactOnCover} and \Cref{cor:GlueFactCstrOnCover}). 

For the rightmost column we have that \(\disjunComp{\cdown{\Ua}}\) is a factorizing basis of \(X\).
Hence, the results from \Cref{subsect:ExtendingUFactAlgs} regarding extending factorization algebras from a factorizing basis apply.
Explicitly, in \Cref{prop:AlgWeissEquivFactBasis} we show that Weiss algebras on a factorizing basis extend, 
and in \Cref{lemma:FactBasisExtends} that this equivalence restricts to the corresponding categories of factorization algebras. 
Combining extension from a factorizing basis with \Cref{thm:constructle-local}, i.e.\ that constructibility is local, 
gives the top rightmost equivalence as proven in \Cref{prop:extension-from-sieve-cover}. 
Note that this step requires that we work with a cover and does \emph{not} hold for any factorizing basis. 

Finally, for the middle column, observe that \(\cdown{\Ua}\) by definition only contains those disjoint unions 
that are subordinate to some open \(U\) of the cover. Hence, to get a factorizing basis from \(\cdown{\Ua}\) 
we need to add in all disjoint unions. In \Cref{prop:MultAlgToDisjCompletion} of \Cref{subsect:ExtendingToDisjUniCompletion} 
we prove that at the level of multiplicative algebras we get an equivalence. This naturally raises the question:

\begin{question}
  \label{question:FactAlgsDisjComp}
	Does the equivalence on multiplicative algebras in \Cref{prop:MultAlgToDisjCompletion}
  descend to an equivalence on factorization algebras
  \(\Fact{\cdown{\Ua}} \xrightarrow{\simeq} \Fact{\disjunComp{\cdown{\Ua}}}\)?
\end{question}

In \Cref{rem:Weiss-covers-not-sifted} we briefly outline what the problem with the straightforward approach to this questions is, and hence why it is still open. 
As a solution we restrict to the subcategories of constructible factorization algebras and in \Cref{prop:final_step_in_gluing} show that at this level we get an equivalence. 
This relies on \Cref{thm:constructible=>Weiss}, or more precisely \Cref{cor:constructible=>Weiss}, which extends techniques of Ayala--Francis and Ayala--Francis--Tanaka to prove
that the Weiss cosheaf condition is automatically satisfied
when left Kan extending from (a suitable collection of) disks.

\begin{equation}
\rotatebox{90}{
 \begin{adjustbox}{width=130ex}
\begin{tikzcd}[ampersand replacement=\&, cramped, row sep=7.0ex, column sep=1.0ex]
	\& \displaystyle {\lim_{\beta\supset I \supsetneq \emptyset}\FactCstr{U_I}} \&\&\& {\FactCstr{\cdown{\Ua}}} \&\&\& {\FactCstr{\disjunComp{\cdown{\Ua}}}} \&\&\& {\FactCstr{X}} \\
	\\
	\& \displaystyle {\lim_{\beta\supset I \supsetneq \emptyset} \Fact{U_I}} \&\&\& {\Fact{\cdown{\Ua}}} \&\&\& {\Fact{\disjunComp{\cdown{\Ua}}}} \&\&\& {\Fact{X}} \\
	\displaystyle {\lim_{\beta\supset I \supsetneq \emptyset} \AlgM{U_I} } \&\&\& {\AlgM{\cdown{\Ua}}} \&\&\& {\AlgM{\disjunComp{\cdown{\Ua}}}} \&\&\& {\AlgM{X}} \\
	\displaystyle {\lim_{\beta \supset I \supsetneq \emptyset} \AlgM{\cdownR{U_I}}} \&\&\& {\AlgM{\cdownR{\Ua}}} \&\& {\AlgWeiss{\cdown{\Ua}}} \&\&\& {\AlgWeiss{\disjunComp{\cdown{\Ua}}}} \&\&\& {\AlgWeiss{X}} \\
	\& \displaystyle {\lim_{\beta \supset I \supsetneq \emptyset} \Alg{\cdownR{U_I}}} \&\&\& {\Alg{\cdownR{\Ua}}} \&\&\& {\Alg{\disjunComp{\cdown{\Ua}}}} \&\&\& {\Alg{X}}
	\arrow[hook, from=1-2, to=3-2]
	\arrow[from=1-5, to=1-2, "\text{Cor.}~\ref{cor:GlueFactCstrOnCover}", "\simeq"'] %corollary from constructibility is local]
	\arrow[hook, from=1-5, to=3-5]
	\arrow[from=1-8, to=1-5, "\text{Prop.}~\ref{prop:final_step_in_gluing}", "\simeq"'] 
	\arrow[hook, from=1-8, to=3-8]
	\arrow[from=1-11, to=1-8, "\text{Prop.}~\ref{prop:extension-from-sieve-cover}", "\simeq"'] %equivalence of constructible on all of x and on disjuncomp(Ua)
	\arrow[hook, from=1-11, to=3-11]
	\arrow[hook,from=3-2, to=4-1]
	\arrow[from=3-5, to=3-2, "\text{Prop.}~\ref{prop:GlueFactOnCover}", "\simeq"']
	\arrow[hook,from=3-5, to=4-4]
	\arrow[from=3-8, to=3-5, dashed, "\text{\Cref{question:FactAlgsDisjComp}}"'] %the equivalence we wanted but can't prove
	\arrow[hook,from=3-8, to=4-7]
	\arrow[hook,from=3-8, to=5-9]
	\arrow[from=3-11, to=3-8, "\text{Prop.}~\ref{lemma:FactBasisExtends}", "\simeq"'] %fact alg equivalence between X and fact basis
	\arrow[hook,from=3-11, to=4-10]
	\arrow[hook,from=3-11, to=5-12]
	\arrow[from=4-1, to=5-1,"\simeq"',"\text{Lemma}~\ref{lem:remove-empty-set-multiplicative}"]
	\arrow[from=4-4, to=4-1, "\text{Cor.}~\ref{cor:equi-multiplicative-limit-with-emptyset}", "\simeq"']
	\arrow[from=4-4, to=5-4,"\simeq"',"\text{Lemma}~\ref{lem:remove-empty-set-multiplicative}"]
	\arrow[from=4-10, to=4-7, crossing over]
	\arrow[hook,from=5-1, to=6-2]
	\arrow[from=5-4, to=5-1, "\text{Lemma}~\ref{lemma:alg-open-cover-eq-descends-to-mult}", "\simeq"'] %new ref 
	\arrow[hook,from=5-4, to=6-5]
	\arrow[from=5-6, to=6-5]
	\arrow[from=5-9, to=5-6]
	\arrow[hook,from=5-9, to=6-8]
	\arrow[from=5-12, to=5-9,  "\simeq" ' near start,  "\text{Prop.}~\ref{prop:AlgWeissEquivFactBasis}" near start] %weiss arrow from all of X to disjuncomp(U) 
	\arrow[hook,from=5-12, to=6-11]
	\arrow[from=6-5, to=6-2, "\text{Cor.}~\ref{cor:GluingAlgebras}", "\simeq"']
	\arrow[from=6-8, to=6-5]
	\arrow[from=6-11, to=6-8]
	\arrow[hook,from=3-5, to=5-6]
	\arrow[hook,from=4-10, to=6-11, crossing over]
	\arrow[hook,from=4-7, to=6-8, crossing over]
	\arrow[from=4-7, to=4-4,  crossing over, "\simeq"' near start, "\text{Prop.}~\ref{prop:MultAlgToDisjCompletion}" near start] %equivalence on multiplicative for disjuncomp
\end{tikzcd}
\end{adjustbox}  }
\label{fig:diagram-outlining-proof-of-gluing}
\end{equation}

\newpage
\subsection{A colimit of \infy-operads}
\label{subsect:FirstStepGluing}
In this subsection we prove the first step towards gluing of factorization algebras. 
Given an open cover \(\Ua\) of \(X\), we will, roughly speaking, ``glue together" multiplicative prefactorization algebras on opens \(U\in \Ua\) 
of the cover to give a multiplicative prefactorization algebra defined on \(\cdown{\Ua}\). 
In particular, the key step is exhibiting the \infy-operad corresponding to a restricted version of \(\cdown{\Ua}\) as a suitable (homotopy) colimit cube.

It turns out that to understand the colimit cube of \infy-operads it is helpful to work in the model category of dendroidal sets. 
To that end we briefly recall the basics needed for us, and refer the reader to e.g.\ \cite{CM11} for further details. 

Recall that \(\strictOpd\) is the \((1,1)\)-category of operads, in contrast to the \((2,1)\)-category \(\Opd\).

\begin{defn}[\cite{MW07}, Definition~4.1]
	Let \(\Omega\) denote the \emph{symmetric tree category},
  i.e.\ the full subcategory \(\Omega \hookrightarrow \strictOpd\)
  of those symmetric operads that are free on finite rooted non-planar trees.
  The category of \emph{dendroidal sets} is the presheaf category
	\begin{equation}
		\dSets= \Fun(\Omega^{\op}, \Set).
	\end{equation}
\end{defn}

\begin{rem}\label{rem:NormalMono}
\begin{enumerate}[label=(\arabic*), ref=(\arabic*)]
	\item[]
	\item The category \(\dSets\) of dendroidal sets has a model structure \cite[Theorem 2.4]{CM11}. 
		The cofibrations are the \emph{normal monomorphisms}. 
	\item \label{item:NormalMonoCofib} A dendroidal set \(B\) is said to be \emph{normal} if and only if for any tree \(T\), 
		the action of the group \(\mathrm{Aut}(T)\) on \(B(T)\) is free. 
		By \cite[Corollary 1.8]{CM11} we know that any monomorphism of dendroidal sets \(\iota\colon A \hookrightarrow B\) 
		where \(B\) is normal is a normal monomorphism. 
		%In particular, all dendroidal sets appearing below are normal, so any monomorphism is automatically normal, and hence a cofibration. 
  	\item The model structure is left proper \cite[Proposition~2.6]{CM11}.
\end{enumerate}
\end{rem}

\begin{defn}
  [\cite{MW07}, Example~4.2]
	The {\em dendroidal nerve functor} \(\ND{}\colon \strictOpd \ra \dSets\) is defined by 
	\begin{equation}
		\ND{\Oa}(T)= \strictOpd\left(T , \Oa \right),
	\end{equation}
	i.e.\ it picks out operations of the shape of the tree \(T\) in the operad \(\Oa\). 
\end{defn}

Since the nerve does \emph{not} preserve homotopy colimits we are extra careful to include the nerve in this section. 
The colimits we want to compute take the shape of (punctured) cubes. 
To avoid finiteness assumptions on the covers we are working with, we will even work with transfinite cubes. 
The details regarding these transfinite cubes are given in \Cref{appendix:transfinite-cubes}. 

One crucial prerequisite for applying the results from \Cref{appendix:transfinite-cubes} is that certain maps related to the (punctured) cubes needs to be \emph{cofibrations}. 
In our setting this follows from all of the dendroidal sets of the cube being normal, and all the maps being monomorphisms (by \Cref{rem:NormalMono}\ref{item:NormalMonoCofib}). 
The following observation explains why we cannot employ these types of arguments directly to the dendroidal sets corresponding to sieves of precovers. 

\begin{obs}\label{obs:usual-nerve-of-sieve-not-normal}
	Let \(\Ua\) be a precover (of \(X\)). 
	The dendroidal set \(\ND{\cdownoperad{\Ua}}\) corresponding to the sieve \(\cdown{\Ua}\) is \emph{not} a normal dendroidal set
  because for any \(U\in \Ua\),
  the dendrex corresponding to the \(2\)-ary operation \((\emptyset,\emptyset)\to \Ua\)
  is a fixed point for the action of the \(2\)-corolla \(\corolla 2\),
  so that the action of \(\corolla 2\) on \(\ND{\cdownoperad{\Ua}}(\corolla 2)\) is not free.
\end{obs}

Hence, to be able to employ the results from \Cref{appendix:transfinite-cubes} we need to work with a version of sieves which excludes the empty set. 
This ensures that we obtain normal dendroidal sets after taking the dendroidal nerve. 

\begin{defn}\label{defn:restricted-sieve}
	Let \(\Ua\) be a precover (of \(X\)). 
	We define the associated \emph{\restrictedsieve{}} to be 
	\begin{equation}
		\cdownR{\Ua} \coloneqq
    \cdown{\Ua}\setminus\{\emptyset\}=
    \{ U' \in \openX \mid U' \neq \emptyset  \text{ and }  \exists U \in \Ua \colon U'\subseteq U  \} \ . 
	\end{equation}
\end{defn}

In words we have that the \restrictedsieve{} associated to a precover corresponds to the usual associated sieve without the empty set.  

\begin{obs}
  Let \(\Vf\) be a \posetofopens{} with \(\emptyset\notin\Vf\).
  Any dendrex of \(\ND{\operad\Vf}(T)\)
  has input colors which are all pairwise disjoint,
  hence pairwise different (because \(\emptyset\notin\Vf\)).
  Since an automorphism of \(T\) is uniquely determined by how it permutes the input edges,
  this means that no tree-automorphism can fix such a dendrex.
  We conclude that the action of \(\Aut(T)\) on \(\ND{\operad\Vf}(T)\) is free,
  which means that the dendroidal set \(\ND{\Vf}\) is normal.
\end{obs}

We now have all the prerequisites to tackle the main results of this section. 
However, we first give a small example before considering the general case to make the arguments more digestible. 

\begin{ex}\label{ex:3CubeVeryCofibrant}
	Let \(\Ua=\{U_i\}_{i\in \{0,1,2\}}\) be an open cover of \(X\) with three opens. 
	For any subset \(J \subseteq \{0, 1, 2\}\), let \(U_J=\bigcap_{j\in J}U_j\). 
	Then we have a 3-cube (see \Cref{ex:3Cube}) coming from inclusions: 
	\begin{equation}\label{eq:PushoutEx3Opens}
   		F \coloneqq\ND{\cdownRoperad{{\Ua}_{-}}}=\quad\quad
    		\cdcubeNA[small]
    		{\ND{\cdownRoperad{U_{012}}}}{\ND{\cdownRoperad{U_{12}}}}
    		{\ND{\cdownRoperad{U_{01}}}}{\ND{\cdownRoperad{U_1}}}
    		{\ND{\cdownRoperad{U_{02}}}}{\ND{\cdownRoperad{U_2}}}
    		{\ND{\cdownRoperad{U_0}}}{\ND{\cdownRoperad{X}} \,.}
  	\end{equation}
	
  	We claim that this is a very cofibrant cube as in \Cref{defn:VeryCofibrantCube}. 
	Recall from Equation \eqref{eq:SLatchingMap} that for every subset \(S\subseteq \{0,1,2\}\) 
	the \(S\)-restricted latching map is given by the transformation \(\partial_S F \ra \tau_S F \ra F \). 
	That is, one first restricts the cube (but still views it as a 3-cube by setting \(F(I)\simeq \emptyset\) for all \(I\nsubseteq S\)), then apply the latching map. 
	We want all of these \(S\)-restricted latching maps to be levelwise cofibrations. More explicitly, let \(S=\{0,2\}\). 
	The corresponding restricted cube \(\tau_S F\) is illustrated in \Cref{ex:SRestrictionAndLatchingMap}. 
	We have for example\footnote{We refrain from writing out all the corresponding copies of \(\emptyset\) since they do not change the colimits.} 
  	\begin{equation}
		\partial_S F(02) \simeq \emptyset \ra F(02) = \ND{\cdownRoperad{U_{02}}} \ , 
	\end{equation}
	where we get \(\emptyset\) from taking the colimit over the empty diagram. 
	Asking for this map to be a cofibration is the same as asking for \(\ND{\cdownRoperad{U_{02}}}\) to be normal, which it is. 
	Now consider
  	\begin{equation}
		\partial_{S} F(2) \simeq F(02) = \ND{\cdownRoperad{U_{02}}} \longrightarrow \ND{\cdownRoperad{U_2}} = F(2) \ . 
	\end{equation}
	It follows from the inclusion \(U_{02} \hookrightarrow U_2\) that this is a monomorphism of dendroidal sets, 
	and since \(\ND{\cdownRoperad{U_2}}\) is normal we get from \Cref{rem:NormalMono}\ref{item:NormalMonoCofib} that this even is a cofibration. 
	Lastly we look at
	\begin{equation} \label{eq:SLatchingMapEx}
		\partial_S F(\emptyset) \simeq \colim \Big( \ND{\cdownRoperad{U_{0}}} \leftarrow \ND{\cdownRoperad{U_{02}}} \ra \ND{\cdownRoperad{U_{2}}} \Big) 
			\longrightarrow \ND{\cdownRoperad{X}} = F(\emptyset) \ .
	\end{equation}
	Recall that colimits of functors are computed pointwise. Let \(T\in \Omega\) be any tree and consider the map 
	\begin{equation}
		\colim \Big( \strictOpd(T, \cdownRoperad{U_0}) \leftarrow \strictOpd(T, \cdownRoperad{U_{02}}) \ra \strictOpd(T, \cdownRoperad{U_2})   \Big)  \longrightarrow \strictOpd(T, \cdownRoperad{\Ua})
	\end{equation}
	in the category of sets. This is a monomorphism for any tree \(T\). 
	Thus it follows, as above, that the map in \eqref{eq:SLatchingMapEx} is a cofibration. 

	Note that all maps considered are cofibrations for essentially the same reason: 
	all the maps in the cube in \eqref{eq:PushoutEx3Opens} are monomorphisms between normal dendroidal sets. 
	The remaining maps to check cofibrancy for are completely analogous. 
\end{ex}

We are now ready to look at the general statement. From now on, let \(\Ua\) be an open cover of \(X\). 
Let \(\beta\) be an ordinal and \(V_\bullet \colon \beta \twoheadrightarrow \Ua\) be a surjection, so that we can write \(\Ua = \{U_i\}_{i < \beta}\). 
For any finite subset \(I \subset \beta\), denote \(U_I \coloneqq\bigcap_{i\in I} U_i\).
This construction yields a \(\beta\)-cube in \(\strictOpd\)
which we can compose with the dendroidal nerve to obtain a cube
\begin{align}
  \label{eq:cube-dendroidal-nerve}
  \ND{\cdownRoperad{U_{-}}} \colon \Pfinop{\beta} & \longrightarrow\dSets
  \\ I &\longmapsto \ND{\cdownRoperad{U_I}}
\end{align}
of dendroidal sets.

\begin{lemma}\label{lemma:GluingOperads}
  The cube \eqref{eq:cube-dendroidal-nerve} induces an isomorphism 
	\begin{equation}
		 \colim_{\beta \supset I\supsetneq \emptyset} \ND{\cdownRoperad{U_I}} \xrightarrow{\cong} \ND{\cdownRoperad{\Ua}}.
	\end{equation}
	of dendroidal sets, and the strict colimit coincides with the homotopy colimit. 
\end{lemma}
\begin{proof}
  	We claim that the cube \(F\coloneqq\ND{\cdownRoperad{U_{-}}}\) is very cofibrant.
  	This is an immediate generalization of the argument from \Cref{ex:3CubeVeryCofibrant}. 
	Recall that we need to show that for every subset \(S\subset \beta\), the corresponding \(S\)-restricted latching map is a levelwise cofibration. 
	The pointwise formula for the \(S\)-restricted latching map is given in \eqref{eq:SLatchingMapPointwise}. That is, let \(I\) be a finite subset of \(S\) and consider
	\begin{equation}
		\partial_S \ND{\cdownRoperad{U_I}} = \colim_{S \supseteq J \supsetneq I} \ND{\cdownRoperad{U_J}} \longrightarrow \ND{\cdownRoperad{U_I}}. 
	\end{equation}
	The lefthand colimit is a colimit in the presheaf category \(\dSets\),
	and is hence computed pointwise. Let \(T\in \Omega\) be any tree, and consider
	\begin{equation}
		\partial_S \ND{\cdownRoperad{U_I}}(T) = \colim_{S \supseteq J \supsetneq I}\strictOpd(T, \cdownR{U_J}) \longrightarrow \strictOpd(T, \cdownR{U_I}). 
	\end{equation}
	For all the opens \(U_J\) appearing in the colimit above we have inclusions \(U_J \subset U_I\),
	hence \(\cdownR{U_J}\subset\cdownR{U_I}\).
	As a result all maps appearing in the colimit are inclusions of subsets into the codomain \(\strictOpd(T, \cdownR{U_I})\). 
	Hence the colimit, which is now computed in the category of sets, is simply the union of the sets. 
	From the inclusions it also follows that the map always is a monomorphism. 
	By \Cref{rem:NormalMono}\ref{item:NormalMonoCofib} it is even a normal monomorphism 
	(since we are working with \restrictedsieves{}), hence a cofibration. 
	We conclude that the cube is indeed very cofibrant. 

	Since the cube \(\ND{\cdownRoperad{U_{-}}}\)
	is very cofibrant (and the model structure on \(\dSets\) is left proper) we can employ \Cref{cor:VeryCofibrantCubeHomotopyColimit}.
	Hence the strict colimit of nerves,
	which, as just explained, we can compute as
	\begin{equation}
		\colim_{I\supsetneq \emptyset} \ND{\cdownRoperad{U_I}}
    		= \bigcup_{i\in\beta}\ND{\cdownRoperad{U_i}}
   		= \ND{\cdownRoperad{\Ua}},
	\end{equation}
	coincides with the homotopy colimit in the model category \(\dSets\). 
\end{proof}

We used the model category of dendroidal sets to be able to work explicitly with the above (homotopy) colimit, but now we move away from this model. 
Recall the following comparison results between different models of \infy-operads:
\begin{itemize}
	\item Barwick shows that Lurie's model for \infy-operads is Quillen equivalent to Barwick's model in \cite{Bar13}. 
	\item Chu, Haugseng and Heuts show that Barwick's model is Quillen equivalent to that of complete dendroidal Segal spaces in \cite{CHH}. 
	\item Finally, Cisinski and Moerdijk compare the three dendroidal models amongst each other,
  		and also dendroidal sets to simplicial operads in \cite{CM13a, CM13b}.
\end{itemize}

In total this means that the homotopy theory of Lurie's \infy-operads is equivalent to that of dendroidal sets and to that of simplicial operads \cite[Corollary 1.2]{CHH}.

\begin{cor}
	The inclusions $U_I \hookrightarrow \cdownR{\Ua}$ 
	induce an equivalence of \infy-operads
	\begin{equation} \label{eq:colimitInfyOperadStatement}
		\colim_{I\supsetneq \emptyset} \cdownRoperad{U_I}  \xrightarrow{\cong} \cdownRoperad{\Ua} \ . 
	\end{equation}
\end{cor}

The above corollary has the following consequence at the level of algebras. 

\begin{cor}\label{cor:GluingAlgebras}
	Restriction induces an equivalence of \infy-categories
	 \begin{equation}\label{eq:GluingAlgebras}
	 	\Alg{\cdownR{\Ua}} \xrightarrow{\simeq}\lim_{\beta \supset I\supsetneq \emptyset} \Alg{\cdownR{U_I}} \,.
	 \end{equation}
\end{cor}
\begin{proof}
  	Taking maps of \infy-operads out of \eqref{eq:colimitInfyOperadStatement},
  	into some fixed symmetric monoidal \(\infty\)-category \(\targetcatOT\) gives the equivalences: 
  	\begin{equation} \label{eq:Alg-lim-on-spaces}
    		(\Alg[\targetcat]{\cdownR{U}})^\simeq
    		\simeq
    		\Opdinfty\left(\colim\limits_{\beta\supset I\supsetneq \emptyset}  \cdownRoperad{U_I} , \targetcatOT\right)
    		\simeq
    		\lim\limits_{\beta\supset I\supsetneq \emptyset}
    		\Opdinfty\left(\cdownRoperad{U_I} , \targetcatOT\right)
    		\simeq
    		\lim_{\beta\supset I\supsetneq \emptyset}
    		(\Alg[\targetcat]{\cdownR{U_I}})^\simeq \ .
  	\end{equation}
  	This establishes the desired equivalence \eqref{eq:GluingAlgebras}
  	on the underlying \(\infty\)-groupoids.

  	To get an actual equivalence of \(\infty\)-categories
  	we employ a standard argument by bootstrapping to all simplices
  	using the functoriality of the equivalence \eqref{eq:Alg-lim-on-spaces}.
  	More explicitly we observe that
  	functorially in \([n]\in\Delta^\op\) and any \infy-operad \(\Oa\) we have
  	\begin{equation}
    		\Fun\left([n],\Alg[\targetcat]{\Oa}\right)
    		\simeq
    		\Alg[{\Alg[\targetcat]{\Oa}}]{[n]}
    		\simeq
    		\Alg[{\Alg[\targetcat]{[n]}}]{\Oa}
  	\end{equation}
  	(where we view \([n]=\{0\to\dots\to n\}\) as an operad with only \(1\)-ary operations)
  	so that we can replace the symmetric monoidal \(\infty\)-category \(\targetcat\)
  	by all of the symmetric monoidal \(\infty\)-categories
  	\(\Alg[\targetcat]{[n]}\)
  	to promote the equivalence \eqref{eq:Alg-lim-on-spaces}
  	to an equivalence of complete Segal spaces
  	\begin{equation}
    		\Fun\left( [\bullet],\Alg{\cdownR{\Ua}} \right)^\simeq
    		\xrightarrow{\simeq}
    		\lim_{\beta\supset I \supsetneq \emptyset}\Fun\left( [\bullet] , \Alg{\cdownR{U_I}} \right)^\simeq
    		\simeq
    		\Fun\left( [\bullet] , \lim_{\beta\supset I \supsetneq \emptyset}\Alg{\cdownR{U_I}} \right)^\simeq
  	\end{equation}
  	presenting the desired equivalence of \(\infty\)-categories \eqref{eq:GluingAlgebras}.
\end{proof}

So far in this section we have worked with the \restrictedsieves{} $\cdownR{U_I}$ and $\cdownR{\Ua}$,
which is \emph{not} what one wants for factorization algebras on $U_I$ and \(\cdown{\Ua}\), respectively.
Hence, we will now restrict our attention to multiplicative algebras and explain how this helps us add the empty set back in. 

\begin{lemma}\label{lemma:alg-open-cover-eq-descends-to-mult}
	The equivalence in \Cref{cor:GluingAlgebras} restricts to an equivalence 
	\begin{equation}
		\AlgM{\cdownR{\Ua}} \xrightarrow{\simeq} \lim_{\beta \supset I \supsetneq \emptyset} \AlgM{\cdownR{U_I}} \ . 
	\end{equation}
\end{lemma}
\begin{proof}
	Restricting a multiplicative algebra on \(\cdownR{\Ua}\) gives multiplicative algebras on all intersections \(\cdownR{U_I}\). 
	
	Conversely, let \(\Aa \in \Alg{\cdownR{\Ua}}\) be an algebra such that each restriction lies in \(\AlgM{\cdownR{U_I}}\). 
	Since any open in \(\cdownR{\Ua}\) is subordinate to some \(U\in \Ua\) it follows that \(\Aa\) is multiplicative from 
	each restriction being multiplicative. 
\end{proof}

We now want an analogous statement to the above where the \infy-operads comes from the honest (and not restricted) sieves. 
For this we employ the following lemma that tells us that
removing the empty set was immaterial when working with multiplicative prefactorization algebras.
Heuristically, this can be explained as follows:
For a general prefactorization algebra \(\Aa\),
the value \(\Aa(\emptyset)\) on the empty set
is a commutative algebra which acts on all other values \(\Aa(U)\).
But in the multiplicative case, \(\Aa(\emptyset)\) is the monoidal unit,
making both the algebra and module structures unique, hence redundant.
The actual proof makes use of operadic left Kan extensions, see \Cref{app:OLKE}.

\begin{lemma}
  \label{lem:remove-empty-set-multiplicative}
  Let \(\PS\) be a \posetofopens{}.
  Then restriction induces an equivalence
  \begin{equation}
    \AlgM{\PS}\xrightarrow{\simeq}\AlgM{\PS\setminus\{\emptyset\}}.
  \end{equation}
\end{lemma}
\begin{proof}
  It suffices to show that operadic left Kan extension along
  \(\PS'\coloneqq \PS\setminus\{\emptyset\}\hookrightarrow\PS\)
  yields an equivalence on multiplicative prefactorization algebras.
  For this it suffices to show that for each \(\Aa\in\Alg{\PS}\)
  which is multiplicative on \(\operad{\PS'}\)
  we have that \(\Aa\) is an operadic left Kan extension from \(\PS'\)
  if and only if it is multiplicative on all of \(\operad{\PS}\).

  The only object where the pointwise left Kan extension condition can possibly fail
  is \(\emptyset\in\PS\),
  where the relevant slice \(\overcatAct{\PS'}{\emptyset}\)
  is a singleton only containing the \pcocartesian{} map
  \(()\xrightarrow{!} \emptyset\)
  from the empty tuple\footnote{
    Recall that \(\emptyset\) is an object of the \posetofopens{} \(\PS\)
    which we view as an object of \(\operad{\PS}\)
    over \(1_+\in\Fin\) with the abbreviation \(\emptyset:=(1_+,(\emptyset))\).
    Contrast this with \(()\) which is our notation for the unique object of \(\operad{\PS}\)
    over \(0_+\in\Fin\).
    
  }.
  Therefore \(\Aa\) is an operadic left Kan extension from \(\PS'\)
  if and only if the map \(\Aa(()\xrightarrow{!} \emptyset)\) is cocartesian;
  this is trivially satisfied if \(\Aa\) is multiplicative.

  To prove the converse, assume that \(\Aa\) is multiplicative on \(\PS'\)
  and that the map \(\Aa(()\xrightarrow{!}\emptyset)\)
  is cocartesian.
  Let \(\bar{t}\xrightarrow{!} \bar{r}\)
  be an arbitrary active-\pcocartesian{} map of \(\operad{\PS}\);
  we have to show that \(\Aa(\bar{t}\xrightarrow{!}\bar{r})\) is cocartesian.
  Denote by \(\bar{b}\in\operad{\PS'}\) the tuple obtained from \(\bar{t}\)
  by omitting all components that contain the empty set;
  in formulas:
  \begin{equation}
    \pi(\bar{b})\coloneqq\{j\in\pi(\bar{t})\mid t_j\neq \emptyset\}
    \quad
    \text{and}
    \quad
    b_j\coloneqq t_j
    \text{ for all } j\in \pi(\bar{b}).
  \end{equation}
  We have a unique \pcocartesian{} arrow \(\bar{b}\to\bar{t}\)
  lifting the inclusion \(\pi(\bar{b})\hookrightarrow \pi(\bar{t})\)
  and an identification
  \begin{equation}
    (\bar{b}\xrightarrow{!}\bar{t})
    \,
    \cong
    \,
    \left(
      (\bar{b}\xrightarrow{=}\bar{b})
      \oplus
      \bigoplus_{j\in \pi(\bar{t})\setminus \pi(\bar{b})} (()\xrightarrow{!}\emptyset)
    \right),
  \end{equation}
  where \(\oplus\) denotes the concatenation of tuples.
  Since we assume that \(\Aa(()\xrightarrow{!}\emptyset)\) is cocartesian,
  the same is also true for \(\Aa(\bar{b}\xrightarrow{!}\bar{t})\).
  Thus it suffices to show that \(\Aa(\bar{b}\xrightarrow{!}\bar{r})\)
  is cocartesian.
  If \(\bar{r}=\emptyset\), then \(\bar{b}\xrightarrow{!}\bar{r}\) must just be
  \(()\xrightarrow{!}\emptyset\)
  because the tuple \(\bar{b}\) does not contain the empty set;
  if \(\bar{r}\neq \emptyset\), then \(\bar{b}\xrightarrow{!}\bar{r}\)
  is fully contained in \(\PS'\).
  In either case we have that \(\Aa(\bar{b}\xrightarrow{!}\bar{r})\)
  is cocartesian by assumption.
\end{proof}

\begin{cor}\label{cor:equi-multiplicative-limit-with-emptyset}
	Restriction induces an equivalence of \infy-categories
	\begin{equation}
		\AlgM{\cdown{\Ua}} \xrightarrow{\simeq} \lim_{\beta \supset I \supsetneq \emptyset} \AlgM{\cdown{U_I}} \ . 
	\end{equation}
\end{cor}
\begin{proof}
	Consider the commuting square where all functors are given by restriction
	\begin{equation}\label{eq:composite-eq-to-get-honest-sieve-mult-version}
	\begin{tikzcd}
    \AlgM{\cdown{\Ua}}  \arrow[r] \arrow[d, "\simeq"']
    &
    \displaystyle  \lim_{\beta \supset I \supsetneq \emptyset} \AlgM{\cdown{U_I}}
    \arrow[d, "\simeq"]
    \\
    \AlgM{\cdownR{\Ua}} \arrow[r, "\simeq"]
    &
    \displaystyle \lim_{\beta \supset I \supsetneq \emptyset} \AlgM{\cdownR{U_I}} 
	\end{tikzcd} \ . 
	\end{equation}
	The vertical equivalences follow from \Cref{lem:remove-empty-set-multiplicative} . 
	The lower horizontal equivalence is that of \Cref{lemma:alg-open-cover-eq-descends-to-mult},
  and thus it follows that the top horizontal map is also an equivalence,
  as desired.
\end{proof}

We now restrict the above result to the subcategories of factorization algebras, respectively constructible factorization algebras. 

\begin{prop}\label{prop:GlueFactOnCover}
  	The equivalence in \Cref{cor:equi-multiplicative-limit-with-emptyset} restricts to an equivalence
  	\begin{equation}\label{eq:GlueFactCoverStep1} 
		\Fact{\cdown{\Ua}}
    		\xrightarrow{\simeq} \lim_{\beta \supset I\supsetneq \emptyset}
    		\Fact{\cdown{U_I}} \ . 
  	\end{equation}
\end{prop}
\begin{proof}
  	Note that the \posetsofopens{} here are all presieves (even sieves), so by \Cref{not:indep_of_poo} the Weiss covers evaluate the same. 
  	Hence, restricting a factorization algebra on \(\operad{\cdown{\Ua}}\) gives
  	factorization algebras on all intersections \(\operad{\cdown{U_I}}\).

  	Conversely, let \(\Aa\in \Alg{\cdown{\Ua}}\) be an algebra such that each
  	restriction \(\Aa\restrict{U_I}\) lies in \(\Fact{\cdown{U_I}}\). Recall that
  	any open in \(\cdown{\Ua}\) is subordinate to some \(U\in \Ua\). Hence, checking
  	the Weiss cosheaf condition for \(\Aa\) boils down to
  	checking it for opens that are all subordinate to the same \(U\in \Ua\). This
  	follows from the restrictions of \(\Aa\) being in \(\Fact{\cdown{U_I}}\), so we
  	are done.
\end{proof}

\begin{cor}\label{cor:GlueFactCstrOnCover}
 	The equivalence in \Cref{prop:GlueFactOnCover} restricts to an equivalence
  	\begin{equation} 
		\FactCstr{\cdown{\Ua}} \xrightarrow{\simeq}\lim_{\beta \supset I\supsetneq \emptyset} \FactCstr{\cdown{U_I}} \ . 
 	\end{equation}
\end{cor}
\begin{proof} 
	As before, restricting a constructible factorization algebra on \(\operad{\cdown{\Ua}}\) gives 
	constructible factorization algebras on all intersections.

  	Conversely, let \(\Aa \in \Fact{\cdown{\Ua}}\) be an algebra such that each
	restriction \(\Aa\restrict{U_I}\) lies in \(\FactCstr{\cdown{U_I}}\). By
	definition, we know that any inclusion of abstractly isomorphic disks in
	\(\cdown{\Ua}\) is subordinate to some open \(U \in \Ua\) of the cover. 
	Using constructibility of \(\Aa\restrict{U}\), we see that this inclusion is $\Aa$-local. Thus it
	follows that \(\Aa\) is constructible.
\end{proof}

\subsection{Extending from factorizing bases} \label{subsect:ExtendingUFactAlgs}
In this subsection we prove that whenever \(\Bf\) is a factorizing basis of
\(\Uf\) we can extend \(\Bf\)-factorization algebras to \(\Uf\)-factorization
algebras. Moreover, this gives rise to an equivalence of the corresponding
\infy-categories. In the special case where \(\Bf\) is a factorizing basis of
\(X\), this produces factorization algebras on all of \(X\).

For the reader unfamiliar with operadic left Kan extensions we refer to
\Cref{app:OLKE} where the definitions and results relevant to our situation are
summarized. In particular, any map \(\iota \colon \PS \ra \Uf\) gives rise to an
adjunction
\begin{equation}
  \label{eq:OLKEAdjInText}
  % \label{eq:OLKEAdjunction} in appendix
\iota_!\colon \Alg[\targetcat]{\PS} \leftrightarrows \Alg[\targetcat]{\Uf} :
\iota^*.
  \end{equation} The right adjoint \(\iota^*\) is restriction, and the left adjoint
\(\iota_!\) is given by operadic left Kan extension, which is computed by the
pointwise operadic colimit formula given in Equation \eqref{eq:pointwise-oLKE}.

We now consider the situation where we operadic left Kan extend from a
factorizing basis, and how this interacts with the Weiss condition.

\begin{prop}
  \label{prop:AlgWeissEquivFactBasis}
  Let \(\Uf\) be a presieve
  and let \(\Bf\) be a factorizing basis of \(\Uf\). Then
restriction induces an equivalence
	\begin{equation} \iota^\ast \colon \AlgWeiss{\Uf} \xrightarrow{\simeq}
\AlgWeiss{\Bf}.
	\end{equation}
\end{prop}
\begin{proof}
  First we need to check that restricting Weiss algebras on \(\Uf\)
indeed gives Weiss algebras on \(\Bf\). Assume that \(\Aa\) is a Weiss
  algebra on \(\Uf\). 
  Since \(\Bf\) is a factorizing basis it follows from \Cref{prop:UniqueExtFromBasis} \ref{it:Weiss-is-LKE-from-Weiss} that 
  \(\iota^\ast(\Aa) = \Aa\restrict{\Bf}\) is also a Weiss algebra.  
  
Conversely, consider the diagram
	\begin{equation}
\begin{tikzcd}[column sep=2.5ex, row sep=2.2ex] &&&& \Alg{\Bf} &&&&&
\Alg{\Uf} \\ &&&&&&& \\ \AlgWeiss{\Bf} &&&&& \AlgWeiss{\Uf} \\ && &&
\Fun(\Bf, \targetcat) &&&&& \Fun(\Uf, \targetcat)\,. \\ &&&&&&& \\ \cShvWeiss{\Bf}
&&&&& \cShvWeiss{\Uf} \arrow[from=1-5, to=4-5] \arrow[from=6-1, to=6-6,
"\simeq", "\ref{prop:UniqueExtFromBasis}"'] \arrow[from=6-1, to=4-5]
\arrow[from=4-5, to=4-10, "\textrm{LKE}"] \arrow[from=1-10, to=4-10]
\arrow[from=3-1, to=3-6, dashed, crossing over] \arrow[from=3-1, to=1-5, hook]
\arrow[from=1-5, to=1-10, "\textrm{oLKE}"] \arrow[from=3-6, to=1-10, hook]
\arrow[from=3-6, to=6-6, crossing over] \arrow[from=6-6, to=4-10, hook]
\arrow[from=3-1, to=6-1]
\end{tikzcd}
\end{equation} The left and right side squares are pullbacks (by
\Cref{defn:WeissAlgebrasCat}).
The backmost square commutes and is a pullback by
\Cref{lem:oLKE-LKE},
which applies because
\(\Bf\) is closed under disjoint unions subordinate to \(\Uf\).
Since the cube commutes it
follows from the pasting law of pullbacks that the front square is also a
pullback square.

Finally, since \(\Bf\) is a factorizing basis and \(\Uf\) is a presieve, the bottom horizontal map between the categories of Weiss cosheaves on \(\Bf\) and \(\Uf\) is an equivalence by \Cref{prop:UniqueExtFromBasis}. But since the front square is a pullback and pullbacks preserves equivalences it follows that the top horizontal dashed map is also an equivalence. 
\end{proof}

\newcommand{\WeissOne}{\Wa_1}
\newcommand{\WeissTwo}{\Wa_2}

Recall from \Cref{rem:Weiss-cosheaf-Weiss-presieves} that whenever the ambient \posetofopens{} is a presieve we have two ways of checking the Weiss condition. Namely, either one checks the Weiss cosheaf condition for all Weiss covers, or equivalently one checks the Weiss hypercosheaf condition for all Weiss presieves. We will often use the later approach in this section.

\begin{prop}
  \label{lemma:FactBasisExtends}
	Let \(\Uf\) be a presieve. For any factorizing basis \(\Bf\) of \(\Uf\), the equivalence in \Cref{prop:AlgWeissEquivFactBasis} restricts to an equivalence
	\begin{equation}
		\iota^\ast \colon \Fact{\Uf} \xrightarrow{\simeq}\Fact{\Bf}.
	\end{equation}
\end{prop}
\begin{proof}
  We need to check that the equivalence from \Cref{prop:AlgWeissEquivFactBasis} restricts to multiplicative algebras. It is clear that restricting a multiplicative algebra on \(\operad{\Uf}\) gives a multiplicative algebra on \(\operad{\Bf}\). 

  Conversely, let \(\Aa\in \AlgWeiss{\Uf}\) be a Weiss algebra on \(\Uf\) such that \(\Aa\restrict{\Bf} \in \Fact{\Bf}\).
  We want to show that \(\Aa\) is multiplicative on \(\Uf\), and thus a \(\Uf\)-factorization algebra.
  For notational simplicity we only consider a \pcocartesian{} map of the form
  \((2_+, (U_1, U_2)) \xrightarrow{!} (1_+, U_1\disjun U_2)\) in \(\Uf\);
  the general case is completely analogous.
  
  A factorizing basis is in particular a multiplicative basis so by \Cref{lem:basis-to-cover} 
  we have Weiss hypercovers \(\Bf \cap \cdown{U_i} \hyprefines U_i\), for \(i\in \{1,2\}\). 
  Moreover, since \(\Bf\) is even a factorizing basis these Weiss hypercovers are Weiss presieves. 
  Hence the top horizontal map in the following commutative diagram is an equivalence
  since \(\Aa\) is a Weiss cosheaf.
  \begin{equation}\label{eq:ComDiagramMultWeiss}
    \begin{tikzcd}
      \displaystyle \colim_{B_1 \in \Bf \cap \cdown{U_1}} \Aa(B_1) \otimes \colim_{B_2 \in \Bf \cap \cdown{U_2}} \Aa(B_2) \arrow[d, swap, "\simeq"]	\arrow[r, "\simeq"]
      & \Aa(U_1)\otimes \Aa(U_2)  \arrow[dd ]
      \\\displaystyle  \colim_{(B_1, B_2) \in \Bf \cap \cdown{U_1} \times \Bf\cap \cdown{U_2}} \Aa(B_1) \otimes \Aa(B_2) \arrow[d, swap, "\simeq"]
      \\ \displaystyle \colim_{(B_1, B_2) \in \Bf \cap \cdown{U_1} \times \Bf\cap \cdown{U_2}} \Aa(B_1\disjun B_2) 	\arrow[r, "\simeq"] 	& \Aa(U_1\disjun U_2)
    \end{tikzcd}.
  \end{equation}
  The top left vertical map is an equivalence because our target category is \(\otimes\)-presentable.
  The bottom left vertical map is an equivalence because \(\Aa\) is multiplicative when restricted to \(\Bf\).
  Note that we have a Weiss presieve
  \begin{equation}
    \label{eq:disjun-Weiss-presieve}
    \{B_1\disjun B_2\mid B_i\in \Bf\cap\cdown{U_i}\}\hyprefines U_1\disjun U_2
  \end{equation}
  because each \(\Bf\cap\cdown{U_i}\hyprefines U_i\) is a Weiss presieve (for \(i=1,2\))
  and \(\Bf\) has all of the disjoint unions appearing in \eqref{eq:disjun-Weiss-presieve}
  since they are subordinate to \(U_1\disjun U_2\in \Uf\).
  Hence we can, once again, use that \(\Aa\) is a Weiss algebra to get the bottom horizontal equivalence.

  We conclude that the rightmost vertical map is also an equivalence;
  thus \(\Aa\) is multiplicative on \(\Uf\). 
\end{proof}

\begin{ex}\label{extend_algebra_from_basis}
  Recall \Cref{ex:algebra_on_disjoints}. 
  In that case, $\disjunComp{\Bf}$ is a factorizing basis of $\open{\reals}$ and $\open{\sphere{1}}$, respectively,
  so using \Cref{lemma:FactBasisExtends} we can extend $\Aa$, $\Aa_{\sphere{1}}$ and $\Ma$ to  factorization algebras on all of~$\reals$ and $\sphere{1}$, respectively.
Hence, we can now evaluate $\Aa_{\sphere{1}}$ on $\sphere{1}$.
The standard trick for this computation is to pushforward $\Aa_{\sphere{1}}$ along the projection $p\colon \sphere{1} \to \reals$ as indicated in the following picture:
\begin{center}
\begin{tikzpicture}[scale=0.5]
\draw (0,0) circle (1cm);
\draw[->] (1.8,0) --node[anchor=south]  {$p$} (3.2,0); 
\draw (4,2) -- (4,-2);
\fill (4, 1) node[anchor=west] {$s_1$} circle (0.1);
\fill (4, -1) node[anchor=west] {$s_2$} circle (0.1);
\end{tikzpicture}
\end{center}
Then $\pf{(\Aa_{\sphere{1}})}$ is a constructible factorization algebra on $\reals$ stratified by the two points indicated.
From the pushforward formula, we immediately see that for any interval $I_1$ containing $s_1$ but not $s_2$ we have that $\pf{p}(\Aa_{\sphere{1}}) (I_1) = \algebra$, and similarly for any interval $I_2$ containing $s_2$ but not $s_1$ we have $\pf{p}(\Aa_{\sphere{1}}) (I_2) = \algebra$. Moreover, any interval between the points is sent to $\algebra \otimes \algebra^{\mathrm{op}}$, and any interval outside $[s_1,s_2]$ is sent to $\field$.

We will see in \Cref{ex:GlueFactCstr} that then we have that 
\begin{equation}
\Aa_{\sphere{1}} (\sphere{1}) = \pf{p}(\Aa_{\sphere{1}}) (\reals) = \algebra \otimes_{\algebra\otimes \algebra^{\mathrm{op}}} \algebra \,.
\end{equation}
\end{ex}

Without assuming the Weiss cosheaf condition above one quickly runs into problems. An explicit example of this is: 
\begin{ex} %Tashi's counterexample
  Consider the factorizing basis \(\Bf = \cdown{(X \amalg X)} \setminus \{X\amalg \emptyset, \emptyset \amalg X\}\) of the space \(X\amalg X\).
  Let \(\Aa\) be a multiplicative algebra on \(\operad{\Bf}\) such that 
  \begin{align}
    \Aa \colon \operad{\Bf} &\longrightarrow \targetcatOT
    \\ X\amalg X &\mapsto L \neq \Aa(X)\otimes \Aa(X).
  \end{align}
  Even though \(\Aa\) is multiplicative we can choose \(L  \neq \Aa(X) \otimes \Aa(X)\)
  because the \pcocartesian{} morphism
  \((2_+, (X\amalg \emptyset, \emptyset\amalg X)) \xrightarrow{!} (1_+, (X\disjun X)) \) is not actually in \(\Bf\).
  Left Kan extending \(\Aa\) along the inclusion \(\Bf \hookrightarrow \cdown{(X\amalg X)}\) will, by construction, not give a multiplicative algebra. 
\end{ex}

\subsection{Extending multiplicative algebras to disjoint union completion}
\label{subsect:ExtendingToDisjUniCompletion}
\newcommand{\EqClass}[1]{\Ja_{#1}} 
\newcommand{\cCFactor}{\text{\textinterrobang}} %Notation for edges that can be factored into active followed by coCart
\newcommand{\CanFactorCat}{\Ia} 
\newcommand{\ActOverQ}{\Ja} 
\newcommand{\FirstMap}{g} 
\newcommand{\SecondMap}{h} 
\newcommand{\FinUnpointed}{\mathrm{Fin}} 

\newcommand\quotrel[2]{{#1}_{/#2}}%quotient of #1 by equivalence relation #2

In this subsection we will explain how one can extend multiplicative algebras on some decomposable \posetofopens{} \(\PS\) 
to multiplicative algebras on the corresponding disjoint union completion \(\disjunComp{\PS}\).
Explicitly:

\begin{prop}
  \label{prop:MultAlgToDisjCompletion}
	Let \(\PS \hookrightarrow \disjunComp{\PS}\)
  be the inclusion of a decomposable \posetofopens{} into its disjoint union completion.
  Then the adjunction \eqref{eq:OLKEAdjInText} restricts to an equivalence of \infy-categories
	\begin{equation}
    \label{eq:AlgMdisjunComp}
		\AlgM{\disjunComp{\PS}} \xrightarrow{\simeq}\AlgM{\PS}.  
	\end{equation}
\end{prop}

We first record some easy consequences of this result. 

\begin{cor}
  \label{cor:AlgM-to-Alg}
  Let \(\PS\) be a \posetofopens{} all of whose elements are connected and non-empty.
  Then restriction induces an equivalence
  \begin{equation}
		\AlgM{\disjunComp{\PS}} \xrightarrow{\simeq}\Alg{\PS}.  
  \end{equation}
\end{cor}
\begin{proof}
  Since all elements of \(\PS\) are connected, it is trivially decomposable;
  and multiplicativity is a vacuous condition
  for prefactorization algebras on \(\PS\) because it contains
  no non-trivial disjoint union (not even nullary ones, which would yield the empty set).
  Hence the equivalence \eqref{eq:AlgMdisjunComp} simplifies
  to the desired equivalence.
\end{proof}

\begin{rem}
  The analog of \Cref{cor:AlgM-to-Alg} is also obtained in Carmona's thesis~\cite{Carmona-thesis}
  in the special case where \(X\) is a smooth manifold
  and \(\PS\coloneqq \cdiskX\) is the poset of contractible disks;
  whose completion under disjoint unions is
  \(\cdisj{\cdiskX}=\diskX\).
\end{rem}

\begin{cor}
  \label{cor:add-disjun-constr}
  Assume that the space \(X\) is a conical manifold.
  Let \(\Bf\) be a decomposable \posetofopens{}.
  Then the restriction map
  \begin{equation}
		\AlgMCstr{\disjunComp{\PS}} \xrightarrow{\simeq}\AlgMCstr{\PS}
  \end{equation}
  is an equivalence.
\end{cor}
\begin{proof}
  This follows directly from the equivalence \eqref{eq:AlgMdisjunComp}
  by restricting to the constructible algebras on both sides.
  Indeed, the definition of constructibility
  only refers to conical disks;
  since these are contractible, no new conical disks lie in \(\cdisj{\PS}\)
  that did not already lie in~\(\PS\).
\end{proof}

As explained in \Cref{question:FactAlgsDisjComp},
we do not know whether the equivalence in \Cref{prop:MultAlgToDisjCompletion}
restricts to an equivalence on the \infy-categories of factorization algebras.
The key missing ingredient is the following seemingly harmless statement:

\begin{question} 
  \label{question:Weiss-on-coproduct}
  Let \(X=X_1\disjun X_2\) be a disjoint decomposition into open subsets.
  Let \(\Aa\) be a multiplicative prefactorization algebra on \(X\)
  which is a Weiss cosheaf when restricted to \(X_1\) and \(X_2\).
  Does \(\Aa\) satisfy descent for each Weiss presieve \(\Wa\hyprefines X\)? 
\end{question}

\begin{rem} \label{rem:Weiss-covers-not-sifted}
  If \(\Wa_1\hyprefines X_1\) and \(\Wa_2\hyprefines X_2\) are Weiss presieves,
  then we obtain a Weiss presieve
  \begin{equation}
    \Wa_1\boxtimes\Wa_2\coloneqq \{W_1\disjun W_2\mid W_i\in\Wa_i\}\hyprefines X
  \end{equation}
  which is \(\Aa\)-local by multiplicativity
  and the Weiss cosheaf condition on each \(X_i\): 
  \begin{align}
    \Aa(\Wa_1\boxtimes\Wa_2)
    &\simeq
      \colim_{W_1,W_2}\Aa(W_1\disjun W_2)
      \simeq
      \colim_{W_1}\Aa(W_1) \otimes \colim_{W_2}\Aa(W_2)
    \\
    &=
      \Aa(\Wa_1)\otimes\Aa(\Wa_2)
      \xrightarrow{\simeq}
      \Aa(X_1)\otimes\Aa(X_2)\simeq\Aa(X).
  \end{align}

  Of course not every Weiss presieve of \(X\) is of this form;
  a general Weiss presieve \(\Wa\hyprefines X\) can always be factored as
  \begin{equation}
    \Wa\hookrightarrow \Wa_1\boxtimes\Wa_2 \hyprefines X,
  \end{equation}
  where \(\Wa_i\coloneqq \{W\cap X_i\mid W\in \Wa\}\) for \(i=1,2\).
  By the previous observation, we are thus reduced to the question
  whether the left inclusion is \(\Aa\)-local, i.e.,
  whether the canonical map
  \begin{equation}
    \colim_{W\in\Wa}\Aa(W\cap X_1)\otimes\Aa(W\cap X_2)
    \longrightarrow
    \colim_{W,W'\in \Wa}\Aa(W\cap X_1)\otimes \Aa(W'\cap X_2)
  \end{equation}
  induced by the diagonal \(\Wa\to \Wa\times\Wa\) is an equivalence.

  Note that in general the diagonal \(\Wa\to \Wa\times \Wa\) is not cofinal
  since there might be \(W',W''\in \Wa\)
  for which there is no \(W \in \Wa\) such that \(W'\cup W'' \subseteq W\). 
  For example, if \(\Wa=\{X\setminus S\mid S\subset Y \text{ finite}\}\)
  and \(X\) is not finite.
\end{rem}

We now start working towards the proof of
\Cref{prop:MultAlgToDisjCompletion}.
The main work consists in simplifying the (operadic) colimits which appear in the
formula for the operadic left Kan extension along
\(\PS\hookrightarrow\cdisj{\PS}\).
For \(\bar{r} \in \DisjunComp{\PS}\),
a priori this is given as a colimit over the category
$\left\{ \bar{b} \xrightarrow{\act} \bar{r} \right\}$ of active maps into it;
we will reduce this to the poset
\(\left\{ \bar{t}\xrightarrow{!}\bar{r} \right\}\)
of active-\pcocartesian{} maps
by a sequence of cofinality and Kan extension arguments.

For many statements below we will only consider the case where
the \posetofopens{} \(\PS\) does \emph{not} contain the empty set
to avoid having to deal with unpleasant edge cases of our constructions.
We will then invoke \Cref{lem:remove-empty-set-multiplicative}
to reduce to this case.

We now introduce some notation that will be convenient to use,
before embarking on the technical preliminary results we need. 

\begin{notation}\label{not:ObjectsMorphismsPresieveOp} Let \(\PS\) be a \posetofopens{} with corresponding \infy-operad \(\operad{\PS}\), and let \(\pi_+ \colon\operad{\PS} \ra \Fin\) denote the forgetful functor down to finite pointed sets. For any object \(\bar{b}\in \operad{\PS}\), let \(\pi(\bar{b})\) denote the subset of \(\pi_+(\bar{b})\) where the base point is discarded. 
  \begin{itemize}
  \item Objects of \(\operad{\PS}\) will be abbreviated by only the finite tuple, i.e.\ \(\bar{b}\coloneqq(\pi_+(\bar{b}), (\bar{b}))\),
    where \(\bar{b} = (b_i)_{i\in \pi(\bar{b})}\).
    A 1-tuple is simply \(b\coloneqq  (1_+, (b))\). 
  \item If an edge is \pcocartesian{} we label the arrow \(\xrightarrow{!}\). 
  \item If a morphism admits a factorization into an active morphism followed by a cocartesian morphism we label it by \(\xrightarrow{\cCFactor}\). 
  \end{itemize}
\end{notation}

For many of the below arguments we need to understand the data of active and active-\pcocartesian{} maps in \infy-operads of the form \(\operad{\PS}\). 
It turns out to be convenient to think of active and active-\pcocartesian{} maps by the equivalence classes that they induce on the opens of the target.  

\begin{constr}
  \label{obs:ActiveMapsAsEqClasses}
  Let \(\PS\) be a \posetofopens{}. All objects and morphisms below are in~\(\operad{\PS}\). 
\begin{enumerate}[label=(\arabic*), ref=(\arabic*)]
	\item\label{item:FirstObsActiveEqClass} An active map \(\alpha\colon \bar{r} \ra b\) defines an equivalence relation \(\sim_\alpha\) on \(\pi_0(b)\) as follows: Identify two connected components \(b', b'' \in \pi_0(b)\) if \(r_i \cap b' \neq \emptyset\) and \(r_i \cap b'' \neq \emptyset\) for some \(i\in \pi(\bar{r})\), and then transitively complete to get an actual equivalence relation. In equations,
	\begin{equation} \label{eq:EqClassFromActiveMap}
		\alpha\colon \bar{r} \ra b \,\,\text{active} \qquad \rightsquigarrow \qquad \pi_0(b) \twoheadrightarrow \EqClass{\alpha}\coloneqq  {\pi_0(b)}_{/_{\sim_\alpha}}.
	\end{equation}
\item\label{item:SecondObsCocartesianEqClass}
  For active-\pcocartesian{} maps the above procedure automatically gives an equivalence relation of a particularly simple form.
	\begin{equation}\label{eq:coCartEqClasses}
		\gamma \colon \bar{s} \xrightarrow{!} b \,\, \text{active-\pcocartesian{}} \qquad
    \rightsquigarrow \qquad  \displaystyle \pi_0(b) = \coprod_{i\in \pi(\bar{s})} \pi_0(s_i) \twoheadrightarrow \{i \in \pi(\bar{s})\mid s_i\neq \emptyset\} = \EqClass{\gamma}.
	\end{equation}
	In words, the active-\pcocartesian{} map yields a partitioning of
  the connected components of \(b\) into at most \(\pi(\bar{s})\) parts.
  If all entries of the tuple \(\bar{s}\) are non-empty,
  then this partition fully recovers the active-\pcocartesian{} arrow \(\gamma\).
\end{enumerate}	
\end{constr}

\begin{rem} \label{item:ActiveMapIntoOneuple}
	Whenever we need to produce active maps into some fixed tuple \(\bar{b}\) we will 
	always implicitly use the identification \(\overcatAct{\operad{\PS}}{\bar{b}} \cong \prod_{i\in \pi(\bar{b})} \overcatAct{\operad{\PS}}{b_i}\). 
	In words, it is sufficient to produce active maps into each object \(b_i\), with \(i\in \pi(\bar{b})\), of the tuple \(\bar{b}\).
\end{rem}

\begin{ex}
  We provide an example of \Cref{obs:ActiveMapsAsEqClasses}\ref{item:FirstObsActiveEqClass}.
  Let
  \begin{equation}
    \bar{r} = (4_+, (r_1, r_{2,1}\disjun r_{2,2}, r_{3,1}\disjun r_{3,2}, r_{4,1}\disjun r_{4,2}))
  \end{equation}
  be the object whose opens are as illustrated in \Cref{fig:ExampleActiveMapEqClasses} below.
  Similarly, \(b= (1_+, (b_1 \disjun b_2 \disjun b_3 \disjun b_4 \disjun b_5))\)
  has connected components as depicted below. Let \(\alpha \colon \bar{r} \ra b\)
  be the unique active map \(4_+\ra 1_+\), with corresponding inclusions.
  For example, \(r_2\) maps into both the first and third connected component of \(b\),
  which means \(b_1\) and \(b_3\) are in the same equivalence class.
  Also note that \(r_3\) maps into the third and fourth connected component of \(b\),
  hence \(b_3\) and \(b_4\) are also in the same equivalence class.
  Upon transitively completing, we get that \(\alpha\) induces
  the three equivalence classes \([b_1] = [b_3] = [b_4]\), \([b_2]\) and \([b_5]\). 
  \begin{figure}[H]
    \begin{overpic}[scale=1.0, tics=10]{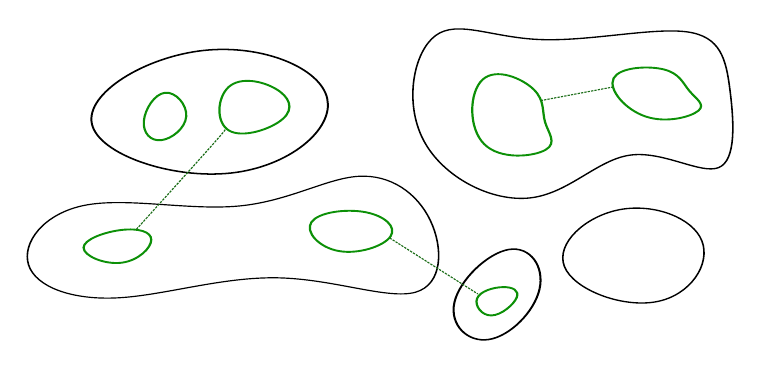}
      %% opens of \bar{r}
      \put(14,16){\small  $r_{2,2}$}
      \put(44,18){\small $r_{3,1}$}
      \put(65,12){\small $r_{3,2}$}
      \put(21,33){\small $r_1$}
      \put(32,34) {\small $r_{2,1}$}
      \put(65,33) {\small $r_{4,1}$}
      \put(85,36) {\small $r_{4,2}$}
      %% connected components of b
      \put(21, 43){\small $b_1$}
      \put(94, 44){\small $b_2$}
      \put(32,9){\small $b_3$}
      \put(70,5){\small $b_4$}
      \put(93,13) {\small $b_5$}
    \end{overpic}  
    \captionsetup{format=plain}
    \caption{
      An example of an active map \(\alpha \colon \bar{r} \ra b\) giving rise to three equivalence classes on \(\pi_0(b)\).
      The dashed lines connect the different connected components of the same \(r_i\). }
    \label{fig:ExampleActiveMapEqClasses}
  \end{figure}
\end{ex}

Just knowing the equivalence relation \(\sim_\alpha\) on \(\pi_0(b)\)
does not suffice to recover an active map \(\alpha\colon \bar{r}\to b\)
because we only know which connected components should be equivalent,
but not by which opens this should be realised.

The following lemma shows that one
can recover an active-\pcocartesian{} map \(\gamma \colon \bar{s} \xrightarrow{!} b\)
from its equivalence relation \(\sim_\gamma\).
Recall from \Cref{defn:DecomposablePresieve} that a \posetofopens{}
\(\PS\) is \emph{decomposable} if for any open \(b \disjun b'\in \PS\)
with \(b\neq \emptyset\neq b'\)
it follows that both \(b\in \PS\) and \(b'\in \PS\).

\begin{lemma}
  \label{lem:cocartesian-eq-rel}
  Let \(\PS\) be a poset of opens with \(\emptyset \notin\PS\) and \(b\in \PS\).
  The assignment \(\gamma\mapsto {\sim}_\gamma\) of
  \Cref{obs:ActiveMapsAsEqClasses}
  yields a fully faithful embedding
  \begin{equation}
    \label{eq:embedding-cocart-eq-rel}
    \left\{\bar{s}\xrightarrow{!}b\right\}
    \hookrightarrow
    \left\{
      {\sim} \text{ on } \pi_0(b)
      \,
      \middle |
      \,
      |\quotrel{\pi_0(b)}{\sim}| < \infty 
    \right\}
  \end{equation}
  of the category of active-\pcocartesian{} arrows in \(\PS\) over \(b\)
  (and active-\pcocartesian{} arrows between them)
  into the poset of equivalence relations on \(\pi_0(b)\)
  with finitely many equivalence classes
  (and inclusions between them).

  If the poset \(\PS\) is decomposable, then this embedding is an equivalence.
\end{lemma}
\begin{proof}
  The inverse functor sends an equivalence relation \(\sim\)
  on \(\pi_0(b)\)
  to the active-\pcocartesian{} arrow
  \begin{equation}
    \left(I_+, \left(\bigdisjun V\right)_{V\in I}\right)\to (1_+, b)
  \end{equation}
  where \(I\coloneqq \quotrel{\pi_0(b)}{\sim}\)
  is the finite set of equivalence classes of \(\sim\).
  This is well defined on the full subposet
  of those equivalence relations \(\sim\)
  such that for each equivalence class \(V\in \quotrel{\pi_0(b)}{\sim}\)
  the disjoint unions \(\bigdisjun V\) lies in \(\PS\). 
  Moreover, each active-\pcocartesian{} arrow \(\bar{s}\xrightarrow{!}b\)
  is recovered from the equivalence relation \({\sim}_{\gamma}\)
  in this way because the empty set cannot appear as one of the \(s_i\).

  If \(\PS\) is decomposable then every equivalence relation \(\sim\) lies in the image
  of the embedding~\eqref{eq:embedding-cocart-eq-rel}:
  in general we have a decomposition  \(b=\bigdisjun V \disjun b'\),
  where \(b'=\bigdisjun_{W\neq V}\bigdisjun W\)
  is the disjoint union of all components of \(b\) that do not lie in the equivalence class \(V\).
  When \(\PS\) is decomposable it follows that \(\bigdisjun V \in \PS\).
  (Note that the case \(\bigdisjun V=\emptyset\) cannot happen and in the case \(b'=\emptyset\)
  we have \(\bigdisjun V=b\in \PS\) anyway.)
\end{proof}

We now consider some different ways of capturing the data of an active-\pcocartesian{} factorization of an active map. 

\begin{lemma}
  \label{lemma:TFAEActiveCocartesianFactorizations}
	Let \(\PSTwo\) be a \posetofopens{} and fix an object \(r\in \PSTwo\).
  Let \(\alpha \colon \bar{b} \ra r\) and \(\gamma \colon \bar{t}\xrightarrow{!} r\)
  be an active and an active-\pcocartesian{} map into \(r\), respectively.
  Assume that the tuple \(\bar{b}\) does not contain the empty set as an element.
  Then the following are equivalent pieces of data:
	\begin{enumerate}[label=(\arabic*), ref=(\arabic*)]
  \item \label{item:TFAE_Lemma_activeMap}
    an active map \(\beta\colon \bar{b} \to \bar{t}\) such that the diagram commutes
		\begin{equation} \label{eq:CommDiagramToSurjMaps}
		\begin{tikzcd}
			\bar{b} \arrow[rr, "\alpha"] \arrow[dr,swap, dashed, "\beta"] && r
			\\ & \bar{t} \arrow[ur, "!"] \arrow[ur, swap, "\gamma"]
		\end{tikzcd}
		\end{equation}
  \item \label{item:TFAE_Lemma_surjectiveMap}
    a surjective map \(\EqClass{\alpha}\twoheadrightarrow \pi(\bar{t})\) such that the diagram commutes
		\begin{equation}
		\begin{tikzcd}
			& \EqClass{\alpha} \arrow[rd, dashed, twoheadrightarrow]
			\\ \pi_0(r) \arrow[ru, twoheadrightarrow] \arrow[rr, twoheadrightarrow] && \pi(\bar{t}) = \EqClass{\gamma}
		\end{tikzcd}
		\end{equation}
  \item\label{item:TFAE_Lemma_subseteq}
    an inclusion \(\sim_{\alpha}\, \subseteq\, \sim_{\gamma}\).
	\end{enumerate}
  If any of these pieces of data exists
  then so do the others and each of them is unique.
\end{lemma}
\begin{proof}
	Follows directly by unravelling the definitions. 
\end{proof}

We now construct a category that assembles all active-\pcocartesian{} factorizations of a fixed active map:
	
\newcommand{\Kalpha}{K_\alpha}

\begin{constr}\label{cstr:KalphaCategory}
	Let \( \PS \hookrightarrow \PSTwo\) be an inclusion of \posetsofopens{}.
  Let \(\bar{r}\in \operad{\PSTwo}\) and fix an object \(\alpha \colon \bar{b} \ra \bar{r} \) in \(\overcatAct{\operad{\PS}}{\bar{r}}\).
  We define the \emph{category of factorizations of \(\alpha\)}, denoted \(\Kalpha\), as follows: 
		\begin{itemize}
		\item
      objects are factorizations \(\bar{b} \xrightarrow{\beta} \bar{t} \xrightarrow{\gamma} \bar{r}\),
      where \(\beta\) is an active morphism in \(\operad{\PS}\) and \(\gamma\) is an active-\pcocartesian{} morphism in \(\operad{\PSTwo}\) 
      such that \(\gamma \circ \beta = \alpha\).
		\item
      Let \(\bar{b} \xrightarrow{\beta'} \bar{t}' \xrightarrow{\gamma'} \bar{r}\) be a second object.
      A morphism from \((\gamma, \beta)\) to \((\gamma',\beta')\)
      is an active-\pcocartesian{} morphism \(\lambda\) such that the diagram 
			\begin{equation}
        \label{eq:map-in-Kalpha}
				\begin{tikzcd}[row sep=1.5ex, column sep=9.0ex, ampersand replacement=\&]
				 \& \bar{t} \arrow[dr, "\gamma"', "!"] \arrow[dd, swap, "\lambda", "!"'] 
		 		\\ \bar{b} \arrow[ur, "\beta"] \arrow[dr, swap, "\beta'"] \& \& \bar{r}
				 \\ \& \bar{t}' \arrow[ur, swap, "\gamma' ", "!"']
			\end{tikzcd}
			\end{equation}
			commutes. 
	\end{itemize}
\end{constr}

\begin{lemma}
  \label{lem:Kalpha-eq-rel}
	Let \( \PS \hookrightarrow \PSTwo\) be an inclusion of \posetsofopens{}
  and assume that \(\emptyset\notin\PS\).
  If \(r=\bar{r}\) lies over \(1_+\) then
  the category \(\Kalpha\) embeds fully faithfully into
  the poset of equivalence relations \(\sim\) on \(\pi_0(r)\) (under inclusion)
  such that
  \begin{enumerate}
  \item \label{enumerate:condition1}
    there are finitely many equivalence classes, 
  \item \label{enumerate:condition2}
    for each equivalence class \(V\), the open set \(\bigdisjun V\) lies in \(\PS\), 
  \item \label{enumerate:condition3}
    and \(\sim_\alpha\ \subseteq\ \sim \). 
  \end{enumerate}
  If the \posetofopens{} \(\PS\) is decomposable,
  then this embedding is an equivalence.
\end{lemma}
\begin{proof}
  Since for a given \(\gamma\colon \bar{t}\to r\)
  there is at most one \(\beta\colon \bar{b}\to\bar{t}\)
  with \((\gamma,\beta)\in \Kalpha\),
  we have that \(\Kalpha\) includes as a subposet into
  the poset of \pcocartesian{} arrows over \(r\)
  with domain in \(\operad{\PS}\)
  (which is indeed a poset by \Cref{lem:cocartesian-eq-rel}).
  It is actually a full subposet since in a diagram
  of the form \eqref{eq:map-in-Kalpha}
  the left triangle automatically commutes
  because there is at most one active map \(\beta''\colon\bar{b}\to\bar{t}'\)
  with \(\gamma'\circ\beta''=\alpha\).

  Thus \(\Kalpha\) is equivalent to the essential image of this embedding
  which consists of those active-\pcocartesian{} arrows \(\gamma\colon\bar{t}\xrightarrow{!}r\)
  where \(\bar{t}\) lies in \(\operad{\PS}\)
  and where \(\alpha\) factors through \(\gamma\).
  Under the embedding of \Cref{lem:cocartesian-eq-rel},
  these correspond precisely to the equivalence relations which in addition to
  having finitely many equivalence classes
  satisfy the conditions
  \ref{enumerate:condition2} and \ref{enumerate:condition3}.
\end{proof}
 
\begin{lemma} \label{lemma:FactCoCartContractible}
	Let \( \PS \hookrightarrow \PSTwo\) be an inclusion of \posetsofopens{}.
  Assume that \(\emptyset\notin\PS\) and that \(\PS\) is decomposable or a presieve. 
  Then the category \(\Kalpha\) is either empty or weakly contractible
\end{lemma}
\begin{proof}
	Without loss of generality, assume that \(\bar{r} = r\).
  We use the description of \Cref{lem:Kalpha-eq-rel} and show that \(\Kalpha\)
  is codirected as long as it is not empty.
  Indeed for two equivalence relations \(\sim\) and \(\sim'\) on \(\pi_0(r)\)
  satisfying conditions \ref{enumerate:condition1}--\ref{enumerate:condition3},
  the same is true for the equivalence relation \(\sim\cap \sim'\).
  This is clear for \ref{enumerate:condition1} and \ref{enumerate:condition3}.
  For \ref{enumerate:condition2} we use that \(\PS\) is either closed under intersections (i.e.\ a presieve) or decomposable. 
\end{proof}

\begin{rem}
	While conditions \ref{enumerate:condition1} and \ref{enumerate:condition3} are closed upward,
  condition \ref{enumerate:condition2} is closed downward,
  so it can happen that there are no equivalence relations satisfying all three conditions at the same time.
\end{rem}

\begin{cor}\label{cor:CoCartMapsIntoFixedTarContractible}
	Let \(\PS \hookrightarrow \PSTwo\) be an inclusion of \posetsofopens{},
  where \(\PS\) is decomposable and \(\emptyset\notin\PS\).
  Fix some object \(\bar{r}\in \operad{\PSTwo}\).
  Then the subcategory of active-\pcocartesian{} maps from \(\operad{\PS}\) into \(\bar{r}\), i.e.
	\begin{equation}
		\left\{	\bar{t} \xrightarrow{!} \bar{r}	\right\} \subset \overcatAct{\operad{\PS}}{\bar{r}},
	\end{equation}
	is either empty or weakly contractible. 
\end{cor}
\begin{proof}
  Note that each object of \(\PS\) admits a unique nullary operation,
  so that the empty tuple is an initial object of \(\operad{\PS}\).
  Hence we obtain the claim as a special case
  of \Cref{lemma:FactCoCartContractible}
  corresponding to setting \(\bar{b}\) to the empty tuple.
\end{proof}

\newcommand{\CaseOneMap}{g} %the map below
\newcommand{\TarCat}{\CanFactorCat } % Target category
\newcommand{\SCat}{E} %Source category
\newcommand{\Functor}{F}
\newcommand{\KalphaTwo}{K_{\alpha'}}

\begin{constr}\label{constr:CartesianFibUnstraightening}
	Let \( \PS \hookrightarrow \PSTwo\) be an inclusion of \posetsofopens{}
  with \(\emptyset\notin\PS\).
  Fix an object \(\bar{r}\in \PSTwo\).  Let \(\TarCat \coloneqq\{\bar{b} \xrightarrow{\cCFactor} \bar{r} \} \) denote the subcategory 
  of \( \overcatAct{\operad{\PS}}{\bar{r}}\) where the active map admits at least one factorization into an active morphism followed by 
  an active-\pcocartesian{} morphism. We define a functor 
	\begin{equation}
	 \Functor \colon \TarCat^{\op}  \longrightarrow \Cat
	\end{equation}
	as follows: 
	\begin{itemize}
	\item Objects \(\alpha\) are sent to the category \(\Kalpha\) of factorizations by an active and active-\pcocartesian{} morphism, as defined in \Cref{cstr:KalphaCategory}. 
	\item
    Let
    \(\alpha\colon \bar{b}\to \bar{r}\) and \(\alpha'\colon \bar{b'}\to\bar{r}\)
    be two objects of \(\TarCat\), and \(\mu\colon \bar{b}\to\bar{b'}\)
    a morphism between these two objects,
    i.e.\ such that \(\alpha = \alpha' \circ \mu\).
    The functor \(\Functor\) sends the morphism \(\mu\)
    to the functor \(\mu^\ast \colon \KalphaTwo \ra \Kalpha\)
    defined by precomposing with \(\mu\).
    Explicitly, on objects \(\mu^\ast\) is defined by 
    
    \begin{equation}
      \begin{tikzcd}[row sep=1.0ex, column sep=3.5ex]
        \bar{b}' \arrow[rr, "\alpha'"] \arrow[ddr, "\beta' "']  && \bar{r}   &&&& \bar{b} \arrow[rr, "\alpha' \circ \mu"] \arrow[ddr, "\beta' \circ \mu" ']&& \bar{r}
        \\ &&& \arrow[rr, mapsto] && {}
        \\ &\bar{t}' \arrow[uur, swap, "\gamma'"] &&  &&&&  \bar{t}' \arrow[uur, "\gamma'"']
      \end{tikzcd},
    \end{equation}
	and analogously on morphisms. 
	\end{itemize}
Unstraightening the functor \(\Functor\) gives the cartesian fibration
	\begin{equation}\label{eq:CartFibLemma}
		 \left\{ \bar{b} \ra \bar{t} \xrightarrow{!} \bar{r} \right\} \stackrel{\CaseOneMap}{\longrightarrow} \left\{ \bar{b} \xrightarrow{\cCFactor} \bar{r} \right\}. 
	\end{equation}	
	All objects of \(\TarCat\) admits at least one factorization by construction, so from \Cref{lemma:FactCoCartContractible} we know that all the fibers of \(\CaseOneMap\) are weakly contractible. Thus, the functor \(\CaseOneMap\) is colimit cofinal because it is a cartesian fibration with weakly contractible fibers. 
\end{constr}

\newcommand{\MapName}{f}
\newcommand{\SCatTwo}{E'}
\begin{cor} \label{cor:CartesianFibrationTrivFibersArg}
	Let \(\PS \hookrightarrow \PSTwo\) be an inclusion of \posetsofopens{}, where \(\PS\) is decomposable and \(\emptyset\notin\PS\).
  Let \(\bar{r}\in \operad{\PSTwo}\) and fix an active morphism \(\alpha\colon \bar{s} \ra \bar{r}\) 
	in \(\overcatAct{\operad{\PS}}{\bar{r}}\). The functor 
  \begin{equation} 
    \SCatTwo \coloneqq
    \left\{
      \begin{tikzcd}[ampersand replacement=\&]
        \bar{b} \arrow[r] \arrow[d] \& \bar{s}  \arrow[d, "\alpha"]
        \\ \bar{t} \arrow[r, "!"] \& \bar{r}
      \end{tikzcd}
    \right\} 
    \stackrel{\MapName}{\longrightarrow} 
    \left\{
      \begin{tikzcd}[ampersand replacement=\&, column sep=1.2em]
        \bar{b} \arrow[rr, ] \arrow[rd, swap, "\cCFactor"] \&\& \bar{s} \arrow[dl, "\alpha"] 
        \\ \& \bar{r}	 
      \end{tikzcd}
    \right\}
    \eqqcolon \overcat{\TarCat}{\alpha},
  \end{equation}
  where \(\bar{b} \in \operad{\PS}\), is a cartesian fibration with weakly contractible fibers. In particular, the functor \(\MapName\) is colimit cofinal. 
\end{cor}
\begin{proof}
	Observe that the category \(\SCatTwo\) is the pullback of the functor \(\CaseOneMap\) from \Cref{constr:CartesianFibUnstraightening} along
  the forgetful functor \(\overcat{\TarCat}{\alpha} \ra \TarCat\),
  and recall that the pullback of a cartesian fibration is again a cartesian fibration. Since the fibers of \(\CaseOneMap\) are weakly contractible, it follows that the fibers of \(\MapName\) are also weakly contractible, and hence the functor \(\MapName\) is also colimit cofinal. 
\end{proof}

Note that one important consequence of decomposability is that for any object \(\bar{s} \in \operad{\PS}\)
the \infy-operad \(\operad{\PS}\) contains all \pcocartesian{} maps of the form
\(\bar{k} \xrightarrow{!} \bar{s}\) in \(\operad{\PSTwo}\)
where the tuple \(\bar{k}\) does not contain the empty set as an element.

\newcommand{\RightAdj}{R} %the right adjoint below
\begin{lemma} \label{lemma:RightAdjToProj}
  Let \(\PS\) be a decomposable \posetofopens{} with \(\emptyset\notin\PS\)
  and let \(\PS \hookrightarrow \PSTwo\coloneqq\disjunComp{\PS}\) be the inclusion into its disjoint union completion.
  Fix some object \(\bar{r} \in \DisjunComp{\PS}\), and fix an active morphism \(\alpha\colon \bar{s} \ra \bar{r}\) in \(\overcatAct{\operad{\PS}}{\bar{r}}\).
  For each active-\pcocartesian{} arrow \(\bar{t}\xrightarrow{!}\bar{r}\)
  with \(\bar{t}\in \operad{\PS}\)
  the pullback
  \begin{equation}
    \cdsquare[pb]
    {\bar{t}\times_{\bar{r}}\bar{s}}
    {\bar{s}}
    {\bar{t}}
    {\bar{r}}
    {!}
    {}{}
    {!}
  \end{equation}
  in \(\DisjunComp{\PS}\) exists
  and the top structure map is active-\pcocartesian{} as indicated.
  Moreover the object \(\bar{t}\times_{\bar{r}}\bar{s}\) again lies in
  \(\operad{\PS}\)

  These pullbacks assemble to yield a right adjoint \(R\)
  to the projection functor
	\begin{equation} \label{eq:ReducingColimitRLKE}
		\left\{ \begin{array}{c}
		\begin{tikzcd}[ampersand replacement=\&]
			\operad{\PS}\ni\bar{b} \arrow[r] \arrow[d] \& \bar{s}  \arrow[d, "\alpha"]
			\\ \bar{t} \arrow[r, "!"] \& \bar{r}
		\end{tikzcd}
		\end{array} \right\} 
		\begin{array}{c}
		\begin{tikzcd}[ampersand replacement=\&, column sep=1.8em]
			{}  \arrow[r,  "\mathrm{proj}"] \& {} \arrow[l, bend left=30, yshift=-1.5ex, "\RightAdj"]
		\end{tikzcd}
		\end{array} 
		\left\{ \begin{array}{c}
		\begin{tikzcd}[ampersand replacement=\&, column sep=1.2em]
			  \& \bar{s} \arrow[d, "\alpha"] 
			\\ \bar{t} \arrow[r, "!"]\& \bar{r}	 
		\end{tikzcd}
		\end{array} \right\}  
	\end{equation}
\end{lemma}
\begin{proof}
  We only have to show that the indicated pullbacks exist;
  the existence of the right adjoint \(R\) is then a formal consequence.

  To construct the pullbacks,
  we start with the solid arrows on the left-hand side of
  \eqref{eq:DiagramDefineK} below.
  Using \Cref{obs:ActiveMapsAsEqClasses} we translate this to the solid arrows on the right-hand side.
  Let \(\bar{k} \ra s\) be the \pcocartesian{} map corresponding to the epi-mono factorization of the composite map \(\pi_0(s) \ra \pi(\bar{t})\), 
  i.e.\ \(\pi_0(s) \twoheadrightarrow \pi(\bar{k}) \hookrightarrow \pi(\bar{t})\). This is well-defined because \(\PS\) is decomposable. 
	\begin{equation} \label{eq:DiagramDefineK}
		 \begin{tikzcd}[column sep=1.2em, row sep=1.2em]
			\bar{k} \arrow[rrr, dashed, "!"] \arrow[dd, dashed] & & & s \arrow[dd, "\alpha"] & & & & \pi_0(s) \arrow[dd, swap, "\pi_0(\alpha)"]  \arrow[rrrr, dashed, twoheadrightarrow]& & & & \pi(\bar{k}) \arrow[dd, hook', dashed]
			\\ &   		&&& 		\arrow[rr, rightsquigarrow] & & {} 
			 \\ \bar{t} \arrow[rrr, "!"]&& & r & & & & 			 \pi_0(r) \arrow[rrrr, twoheadrightarrow]& & & & \pi(\bar{t})
		\end{tikzcd}
	\end{equation}

  We now want to show that \(\bar{k}\) as defined above indeed is a pullback.
  Let \(\bar{b}\) together with the maps on the left-hand side of \eqref{eq:ComDiagramToSurjectiveMapsOne} below be an arbitrary object of the target category of \(\RightAdj\).
  By \Cref{obs:ActiveMapsAsEqClasses} and \Cref{lemma:TFAEActiveCocartesianFactorizations} this translates to the commuting diagram on the right-hand side below. 

	\begin{equation} \label{eq:ComDiagramToSurjectiveMapsOne}
		 \begin{array}{c}
		 \begin{tikzcd}[ampersand replacement=\&, column sep=1.2em, row sep=1.2em]
			\bar{b} \arrow[rrr, "\beta"] \arrow[dd]  \& \& \& s \arrow[dd, "\alpha"]
			\\ \& {}  
			 \\ \bar{t} \arrow[rrr, "!"]\&\& \& r
		\end{tikzcd}
		\end{array} 
		\rightsquigarrow
		 \begin{array}{c}
		 \begin{tikzcd}[ampersand replacement=\&, column sep=1.0em, row sep=0.9em]
			 \& \& \& \EqClass{\beta} \arrow[dd]
			 \\ \pi_0(s) \arrow[dd, swap, "\pi_0(\alpha)"] \arrow[rrru, twoheadrightarrow] \& \& \& \& {}
			 \\ \&\&\& \EqClass{\alpha\circ\beta} \arrow[rd, twoheadrightarrow]
			 \\ \pi_0(r) \arrow[rrru, twoheadrightarrow] \arrow[rrrr, twoheadrightarrow]\& \& \& \& \pi(\bar{t})
		\end{tikzcd}
		\end{array} 
	\end{equation}

  To say that \(\bar{k}\) is a pullback we need to produce
  a (necessarily unique) active map \(\delta\)
  as indicated by the dashed arrow in \eqref{eq:SliceCatSomeObjectSecond} below,
  or equivalently a surjective map \(\EqClass{\beta} \twoheadrightarrow \pi(\bar{k})\).
  By the definition of \(\bar{k}\)
  (using the epi-mono factorization of \(\pi_0(s) \ra \pi(\bar{t})\)),
  we get a dashed surjective map \(\EqClass{\beta} \twoheadrightarrow \pi(\bar{k})\) making the right-hand side diagram below commute, which ensures the existence of \(\delta\).

	\begin{equation} \label{eq:SliceCatSomeObjectSecond}
		 \begin{array}{c}
		 \begin{tikzcd}[ampersand replacement=\&, column sep=1.2em, row sep=1.2em]
			\bar{b} \arrow[rrr, "\beta"] \arrow[dd]  \arrow[dr, dashed, "\delta"]\& \& \& s \arrow[dd, "\alpha"]
			\\ \& \bar{k} \arrow[rru, "!", " "'] \arrow[dl]
			 \\ \bar{t} \arrow[rrr, "!"]\&\& \& r
		\end{tikzcd}
		\end{array} 
		\leftrightsquigarrow
		 \begin{array}{c}
		 \begin{tikzcd}[ampersand replacement=\&, column sep=1.0em, row sep=0.9em]
			 \& \& \& \EqClass{\beta} \arrow[rd, dashed, twoheadrightarrow] \arrow[dd]
			 \\ \pi_0(s) \arrow[dd, swap, "\pi_0(\alpha)"] \arrow[rrru, twoheadrightarrow] \arrow[rrrr, twoheadrightarrow, crossing over]\& \& \& \& \pi(\bar{k}) \arrow[dd]
			 \\ \&\&\& \EqClass{\alpha\beta} \arrow[rd, twoheadrightarrow]
			 \\ \pi_0(r) \arrow[rrru, twoheadrightarrow] \arrow[rrrr, twoheadrightarrow]\& \& \& \& \pi(\bar{t})
		\end{tikzcd}
		\end{array}
	\end{equation} 
	This completes the proof. 
\end{proof}

\newcommand{\LastCat}{\Da}

\begin{lemma} \label{lemma:RelLKE}
  Let \(\PS\) be a decomposable \posetofopens{} with \(\emptyset\notin\PS\).
  Fix some object \(\bar{r} \in \DisjunComp{\PS}\), and let \(\Aa \in \Alg{\disjunComp{\PS}}\) be an algebra such that \(\Aa\) is multiplicative when restricted to \(\operad{\PS}\). The commutative diagram 
\begin{equation}
	\begin{tikzcd}
		\CanFactorCat \coloneqq\{ \bar{b} \xrightarrow{\cCFactor} \bar{r} \}  \arrow[d, hook] \arrow[r] & \operad{\PS}_{\act} \arrow[r, "\Aa_{\act}"] & \targetcatOT_{\act} \arrow[d, ""]
		\\ \ActOverQ \coloneqq\{ \bar{b} \xrightarrow{\act} \bar{r}\} \arrow[ru] \arrow[rr] & & \Fin
	\end{tikzcd}
\end{equation}
exhibits the diagonal composite \(\Aa\restrict{\Ja}\) as an operadic left Kan extension of its restriction to \(\Ia\). 
\end{lemma}
\begin{proof}
	Fix some arbitrary object \(\bar{s} \xrightarrow{\alpha} \bar{r} \in \ActOverQ\). We need to show that
	\begin{equation}\label{eq:RelLKEWant}
    \ocolim_{{\bar{b}}\,\in\,{\overcat{\CanFactorCat}{\alpha}}} \Aa (\bar{b})
    \xrightarrow{\simeq} 
		\Aa(\bar{s}).
	\end{equation}
	We first simplify the operadic colimit. Consider the functors
	\begin{equation} \label{eq:SliceCatSomeObject}
		\LastCat \coloneqq\left\{ \begin{array}{c}
		\begin{tikzcd}[ampersand replacement=\&]
			{} \& \bar{s}  \arrow[d, "\alpha"]
			\\ \bar{t} \arrow[r, "!"] \& \bar{r}
		\end{tikzcd}
		\end{array} \right\} 
		\xrightarrow{\RightAdj} 
		\left\{ \begin{array}{c}
		\begin{tikzcd}[ampersand replacement=\&]
			\bar{b} \arrow[r] \arrow[d] \& \bar{s}  \arrow[d, "\alpha"]
			\\ \bar{t} \arrow[r, "!"] \& \bar{r}
		\end{tikzcd}
		\end{array} \right\} 
		\xrightarrow{\MapName}
		\left\{ \begin{array}{c}
		\begin{tikzcd}[ampersand replacement=\&, column sep=1.2em]
			 \bar{b} \arrow[rr, ] \arrow[rd, swap, "\cCFactor"] \&\& \bar{s} \arrow[dl, "\alpha"] 
			\\ \& \bar{r}	 
		\end{tikzcd}
		\end{array} \right\}  = \overcat{\CanFactorCat}{\alpha}.
	\end{equation}
By \Cref{cor:CartesianFibrationTrivFibersArg} the rightmost functor \(\MapName\) is colimit cofinal. The leftmost functor \(\RightAdj\) is a right adjoint by \Cref{lemma:RightAdjToProj}, so in particular it is also colimit cofinal.  Since the composition of colimit cofinal functors is again colimit cofinal we have 
\begin{equation}\label{eq:ColimCofinalRelColim}
	\ocolim_{\bar{t}\, \in\, \LastCat}\big( \Aa \circ \MapName \circ \RightAdj \big)(\bar{t}) = \ocolim_{\bar{t}\, \in\, \LastCat} \Aa(\bar{t} \times_{\bar{r}} \bar{s}) \xrightarrow{\simeq}	\ocolim_{\bar{b}\, \in\, \overcat{\CanFactorCat}{\alpha}} \Aa (\bar{b}).
\end{equation}

Observe that the category \(\LastCat\) is non-empty because \(\bar{r} \in \DisjunComp{\PS}\). Since \(\alpha\) is fixed, \(\LastCat\) is weakly contractible by \Cref{cor:CoCartMapsIntoFixedTarContractible}. 
Additionally, all the morphisms of the diagram \(\LastCat\) are active-\pcocartesian{} morphisms in \(\operad{\PS}\), 
so, by assumption, \(\Aa\) sends them to cocartesian morphisms of \(\targetcatOT_{\act}\). 
This puts us in the situation of  \Cref{lemma:CocartContractibleDiagramCatOColim}, which says that the operadic colimit exists 
and all the structure maps are cocartesian. Hence, for any fixed object \(\bar{t}_0\xrightarrow{!} \bar{r} \in \LastCat\) we have 
\begin{equation}
	\ocolim_{\overcat{\CanFactorCat}{\alpha}} \Aa \stackrel{\eqref{eq:ColimCofinalRelColim}}{\simeq} \ocolim_{\bar{t}\, \in\, \LastCat} \Aa(\bar{t} \times_{\bar{r}} \bar{s}) 						\stackrel{\ref{lemma:CocartContractibleDiagramCatOColim}}{\simeq} \bigotimes \Aa(\bar{t}_0 \times_{\bar{r}} \bar{s}) \simeq \Aa(\bar{s}),
\end{equation}
where the last equivalence follows from \(\Aa\) being multiplicative on \(\operad{\PS}\) since the (fixed) pullback comes with 
a \pcocartesian{} map into \(\bar{s}\), i.e.\ \(\bar{t}_0 \times_{\bar{r}} \bar{s} \xrightarrow{!} \bar{s}\). 
\end{proof}

\newcommand{\ThirdMap}{f}
\newcommand{\RAdj}{R}
\newcommand{\LAdj}{L}

Before we are ready to prove the main result of this subsection,
we have a final preliminary lemma.

\begin{lemma} \label{obs:MultiplicativeReducedDisjunComp}
  Let \(\PS\) be a \posetofopens{} and let \(\Aa \in \Alg{\disjunComp{\PS}}\) be an algebra.
  The algebra \(\Aa\) is multiplicative on \(\DisjunComp{\PS}\)
  if and only if \(\Aa\) sends \pcocartesian{} morphisms \(\bar{t} \xrightarrow{!} \bar{r}\)
  with \(\bar{t}\in\operad{\PS}\) to cocartesian morphisms in the target category \(\targetcatOT\). 
\end{lemma}
\begin{proof}
  The forward implication is trivial.
  For the converse, let \(\gamma\colon\bar{k} \xrightarrow{!} \bar{r}\)
  be any \pcocartesian{} morphism in \(\DisjunComp{\PS}\).
  Since \(\disjunComp{\PS}\) is the disjoint union completion of \(\PS\)
  there exists some \(\bar{t}_0 \in \operad{\PS}\) and a \pcocartesian{} arrow
  \(\lambda\colon \bar{t}_0 \to \bar{k}\).
  Then both \(\lambda\) and \(\gamma\lambda\) are \pcocartesian{} morphisms with domain \(\bar{t}_0\in\operad{\PS}\),
  hence are sent to cocartesian arrows by \(\Aa\).
  By the 2-out-of-3 property that holds for cocartesian morphisms,
  it follows that \(\Aa(\gamma)\) is also cocartesian, which is what we had to show.
\end{proof}

We are finally ready to prove the main result of this subsection.

\begin{proof}[Proof of Proposition \ref{prop:MultAlgToDisjCompletion}] 
\label{proof:DisjCompletion}
We start by proving the case where \(\emptyset\notin\PS\).
In this case we have built up various preliminary results that allow us to
easily manipulate the relevant operadic Kan extensions.

		It is clear that restricting a multiplicative algebra on \(\DisjunComp{\PS}\) gives a multiplicative algebra on \(\operad{\PS}\). Conversely, since we are operadic left Kan extending along a fully faithful map, we know that there is an equivalence between \(\AlgM{\PS}\) and the essential image of \(\iota_!\). Thus, we need to compare the essential image to \(\AlgM{\disjunComp{\PS}}\). Explicitly, given \(\Aa \in \Alg{\disjunComp{\PS}}\) such that \(\Aa\restrict{\PS}\) is multiplicative we need to show that 
\begin{equation}\label{eq:IFFEssImEqualsMultiplicative}
	\Aa\simeq \iota_!\Aa\restrict{\PS} \quad	\Longleftrightarrow	\quad \Aa \in \AlgM{\disjunComp{\PS}}.
\end{equation}

Before addressing the above implications, let us once again simplify the colimit of the operadic left Kan extension. Fix some arbitrary object \(\bar{r} \in \DisjunComp{\PS}\), and consider the diagram
	
	\begin{equation}\label{eq:ReducingColimMultEq}
	\begin{tikzcd}
		\left\{ \bar{t} \xrightarrow{!} \bar{r} \right\} \arrow[r, swap, "\RAdj"] &\arrow[l, bend right=30, start anchor={[xshift=-2.2ex]}, swap, "\LAdj"] \left\{ \bar{b} \ra \bar{t} \xrightarrow{!} \bar{r} \right\} \arrow[r, "\CaseOneMap"]& \left\{\bar{b} \xrightarrow{\cCFactor} \bar{r} \right\} \arrow[r, hook]& \left\{ \bar{b} \xrightarrow{\act} \bar{r} \right\},
	\end{tikzcd}
	\end{equation}
  where the objects \(\bar{b}\) and \(\bar{t}\) are always assumed to be in \(\operad{\PS}\).
  Note that \(\iota_! \Aa(\bar{r})\) is computed by the operadic colimit over the rightmost category.
  By \Cref{lemma:RelLKE}, \(\Aa\) restricted to the rightmost category \(\{\bar{b} \xrightarrow{\act} \bar{r}\}\) is an operadic left Kan extension of \(\Aa\) restricted to \(\{\bar{b} \xrightarrow{\cCFactor} \bar{r}\}\). The middle functor \(\CaseOneMap\) is colimit cofinal by \Cref{constr:CartesianFibUnstraightening}. The functors \(\LAdj \) and \(\RAdj\) in \eqref{eq:ReducingColimMultEq} are defined as follows on objects:
\begin{equation}
	\LAdj(\bar{b} \ra \bar{t} \xrightarrow{!} \bar{r}) = \bar{t}\xrightarrow{!} \bar{r}, \qquad \qquad \RAdj(\bar{t} \xrightarrow{!} \bar{r})= \bar{t} \xrightarrow{\id} \bar{t} \xrightarrow{!} \bar{r},
\end{equation}
 and analogously on morphisms. This defines an adjunction where \(\RAdj\) is a right adjoint, hence it is in particular colimit cofinal. 
 
Note that the leftmost category in \eqref{eq:ReducingColimMultEq} is non-empty because \(r\in \DisjunComp{\PS}\), and hence it is weakly contractible by \Cref{cor:CoCartMapsIntoFixedTarContractible}. Moreover, every morphism of this category is carried to a cocartesian morphism in the target category because \(\Aa\restrict{\PS}\) is multiplicative. Hence, we are in the situation of \Cref{lemma:CocartContractibleDiagramCatOColim}, which tells us that for any fixed object \(\bar{t}_0 \in \{\bar{t} \xrightarrow{!} \bar{r}\}\) the structure map into the corresponding operadic colimit is cocartesian. 
In summary, for any fixed active-\pcocartesian{} morphism \(\bar{t}_0\xrightarrow{!}\bar{r}\) we have
\begin{equation}\label{eq:proofIFFMultOLKE}
	\iota_! \Aa\restrict{\PS}(\bar{r}) = \ocolim_{\{\bar{b} \xrightarrow{\act} \bar{r}\}} \Aa\restrict{\PS} (\bar{b}) \simeq \ocolim_{\{\bar{t} \xrightarrow{!} \bar{r}\}} \Aa\restrict{\PS}(\bar{t}) \simeq \bigotimes \Aa(\bar{t}_0).
\end{equation}
Recall from \Cref{obs:MultiplicativeReducedDisjunComp} that \(\Aa\) is multiplicative on \(\DisjunComp{\PS}\) if and only if it is multiplicative on the objects of the category \(\{\bar{t}\xrightarrow{!}\bar{r}\}\). Thus the forward implication of \eqref{eq:IFFEssImEqualsMultiplicative} follows from the left-hand side operadic colimit of \eqref{eq:proofIFFMultOLKE} by assumption being equivalent to \(\Aa(\bar{r})\). Conversely, assuming that \(\Aa\) is multiplicative we know that \(\bigotimes \Aa(\bar{t}_0)\simeq \Aa(\bar{r})\), which, by reading \eqref{eq:proofIFFMultOLKE} backwards, gives that \(\Aa\simeq \iota_! \Aa\restrict{\PS}\). 

This concludes the proof in the case \(\emptyset\notin\PS\).
Let us now prove the remaining case by assuming \(\emptyset\in\PS\).
Define \(\PS'\coloneqq \PS\setminus\{\emptyset\}\) obtained from \(\PS\)
by simply removing the empty set.
Since \(\PS\) was decomposable, the same is true for \(\PS'\);
and moreover we have \(\cdisj{\PS'}=\cdisj{\PS}\)
because \(\emptyset\) is the disjoint union of the empty list.

By applying the already proven case to \(\PS'\)
(which satisfies \(\emptyset\notin\PS'\) by construction),
we thus obtain that the composite of the two restriction functors
\begin{equation}\label{eq:composite-restriction-functors-disjun-comp-and-adding-empty-set}
  \AlgM{\cdisj{\PS}}
  \to
  \AlgM{\PS}
  \to
  \AlgM{\PS'}
\end{equation}
is an equivalence.
The result follows, because the second map is an equivalence by
\Cref{lem:remove-empty-set-multiplicative}.
\end{proof}

\newpage

\section{Constructibility tools}
\label{sec:external-tools}

In this section we collect various tools to deal
specifically with \emph{constructible} factorization algebras;
tools which are not available or meaningful in the general situation.

Some versions of these tools were already present in the literature,
mainly in \cite{Ginot}, \cite{AFprimer} and \cite{AFT-fh-stratified}.
Unfortunately the available proofs
are either not general enough or not detailed enough for our purposes,
thus requiring us to expand them and record our own versions.
As an added benefit, it makes this paper mostly self-contained
and better suited as a reference document.

We also extract some of the abstract and reusable intermediate results
that were previously left implicit;
hopefully, this make future generalizations and variations less laborious.
Since these results are of very general category-theoretic nature,
we collect them in \Cref{app:categories} for ease of navigation and reference.

\subsection{Constructibility is local}
\label{sec:constructible}
\newcommand\mixedsector[3]{S_{#1}(#2,#3)} %sectors of the mixed cone

In this section we prove that constructibility of a factorization algebra
can be checked locally.
This is claimed as \cite[Proposition 24]{Ginot}, albeit without proof.
The idea of this proof is adapted from \cite[Proposition 13]{Ginot},
but allows for general targets,
streamlines the argument, and works in the stratified setting.

\begin{defn}\label{defn:LocallyConstructibleFA}
	Let \(\Aa\) be a factorization algebra on \(X\).
  We say that \(\Aa\) is \emph{locally constructible}
  if there is an open cover \(\Ua\) of \(X\)
  such that for each \(U\in \Ua\), the restriction
  \(\Aa\restrict{U}\) is a constructible factorization algebra on \(U\).
\end{defn}

\begin{thm}[Constructibility is local]
  \label{thm:constructle-local}
  Let \(M\) be a conical manifold with enough good disks
  and \(\Aa\) a factorization algebra on \(M\).
  Then \(\Aa\) is constructible if and only if it is locally constructible.
\end{thm}

\begin{cor}
  \label{cor:locally-locally-constant}
  Every locally locally constant\footnote{
    Usually in topology,
    ``locally P'' means ``P for every sufficiently fine open cover''.
    However it does not really make sense to talk
    about ``constant'' factorization algebras,
    unless one means the trivial one sending everything
    to the monoidal unit of \(\targetcat\).
    Therefore one might argue that the established terminology
    ``locally constant'' is a bit of a misnomer since it does not mean ``constant on any sufficiently small open set''.
    \Cref{cor:locally-locally-constant} shows that at
    least being ``locally constant'' is a local property, as the name suggests.
  } factorization algebra on a smooth manifold
  is locally constant.
\end{cor}

We collect some intermediate lemmas before we prove \Cref{thm:constructle-local}.

\begin{lemma}
  \label{lem:ball-in-one-direction}
  Let \(X\) be a conical manifold with enough good disks and
  \(\Aa\) a constructible factorization algebra
  on \(\reals^n\times X\).
  Then the pushforward $\pf{p}\Aa$ along the projection $p\colon \reals^n\times X \to \reals^n$ of $\Aa$ is constructible,
  that is, for each inclusion \(U\subseteq V\) of disks in \(\reals^n\)
  the inclusion \(U\times X\hookrightarrow V\times X\) is \(\Aa\)-local.
\end{lemma}

\begin{proof}
  Let \(\Bf\subset\diskX\) be a factorizing basis of \(X\).
  Then for each \(W\in \Bf\) we have an inclusion
  \(U\times W\hookrightarrow V\times W\)
  of abstractly isomorphic conical disks,
  which is hence \(\Aa\)-local.
  This yields the top horizontal equivalence in the following commutative square
  \begin{equation}
    \begin{tikzcd}
      \colim\limits_{W\in \Bf}\Aa(U\times W)
      \ar[d,"\simeq"]
      \ar[r,"\simeq"]
      &
      \colim\limits_{W\in \Bf}\Aa(V\times W)
      \ar[d,"\simeq"]
      \\
      \Aa(U\times X)
      \ar[r]
      &
      \Aa(V\times X)
    \end{tikzcd}
  \end{equation}
  Moreover, the vertical maps are also equivalences since they are induced
  by the Weiss presieves
  \begin{equation}
    \{U\times W\mid W \in \Bf\}\hyprefines U\times X
    \quad
    \text{and}
    \quad
    \{V\times W\mid W \in \Bf\}\hyprefines V\times X.
  \end{equation}
  Thus we conclude that the lower map is an equivalence, as desired.
\end{proof}

\begin{figure}[H]
  \includegraphics*[width=.45\textwidth]{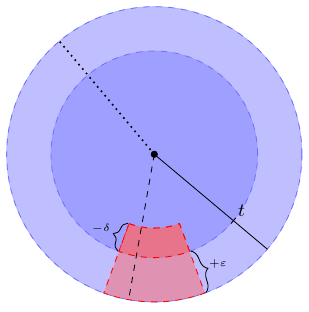}
  \caption{
    The key construction in the proof of
    \Cref{lem:locally-constr-cones}.
    The whole disk is the cone (where we omit the \(\reals^0\)-component)
    on the circle \(X=\sphere{1}\) with three points,
    see \Cref{fig:StratCircleCone}.
    The whole red region is one of the sectors
    \(\mixedsector{(t-\delta,t+\epsilon)}{U}{\reals^0}\)
    with \(U\) being an interval on the circle.
    The darker inner red region is the smaller sector
    \(\mixedsector{(t-\delta,t)}{U}{\reals^0}\)
  yielding the inclusion \eqref{eq:sector-inclusion-in-proof}.
  }
    \label{fig:sectors}
\end{figure}

\begin{lemma}
  \label{lem:locally-constr-cones}
  Fix a compact conical manifold \(X\) with enough good disks,
  a natural number \(n \in \naturals\) and
  a factorization algebra \(\Aa\) on
  the stratified space
  \(\mixedcone[\infty]\coloneqq\topcone{X}\times \reals^n\). 
  Assume that \(\Aa\) is locally constructible away from the cone point
  (i.e., locally constructible on \(\mixedcone[\infty]\setminus\{0\}\)).
  Then for every \(0<t\leq \infty\),
  the inclusion
  \begin{equation}
    \mixedcone[t]\coloneqq \topcone[t]{X}\times \ball{t}{0}
    \hookrightarrow \topcone{X}\times \reals^n=\stratdisk{X}{n}=:\mixedcone[\infty]
  \end{equation}
  is \(\Aa\)-local.
\end{lemma}

\begin{proof}
  We wish to apply \Cref{lem:transfinite-extension}
  to the diagram 
  \(F\colon (0,\infty]\xrightarrow{\mixedcone[-]} \open{\mixedcone[\infty]}\xrightarrow{\Aa} \targetcat\),
  which will immediately yield the desired result
  if we can establish its assumptions.
  Since \(\Aa\) is a Weiss cosheaf,
  assumption \ref{it:equivalence-to-supremum} follows directly
  from the fact that for all subsets \(S\subseteq (0,\infty]\),
  the presieve \(\{\mixedcone[s]\mid s\in S\}\) is a Weiss cover of
  \(\mixedcone[\sup S]\).

  The rest of the proof is devoted to proving assumption
  \ref{it:can-extend-equivalence}.
  Fix \(0<t\leq \infty\).
  We have to find an \(\epsilon>0\) such that the inclusion
  \(\mixedcone[t]\hookrightarrow\mixedcone[t+\epsilon]\)
  is \(\Aa\)-local.
  Using the identification of \Cref{rem:cone-sphere},
  we introduce the following notation:
  Let \(I\subseteq [0,\infty)\) be an interval
  and \(U\subseteq X\) and \(V\subseteq \sphere{n-1}\) open subsets.
  We write \(\mixedsector{I}{U}{V}\subseteq \mixedcone[\infty]\) for the image of 
  \begin{equation}
    I\times U\times I\times V\subseteq
    [0,\infty)\times X\times [0,\infty)\times \sphere{n-1}
    \twoheadrightarrow \topcone{X}\times\topcone{\sphere{n-1}}\cong\mixedcone[\infty];
  \end{equation}
  see \Cref{fig:sectors} for an example illustration.
  Note that the projection  
  \(I^2\times U\times V\twoheadrightarrow \mixedsector{I}{U}{V}\)
  is an isomorphism as long as \(0\notin I\).
  
  Since \(\Aa\) is locally constructible away from the cone point
  and \(X\) and \(\sphere{n-1}\) are compact,
  we can find finite open covers
  \(\Ua\) of \(X\) and \(\Va\) of \(\sphere{n-1}\)
  as well as an \(\frac{t}{2}>\epsilon>0\)
  such that \(\Aa\) is constructible on each
  \(\mixedsector{(t-\epsilon,t+\epsilon)}{U}{V}\)
  with \(U\in\cdown{\Ua}\) and \(V\in\cdown{\Va}\).

  Consider the cover 
  \begin{equation}
    \Wa' \coloneqq 
      \open{\mixedcone[t]} \disjun
    \{\mixedsector{(t-\delta,t+\epsilon)}{U}{V}\mid 0<\delta<\epsilon, U\in \cdown{\Ua}, V\in \cdown{\Va}\}
  \end{equation}
  of \(\mixedcone[t+\epsilon]\).
  Observe that \(\Wa'\) is closed under intersections since
  for all \(O\in \open{\mixedcone[t]}\)
  and \(U,U'\in\cdown{U}\) and \(V,V'\in \cdown{V}\)
  we have
  \begin{equation}
    O \cap \mixedsector{(t-\delta,t+\epsilon)}{U}{V} \in \open{\mixedcone[t]}
  \end{equation}
  \begin{equation}
    \mixedsector{(t-\delta,t+\epsilon)}{U}{V}
    \cap
    \mixedsector{(t-\delta',t+\epsilon)}{U'}{V'}
    =
    \mixedsector{(t-\min(\delta,\delta'),t+\epsilon)}{U\cap U'}{V\cap V'}.
  \end{equation}
  Write \(\Wa\coloneqq \disjunComp{\Wa'} \)
  for the cover obtained from \(\Wa'\) by
  closing under finite disjoint unions;
  it is again a presieve because \(\cap\) distributes over \(\disjun\).
  Note that for each \(W'\in\Wa'\),
  the inclusion \(W'\cap \mixedcone[t]\hookrightarrow W'\)
  is \(\Aa\)-local:
  this is trivial in the case \(W'\in \open{\mixedcone[t]}\);
  in the case \(W'=\mixedsector{(t-\delta,t+\epsilon)}{U}{V}\)
  we can identify this inclusion with
  \begin{equation}
    \label{eq:sector-inclusion-in-proof}
    W'\cap \mixedcone[t]
    \cong
    (t-\delta,t)^2\times U\times V
    \hookrightarrow 
    (t-\delta,t+\epsilon)^2\times U\times V
    \cong W'
  \end{equation}
  which is \(\Aa\)-local by \Cref{lem:ball-in-one-direction}
  because \(\Aa\) is constructible on \(W'\) by the choice of \(\Ua\) and \(\Va\).
  Since \(\Aa\) is multiplicative,
  the same is also true for each \(W\in \Wa\).

  It follows that in the following commutative square of refinements
  \begin{equation}
    \begin{tikzcd}
      \Wa\wedge\{\mixedcone[t]\}
      \ar[d]\ar[r]
      &
      \Wa
      \ar[d]
      \\
      \mixedcone[t]
      \ar[r]
      &
      \mixedcone[t+\epsilon]
    \end{tikzcd}
  \end{equation}
  the top horizontal refinement is \(\Aa\)-local by \Cref{lem:coll-restrict-local}.
  On the left we have a degenerate cover which is trivially \(\Aa\)-local.
  Finally, we claim that on the right side we have a Weiss cover,
  which then is also \(\Aa\)-local because \(\Aa\) is a Weiss cosheaf.
  Indeed, for every finite subset \(T\subseteq \mixedcone[t+\epsilon]\)
  we can choose \(0<\delta<\epsilon\) such that
  \(T\cap \mixedcone[t]\subset \mixedcone[t-\delta]\);
  then \(T\) can be disjointly covered by \(\mixedcone[t-\delta]\)
  and opens of the form \(W_j\coloneqq\mixedsector{(t-\delta,t+\epsilon)}{U_j}{V_j}\)
  for \(U_j\in \cdown{\Ua}\) and \(V_j\in\cdown{\Va}\) as above;
  taking the disjoint union of these then yields the desired open
  \(T\subseteq W\coloneqq \mixedcone[t-\delta]\disjun\bigdisjun_j W_j\in \Wa\).

  We conclude that the bottom horizontal map is \(\Aa\)-local,
  which is exactly what we set out to prove.
\end{proof}

\begin{lemma}
  \label{lem:strat-cones-zero-stratum}
  Let \(U\subset V\) be an inclusion of conical disks
  such that \(U\) and \(V\) are (abstractly) isomorphic.
  Then there exists an isomorphism \(\phi\colon V\cong\stratdisk{X}{n}\)
  and inclusions
  \begin{equation}
    \mixedcone[t] \subset \phi(U)\subseteq \stratdisk{X}{n}.
  \end{equation}
  for all sufficiently small \(t>0\).
\end{lemma}

\begin{proof}
  Without loss of generality we may assume that \(V=\stratdisk{X}{n}\).
  Choose a point \(u\in U_0\) in the minimal stratum \(U_0\) of \(U\).
  Note that \(U_0\subseteq V_0=\{0\}\times \reals^n\),
  hence \(u=(0,t)\) for some \(t\in \reals^n\).
  Hence with the isomorphism 
  \begin{equation}
    \phi\colon V=
    \stratdisk{X}{n}\xrightarrow[\cong]{\id_{\topcone{X}}\times (x\mapsto x-t)}
    \stratdisk{X}{n}
  \end{equation}
  we have \(0\in \phi(U)\).
  Since \(\phi(U)\subset\stratdisk{X}{n}\) is open,
  we also have
  \begin{equation}
    0\in \mixedcone[t]\subset \phi(U)
  \end{equation}
  for every sufficiently small \(t>0\).
\end{proof}

We can now finally prove that constructibility is local.

\begin{proof}[Proof of \Cref{thm:constructle-local}]
  Every constructible factorization algebra is automatically locally constructible.

  To prove the converse direction, assume that \(\Aa\) is locally constructible
  and let \(U\subset V\subseteq M\) be an inclusion of opens,
  both isomorphic to the same standard conical disk
  \(\stratdisk{X}{n}\). We aim to show that the inclusion
  \(U\hookrightarrow V\) is \(\Aa\)-local.

  We apply \Cref{lem:strat-cones-zero-stratum} twice to produce the following diagram of inclusions and identifications
  \begin{equation}
    \begin{tikzcd}
      T\ar[r,hookrightarrow]
      \ar[dd,"{\cong}"]
      &W
      \ar[d,"{\cong}"]
      \ar[r,hookrightarrow]
      &U
      \ar[d,"{\cong}"]
      \ar[r,hookrightarrow]
      &V
      \ar[d,"{\cong}"]
      \\
      &\mixedcone[t]
      \ar[r,hookrightarrow]
      \ar[d,"{\cong}"]
      & \phi(U)
      \ar[d,"{\cong}"]
      \ar[r,hookrightarrow]
      &\stratdisk{X}{n}
      \\
      \mixedcone[t']
      \ar[r,hookrightarrow]
      &\varphi(W)
      \ar[r,hookrightarrow]
      &\stratdisk{X}{n}
    \end{tikzcd} \ . 
  \end{equation}
  Explicitly, we first apply \Cref{lem:strat-cones-zero-stratum} to the inclusion \(U\hookrightarrow V\) 
  producing an isomorphism \(\phi \colon V \cong \stratdisk{X}{n}\),
  from which we obtain \(W \coloneqq \phi^{-1}(\mixedcone[t])\). 
  Secondly, we apply \Cref{lem:strat-cones-zero-stratum} to the inclusion \(W \hookrightarrow U\) 
  producing an isomorphism \(\varphi \colon U \cong \stratdisk{X}{n}\),
  from which we obtain \(T\coloneqq\varphi^{-1}(\mixedcone[t'])\). 
%  by first applying \Cref{lem:strat-cones-zero-stratum}
%  to \(U\subseteq V\) (to obtain $W=\phi^{-1}(\mixedcone[t])$),
 % and then again to \(W\subseteq U\).
  Then the inclusions \(W\hookrightarrow V\)
  and \(T\hookrightarrow U\)
  are \(\Aa\)-local by \Cref{lem:locally-constr-cones}
  which explicitly means that \(\Aa(W)\to \Aa(U)\)
  has both a left and a right inverse, hence is an equivalence.
  But then \(\Aa(U)\to\Aa(V)\) is also an equivalence as desired,
  since it is the retraction of an equivalence.
\end{proof}

Let us also record the following consequence of \Cref{thm:constructle-local}
which says that constructibility is only a non-vacuous condition
away from the \(0\)-dimensional strata.

\begin{cor}
  \label{cor:constructible-away-from-0}
  Let \(M\) be a conical manifold with enough good disks
  and let \(\Aa\) be a factorization algebra on \(M\).
  Let \(S\) be the union of all \(0\)-dimensional strata of \(M\)
  and assume that \(\Aa\) is constructible on \(M\setminus S\).
  Then \(\Aa\) is constructible on all of \(M\).
\end{cor}

\begin{proof}
  By \Cref{thm:constructle-local}, it suffices to show that \(\Aa\)
  is locally constructible at each point \(s\in S\).
  For each \(s\in S\) choose a conical disk
  \(\stratdisk{X}{n}\cong O\subset M\) with \(s\) on its deepest stratum
  and \(O\cap S =\{s\}\).
  Because \(S\) is discrete, such a conical disk exists and we have \(n=0\),
  so that \(O\cong \topcone X\) is a pure cone and \(s\) its cone point.
  We claim that \(\Aa\) is constructible on \(O\).

  Let \(U\hookrightarrow V\)
  be an inclusion of abstractly isomorphic conical disks in \(O\).
  Since they are abstractly isomorphic,
  they either both contain the cone point \(s\) of \(O\)
  or none of them do.
  In the second case, they are both contained in \(M\setminus S\),
  where \(\Aa\) is constructible.
  In the first case we may proceed as in the proof of \Cref{thm:constructle-local}:
  we invoke \Cref{lem:strat-cones-zero-stratum}
  and \Cref{lem:locally-constr-cones}
  (which only requires constructibility away from the cone point) twice
  to construct \(T\hookrightarrow W\hookrightarrow U\)
  such that both \(T\hookrightarrow U\) and \(W\hookrightarrow V\) are \(\Aa\)-local
  and conclude that \(U\hookrightarrow V\) is \(\Aa\)-local as well.
\end{proof}

\subsection{Localizing at isotopy equivalences}
\label{sec:localizing-disk}

Let \(X\) be a smooth conical manifold.
In \cite{AFT-fh-stratified} Ayala--Francis--Tanaka computed
the localization of \(\diskX\) with respect to all isotopy equivalences of disks.
In this section we consider a slight generalization where we localize a sufficiently nice \posetofopens{}
\(\Bf \subset \diskX\) with respect to its subordinate isotopy equivalences and get the same result.

We first define what we mean by the subposet of isotopy equivalences.
\begin{defn}
  Let \(\Bf \subset \diskX\) be a poset of open disks of \(X\).
  Denote by \(\Jmaps[\Bf]\subset \Bf\) the wide subposet of those inclusions
  \(U \hookrightarrow V\) of opens in \(\Bf\)
  which induce a bijection \(\pi_0(U)\xrightarrow{\cong}\pi_0(V)\)
  and such that for each component \(V^j\) of $V$ the conical disks
  \(U^j=V^j\cap U\) and \(V^j\) are abstractly isomorphic.
\end{defn}
By definition, a factorization algebra \(\Aa\) on \(\Bf\) is constructible
precisely if the underlying functor
\(\Aa\colon\Bf \to \targetcat\)
sends all maps in \(\Jmaps[\Bf]\) to equivalences, i.e.,
if it factors through
the localization \(\Bf\to\localize{\Bf}{\Jmaps[\Bf]}\).

As mentioned, when \(\Bf=\diskX\) this localization was computed in \cite{AFT-fh-stratified}:
The result is the slice \(\infty\)-category \(\ddiskofX\),
where \(\ddisk\subseteq \SSnglr\)
is the full subcategory spanned by those smooth conical manifolds which
are disjoint unions of conical disks.

We want to use this localization-result in the context of operadic left Kan extending
from certain presieves later on,
but since arbitrary disks are not closed under intersections we cannot use all of \(\diskX\).
However, a version of their proof can be made to work
as long as we work with a \posetofopens{} with sufficiently many disks;
the authors thank Pelle Steffens for this observation. 

\begin{lemma}[\cite{AFT-local-structures}, Theorem~4.3.1]
  \label{lem:char-iso-disk}
  Let \(U \hookrightarrow V\) be an inclusion of conical disks.
  The following are equivalent:
  \begin{itemize}
  \item
    \label{it:U-V-abstract-iso}
    the disks \(U\) and \(V\) are abstractly isomorphic;
  \item
    \label{it:U-V-touch-lowest-stratum}
    the disk \(U\) contains a point in the deepest stratum of \(V\);
  \item
    \label{it:U-V-isotopy-equiv}
    the inclusion \(U\hookrightarrow V\) becomes an equivalence in \(\ddisk\).
  \end{itemize}
\end{lemma}
\Cref{lem:char-iso-disk} implies that the canonical functor
\begin{equation}
  {\Bf}\xrightarrow{} \ddiskX
\end{equation}
sends all maps of \(\Jmaps[\Bf]\) to equivalences.
We aim to show that this functor is an \(\infty\)-categorical localization
at these maps.

\begin{thm}
  [Variation of {\cite[Proposition~2.22]{AFT-fh-stratified}}]
  \label{thm:Disk-localized}
  Let \(X\) be a smooth conical manifold and
  \(\Bf\subseteq\diskX\) be a decomposable multiplicative disk-basis of \(X\).
  Then the induced map
  \begin{equation}
    \localize{\Bf}{\Jmaps[\Bf]} \xrightarrow{\simeq}
    \ddiskX
  \end{equation}
  is an equivalence of \(\infty\)-categories.
\end{thm}

\begin{rem}\label{rem:localizing_gap}
  As explained to us by David Ayala in private communication,
  the original proof of \cite[Proposition~2.22]{AFT-fh-stratified}
  (which \Cref{thm:Disk-localized} generalizes)
  contains a major gap,
  which was only later closed using a result of Mazel-Gee~\cite{Mazel-Gee}.
  Unfortunately, a revised version has yet to appear in print.
  
  We learned about a version of the missing argument
  from Berry's thesis~\cite{Eric};
  a version of this argument also appeared in \cite{Cepek}.
  
  To make this argument more transparent
  and reusable for later variations on the common localization theme,
  we have abstracted it as a separate, purely category-theoretic statement,
  which we prove as \Cref{lem:Berry-localization} in the appendix.

  After the completion of this paper we were made aware
  of the independent and contemporaneous work by Arakawa,
  who in the unstratified setting also gave a complete proof of
  this localization result for \(\Bf=\diskX\)
  following the same outline and using the same localization lemma;
  see \cite[Theorem 2.24 and Proposition~2.27]{Arakawa}.
\end{rem}

Before we turn to the proof of \Cref{thm:Disk-localized}
we need some notation and preliminary results.

\begin{figure}[H]
\begin{subfigure}{0.4\linewidth}
\begin{minipage}{0.9\linewidth} \begin{center}
\begin{overpic}[scale=0.6, tics=10]{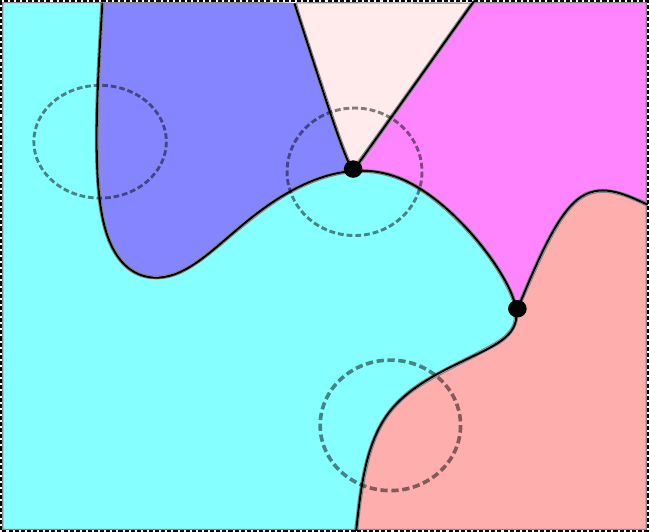}
%% labels for the opens
	\put(72,16){ $U_2^2$}
	\put(63,62){ $U_1$}
	\put (22, 69) {$U_2^1$}
\end{overpic}  
\end{center}
  %\caption{}
 % \label{}
\end{minipage} 
\end{subfigure} \qquad \qquad %\hfill
\begin{subfigure}{0.4\linewidth}
\begin{minipage}{0.9\linewidth} \begin{center}
  \begin{overpic}[scale=0.6, tics=10]{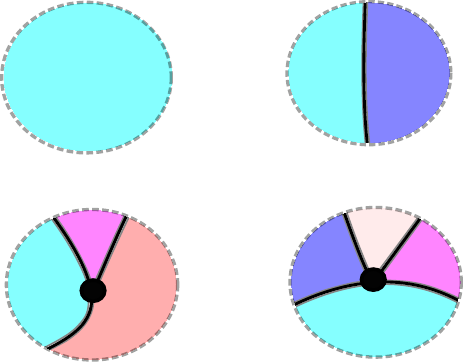} 
  \end{overpic}
  \end{center}
  %\caption{}
  %\label{}
\end{minipage} 
\end{subfigure}
\captionsetup{format=plain}
\caption{ With \(X=\reals^2\), an example of an object \(U \in \disk{X}\)
  and its decomposition \(U=U_1\disjun U_2^1 \disjun U_2^2\)
  from \Cref{not:U-types}~\ref{item:not-U-typesIntDisjUnion}.
  The space \(X\) decomposes into four different pieces \(X_{[D]}\),
  where \(D\) ranges over the types
  indicated on the right.}
\label{fig:NotationDisjoinConicalDisksEmb}
\end{figure}

\begin{notation}
  \label{not:U-types}
  \begin{enumerate}[label=(\arabic*), ref=(\arabic*)]
  \item[]
  \item
    Every object \(U\in\ddisk\) can be written as an (external) disjoint union
    \begin{equation}
      U=\coprod_{i\in I}k_i \times U_i
    \end{equation}
    of pairwise non-isomorphic conical disks \(U_i\), each appearing \(k_i\) times.
  \item \label{item:not-U-typesIntDisjUnion}
    Similarly, we can decompose each \(U\in\diskX\) as an internal disjoint union
    \begin{equation}
      U= \bigdisjun_{i\in I}\bigdisjun_{j=1}^{k_i}U_i^j,
    \end{equation}
    where the \(U_i^j\subset X\) are conical disks,
    two of which are abstractly isomorphic
    if and only if they have the same lower index ``\(i\)''. An example of this is depicted in \Cref{fig:NotationDisjoinConicalDisksEmb}. 
  \item
    In either case, we write \([U]\) for the isomorphism type
    of such a disjoint union of conical disks.
  \item
    For each contractible conical disk \(D\) we write \(X_{[D]}\subset X\)
    for the subspace of those points whose local neighborhoods have type \(D\).
    Each \(X_{[D]}\) is an (unstratified) manifold that is locally closed in \(X\)
    and we have a decomposition \(X=\bigdisjun_{[D]}X_{[D]}\)
    since local neighborhood types are unique
    by \Cref{lem:char-iso-disk}.
  \item
    For a topological space \(Y\), we write \(\Conf[k]{Y}\)
    for the \emph{unordered} configuration space of \(k\) points in \(Y\).
  \item
    Given an isomorphism type \([U]\),
    we write \(\Jmaps[\Bf]^{[U]}\subseteq \Jmaps[\Bf]\)
    for the full poset consisting only of those disjoint unions of disks
    which are of that prescribed type.
  \end{enumerate}
\end{notation}

\begin{lemma}[\cite{AFT-fh-stratified}, Lemma~2.21]
  \label{lem:core-DDisk-is-Conf}
  The underlying \(\infty\)-groupoid of \(\ddiskX\)
  is canonically identified with the space
  \begin{equation}
    \label{eq:core-DDisk-is-Conf}
    (\ddiskX)^\simeq \xrightarrow{\simeq} \coprod_{[U]}\prod_{i\in I}\Conf[k_i]{X_{[U_i]}},
  \end{equation}
  where the coproduct ranges over the isomorphism classes of disjoint unions of conical disks.
\end{lemma}

\begin{rem}
  \label{rem:explicit-core-ddisk-is-Conf}
  From the proof of \cite[Lemma~2.21]{AFT-fh-stratified}
  we can extract an explicit description of the equivalence~\eqref{eq:core-DDisk-is-Conf}:
  an embedding \((\phi\colon U=\coprod_ik_i\times U_i\hookrightarrow X)\)
  is sent to the tuple of configurations \((\phi(S_i))_i\),
  where each \(S_i=\{x_1,\dots,x_{k_i}\}\) contains precisely
  one point \(x_j\in (U_i)_{[U_i]}\) in the deepest stratum of
  each of the \(k_i\) copies of \(U_i\);
  this is well defined because the space of such \(S_i\) is contractible.
\end{rem}

\begin{lemma}
  \label{lem:class-J-conical-manifold}
  Let \(\Bf\subseteq \diskX\) be a decomposable multiplicative basis of \(X\).
  Fix an isomorphism type \([U]=[\coprod_{i\in I}k_i\times U_i]\).
  Then we have a canonical equivalence of \(\infty\)-groupoids
  \begin{equation}
    \classspace{\Jmaps[\Bf]^{[U]}} \simeq \prod_{i\in I}\Conf[k_i]{X_{[U_i]}}.
  \end{equation}
\end{lemma}

\begin{proof}
  For each element \(V=\bigdisjun_{i,j}V_i^j \in \Jmaps[\Bf]^{[U]}\)
  we have the open subset 
  \begin{equation}
    \smallconf{V}\coloneqq
    \left\{\left(S_i \in \Conf[k_i]{X_{[V_i]}}\right)_{i\in I}
      \mid
      \forall_{i\in I}\,\forall_{1\leq j\leq k_i}\colon
      \# (S_{i}\cap V_i^j)= 1\right\}
    \subset
    \prod_{i\in I}\Conf[k_i]{X_{[U_i]}}
  \end{equation}
  consisting of those tuples of configurations \((S_i)_{i}\) of points
  where all points of \(S_i\) have local neighborhood type \([U_i]\)
  (in particular, the \(S_i\) are pairwise disjoint)
  and which induces a bijection \(\pi_0(S_i)\xrightarrow{\cong}\pi_0(V_i)\).
  This assignment assembles into a functor
  \begin{equation}
    \smallconf{-}\colon \Jmaps[\Bf]^{[U]} \to \open{\prod_{i\in I}\Conf[k_i]{X_{[U_i]}}};
    \quad V \mapsto \smallconf{V}.
  \end{equation}

  We now want to apply Lurie's Seifert--Van Kampen
  theorem~\cite[Theorem~A.3.1]{LurHA}
  to \(\smallconf{-}\).
  For each point \(S=(S_i=\{s_i^1,\dots,s_i^{k_i}\})_i\in \prod_{i\in I}\Conf[k_i]{X_{[U_i]}}\)
  we have to check that the full subposet
  \begin{equation}
    K_S\coloneqq\{V \in \Jmaps[\Bf]^{[U]}\mid S \in \smallconf{V}\}\subseteq \Jmaps[\Bf]^{[U]}
  \end{equation}
  is weakly contractible.
  \begin{itemize}
  \item[Claim:] The poset \(K_S\) is cofiltered.
  \item
    Let \(W\subseteq X\) be an arbitary open subset containing \(S\).
    Since \(\Bf\) is a basis for \(X\)
    we can find, for each \(s_i^j\in S_i\),
    an \(V_i^j\in \Bf\) with \(s_i^j\in V_i^j\).
    By possibly shrinking them further, we may assume
    that all the \(V_i^j\) are pairwise disjoint
    (because the sets \(S_i\) are pairwise disjoint and \(X\) is Hausdorff),
    and that each \(V_i^j\) is connected (because \(\Bf\) is decomposable)
    and of type \([U_i]\) (because the local type around \(s_i^j\) is \([U_i]\)
    by definition).
    But this means precisely that \(V\coloneqq\bigdisjun_{i,j}V_i^j \subset W\)
    lies in \(K_S\)
    (because \(\Bf\) is multiplicative, i.e.\ closed under disjoint unions)
    showing that this poset is non-empty.
  \item
    Given \(V,V' \in K_S\) we can simply apply the previous argument to
    the open set \(W\coloneqq V\cap V'\)
    to get a \(V''\in K_S\) with \(V''\subset V'\) and \(V''\subset V\).
    Recall that \(K_S\subset \Bf\) is not a full subposet,
    hence it remains to show that these inclusions lie in \(K_S\).
    By construction,
    the inclusions induce bijections
    \(\bigdisjun_i\pi_0(S_i)\xrightarrow{\cong}\pi_0(V'')\xrightarrow{\cong}\pi_0(V')\)
    and for each \(i,j\), the corresponding connected components
    \({V''}_i^j\) and \({V'}_i^j\) have a common point
    \(s_i^j\) in their deepest strata, which by \Cref{lem:char-iso-disk}
    implies that they are abstractly isomorphic.
    Therefore the inclusion \(V''\subset V'\) lies in \(\Jmaps[\Bf]\),
    hence in \(K_S\) which is a full subposet thereof.
    The same argument applies to \(V''\subset V\).
  \end{itemize}

  Since cofiltered posets are weakly contractible,
  we can apply the Seifert--Van Kampen theorem and obtain the desired equivalence
  \begin{equation}
    \label{eq:SvK-via-small-conf}
    \classspace{\Jmaps[\Bf]^{[U]}}
    =
    \colim_{V\in\Jmaps[\Bf]^{[U]}}\{\star\}
    \xleftarrow{\simeq}
    \colim_{V\in\Jmaps[\Bf]^{[U]}}{\smallconf{V}}
    \xrightarrow{\simeq}
    \prod_{i\in I}\Conf[k_i]{X_{[U_i]}},
  \end{equation}
  where on the left we use that
  each \(\smallconf{V}\) is homeomorphic to \(\prod_{i,j}(V_i^j)_{[U_i]}\),
  hence contractible
  (because each \(V_i^j\) is a disk of type \(U_i\),
  so that \((V_i^j)_{[U_i]}\) is just the lowest stratum of a cone,
  which is some \(\reals^{n_i}\)).
\end{proof}

Finally, we arrive at the proof of \Cref{thm:Disk-localized}.
The heavy category-theoretic lifting in this proof
is done by \Cref{lem:Berry-localization},
which is a general criterion for detecting localizations of \(\infty\)-categories.

\begin{proof}[Proof of \Cref{thm:Disk-localized}]
  \Cref{lem:char-iso-disk} tells us that \(\Jmaps[\Bf]\)
  is exactly the collection of those arrows in \(\Bf\) which are sent
  to equivalences in \(\ddiskX\).
  Hence we only have to
  invoke \Cref{lem:Berry-localization} for the functor
  \(D\coloneqq\Bf\to E\coloneqq \ddiskX\) with \(W\coloneqq\Jmaps[\Bf]\).
  We check its two assumptions:
  \begin{enumerate}
  \item
    \label{it:class-JX=ddiskX}
    We have to check that the induced map
    \(\classspace{\Jmaps[\Bf]}\to\left( \ddiskof{X} \right)^\simeq\)
    is an equivalence.
    Consider the commutative square
    \begin{equation}
      \begin{tikzcd}[cells={font=\everymath\expandafter{\the\everymath\displaystyle}}] % puts every cell in displaystyle instead of manually writing it everywhere
        &\Jmaps[\Bf]\ar[dr]\ar[dl,dashed]
        \\
        \colim_{V\in \Jmaps[\Bf]}
        \left(\overcat{\ddisk}{V}\right)_{\mathrm{term}}^\simeq
        \ar[r,hookrightarrow]
        \ar[d,"\simeq",dashed]
        &
        \colim_{V\in \Jmaps[\Bf]}
        \left(\overcat{\ddisk}{V} \right)^\simeq
        \ar[r]
        \ar[d,"\simeq"]
        &
        \left(\overcat{\ddisk}{X} \right)^\simeq
        \ar[d,"\simeq"]
        \\
        \colim_{V\in \Jmaps[\Bf]}\smallconf{V}
        \ar[r,hookrightarrow]
        \ar[r]
        \ar[rr,bend right=17,"\simeq"']
        &
        \colim_{V\in \Jmaps[\Bf]}\coprod_{[U]}\prod_{i\in I}\Conf[k_i]{V_{[U_i]}}
        \ar[r]
        &
        \coprod_{[U]}\prod_{i\in I}\Conf[k_i]{X_{[U_i]}} 
      \end{tikzcd}
    \end{equation}
    where
    \begin{itemize}

        \item
      \((\overcat{\ddisk}{V})_{\mathrm{term}}\)
      is the full subcategory of the slice consisting only of those
      objects where the structure map \(U\hookrightarrow V\) is an equivalence,
      i.e.\ the subgroupoid of terminal objects,
     \item
      the right diagonal map is given explicitly as
      \(V\mapsto (V\hookrightarrow X)\),
      hence tautologically factors as the dashed diagonal map
      \(V \mapsto (V,V\xrightarrow{\id}V)\),
    \item
      the lower left horizontal map is given for each \(V\) by the inclusion
      \(\smallconf{V}\hookrightarrow \prod_{i\in I}\Conf[k_i]{V}\)
      into the summand corresponding to the type \([U]=[V]\),
    \item
      the middle and right vertical equivalences are those of \Cref{lem:core-DDisk-is-Conf},
    \item
      the lower horizontal composite equivalence is that of Equation \eqref{eq:SvK-via-small-conf},
    \item
      the left vertical dashed map exists because
      if \(U\hookrightarrow V\) is an equivalence
      then every choice of points \(S_i\) as in \Cref{rem:explicit-core-ddisk-is-Conf}
      will yield a tuple of configurations that lies in~\(\smallconf{V}\),
    \item
      the dashed vertical map is an equivalence since both
      \((\overcat{\ddisk}{V})_{\mathrm{term}}\)
      and
      \(\smallconf{V}\)
      are contractible.
    \end{itemize}
    Finally we observe that since every \(\left(\overcat{\ddisk}{V}\right)_{\mathrm{term}}^\simeq\)
    is contractible,  we have that
    \[\classspace{\Jmaps[\Bf]} \simeq 
    \colim_{V\in \Jmaps[\Bf]}
        \ast 
        \simeq \colim_{V\in \Jmaps[\Bf]}
        \left(\overcat{\ddisk}{V}\right)_{\mathrm{term}}^\simeq\,.\]
    Hence, the dashed diagonal map presents the localization
    \(\Jmaps[\Bf]\to \classspace{\Jmaps[\Bf]}\) to the classifying space.
    We conclude that the desired map
    \(\classspace{\Jmaps[\Bf]}\to\left( \ddiskof{X} \right)^\simeq\)
    is presented as the upper horizontal composite,
    which we have shown is a composite of equivalences.
  \item
    For each \(U\in \Bf\) we have to check that the map
    \begin{equation}
      \classspace{\Jmaps[\Bf]\times_{\Bf}(\Bf\cap \cdown{U})}
      \to
      (\overcat{\left(\ddiskX\right)}{U\hookrightarrow X})^\simeq
    \end{equation}
    is an equivalence.
    But this map is canonically equivalent to the map
    \begin{equation}
      \classspace{\Jmaps[\Bf\cap\cdown{U}]}
      \to
      (\overcat{\ddisk}{U})^\simeq
    \end{equation} 
    which is an equivalence by applying 
    the previously proved part~\ref{it:class-JX=ddiskX}
    to the decomposable multiplicative basis \(\Bf\cap\cdown{U}\) of the
    smooth conical manifold \(U\).
    \qedhere
  \end{enumerate}
\end{proof}

We end this section with some easy consequences of the localization result.

\begin{prop}
  \label{obs:oKLE-from-disks-on-disks}
  \label{rem:oLKE-from-basis-on-disks}
  Let \(\Uf\) be a \posetofopens{} and \(\iota\colon\Bf\subseteq \Uf\) a decomposable multiplicative disk-basis.
  Let \(\Aa\) be a constructible multiplicative prefactorization algebra on \(\Bf\).
  Let \(U=U_1\disjun\cdots \disjun U_n\) be a disjoint union of contractible disks.
  Since \(\Bf\) is a decomposable multiplicative disk-basis of \(\Uf\)
  we can choose
  \(B=B_1\disjun\cdots\disjun B_n\in \Bf\),
  where each \(B_i\subseteq U_i\) is an abstractly isomorphic subdisk.
  Then for any such choice, the inclusion \(B\hookrightarrow U\) induces an equivalence
  \begin{equation}
    \Aa(B)=\iota_!\Aa(B)\xrightarrow{\simeq} \iota_!\Aa(U).
  \end{equation}
\end{prop}

\begin{proof}
  By assumption, \(\Bf\) is closed under disjoint unions subordinate to \(\Uf\),
  so that by \Cref{lem:oLKE-LKE} we may simply work with
  the underlying copresheaf and its ordinary left Kan extension.
  Observe that \(\Bf\cap \cdown{U}\) is a decomposable multiplicative disk-basis
  of \(U\),
  so that \Cref{thm:Disk-localized} yields that the functor
  \begin{equation}
    \Bf\cap\cdown{U} \to \ddiskof{U}
  \end{equation}
  is a localization, hence in particular colimit cofinal.
  Since \(\Aa\) is multiplicative and constructible
  it factors through this localization,
  yielding the right equivalence in the composite
  \begin{equation}
    \label{eq:composite-oLKE-tru-curly}
    \Aa(B)
    \to
    \iota_!\Aa(U)
    =
    \colim \Aa\restrict{\Bf\cap\cdown{U}}
    \xrightarrow{\simeq}
    \colim \Aa\restrict{\ddiskof U}.
  \end{equation}
  By \Cref{lem:char-iso-disk},
  the inclusion \(B\hookrightarrow U\) becomes an equivalence in \(\ddisk\),
  hence a terminal object in \(\ddiskof{U}\);
  thus the composite map
  \eqref{eq:composite-oLKE-tru-curly}
  is an equivalence, which concludes the proof.
\end{proof}

\begin{cor}
  \label{cor:oLKE-from-basis-constr}
  Let \(\Uf\) be a \posetofopens{} and
  \(\iota\colon \Bf\subset \Uf\) be a decomposable multiplicative disk-basis.
  Then the operadic left Kan extension \(\iota_!\Aa\)
  of a constructible multiplicative prefactorization algebra \(\Aa\) on \(\Bf\)
  is constructible on \(\Uf\).
\end{cor}

\begin{proof}
  Consider an inclusion \(U\hookrightarrow V\)
  of abstractly isomorphic conical disks.
  Then for any choice of \(B\hookrightarrow U\) as in
  \Cref{obs:oKLE-from-disks-on-disks},
  both the lower horizontal and the left diagonal map are equivalences
  \begin{equation}
    \begin{tikzcd}
      &\iota_!\Aa(U)\ar[dr]
      \\
      \Aa(B)\ar[rr,"\simeq"]\ar[ru,"\simeq"]&&\iota_!\Aa(V)
    \end{tikzcd}
  \end{equation}
  hence also the right diagonal one, as desired.
\end{proof}

\begin{cor}
  \label{cor:cstr-B=cstr-disk}
  Let \(\iota\colon\Bf\subseteq\diskX\) be a decomposable multiplicative disk-basis.
  Then restriction along \(\iota\) yields the equivalence
  \begin{equation}
    {\AlgMCstr{\diskX}}
    \xrightarrow{\simeq}
    {\AlgMCstr{\Bf}}
  \end{equation}
  between multiplicative constructible prefactorization algebras
  on \(\Bf\) and \(\diskX\), respectively.
\end{cor}

\begin{proof}
  Starting from the usual operadic left Kan extension adjunction
  \begin{equation}
    \iota_!\colon
    \Alg{\Bf}
    \leftrightarrows
    \Alg{\diskX}
    : \iota^*
  \end{equation}
  with fully faithful left adjoint \(\iota_!\),
  we only have to prove the following two things:
  \begin{itemize}
  \item
    Every multiplicative constructible prefactorization algebra \(\Aa\) on \(\diskX\)
    is an operadic left Kan extension of its restriction to \(\Bf\).
    Indeed, for each \(U\in \diskX\) we have a commutative diagram
    with \(B\) as in \Cref{rem:oLKE-from-basis-on-disks}:
    \begin{equation}
      \begin{tikzcd}
        & \iota_!\iota^*\Aa (U)
        \ar[dr]
        \\
        \Aa(B)\ar[ur,"\simeq"]
        &&
        \Aa(U)
        \\
        \Aa(B_1)\otimes\cdots\otimes\Aa(B_n)
        \ar[u,"\simeq"]
        \ar[rr,"\simeq"]
        &&
        \ar[u,"\simeq"]
        \Aa(U_1)\otimes\cdots\otimes\Aa(U_n)
      \end{tikzcd}
    \end{equation}
    where the upper left diagonal map is an equivalence by 
    \Cref{rem:oLKE-from-basis-on-disks},
    the lower horizontal map is an equivalence by constructibility,
    and the lower vertical maps are equivalences by multiplicativity.
  \item
    If \(\Aa\) is a constructible multiplicative prefactorization algebra on \(\Bf\),
    then \(\iota_!\Aa\) is also multiplicative and constructible.
    Constructibility is precisely \Cref{cor:oLKE-from-basis-constr} for \(\Uf\coloneqq \diskX\).
    To prove multiplicativity, let \(U,U'\in \diskX\) be disjoint.
    We can choose corresponding \(B,B'\) as in
    \Cref{rem:oLKE-from-basis-on-disks}
    and obtain a commutative square
    \begin{equation}
      \cdsquare
      {\Aa(B)\otimes\Aa(B')}
      {\iota_!\Aa(U)\otimes\iota_!\Aa(U')}
      {\Aa(B\disjun B')}
      {\iota_!\Aa(U\disjun U')}
      {\simeq}
      {\simeq}
      {}
      {\simeq}
    \end{equation}
    where the two horizontal maps are equivalences by
    \Cref{rem:oLKE-from-basis-on-disks}
    (note that \(B\disjun B'\in\Bf\) is a valid choice for \(U\disjun U'\in\diskX\))
    and the left vertical map is an equivalence by multiplicativity of \(\Aa\).
    Thus the right vertical map is an equivalence as desired.
    \qedhere
  \end{itemize}
\end{proof}

We also record the following variant of \Cref{thm:Disk-localized}
regarding the \emph{operads} of disks, as opposed to the mere categories.

Recall that the \infy-category \(\SSnglr\) has a symmetric monoidal structure
\(\amalg\),
given by disjoint union of smooth conical manifolds.
Since the unit \(\emptyset\) of this symmetric monoidal structure is an initial object,
we then have an induced operad structure \(\DDiskofX\)
on each slice \(\ddiskofX\) which is pre-cocartesian
in the sense of \Cref{def:pre-cocartesian};
see \Cref{ex:slice-operad}.
Note that despite the chosen symbol,
\(\amalg\) is \emph{not} the coproduct in \(\SSnglr\),
and the \infy-operad \(\DDiskofX\) is not monoidal or even partially monoidal.

\begin{prop}
  \label{prop:Disk-operad-localization}
  Let \(X\) be a smooth conical manifold and \(\Bf\subseteq\diskX\)
  a decomposable multiplicative disk-basis of X.
  Then the operad map
  \begin{equation}
    \operad{\Bf}\to \DDiskofX
  \end{equation}
  is a localization at the arrows in \(\Jmaps[\Bf]\)
  so that the restriction map
  \begin{equation}
    \Alg{\ddiskofX}\longrightarrow \Alg \Bf
  \end{equation}
  is fully faithful with essential image consisting of those algebras
  which send all operations in \(\Jmaps[\Bf]\) to equivalences.
\end{prop}

\begin{proof}
  Since \(\Bf\) is multiplicative,
  the operad \(\operad{\Bf}\) inherits its partial monoidal structure from
  \(\Open{X}\),
  hence is pre-cocartesian;
  the pre-cocartesian arrows are just the \pcocartesian{} ones
  (which also agree with the cocartesian ones due to multiplicative).
  The operad map
  \(\operad{\Bf}\to \DDiskofX\)
  is pre-cocartesian since the composite
  \(\operad{\Bf}\to \DDiskofX\to \DDisk\) is a map of (partial)
  symmetric monoidal \infy-categories, i.e., sends \(\disjun\) to \(\amalg\).
  
  Hence we can apply \Cref{lem:operadic-Berry-localization}
  with \(W=\Jmaps[\Bf]\)
  by checking its conditions:

  \begin{itemize}
  \item
    [\ref{it:op-localization-underlying}]
    We already verified the localization conditions on the underlying categories
    in \Cref{thm:Disk-localized}.
  \item [\ref{it:op-localization-left-fib}]
    Given an inclusion \(D_1\disjun\cdots\disjun D_n=O \hookrightarrow O'\)
    in \(\Jmaps[\Bf]\), we obtain the unique factorization
    \eqref{eq:op-localization-W-factorize}
    by letting each
    \begin{equation}
      D_i\hookrightarrow D'_i\coloneqq\bigdisjun \{V\in \pi_0(O')\mid V\cap D_i \neq \emptyset\}\hookrightarrow O'
    \end{equation}
    be the union of those connected components of \(O'\)
    which are hit by \(D_i\).
  \item [\ref{it:op-localization-fibers}]
    Both sides of \eqref{eq:op-localization-fibers}
    are just the set of equivalence relations on the finite set
    \(\pi_0(O)=\pi_0(f(O))\).
  \end{itemize}
  Note that we made repeated use of the fact that \(\Bf\) is decomposable
  so that every subset \(I\subseteq\pi_0(O)\)
  does give a well defined object \(\bigdisjun I \in \Bf\cap\cdown{O}\).
\end{proof}

An algebra for the \infy-operad \(\DDiskofX\)
is called \emph{multiplicative}
if it sends pre-cocartesian arrows to cocartesian ones.\footnote{%
  The same definition would still make sense for any pre-cocartesian \infy-operad.
  We avoid this generality, because it would cause an unfortunate clash
  of terminology in the case where the operad arises from
  a \posetofopens{} \(\PS\) which has all finite disjoint joins but
  is not closed under disjoint unions
  (e.g., the poset of open intervals in \(\reals\),
  which even has all finite joints, not just disjoint ones).
  The disjoint join would make such an operad into a partial monoidal category,
  hence into a pre-cocartesian operad with the pre-cocartesian morphisms
  being exactly the cocartesian ones
  (see \Cref{ex:partially-monoidal-pre-cocartesian}).
  But in \Cref{defn:multiplicative-prefact-alg}
  we purposefully only require a multiplicative algebra to send \pcocartesian{}
  morphisms to cocartesian ones.

  Recall \Cref{def:pcocartesian} for the distinction,
  but observe that this distinction is immaterial in all \posetsofopens{}
  that actually appear in this section,
  because they are all closed under subordinate disjoint unions.
}

As with pre-factorization algebras, we use the superscript ``\(\mult\)''
to denote categories of multiplicative algebras.

\begin{cor}
  \label{cor:AlgM-localization}
  Let \(X\) be a smooth conical manifold and \(\Bf\subseteq\diskX\)
  a decomposable multiplicative disk-basis of X.
  Restriction yields an equivalence
  \begin{equation}
    \AlgM{\ddiskofX} \xrightarrow{\simeq}
    \AlgMCstr{\Bf}.
  \end{equation}
\end{cor}

\begin{proof}
  The multiplicativity condition on \(\Bf\) corresponds to that on \(\ddiskof X\)
  since every pre-cocartesian arrow in \(\DDiskofX\)
  is equivalent to one that comes from \(\operad{\Bf}\).
  For multiplicative prefactorization algebras on \(\Bf\)
  the constructibility condition amounts precisely
  to inverting the operations in \(\Jmaps[\Bf]\).
  Hence the result follows directly from 
  \Cref{prop:Disk-operad-localization}.
\end{proof}

\subsection{Constructible implies Weiss on disks}
\label{sec:constr->Weiss}

It has long been folklore that for constructible factorization algebras
extended from disks the Weiss condition is redundant.
A version of this statement first appeared in \cite{AFprimer}
in the unstratified setting.
In contrast, we work in the stratified setting,
and also allow for the additional freedom
of restricting to a suitable collection of disks.

\begin{thm}
  \label{thm:constructible=>Weiss}
  Let \(X\) be a smooth conical manifold and let
  \(\Bf\subseteq\diskX\) be a decomposable multiplicative disk-basis of \(X\).
  Every left Kan extension of a constructible copresheaf on \(\Bf\)
  to any \posetofopens{} \(\Uf\supseteq \Bf\) is automatically a Weiss hypercosheaf.
\end{thm}

We will prove \Cref{thm:constructible=>Weiss} below.
First, we deduce some corollaries.

\begin{cor}
  \label{cor:constructible=>Weiss}
  Let \(X\) be a smooth conical manifold,
  \(\Bf\subseteq\diskX\) a decomposable multiplicative disk-basis of \(X\)
  and \(\Uf\supseteq\Bf\) a presieve.
  Then for every constructible copresheaf on \(\Bf\),
  the left Kan extension to \(\Uf\) is a Weiss cosheaf.
\end{cor}

\begin{proof}
  This is a direct consequence of \Cref{thm:constructible=>Weiss}
  because every Weiss hypercosheaf on a presieve is a Weiss cosheaf
  (see \Cref{rem:Weiss-cosheaf-Weiss-presieves}).
\end{proof}

\begin{cor}
  \label{cor:extending-from-disk-basis}
  Let \(X\) be a smooth conical manifold,
  \(\Uf\subseteq \open{X}\) a presieve
  and \(\Bf\subseteq\Uf\) a decomposable factorizing disk-basis of \(\Uf\).
  Then restriction yields an equivalence
  \begin{equation}
    \FactCstr{\Uf}\xrightarrow{\simeq}\FactCstr{\Bf}=\AlgMCstr{\Bf}
  \end{equation}
  of \(\infty\)-categories.
\end{cor}

\begin{proof}%[Proof of \Cref{cor:extending-from-disk-basis}] 
  We start with the equivalence
  \begin{equation}
    \iota^\ast\colon\Fact{\Uf}\xrightarrow{\simeq}\Fact{\Bf}
  \end{equation}
  of \Cref{lemma:FactBasisExtends}
  which restricts to an equivalence
  \begin{equation}
    \iota^\ast\colon\FactCstr{\Uf}\xrightarrow{\simeq}\FactCstr{\Bf}
  \end{equation}
  by \Cref{cor:oLKE-from-basis-constr}.
  Finally, we observe that the inclusion
  \(\AlgMCstr{\Bf}\subseteq \FactCstr{\Bf}\)
  is an equality by \Cref{cor:constructible=>Weiss}.
\end{proof}

\begin{cor}
  \label{cor:extending-from-disks}
  Let \(X\) be a smooth conical manifold with enough good disks.
  Then restriction yields an equivalence
  \begin{equation}
    \FactCstr{X}\xrightarrow{\simeq}\AlgMCstr{\diskX}
  \end{equation}
  of \(\infty\)-categories.
\end{cor}

\begin{proof}
  Choose a decomposable factorizing disk-basis \(\Bf\) of \(X\).
  Then the claim follows by 2-out-of-3 for the composition
  \begin{equation}
    \FactCstr{X}\to\AlgMCstr{\diskX}\xrightarrow{\simeq}\AlgMCstr{\Bf}
  \end{equation}
  where the composite is the equivalence of \Cref{cor:extending-from-disk-basis}
  and the right restriction functor is the equivalence of
  \Cref{cor:cstr-B=cstr-disk}.
\end{proof}

\begin{rem}
  \label{rem:AF_analog_extension}
  Ayala--Francis state the unstratified analog of
  \Cref{cor:extending-from-disks}
  namely \cite[Proposition~2.22]{AFprimer}
  without using  the existence of enough good disks.
  Their proof seems to go along the same lines as our proof of
  \Cref{cor:extending-from-disk-basis}
  (and in fact we learned some of the key ideas from them)
  but directly applied to the \posetofopens{} \(\diskX\) of \emph{all} disks.
  However, since this \posetofopens{} is \emph{not} closed under intersections,
  we were not able to reproduce their proof.
  For this reason we have to go through the auxiliary factorizing disk-basis \(\Bf\)
  whose existence we have to postulate.

  Still in the unstratified setting,
  essentially the same result was also obtained in Carmona's thesis~\cite{Carmona-thesis},
  who like us requires the existence of enough good disks.

  At this point it is still unclear to us
  whether this statement can be proved without the existence of enough good disks;
  if so, it might require some more high-powered machinery.
\end{rem}

The final corollary is the full comparison of factorization algebras
in our sense
with those defined as multiplicative \infy-categorical disk algebras.

\begin{cor}
  \label{cor:FactX-AlgDiskX}
  Let \(X\) be a smooth conical manifold with enough good disks.
  We have a canonical equivalence
  \begin{equation}
    \FactCstr{X}\simeq \AlgM{\ddiskofX}\simeq \Alg{\cddiskof{X}}
  \end{equation}
\end{cor}
\begin{proof}
  For the first equivalence
  combine \Cref{cor:extending-from-disks}
  with \Cref{cor:AlgM-localization}.
  The second one follows analogously using 
  \Cref{cor:add-disjun-constr}.
\end{proof}

\begin{ex}[Associative algebras and modules]
  \label{ex:ass-and-rmod}
  Using the identification of
  \Cref{cor:FactX-AlgDiskX},
  we can now explain
  \Cref{ex:algebra}
  and
  \Cref{ex:bimod}
  systematically:
  constructible factorization algebras on the open interval \((0,\infty)\)
  are given precisely by associative algebras \(\Aa\);
  constructible factorization algebras on the half-open interval \([0,\infty)\)
  additionally encode the data of a pointed \(\Aa\)-module
  (see also \cite[Section~2.6]{AFT-fh-stratified}).

  Explicitly, let \(\RModp\) be the discrete operad encoding a pointed right action:
  \begin{itemize}
  \item
    It has two colors \(a\) and \(m\)
    (corresponding to the algebra and the module).
  \item
    It has operations \(\sigma\colon (a,\dots, a)\to a\)
    and \(\tau \colon (m,a,\dots,a)\to m\)
    (with \(n\geq 0\) many \(a\)'s)
    where \(\sigma\) and \(\tau\) are linear orders of the multiset
    \(\{a,\dots, a\}\) or \(\{m,a,\dots,a\}\),
    respectively, with \(m\) always being required to be the first element.
  \item
    Unlike the usual operad \(\RMod\) encoding an unpointed right action,
    there is an additional \(0\)-ary operation \(()\to m\)
    that yields a pointing of the module.
  \item
    Composition is evident.
  \end{itemize}

  Moreover, there is an explicit operad map
  \begin{equation}
    \label{eq:discretization-map}
    \cDDiskof{[0,\infty)}\to \RModp
  \end{equation}
  that sends the disks of type \([0,\infty)\) to \(m\)
  and the disks of type \((0,\infty)\) to \(a\).
  If \(I_1,\dots,I_n, I\) are disks with embeddings \(\phi_1,\dots,\phi_n, \phi\)
  into \([0,\infty)\),
  then the space of operations
  \begin{equation}
    ((I_1,\phi_1),\dots,(I_n,\phi_n))\to (I,\phi)
  \end{equation}
  is the space of
  \begin{enumerate}[label=($\dagger$), ref=($\dagger$)]
  \item
    \label{it:space-of-interval-embeddings}
    embeddings \(\psi\colon I_1\amalg\dots\amalg I_n\hookrightarrow I\)
    with isotopies \(h_i\colon \psi\restrict{I_i}\simeq \phi_i\).
  \end{enumerate}
  The map \eqref{eq:discretization-map}
  sends such a pair \((\psi,(h_i))\) to the linear order
  of the set \(v_1,\dots,v_n\) induced by \(\phi\),
  thus yielding the operation \((v_1,\dots,v_n)\to v\);
  here the symbol \(v\) and each \(v_i\) is either ``\(m\)'' or ``\(a\)''
  depending on the type of the disk \(I\) and \(I_i\).

  Finally, we observe that each space
  \ref{it:space-of-interval-embeddings}
  has contractible connected components  and that the operad map \eqref{eq:discretization-map}
  is actually an equivalence\footnote{
    This boils down to the fact that the space of smooth embeddings
    \(\reals^{\amalg n}\hookrightarrow\reals\) is discrete (up to equivalence)
    with \(n! \cdot 2^n\) connected components,
    corresponding to the linear orders of the \(n\) copies
    and a choice of orientation for each.
    The identifications \(h_i\) fix the orientations,
    leaving precisely the choice of linear order.
  }.
  In particular the equivalence of
  \Cref{cor:FactX-AlgDiskX}
  yields the desired identification between constructible factorization algebras
  on \([0,\infty)\) and \(\RModp\)-algebras, i.e., pointed right modules.
\end{ex}

\begin{rem}
  As explained in \Cref{ex:ass-and-rmod} the identification
  with associative algebras and pointed modules thereover is very straightforward
  when working with the topologized version of factorization algebras
  (i.e., the \infy-operad \(\cDDiskof{[0,\infty)}\)).
  For this reason, this identification is often accepted as a basic fact
  of the theory of constructible factorization algebras
  with the implicit understanding that ``factorization algebra''
  is to be read as \(\cddiskof{[0,\infty)}\)-algebra.

  From our definition of (constructible) factorization algebras
  (which, for example, is also that of Costello-Gwilliam and Ginot),
  this is not at all obvious
  and relies on the highly non-trivial localization result
  that underlies the comparison \Cref{cor:FactX-AlgDiskX}.
  
  For associative algebras without a module
  (and more generally \Ealgebra{n}s for each \(n\)),
  this identification goes back to Lurie
  (modulo some results about extensions from a factorizing basis,
  since Lurie only deals with disks and does not use the Weiss-condition\footnote{
    More precisely:
    First we view associative algebras as \Ealgebra{1}s,
    as in \cite[Example~5.1.0.7]{LurHA}.
    Then \cite[Theorem~5.4.5.9]{LurHA} specializes
    under \cite[Example~5.4.5.3]{LurHA} to
    an equivalence between \Ealgebra{1}s
    and constructible multiplicative prefactorization algebras
    on \(\disk{\reals}\).
    These in turn are equivalent to locally constant factorization algebras
    on \(\reals\cong(0,\infty)\)
    by \Cref{cor:extending-from-disks}.
  }).
  While an extension of his techniques to the stratified setting
  (which means taking into account an additional module)
  seems relatively straightforward,
  it has not appeared in print as far as we are aware.
\end{rem}

After having recorded these immediate corollaries,
we start working towards the proof of
\Cref{thm:constructible=>Weiss} for which we need some preliminaries.
One key ingredient is the well-known fact that the colimit over a hypercover
of a topological space computes the homotopy type of that space.
This is another instance of the Seifert--Van Kampen theorem.

\begin{lemma}[\cite{DuggerIsaksen}, Proposition~4.6~(c)]
  \label{lem:DI-hypercover}
  Let \(X\) be a topological space.
  The fundamental \infy-groupoid functor
  \begin{equation}
    \openX \to \Gpdinfty
  \end{equation}
  which sends an open set to its homotopy type
  is a hypercosheaf.
  Explicitly this means that for each hypercover\footnote{
    What we call ``hypercover'' is called ``complete cover''
    by Dugger--Isaksen.
    See also \Cref{rem:hypercover-nonstandard}.
  } \(\Wf\hyprefines X\)
  we have the equivalence
  \begin{equation}
    \colim_{W\in\Wf}W \xrightarrow{\simeq} X
  \end{equation}
  in the \(\infty\)-category of spaces. 
\end{lemma}

We now come to the proof of \Cref{thm:constructible=>Weiss}.
We make use of a general category theoretic fact about cofinality,
which for completeness we record as \Cref{lem:cofinality-right-fibrations} in the appendix.

\begin{proof}[Proof of \Cref{thm:constructible=>Weiss}]
  Let \(\Aa\colon \Uf\to \targetcat\) be a copresheaf left Kan extended from \(\Bf\)
  and assume that \(\Aa\restrict{\Bf}\) is constructible.
  For each \(U\in \Uf\) we can apply
  \Cref{thm:Disk-localized}
  to the decomposable multiplicative basis \(\Bf\cap\cdown{U}\) of \(U\)
  to get a factorization
  \begin{equation}
    \label{eq:in-proof-B-U-localization}
    \Aa\restrict{\Bf\cap\cdown{U}}\colon
    \Bf\cap\cdown{U}\to \ddiskof{U}\to \targetcat
  \end{equation}
  where the first functor is a localization.
  
  Now, let \(\Wf\hyprefines{M}\) be a Weiss hypercover of an open set $M$ in \(\Uf\).
  We aim to show that it is \(\Aa\)-local.
  We apply \Cref{lem:cofinality-right-fibrations}
  (specifically \ref{it:colimit-core}$\implies$\ref{it:colimit-right-fib})
  to
  \begin{itemize}
  \item
    the \(\infty\)-category \(\Da\coloneqq\ddisk\subseteq\SSnglr\)
    of conical disks
  \item
    and the diagram
    \begin{equation}
      F\colon \Wf^\triangleright\xrightarrow{\infty\mapsto M}\Wf\cup\{M\}\to \RFib(\ddisk),
      \quad
      W \mapsto (\ddiskof{W}\to \ddisk)
    \end{equation}
    of right fibrations.
  \end{itemize}
  To establish the condition~\ref{it:colimit-core}
  of \Cref{lem:cofinality-right-fibrations}
  we observe that by \Cref{lem:core-DDisk-is-Conf}
  the map \eqref{eq:colimit-core}
  is identified with the map of \infy-groupoids
  \begin{equation}
    \label{eq:colim-conf-hyper}
    \coprod_{[U]}\colim_{W\in\Wf}\prod_{i\in I}\Conf[k_i]{W_{[U_i]}}
    \to
    \coprod_{[U]}
    \prod_{i\in I}\Conf[k_i]{M_{[U_i]}}
  \end{equation}
  which is induced by the hyperprecovers
  \begin{equation}
    \label{eq:hypcover-conf-space}
    \left\{
      \prod_{i\in I}\Conf[k_i]{W_{[U_i]}}
      \mid W\in\Wf
    \right\}
    \hyprefines
    \prod_{i\in I}\Conf[k_i]{M_{[U_i]}}.
  \end{equation}
  (for each isomorphism class \(U=\coprod_{i}k_i\times U_i\) as in \Cref{not:U-types}).
  The fact that \(\Wf\hyprefines M\) is a Weiss hypercover
  directly implies that 
  \eqref{eq:hypcover-conf-space}
  is a hypercover,
  so that \Cref{lem:DI-hypercover} yields the desired
  equivalence~\eqref{eq:colim-conf-hyper} of \(\infty\)-groupoids.

  Having established condition~\ref{it:colimit-core}
  of \Cref{lem:cofinality-right-fibrations},
  we are entitled to its conclusion~\ref{it:colimit-right-fib}:
  the right vertical functor in the commutative square
  of \(\infty\)-categories
  \begin{equation}
    \cdsquareNA
    {\colim\limits_{W\in \Wf}\Bf\cap\cdown{W}}
    {\colim\limits_{W\in \Wf}\ddiskof{W}}
    {\Bf\cap\cdown{M}}
    {\ddiskof{M}}
  \end{equation}
  is colimit cofinal.
  By setting \(U\coloneqq M\) and \(U\coloneqq W\)
  in \eqref{eq:in-proof-B-U-localization},
  and passing to colimits of \(\infty\)-categories
  we see that the horizontal functors are localizations
  through which \(\Aa\) factors (uniquely).
  After passing to colimits in \(\targetcat\), this implies that
  the right vertical map in the following commutative
  diagram is an equivalence:
  \begin{equation}
    \begin{tikzcd}
     \colim\limits_{W\in \Wf}\Aa(W)
     \ar[d]
     &\colim\limits_{W\in\Wf}\colim\limits \Aa\restrict{\Bf\cap\cdown{W}}
     \ar[l]
     \ar[r,equal]
     &\colim\limits\Aa\restrict{\colim\limits_{W\in \Wf}\Bf\cap\cdown{W}}
     \ar[d,"\simeq"]
     \\
     \Aa(M)
     &
     &\colim\limits\Aa\restrict{\Bf\cap\cdown{M}}
     \ar[ll]
    \end{tikzcd}
  \end{equation}
  Moreover, the two horizontal arrows are equivalences
  because \(\Aa\) is a left Kan extension of its restriction to \(\Bf\).
  We conclude that the left vertical map is an equivalence
  which was the goal.
\end{proof}

\newpage

\section{Main theorems and applications}
\label{sec:applications}

In this section we give some applications of our gluing techniques.

\subsection{Gluing of constructible factorization algebras}

%%%%%%
% local macro:
\newcommand{\basis}{\Df}
\newcommand{\poset}{\Pf}

Our first application is the main theorem of this article, namely that constructible factorization algebras on an open cover glue. In other words, they form a sheaf of $\infty$-categories. In the locally constant case, this was proven in \cite[Theorem 1.3]{Matsuoka}.

Without further ado,

\begin{thm}
  \label{thm:main-gluing}
  Let \(\targetcatOT\) be a \(\otimes\)-presentable symmetric monoidal \(\infty\)-category.
  Let \(X\) be a smooth conical manifold with enough good disks.
  Let \(\Ua=\{U_i\}_{i<\beta}\) be an open cover of \(X\).
  For any finite subset \(I \subset \beta\),
  denote \(U_I \coloneqq\bigcap_{i\in I} U_i\).

  Then we have an equivalence of \infy-categories
  \begin{equation}
    \label{eq:GluingEquivalence}
    \FactCstr[\targetcat]{X} \xrightarrow{\simeq}\lim_{\beta \supset I \supsetneq \emptyset} \FactCstr[\targetcat]{U_I}.
  \end{equation}
  In other words:
  the assignment $U \mapsto \FactCstr[\targetcat]{U}$ is a sheaf of $\infty$-categories.
\end{thm}

\begin{proof}
  We factor the map \eqref{eq:GluingEquivalence}
  as the following composite of equivalences.
  \begin{equation}
    \label{eq:composite-for-main-thm}
    \begin{tikzcd}[column sep=large]
      \FactCstr{X}
      \arrow[r,"\simeq","{\text{Prop~}\ref{prop:extension-from-sieve-cover}}"']
      &
      \FactCstr{\cdisj{\cdown{\Ua}}} \arrow[r, "\simeq", "{\text{Prop~}\ref{prop:final_step_in_gluing}}"']
      & \FactCstr{\cdown{\Ua}} \arrow[r, "\simeq", "{\text{Cor~}\ref{cor:GlueFactCstrOnCover}}"']
      & \displaystyle \lim_{\beta \supset I \supsetneq \emptyset} \FactCstr{U_I}.
    \end{tikzcd}
  \end{equation}
  The first functor is an equivalence by \Cref{prop:extension-from-sieve-cover} below;
  the last by \Cref{cor:GlueFactCstrOnCover}.
  The middle functor is an equivalence by applying
  \Cref{prop:final_step_in_gluing} below
  to $\poset = \cdown{\Ua}$, which (being a sieve) is a decomposable presieve;
  it has enough good disks because \(X\) has enough good disks.
\end{proof}

\begin{prop}
  \label{prop:extension-from-sieve-cover}
  Assume that \(X\) has enough disks and let \(\Ua\) be a cover of \(X\).
  Then the restriction functor
  \begin{equation}
    \label{eq:equivalence-X-Udown-cstr}
    \iota^\ast\colon \FactCstr{X} \xrightarrow{\simeq}  \FactCstr{\cdisj{\cdown{\Ua}}}
  \end{equation}
  is an equivalence.
\end{prop}

\begin{proof}
  Observe that \(\cdisj{\cdown{\Ua}}\)
  is a decomposable factorizing basis
  of the presieve \(\cdown{X}\),
  so that \Cref{lemma:FactBasisExtends} yields the equivalence
  \begin{equation}
    \label{eq:equivalence-X-Udown}
    \iota^\ast\colon \Fact{X} \xrightarrow{\simeq} \Fact{\cdisj{\cdown{\Ua}}}.
  \end{equation}
  Note that a factorization algebra \(\Aa\) on \(X\)
  is constructible on \(\cdisj{\cdown{\Ua}}\)
  if and only if it is constructible on each \(U\in \Ua\)
  (because each instance of the constructibility condition involves only contractible disks).
  Since constructibility is local (see \Cref{thm:constructle-local})
  this happens if and only if \(\Aa\) is constructible on \(X\).
  Thus the equivalence \eqref{eq:equivalence-X-Udown}
  restricts to the desired equivalence
  \eqref{eq:equivalence-X-Udown-cstr}.
\end{proof}

\begin{prop}\label{prop:final_step_in_gluing}
  Let $\poset \subseteq \openX$ be a decomposable presieve.
  Then the restriction 
  \begin{equation}
    \label{eq:equivalence-P-Pdisj}
    \FactCstr{\disjunComp{\poset}} \xrightarrow{\simeq}	\FactCstr{{\poset}} 
  \end{equation}
  is an equivalence of \infy-categories.
\end{prop}

\begin{proof}
  We start with the equivalence
  \begin{equation}
		\AlgMCstr{\cdisj{\poset}} \xrightarrow{\simeq}\AlgMCstr{\poset}
  \end{equation}
  of \Cref{cor:add-disjun-constr}
  and claim that it restricts to the desired equivalence
  \eqref{eq:equivalence-P-Pdisj},
  i.e., that (the copresheaf underlying)
  any given multiplicative constructible prefactorization algebra \(\Aa\)
  on \(\cdisj{\poset}\)
  is a Weiss cosheaf
  if and only if its restriction to \(\poset\) is a Weiss cosheaf.
 
  ``Only if'' is clear because \(\poset\) is a presieve
  so the value of Weiss covers in \(\poset\)
  do not change when evaluated in \(\cdisj{\poset}\).
 
  For the converse, assume that \(\Aa \in \AlgMCstr{\cdisj{\poset}}\) is a Weiss cosheaf when restricted to \(\poset\).
  Choose a factorizing disk-basis \(\basis\) of \(\poset\),
  witnessing that it has enough disks;
  without loss of generality we may assume that \(\basis\) is decomposable.
  We note two facts:
  \begin{itemize}
  \item
    Since \(\Aa\restrict\poset\) is a Weiss cosheaf
    it is a left Kan extension from the factorizing basis \(\basis\)
    by \Cref{prop:UniqueExtFromBasis}~\ref{it:Weiss-is-LKE-from-Weiss}.
    Moreover,
    \(\basis\)
    is closed under disjoint unions subordinate to \(\poset\),
    so that
    \(\Aa\restrict{\operad{\poset}}\)
    is also an \emph{operadic} left Kan extension from \(\operad{\basis}\)
    by \Cref{lem:oLKE-LKE}.
  \item
    Applying \Cref{prop:MultAlgToDisjCompletion} to
    \(\PS\coloneqq\basis\)
    and to
    \(\PS\coloneqq\poset\)
    (which are both decomposable),
    we get that the multiplicative prefactorization algebras
    \(\Aa\restrict{\DisjunComp{\basis}}\)
    and
    \(\Aa\restrict{\DisjunComp{\poset}}\)
    are necessarily operadic left Kan extensions
    from
    \(\Aa\restrict{\operad{\basis}}\)
    and
    \(\Aa\restrict{\operad{\poset}}\),
    respectively.
  \end{itemize}
  By the transitivity and 2-out-of-3 property for operadic left Kan extensions
  for the commutative square
  \begin{equation}
    \cdsquareOpt
    {\operad{\basis}}
    {\operad{\poset}}
    {\DisjunComp{\basis}}
    {\DisjunComp{\poset}}
    {hookrightarrow}
    {hookrightarrow}
    {hookrightarrow}
    {hookrightarrow}
  \end{equation}
  of operads,
  we conclude that \(\Aa\) (which, recall, is defined on \(\operad{\cdisj{\poset}}\))
  is an operadic left Kan extension from \(\operad{\cdisj{\basis}}\).
  Again invoking \Cref{lem:oLKE-LKE}
  this means that the underlying copresheaf \(\Aa\restrict{\cdisj{\poset}}\)
  is a left Kan extension from \(\Aa\restrict{\cdisj{\basis}}\),
  which is constructible (since \(\Aa\) was constructible).
  Since \(\cdisj{\basis}\) is
  a decomposable multiplicative (even factorizing) disk-basis of \(\cdisj{\poset}\),
  it follows from \Cref{cor:constructible=>Weiss}
  that \(\Aa\) is a Weiss cosheaf on \(\cdisj{\poset}\),
  which is what we were required to prove.
\end{proof}

\begin{ex}\label{ex:GlueFactCstr}
Let \(U=(0, 3)\) and \(V=(2, 5)\) be intervals in $\reals$, and \(X=U\cup V\).
We choose to stratify both intervals with a single point, say, choose $\{1\} \subset U$ and $\{4\} \subset V$.
Consider constructible factorization algebras \(\Fa_U\) on $U$ and \(\Fa_V\) on $V$, together with an isomorphism $t(\Fa_U)\coloneqq (\Fa_U)|_{(2,3)} \cong (\Fa_V)|_{(2,3)} \eqqcolon s(\Fa_V)$.
Then they can be glued (up to equivalence) to a factorization algebra on all of $X$.
\begin{center}
\begin{tikzpicture}[scale=2]
	\begin{scope}[xshift=0cm]
	\draw[draw =blue, ultra thick] (0,0) -- (0.3, 0);
	\draw[draw =red, ultra thick] (0.3, 0) -- (1,0);
	\draw (0.3, 0) node[dot] {};
	\path (0,0) -- node[anchor=north] {\scriptsize $\Fa_U$} (1,0);
	\draw (1.2, 0) node (A) {};
	\end{scope}

	\begin{scope}[xshift=2.5cm]
	\draw[draw =redbg, ultra thick] (0, 0) -- (0.7,0);
	\draw[draw =greenbg, ultra thick] (0.7, 0) -- (1,0);
	\draw (0.7, 0) node[dot] {};
	\path (0,0) -- node[anchor=north] {\scriptsize $\Fa_V$} (1,0);
	\draw (-0.2, 0) node (B) {};
	\end{scope}

	\draw[<->] (A) --node[anchor=south] {\scriptsize $t(\Fa_U)\cong s(\Fa_V)$} (B);

	\draw (4, 0) node {$\rightsquigarrow$};

	\begin{scope}[xshift=4.5cm]
	\draw[draw =bluebg, ultra thick] (0,0) -- (0.3, 0);
	\draw[draw =redbg, ultra thick] (0.3, 0) -- (1,0);
	\draw[draw =greenbg, ultra thick] (1.3, 0) -- (1,0);

	\draw (0.3, 0) node[dot] {};
	\draw (1, 0) node[dot] {};

	\draw [thick, decoration={brace, mirror, raise=0.2cm}, decorate] (0, 0) --  node[anchor=north, yshift= -0.25cm] {\scriptsize $\Fa_U$} (1, 0);

	\draw [thick, decoration={brace, raise=0.2cm}, decorate] (0.3, 0) --  node[anchor=south, yshift= 0.25cm] {\scriptsize $\Fa_V$} (1.3, 0);
	\end{scope}
\end{tikzpicture}
\end{center}
The main example of this situation to keep in mind is when \(\Fa_U\) and \(\Fa_V\) arise,
 similarly to \Cref{ex:bimod}, \Cref{ex:bimod_is_fact}, and \Cref{extend_algebra_from_basis},
  from an $(\algebra_1,\algebra_2)$-bimodule $\bimod_1$ and an $(\algebra_2,\algebra_3)$-bimodule $\bimod_2$,
   respectively\footnote{Together with pointings which we do not write out explicitly in this example.}. 
We shall now compute that the value of the glued factorization algebra \(\Fa\)
on the open set \(X=(0,5)\) is the relative tensor product
$\bimod_1\otimes_{\algebra_2} \bimod_2$.
Moreover, $\Fa$ also encodes the data of $\bimod_1\otimes_{\algebra_2} \bimod_2$
as a \((\algebra_1,\algebra_3)\)-bimodule.

To compute the value of \(\Fa\) on \(X\),
we denote by \(\Ka\) the poset of nonempty compact subsets of \(U\cap V=(2,3)\)
with finitely many connected components
and consider the Weiss presieve
\begin{equation}
  \Wa\coloneqq\{X\setminus K\mid K\in \Ka\}\to X.
\end{equation}
It is a well known result that the canonical functor
\begin{equation}
  \label{eq:localization-K-Delta}
  \Wa\cong \Ka^\op\xrightarrow{\pi_0} \Delta^\op,
\end{equation}
(where each set \(\pi_0(K)\) is endowed with the linear order induced from \(\reals\))
is an \(\infty\)-categorical localization;
for convenience we also provide a short proof as
\Cref{lem:localize-to-Delta} in the appendix.
Since \(\Fa\) is constructible and multiplicative, the induced functor
\(\Fa\colon \Wa\to \targetcat\)
factors through this localization.
For each \(n\in \naturals\),
we choose an arbitrary compact subset \(K=K_0\disjun\dots\disjun K_n\) of \(U\cap V\)
with exactly \(n+1\) components (written in ascendent order)
and let \(U\coloneqq X\setminus K=U_0\disjun\dots\disjun U_{n+1}\)
be its open complement with exactly \(n+2\) components (written in ascendent order).
Then we compute
\begin{align}
  \Fa(X)
  &\simeq
    \colim \Fa\restrict{\Wa}
    \simeq
    \colim_{[n]\in \Delta^\op}
    \Fa(X \setminus (K_0\disjun\cdots\disjun K_n))
  \\
  &
    \simeq
    \colim_{[n]\in \Delta^\op}
    \Fa(U_0)\otimes\cdots\otimes\Fa(U_{n+1})
    \simeq
    \colim_{[n]\in\Delta^\op}
    \bimod_1\otimes \algebra_2^{\otimes n}\otimes \bimod_2 \,.
\end{align}
From left to right in the chain of equivalences, we use that
that \(\Fa\) is a Weiss cosheaf,
that the localization functor \eqref{eq:localization-K-Delta} is cofinal,
\(\Fa\) is multiplicative,
and that the singularity points \(1\) and \(4\) lie in \(U_0\) and \(U_{n+1}\),
respectively.
The final simplicial colimit is precisely the usual definition of the relative tensor product
\(\bimod_1\otimes_{\algebra_2} \bimod_2 \).
\end{ex}

Without going into great detail we will quickly explain why these simple examples are important in the context of \emph{the higher Morita category} \(\mathrm{Alg}_1(\targetcat)\) constructed from constructible factorization algebras: 
\begin{rem}
	The above example is the key step in proving composition of 1-morphisms in the higher Morita category \(\mathrm{Alg}_1(\targetcat)\) of \cite{Scheimbauer, GS}. Very informally, objects of this \((\infty, 2)\)-category are locally constant factorization algebras on \((0,1)\), 1-morphisms are constructible factorization algebras for \((0,1)\) stratified with a marked point, and 2-morphisms are maps of such constructible factorization algebras.
	The source and target of a 1-morphism are given by restricting to the left and right of the chosen point, respectively.
	Given two 1-morphisms, i.e.\ \(\Fa_0\) and \(\Fa_1\) such that the target and source are equivalent, we compose as follows. One first glues together the two constructible factorization algebras to produce a constructible factorization algebra \(\Fa\) on an interval with 2 marked points as in \Cref{ex:GlueFactCstr}. 
Since the gluing is only defined up to equivalence, composition is only defined up to equivalence.
Indeed, \(\mathrm{Alg}_1(\targetcat)\) is actually an $(\infty,1)$-category modelled by a complete Segal space; the situation with several points as the stratification corresponds to proving the Segal condition in general.

Then we push forward along a ``collapse-and-rescale map'' \(\varrho\) that collapses everything between the 2 marked points to a single point and rescales back to get \((0,1)\) as illustrated below. 
	
\begin{center}
%\begin{figure}[H]\label{fig:PCovRescale}
\begin{tikzpicture}[scale=1.3]
% top pic
\draw[draw =bluebg, ultra thick] (0,0) -- (0.3, 0);
\draw[draw =redbg, ultra thick] (0.3, 0) -- (1,0);
\draw[draw =greenbg, ultra thick] (1.3, 0) -- (1,0);

\draw (0.3, 0) node[dot] {};
\draw (1, 0) node[dot] {};

% middle pic
\draw[draw =bluebg, ultra thick] (0.35,-1) -- (0.65, -1);
\draw[draw =greenbg, ultra thick] (0.95,-1) -- (0.65, -1);
\path (0.95,-1) node (B) {} -- (0.35,-1) node (A) {};
\draw (0.65, -1) node[dot] {};

%lower pic
\draw[draw =bluebg, ultra thick] (0.15, -2) -- (0.65, -2);
\draw[draw =greenbg, ultra thick] (0.65, -2) -- (1.15, -2);
\draw (0.65, -2) node[dot] {};

\draw[densely dotted] (0, 0) -- (0.35, -1)
(0.3, 0) -- (0.65, -1) -- (1, 0)
(1.3, 0) -- (0.95,-1)

(0.35,-1) -- (0.15, -2)
(0.95, -1) -- (1.15, -2);

%arrow
\draw[->] (2, 0) -- node[anchor=west] {$\varrho \times \id$} (2, -2);
\end{tikzpicture}
%\end{figure}
\end{center}
Pushing forward along this map produces a {\em constructible} factorization algebra \(\varrho_\ast \Fa\) on \((0,1)\) with 1 marked point, which is again a 1-morphism in \(\mathrm{Alg}_1(\targetcat)\). 
	
\end{rem}

The above example of course does not depend on working with the interval in any way. There is a more general notion of the higher Morita category, denoted \(\mathrm{Alg}_n(\targetcat)\), and generalising the above procedure indeed gives compositions of \(k\)-morphisms, for all \(k \leq n\).

\subsection{The symmetric monoidal \infy-category of constructible factorization algebras}
\label{sec:SymMonStructureFactAlgs}

Let \(X\) be a smooth conical manifold and
\(\targetcatOT\) a \(\otimes\)-presentable symmetric monoidal \(\infty\)-category.
The goal of this section is to construct a symmetric monoidal structure
on the $\infty$-category \(\FactCstr[\targetcat]{X}\)
of constructible factorization algebras.

On prefactorization algebras it is easy to write down a symmetric monoidal structure.
Namely, given prefactorization algebras $\Fa$ and $\Ga$ we set
\begin{equation}
  \label{eq:product-of-factalg}
  (\Fa \otimes \Ga) (U) \coloneqq\Fa(U) \otimes \Ga(U).
\end{equation}
This was already suggested in \cite[\S3 Section~1.5]{CG}.
However, it is not immediately obvious that this satisfies the Weiss condition.

The formula \eqref{eq:product-of-factalg}
can be interpreted as the following two-step process:
\begin{itemize}
\item
  Glue the prefactorization algebras \(\Fa\) and \(\Ga\) on two copies of \(U\)
  to a prefactorization algebra \(\Fa\amalg \Ga\) on \(U\amalg U\).
\item
  Pushforward along the fold map \(\nabla\colon U\amalg U\to U\)
  to define the tensor product
  \begin{equation}
    \Fa\otimes\Ga\coloneqq \pf{\nabla}(\Fa\amalg\Ga).
  \end{equation}
\end{itemize}

We now implement this idea to construct the fully coherent
symmetric monoidal structure on the $\infty$-categories
of {\em constructible} factorization algebras.

\begin{constr}
  \label{constr:symmetric_monoidal}
  Let $X$ be a topological space. 
  Consider the functor of \(1\)-categories
  \begin{equation}
    \Fin \to \Opd^\op,
    \quad I_+ \mapsto \Open{X^{\amalg I}}
  \end{equation}
  where a map $f\colon I_+ \to J_+$ is sent to
  \begin{equation}
    \Open{X^{\amalg J}}  \to \Open{X^{\amalg I}},
    \qquad
    \coprod_{j\in J} {U_j} \mapsto \coprod_{i\in I} {U_{f(i)}}.
  \end{equation}
  Here we use the convention that for the base point $+$ of \(J_+\)
  we have $U_+\coloneqq\emptyset$.

  Then we postcompose with
  \begin{equation}
    \Alg[\targetcat]{-} \colon \Opd^\op \to \Catinfty
  \end{equation}
  to get the functor
  \begin{equation}
    \AlgO[\targetcat]{X}\colon \Fin\to\Catinfty,
    \quad
    I_+\mapsto \Alg[\targetcat]{X^{\amalg I}}.
  \end{equation}
  This functor is the one exhibiting the symmetric monoidal structure
  on the \(\infty\)-category \(\Alg[\targetcat]{X}\) of prefactorization algebras.
\end{constr}

\begin{prop}
  We have a subfunctor
  \begin{equation}
    \FactCstrO{X}\subset\AlgO{X} \colon \Fin\to\Catinfty
  \end{equation}
  by pointwise passing to the full subcategories of constructible factorization algebras.
\end{prop}

\begin{proof}
  We need to check that for each $f\colon I_+ \to J_+$ in \(\Fin\),
  the functor \(\AlgO{X}(f)\colon \Alg{X^{\amalg I}}\to \Alg{X^{\amalg J}}\)
  preserves constructible factorization algebras.
  We may show this separately when \(f\) is active or inert
  because every map in \(\Fin\) is a composite of an active and an inert map.

  If \(f\) is active,
  then $\AlgO{X}(f)$ is given by pushing forward along a disjoint union
  \begin{equation}
    \coprod_{j\in J}
    \left(
      X^{\amalg \preim[j]{f}}
      \to
      X
    \right)
  \end{equation}
  of fold maps.
  This pushforward preserves constructible factorization algebras by
  \Cref{lemma:fold_is_constr} and \Cref{lemma:union_is_constr}.

  If \(f\) is inert, then $\AlgO{X}(f)$ is restriction along an open embedding,
  which preserves constructible factorization algebras by
  \Cref{lem:restriction-preserves}.
\end{proof}

\begin{cor}\label{thm:symmetric_monoidal}
  For every smooth conical manifold $X$ with enough good disks,
  the functor
  \begin{equation}
    \FactCstrO{X}\colon \Fin\to \Catinfty,
    \quad\quad
    I_+\mapsto \FactCstr{X^{\amalg I}}
  \end{equation}
  exhibits a monoid object in \((\Catinfty,\times)\),
  i.e., a symmetric monoidal structure
  on the $\infty$-category \(\FactCstr{X}\) of constructible factorization algebras on $X$.
\end{cor}

\begin{proof}
  We need to show that the functor satisfies the Segal condition.
  That is, that for each \(I_+,J_+\in \Fin\)
  we form the wedge sum $I_+\vee J_+ =(I+J)_+$ and show
  the diagram induced by the projections
  \begin{equation}\label{Segal_square}
    \begin{tikzcd}
      (I+J)_+  \arrow{r}{\iota_I} \arrow{d}{\iota_J} &  I_+  \arrow{d}\\
      J_+  \arrow{r}  & \emptyset_+
    \end{tikzcd}
  \end{equation}
  is a pullback square.
  Applying $\FactCstrO{X}$, this amounts to showing that the square
  \begin{equation}
    \begin{tikzcd}
      \FactCstr{X^{\amalg (I+J)}} \arrow{r}\arrow{d} & \FactCstr{X^{\amalg I}}  \arrow{d}\\
      \FactCstr{X^{\amalg J}}   \arrow{r} & *
    \end{tikzcd}
  \end{equation}
  is a pullback.
  This is an instance of \Cref{thm:main-gluing} for the decomposition
  \(X^{\amalg(I+J)}=X^{\amalg I}\coprod X^{\amalg J}\).
\end{proof}

Finally we record an improvement on our main gluing theorem (\Cref{thm:main-gluing}),
exhibiting constructible factorization algebras not just as a sheaf of \(\infty\)-categories
but of \emph{symmetric monoidal} \(\infty\)-categories.

Recall that \(\Snglr\) is the \(1\)-category of smooth conical manifolds
and open embeddings between them.

\begin{constr}
  \label{cstr:presheaf-FactCstrO}
  Consider the functor of \(1\)-categories
  \begin{equation}
    \Snglr\times \Fin \to \Opd,
    \quad
    (X,I_+)\mapsto \Open{X^{\amalg I}}
  \end{equation}
  extending the one of \Cref{constr:symmetric_monoidal}.
  By postcomposing with \(\Alg[\targetcat]{-}\)
  and then passing to the subfunctor of constructible factorization algebras,
  we obtain the functor
  \begin{equation}
    \Snglr\times \Fin \to \Catinfty,
    \quad
    (X,I_+)\mapsto \FactCstr{X^{\amalg I}}.
  \end{equation}
  Since for fixed \(X\in\Snglr\) this defines a symmetric monoidal \(\infty\)-category
  (= a Segal object in \(\Catinfty\)),
  we can transpose it to obtain a functor
  \begin{equation}
    \SheafFactCstrO[\targetcat]{}\colon
    \Snglr\to \SMCatinfty\subset \Fun(\Fin,\Catinfty),
    \quad
    X\mapsto \FactCstrO[\targetcat]{X}.
  \end{equation}
\end{constr}

\begin{prop}\label{prop:symm_mon_gluing}
  The functor \(X\mapsto \FactCstrO[\targetcat]{X}\)
  of \Cref{cstr:presheaf-FactCstrO}
  is a sheaf of symmetric monoidal \(\infty\)-categories.
\end{prop}

\begin{proof}
  We have to show that for each open cover \(\{U_i\}_{i<\beta}\)
  of a smooth conical manifold \(X\),
  the restriction functors induce an equivalence
  \begin{equation}
    \FactCstrO{X}
    \xrightarrow{\simeq}
    \lim_{\beta\supset I \supsetneq \emptyset}\FactCstrO{U_I}
  \end{equation}
  of symmetric monoidal \(\infty\)-categories.
  When encoded as Segal objects \(\Fin\to \Catinfty\),
  limits of symmetric monoidal \(\infty\)-categories are computed pointwise;
  hence we only have to show that for each \(J_+\in \Fin\),
  the induced map
  \begin{equation}
    \FactCstr{X^{\amalg J}}
    \xrightarrow{\simeq}
    \lim_{\beta\supset I \supsetneq \emptyset}\FactCstr{U_I^{\amalg J}}
  \end{equation}
  is an equivalence of \(\infty\)-categories.
  This is an instance of \Cref{thm:main-gluing}
  for the cover \(\{U_i^{\amalg J}\}_{i<\beta}\) of \(X^{\amalg J}\).
\end{proof}

\subsection{Link decomposition of constructible factorization algebras}
\label{sec:link-decomp}
%local macros
\newcommand{\Mdim}[1]{M_{#1}}
\newcommand{\Mgeq}[1]{M_{\geq#1}}
\newcommand{\UfwithES}{\Uf_{\emptyset}} % (shortcut) version of Uf with the empty set added

We now turn to applying our toolbox to answer the following question, which was
proposed to us by David Ayala.

\begin{question}
  Is the data of a constructible factorization algebra on a cone
  $\topcone{Z}$ the same as the following?
  \begin{enumerate}
  \item\label{it:algebra-in-fact-on-Z}
    An associative algebra \(\Aa\) in constructible factorization algebras on $Z$ and
  \item
    a pointed module for the associative algebra $\Aa(Z)=\int_{Z}\Aa$.
  \end{enumerate}
\end{question}

We reinterpret this question in the following way:

\begin{itemize}
\item
  We interpret associative algebras as
  locally constant factorization algebras on \((0,\infty)\)
  and pointed modules with constructible extensions to \([0,\infty)\);
  see \Cref{ex:ass-and-rmod}.
\item
  Using a version of Dunn's additivity theorem,
  we can thus interpret the data in \eqref{it:algebra-in-fact-on-Z}
  as a constructible factorization algebra \(\Aa\) on \(Z\times (0,\infty)\).
  Such additivity was established in Berry's thesis~\cite{Eric} in the unstratified setting;
  the general result for smooth conical manifolds
  is work in progress by Anja \v{S}vraka\footnote{
    To even state such an additivity result in the stratified situation 
    one needs to construct the symmetric monoidal structure on the \(\infty\)-category
    of constructible factorization algebras.
    We do this in 
    \Cref{sec:SymMonStructureFactAlgs}.
  }.
\end{itemize}

With this reformulation, the following theorem gives a positive answer to Ayala's question:

\begin{thm}
  \label{thm:cones_as_modules}
  Let \(\targetcatOT\) be a \(\otimes\)-presentable
  symmetric monoidal \(\infty\)-category.
  Let $Z$ be a compact smooth conical manifold
  with enough good disks.
  We have a pullback square of \(\infty\)-categories
  \begin{equation}
    \label{eq:cone-pb-square}
    \begin{tikzcd}
      \FactCstr[\targetcat]{ \topcone{Z} }  \arrow{d}{\pf{p}} \arrow{r}
      \isCartesian
      &
      \FactCstr[\targetcat]{Z\times (0,\infty) }   \arrow{d}{\pf{p}}
      &
      \Aa
      \ar[d,mapsto]
      \\
      \FactCstr[\targetcat]{[0,\infty) } \arrow{r}
      &
      \FactCstr[\targetcat]{ (0,\infty) }
      &
      \int_Z\Aa
    \end{tikzcd}
  \end{equation}
  where the vertical functors are given by pushforward
  along the quotient map \(p\colon \topcone{Z}\to [0,\infty)\)
  and the horizontal functors are given by restriction to open subspaces.
\end{thm}

\begin{rem}
  Specifically, the square \eqref{eq:cone-pb-square}
  is induced by applying \(\Alg[\targetcat]{-}\)
  to the following square of (discrete, colored) operads
  \begin{equation}
    \label{eq:operad-square-cone}
    \cdsquareOpt
    {\Open{\topcone{Z}}}
    {\Open{Z\times (0,\infty)}}
    {\Open{[0,\infty)}}
    {\Open{(0,\infty)}}
    {hookleftarrow}
    {leftarrow, "p^{-1}"}
    {leftarrow, "p^{-1}"}
    {hookleftarrow}
  \end{equation}
  and then passing
  to the full subcategories of constructible factorization algebras.
\end{rem}

\begin{rem}
  \label{rem:brav-rozenblyum}
  After the completion of this paper we were made aware of a result of
  Brav-Rozenblyum~\cite[Section~3.2]{brav-rozenblyum}
  who independently proved the analog of \Cref{thm:cones_as_modules}
  but on the other side of the equivalence of \Cref{cor:FactX-AlgDiskX}.
  Their proof techniques rely on monadicity results and appear
  to be completely different from ours.
\end{rem}

Before we go into the proof of \Cref{thm:cones_as_modules},
we explain how this local result can be globalized.
This globalization allows constructible factorization algebras on a smooth conical manifold
to be assembled inductively from algebras and pointed modules in
locally constant factorization algebras
on the individual strata, each of which is an (ordinary) smooth manifold.

To be precise,
this globalization relies on a slightly more general version of the cone decomposition theorem,
which for clarity we now formulate as a separate conjecture.
Observe that
by applying \Cref{thm:cones_as_modules} to \(\targetcat=\FactCstr{\reals^d}\),
this conjecture follows directly
once a version of additivity for \(-\times \reals^d\) is established.
For this reason we do not attempt to prove the conjecture directly,
even though we expect that this should be possible
by modifying the arguments in our proof of
\Cref{thm:cones_as_modules}
to account for the extra copy of \(\reals^d\) everywhere.

\begin{conj}
  \label{conj:cone-decomp-Rd}
  In the setting of \Cref{thm:cones_as_modules} we have a pullback square
  \begin{equation}
    \begin{tikzcd}
      \FactCstr[\targetcat]{ \topcone{Z}\times\reals^d }  \arrow{d}{\pf{p}} \arrow{r}
      \isCartesian
      &
      \FactCstr[\targetcat]{Z\times (0,\infty)\times\reals^d }   \arrow{d}{\pf{p}}
      &
      \Aa
      \ar[d,mapsto]
      \\
      \FactCstr[\targetcat]{[0,\infty)\times\reals^d } \arrow{r}
      &
      \FactCstr[\targetcat]{ (0,\infty)\times\reals^d }
      &
      \int_Z\Aa
    \end{tikzcd}
  \end{equation}
  for each \(d\geq 0\).
\end{conj}

\begin{constr}[Global link decomposition, using \Cref{conj:cone-decomp-Rd}]
  \label{constr:step-unlinking}
  Let \(M\) be a smooth conical manifold
  and let \(d\) be the smallest dimension appearing among the strata of \(M\).
  Let \(S\subset M\) be the union of all strata of \(M\) of dimension \(d\);
  in particular \(S\) is closed in \(M\).

  Around \(S\) there is a tubular neighborhood \(T=T(S)\)
  (see \cite[Proposition~8.2.3]{AFT-local-structures})
  which yields a (usually nontrivial) cone-bundle \(p\colon T\to S\)
  so that all fibers of \(p\)
  are of the form \(\topcone{Z}\) for some compact smooth conical manifold \(Z\)
  (the link at \(s\in S\) in \(M\)).
  The bundle \(p\colon T\to S\) has the cone-point-section
  which corresponds to the inclusion \(i\colon S\subset T\) into the tubular neighborhood.

  Denote by \(r\colon V\to S\)
  the ray bundle obtained from \(p\colon T\to S\)
  by fiberwise passing to the quotient \(q\colon \topcone Z\to [0,\infty)\).
  Note that unlike the cone bundle \(p\colon T\to S\),
  the ray bundle \(r\colon V\to S\) is automatically trivial
  because the space of stratified self-diffeomorphisms of \([0,\infty)\) is contractible;
  hence we can write \(V=S\times [0,\infty)\).

  Thus for each local chart \(\reals^d\cong U\subset S\)
  we have the situation
  \begin{equation}
    \begin{tikzcd}
      {T_U}
      \ar[d,"q_U"]
      &
      {\reals^d\times \topcone{Z}}
      \ar[r,hookleftarrow]
      \ar[d,"\reals^d\times q"]
      \ar[l,"\cong"]
      &
      {\reals^d\times Z\times (0,\infty)}
      \ar[d,"\reals^d\times q"]
      \ar[r,"\cong"]
      &
      {T_U\setminus U}
      \ar[d,"q_U"]
      \\
      U\times [0,\infty)
      \ar[d,"r_U"]
      &
      {\reals^d\times [0,\infty)}
      \ar[r,hookleftarrow]
      \ar[l,"\cong"]
      \ar[d]
      &
      {\reals^d\times (0,\infty)}
      \ar[r,"\cong"]
      &
      U\times (0,\infty)
      \\
      U
      \ar[uu,bend left=80,hookrightarrow,"i_U"]
      &
      {\reals^d}
      \ar[l,"\cong"]
    \end{tikzcd}
  \end{equation}

  Therefore we obtain the following pullback square of \(\infty\)-categories
  \begin{equation}
    \cdsquareNA[pb]
    {\FactCstr{T_U}}
    {\FactCstr{T_U\setminus U}}
    {\FactCstr{U\times[0,\infty)}}
    {\FactCstr{U\times(0,\infty)}}
  \end{equation}
  by applying \Cref{conj:cone-decomp-Rd}.

  Finally we use our main gluing result, \Cref{thm:main-gluing},
  to glue all these pullback squares together,
  where we let \(U\) range over an intersection-closed collection of charts of \(S\),
  for example, the geodesically convex disks with respect to a chosen metric.
  Therefore we obtain the lower of the following two pullback squares;
  the upper one is another instance of gluing,
  namely for the cover \(M=T \cup (M\setminus S)\).
  \begin{equation}
    \label{eq:pasted-pullback-globalized-link}
    \begin{tikzcd}
      {\FactCstr{M}}
      \ar[r]\ar[d]
      \isCartesian
      &
      {\FactCstr{M\setminus S}}
      \ar[d]
      \\
      {\FactCstr{T}}
      \ar[r]\ar[d]
      \isCartesian
      &
      {\FactCstr{T\setminus S}}
      \ar[d]
      \\
      {\FactCstr{S\times [0,\infty)}}
      \ar[r]
      &
      {\FactCstr{S\times (0,\infty)}}
    \end{tikzcd}
  \end{equation}
\end{constr}

\begin{cor} [Using \Cref{conj:cone-decomp-Rd}]
  \label{rem:interpret-unlink}
  Let \(M\) be a smooth conical manifold.
  For each \(n\in\naturals\) we have a pullback square
  \begin{equation}
    \label{eq:unlink-n}
    \cdsquareNA[pb]
    {\FactCstr{\Mgeq{n}}}
    {\FactCstr{\Mgeq{n+1}}}
    {\FactCstr{\Mdim{n}\times[0,\infty)}}
    {\FactCstr{\Mdim{n}\times(0,\infty)}}
  \end{equation}
  of \(\infty\)-categories,
  where \(\Mdim{n}\) or \(\Mgeq{n}\)
  denotes the union of all strata of \(M\)
  of dimension \(n\) or of dimension \(\geq n\), respectively.
\end{cor}

\begin{proof}
  For each \(n\), apply \Cref{constr:step-unlinking}
  to the smooth conical manifold \(\Mgeq{n}\),
  whose lowest dimensional strata comprise \(S=\Mdim{n}\).
\end{proof}

\begin{rem}\label{remark:globalizing_unlinking}
  Using the interpretation established at the beginning of the section,
  the pullback square \eqref{eq:unlink-n}
  can be viewed as a
  decomposition of the datum of a constructible factorization algebra on
  \(\Mgeq{n}\) into
  \begin{itemize}
  \item
    a constructible factorization algebra \(\Aa\)
    on \(\Mgeq{n+1}\)
  \item
    together with pointed right modules for the associative algebra
    \(\int_{\Link{M}{n}}\Aa\)
    in locally constant factorization algebras on \(\Mdim{n}\),
  \end{itemize}
  where \(\Link{M}{n}\) denotes the link bundle of \(\Mdim{n}\) in \(\Mgeq{n}\)
  see \Cref{not:link}.

  Overall we can then say that the descending filtration
  \begin{equation}
    M=\Mgeq{0}\supseteq \Mgeq{1} \supseteq \Mgeq{2} \supseteq\cdots,
  \end{equation}
  induces a sequence of cartesian fibrations
  \begin{equation}
    \begin{tikzcd}
      \FactCstr{M}=\FactCstr{\Mgeq{0}}
      \ar[r]&
      \FactCstr{\Mgeq{1}}
      \ar[r]&
      \FactCstr{\Mgeq{2}}
      \ar[r]&
      \cdots
    \end{tikzcd}
  \end{equation}
  whose fibers are all given by categories of right modules in
  \(
  \FactLC{\Mdim{0}}, 
  \FactLC{\Mdim{1}}, 
  \FactLC{\Mdim{2}}, 
  \dots,
  \)
  respectively,
  each of which is locally constant factorization algebras on an (unstratified) smooth manifold.
\end{rem}

The rest of this section is dedicated to the proof of
\Cref{thm:cones_as_modules}.
We start with some notation. 
Recall from \Cref{defn:restricted-sieve} that \(\cdownR{\Ua}\) denotes 
the \restrictedsieve{} associated to a precover \(\Ua\), 
i.e.\ in words it is the usual associated sieve minus the empty set.

\begin{constr}
  \label{cstr:cone-operads}
  Denote by 
  \begin{equation}
    \Nf\coloneqq\cdownR{(Z\times(0\times\infty))}
    \quad
    \text{and}
    \quad
    \Hf\coloneqq p^{-1}\cdownR{[0,\infty)}
  \end{equation}
  the \posetsofopens{} of non-empty opens that
  \begin{itemize}
  \item
    do \emph{not} contain the cone point, or
  \item
    are \emph{horizontal}, 
    i.e., of the form \(p^{-1}(U)\) for some non-empty open \(U\subseteq [0,\infty)\),
  \end{itemize}
  respectively.
  Write \(\Uf\coloneqq \Nf\cup\Hf\) for their union.
  By construction, the associated operad \(\operad{\Uf}\) fits
  into a pullback square 
  \begin{equation}
    \label{eq:dendroidal-set-pushout}
    \cdsquareOpt[po]
    {\operad{\Uf}}
    {\operad{\Nf}}
    {\operad{\Hf}}
    {\operad{\Nf}\cap\operad{\Hf}}
    {hookleftarrow}
    {hookleftarrow}
    {hookleftarrow}
    {hookleftarrow}
  \end{equation}
  of operads and fully faithful embeddings.
  Moreover, we write \(\UfwithES\coloneqq\Uf\cup\{\emptyset\}\)
  obtained from \(\Uf\) by putting the empty set back in.
\end{constr}

\begin{prop}
  \label{lem:cone-operad-pushout}
  The square \eqref{eq:dendroidal-set-pushout} is a pushout
  of operads and even of \infy-operads.
\end{prop}

\begin{figure}[H]
  \includegraphics*[width=0.4\textwidth]{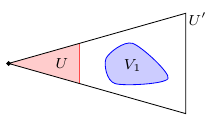} 
  \caption{
    An example of a 2-ary operation \(f\colon U\disjun V_1 \ra U'\) in \(\operad{\Uf}\)
    which does not lie in \(\Hf\) because \(V_1\) is not horizontal,
    but also does not lie in \(\Nf\) because \(U\) and \(U' = \topcone{Z}\)
    contain the cone point.
    The open \(U''\) is the interior of the cone minus the open \(U\). 
    One directly sees that the inclusions \(a\colon U\disjun U'' \hookrightarrow U'\) and \(b_f \colon V_1 \hookrightarrow U''\) compose to the original 2-ary operation~\(f\). }
  \label{fig:cone-inclusions-squiggly}
\end{figure}

\begin{proof}[Key insight]
  Let us just explain the key insight that goes into proving this fact.
  The actual proof is unfortunately rather technical,
  so we delay it until the end of the section;
  see \Cref{prop:operad-pushout} below.
  The main obstacle to overcome is that there are operations in \(\operad{\Uf}\)
  whose output open contains the cone point
  but at least one input open is not horizontal. 
  See \Cref{fig:cone-inclusions-squiggly} for an easy example.
  Such an operation will not lie in \(\Nf\) nor in \(\Hf\);
  nonetheless, we need to show that it lies in the pushout.%
  \footnote{
    The reader might wonder why an analogous obstacle did not arise
    in the setting of \Cref{subsect:FirstStepGluing},
    where it was almost tautological that the cube 
    \eqref{eq:cube-dendroidal-nerve}
    was a strict pushout of dendroidal sets.
    In that setting every operation of the big operad
    \(\cdownR{\Ua}\)
    with output color in some small operad \(\cdownR{U_I}\)
    was already fully contained in that small operad,
    because all \(\cdownR{U_I}\) were \restrictedsieves{}.
  }
  The reason is that such an operation \(f\) can be factored canonically
  into a composition of an operation \(b\in\operad{\Nf}\)
  away from the cone point,
  and a horizontal operation \(a\in\operad{\Hf}\):
  \begin{itemize}
  \item
    Let \(U'\) be the output open of the operation \(f\),
    which we assume contains the cone point
    (otherwise \(f\) lies in \(\Nf\) and there is nothing to do).
  \item
    There is at most one input open of \(f\) which contains the cone point.
    If such an open exists, we denote it by \(U\).
    We denote the other input opens of \(f\) by \(V_1,\dots,V_n\).
  \item
    Consider the open \(U''\) defined as the interior
    of \(U'\cap (Z\times (0,\infty))\setminus U\).
    Note that \(U''\) is again horizontal, since \(U'\) and \(U\) are horizontal.
    By construction, we have an inclusion
    \begin{equation}
      a\colon (U\disjun) U''\hookrightarrow U',
    \end{equation}
    which is a \(2\)- or \(1\)-ary operation of \(\Hf\) depending on whether \(U\) exists or not.
    Note that \(a\) just depends on (\(U\) and) \(U'\), but not on \(f\).
  \item
    Since the \(V_i\) are open, do not contain the cone point and are disjoint to \(U\),
    we have an inclusion
    \begin{equation}
      b_f\colon V_1\disjun \cdots\disjun V_n\hookrightarrow U'',
    \end{equation}
    which lies in \(\Nf\).
  \item
    By construction, the composite \(a\circ b_f\) is the original
    (\((n+1)\)-ary or \(n\)-ary) operation
    \begin{equation}
      f\colon (U\disjun) V_1\disjun \dots \disjun V_n \hookrightarrow U'.
    \end{equation}
  \end{itemize}
  Thus every operation of \(\operad{\Uf}\) can be written as a composite
  of operations in \(\Nf\) and \(\Hf\).
  The technical part is to use the universality of the operations \(a\)
  described above to obtain the desired pushout;
  this is the content of \Cref{prop:operad-pushout}.
  \end{proof}

%%commands for figure
\def\width{3.0} %%the length of the cone
\def\height{1.7} % the height of the cone
\def\hheight{0.5*\height} %half of the height
\def\small{0.4*\width} %small cone length
\def\sheight{0.4*\hheight} % the half-height of the small cone
\newcommand{\point}[3]{\draw[fill=black] (#1, #2) circle (#3);}

\begin{proof}[Proof of \Cref{thm:cones_as_modules}]
  From \Cref{lem:cone-operad-pushout}
  it follows that after passing to algebras in \(\targetcat\),
  we get a pullback square of \(\infty\)-categories
  \begin{equation}
    \cdsquareNA[pb]
    {\Alg{\Uf}}
    {\Alg{\cdownR{Z\times (0,\infty)}}}
    {\Alg{\cdownR{[0,\infty)}}}
    {\Alg{\cdownR{(0,\infty)}}} \ . 
  \end{equation}

Next, we claim that we retain a pullback square when passing to the full subcategories of multiplicative prefactorization algebras. 
That is, we need to show that an algebra \(\Aa\) on \(\Uf\) is multiplicative if its restriction to each of the two \posetofopens{} 
\(\cdownR{(Z\times (0,\infty))}\) and \(p^{-1}\cdownR{[0,\infty)}\) are multiplicative. 
Consider a decomposable open \(U\in\Uf\). 
If \(U\) lies in \(\cdownR{(Z\times(0,\infty))}\), then the same is true for all of its components, so that we can apply the multiplicativity on \(\cdownR{(Z\times(0,\infty))}\).
If \(U\) lies in \(p^{-1}\cdownR{[0,\infty)}\) and contains the cone point, then we can write it as \(U=U'\disjun U''\), 
where \(U'\in p^{-1}\cdownR{[0,\infty)}\) is the connected component of the cone point and \(U''\in \cdownR{(Z\times(0,\infty))}\). 
By multiplicativity on \(p^{-1}\cdownR{[0,\infty)}\), the algebra \(\Aa\) is multiplicative with respect of the decomposition \(U=U'\disjun U''\). 
Any finer decomposition must happen completely within \(U''\) (since \(U'\) is connected), so that we 
can reduce to the case of \(U\in \cdownR{(Z\times (0,\infty))}\) which we already established.
In total we have a pullback square 
  \begin{equation}\label{eq:pullback-after-passing-to-multiplicative-subcats}
    \cdsquareNA[pb]
    {\AlgM{\Uf}}
    {\AlgM{\cdownR{Z\times (0,\infty)}}}
    {\AlgM{\cdownR{[0,\infty)}}}
    {\AlgM{\cdownR{(0,\infty)}}} \ . 
  \end{equation}

Next we want to add the empty set to all of the \posetsofopens{} appearing above. 
For every \posetofopens{} \(\Bf\) we have proven in
\Cref{lem:remove-empty-set-multiplicative}
that the restriction functor \(\AlgM{\Bf} \ra \AlgM{\Bf'}\) is an equivalence. 
Applying this to \(\Bf\) being
\(\UfwithES\coloneqq\Uf\cup\{\emptyset\}\), \(\cdown{(Z\times(0,\infty))}\),
\(\cdown{[0, \infty)}\) and \(\cdown{(0,\infty)}\), respectively,
gives that the pullback square
\eqref{eq:pullback-after-passing-to-multiplicative-subcats}
is equivalent to the square
  \begin{equation}
    \cdsquareNA[pb]
    {\AlgM{\UfwithES}}
    {\AlgM{\cdown{Z\times (0,\infty)}}}
    {\AlgM{\cdown{[0,\infty)}}}
    {\AlgM{\cdown{(0,\infty)}}}
  \end{equation}
  which is therefore also a pullback.
Additionally, observe that \(\UfwithES\) is now a presieve because sets of each type are closed under intersection and 
an intersection \(U\cap V\) with \(U\in \cdown{(Z\times(0,\infty))}\) again lies in \(\cdown{(Z\times(0,\infty))}\). 

We now want to show that we maintain a pullback when passing to the full subcategories of constructible factorization algebras. 
Explicitly, we have to show that a multiplicative algebra \(\Aa\) on \(\UfwithES\) is Weiss and constructible if its restriction to each of the \posetofopens{}  
\(\cdown{(Z\times (0,\infty))}\) and \(p^{-1}\cdown{[0,\infty)}\) has the corresponding property:
\begin{itemize}
  \item
    For constructibility, let \(D\hookrightarrow D'\) be an inclusion of abstractly
    isomorphic conical disks in \(\UfwithES\).
    Note that it cannot happen that one contains the cone point
    and the other one does not,
    because then they would not be abstractly isomorphic (\Cref{lem:char-iso-disk}).
    Therefore they both lie in \(p^{-1}\cdown{[0,\infty)}\)
    or they both lie in \(\cdown{(Z\times(0,\infty))}\);
    in either case we are done.
   \item 
    For the Weiss condition,
    we have to show that every Weiss presieve \(\Wf\hyprefines U\)
    in \(\UfwithES\) is \(\Aa\)-local
    (this makes use of the fact that \(\UfwithES\) is a presieve, see
    \Cref{rem:Weiss-cosheaf-Weiss-presieves}).
    If \(U\) does not contain the cone point, then neither does any open in \(\Wf\),
    so that we can just invoke the Weiss condition on \(\cdown{Z\times (0,\infty)}\).
    It remains to consider the case where \(U\) contains the cone point.
    Let \(\Wf'\subset \Wf\) be the subpresieve consisting of those \(W\in \Wf\)
    which contain the cone point \(0\).
    \begin{itemize}
    \item
      The presieve \(\Wf'\) is still a Weiss cover of \(U\).
      Indeed, for any finite subset \(S\subset U\) there must be a \(W\in \Wf\)
      that contains the finite set \(S\cup \{0\}\)
      (because \(\Wf\hyprefines U\) is a Weiss presieve)
      but then \(W\) lies in \(\Wf'\) by definition and contains \(S\).
    \item
      Next, we claim that the inclusion \(\Wf'\hookrightarrow \Wf\)
      is \(\Aa\)-local.
      To this end, observe that for each \(W\in \UfwithES\cap\cdown{\Wf}\),
      we have a Weiss presieve
      \begin{equation}
        \label{eq:Weiss-presieve-W'-W}
        \UfwithES\cap\cdown{\Wf'}\cap\cdown{W}\hyprefines W
      \end{equation}
      by the same argument as above:
      for each finite \(S\subset W\) there is a \(W'\in \Wf'\)
      with \(W'\supset S\cup\{0\}\),
      hence \(W'\cap W\) is an element of
      \(\Uf\cap\cdown{\Wf'}\cap\cdown{W}\) containing \(S\).
      Moreover, note that the Weiss presieve 
      \eqref{eq:Weiss-presieve-W'-W}
      is degenerate if \(W\) contains the cone point
      (because then \(W\in \cdown{\Wf'}\))
      and fully contained in \(Z\times (0,\infty)\) otherwise.
      In either case, we have that it is \(\Aa\)-local.
      In conclusion, we have proved that the lower horizontal inclusion
      in the following commutative square (by \Cref{rem:SomeHyprefsCompose}) of presieve inclusions
      is strongly \(\Aa\)-local,
      hence \(\Aa\)-local by \Cref{lem:strongly-local-local}.
      \begin{equation}
        \label{eq:W'-W-to-sieves}
        \cdsquareOpt
        {\Wf'}
        {\Wf}
        {\cdown{\UfwithES\cap\Wf'}}
        {\cdown{\UfwithES\cap\Wf}}
        {hookrightarrow}
        {hookrightarrow}
        {hookrightarrow}
        {hookrightarrow}
      \end{equation}      
      The two vertical inclusions are also \(\Aa\)-local
      because \(\Wf'\) and \(\Wf\) are presieves
      so that we have inverse hyperrefinements
      \(\UfwithES\cap\cdown{\Wf'}\hyprefines{\Wf'}\)
      and
      \(\UfwithES\cap\cdown{\Wf}\hyprefines{\Wf}\)
      (as in \Cref{rem:iso-of-refinements}).
      \Cref{lem:hypref-functorial} ensures that the diagram stays commutative after applying \(\Aa\), 
      thus we conclude that the top horizontal inclusion is \(\Aa\)-local,
      as desired.
    \end{itemize}
    In total, we are in the situation where the left inclusion in
    \(\Wf'\hookrightarrow \Wf\hyprefines U \)
    is \(\Aa\)-local.
    The composite is a Weiss presieve that is fully contained in
    \(p^{-1}\cdown{[0,\infty)}\), hence also \(\Aa\)-local.
    We conclude that \(\Wf\hyprefines U\) is \(\Aa\)-local, as desired.
  \end{itemize}

  Summarizing, we have now established the pullback square
  \begin{equation}
    \cdsquare[pb]
    {\FactCstr{\UfwithES}}
    {\FactCstr{Z\times(0,\infty)}}
    {\FactCstr{[0,\infty)}}
    {\FactCstr{(0,\infty)}}
    {} {\pf{p}} {\pf{p}} {}
  \end{equation}

  To complete the proof we have to show that the inclusion
  \(\UfwithES\hookrightarrow \open{\topcone Z}\)
  of posets
  induces an equivalence on the categories of constructible factorization algebras.
  Note that \(\disjunComp{\Uf} = \disjunComp{(\UfwithES)}\) because the disjoint union completion automatically adds in the empty set, 
  and consider the factorization 
  \(\UfwithES \hookrightarrow \cdisj{\Uf}\hookrightarrow \open{\topcone Z}\). 
  We will show separately that both of the induced restriction functors 
  \begin{equation}
    \FactCstr{\topcone Z}
    \longrightarrow
    \FactCstr{\cdisj{\Uf}}
    \longrightarrow
    \FactCstr{\UfwithES}
  \end{equation}
  are equivalences.
  
  \begin{itemize}
  \item
    For the left equivalence we observe that
    \(\cdisj{\Uf}\) is a factorizing basis of \(\topcone Z\)
    so that \Cref{lemma:FactBasisExtends}
    yields the equivalence
    \begin{equation}
      \Fact{\topcone Z}
      \xrightarrow{\simeq}
      \Fact{\cdisj{\Uf}}.
    \end{equation}
    It remains to show that this equivalence restricts to
    an equivalence on constructible factorization algebras,
    i.e., that a factorization algebra \(\Aa\) on \(\topcone Z\)
    is constructible if it is constructible
    on \(\cdisj{\Uf}\), or equivalently,
    on \(Z\times (0,\infty)\) and on \(p^{-1}\cdown{[0,\infty)}\).
    For this we apply \Cref{cor:constructible-away-from-0}
    which in fact says that it suffices
    to check constructibility on
    \(\topcone Z \setminus \{0\} = Z\times (0,\infty)\).
    Note that here we use compactness of \(Z\)
    because otherwise \(M=\topcone Z\) would not be a conical manifold.
  \item
    For the right equivalence we observe that
    the \posetofopens{} \(\UfwithES\) is decomposable so 
    \Cref{cor:add-disjun-constr}
    yields the equivalence
    \begin{equation}
      \AlgMCstr{\cdisj{\Uf}}
      \xrightarrow{\simeq}
      \AlgMCstr{\UfwithES}.
    \end{equation}
    It remains to show that this equivalence restricts to an equivalence on
    constructible factorization algebras, i.e.,
    that a multiplicative constructible algebra \(\Aa\) on \(\cdisj{\Uf}\)
    is a Weiss cosheaf if it is a Weiss cosheaf on \(\UfwithES\).
    Choose a factorizing disk-basis \(\Bf\) of \(Z\),
    witnessing that it has enough good disks;
    without loss of generality we may assume that \(\Bf\) is decomposable.
    Write \(\Df\subseteq \UfwithES\) for the \posetofopens{}
    consisting of
    \begin{itemize}
    \item
      the standard cones \(\topcone[t]Z\), for \(0<t\leq \infty\)
    \item
      disjoint unions of disks \(D\times (a,b)\),
      where \(D\in \Bf\) and \(0<a<b\leq \infty\).
    \end{itemize}
    Observe that \(\cdisj{\Df}\) is a decomposable factorizing disk-basis
    of \(\topcone Z\),
    so that we can invoke \Cref{cor:constructible=>Weiss}:
    to show that \(\Aa\) is a Weiss cosheaf on \(\cdisj{\Uf}\)
    it suffices to prove that it is a left Kan extension from \(\cdisj{\Df}\),
    i.e.,
    that for each \(U\in \cdisj{\Uf}\),
    the presieve \(\cdisj{\Df}\cap \cdown{U} \to U\) is \(\Aa\)-local.
    If we write \(U=U'\disjun U''\), where \(U'\) is the component of the cone point,
    then by multiplicativity of \(\Aa\)
    we can write the desired comparison map as the composite
    \begin{equation}
      \label{eq:cone-LKE-decompose-U'-U''}
      \Aa(\cdisj{\Df}\cap \cdown U)
      \simeq
      \Aa(\cdisj{\Df}\cap \cdown{U'})\otimes \Aa(\cdisj{\Df}\cap\cdown U'')
      \to
      \Aa(U')\otimes\Aa(U'')
      \simeq
      \Aa(U),
    \end{equation}
    where the first equivalence holds because \(\cdisj{\Df}\)
    is decomposable and closed under disjoint unions,
    so that we have the isomorphism
    \({-\disjun-}\colon \overcat{\cdisj{\Df}}{U'}\times \overcat{\cdisj{\Df}}{U''}
    \xrightarrow{\cong} \overcat{\cdisj{\Df}}{U}\).
    Hence it suffices to show that the two \(\otimes\)-components
    of the middle arrow in \eqref{eq:cone-LKE-decompose-U'-U''} are equivalences.
    \begin{itemize}
    \item
      For the \(U''\)-component observe that
      \(\cdisj{\Df}\cap U''\hyprefines U''\) is a Weiss presieve
      fully contained in \(\cdown{(Z\times (0,\infty))} \subset \UfwithES\),
      hence \(\Aa\)-local.
    \item
      For the \(U'\)-component,
      observe that \(U'\) must be a standard cone \(\topcone[t]{Z}\)
      for some \(t\), because these are the only connected opens
      in \(\cdisj{\Uf}\) that contain the cone point.
      Therefore the presieve \(\cdisj{\Df}\cap \cdown{U'}\hyprefines U'\)
      is actually degenerate, hence \(\Aa\)-local.
      \qedhere
    \end{itemize}
  \end{itemize}
\end{proof}

\subsubsection*{Proof of \Cref{lem:cone-operad-pushout}}
% \label{app:operad-pushout}
% local macros for this subsection

\newcommand{\UO}{U^\otimes}
\newcommand{\HO}{H^\otimes}
\newcommand{\NO}{N^\otimes}
\newcommand{\OO}{O^\otimes}

\newcommand{\Qleq}[1][n,d]{Q^\leq{(#1)}}
\newcommand{\Qless}[1][n,d]{Q^{<}{(#1)}}

\newcommand\spawndary[1]{\lambda^<(#1)}

The rest of the section is devoted to the proof of
\Cref{lem:cone-operad-pushout}
using some explicit combinatorics of dendroidal sets.
More precisely, we formulate the following proposition,
which captures the features of the situation
that make \Cref{lem:cone-operad-pushout} work.

The reader is invited to keep in mind the setting of \Cref{cstr:cone-operads},
where we consider the following full suboperads of opens in a cone:

\begin{itemize}
  \item The colors in \(\OO\) are those \emph{horizontal} opens that contain the cone point \(0\).
  \item The colors of \(\HO\) are the \emph{horizontal} opens.
  \item The colors of \(\NO\) are the opens that do \emph{not} contain the cone point.
  \item The colors of \(\UO\) are those that are horizontal or do not contain the cone point.
\end{itemize}
In all cases we exclude the empty set as a color to ensure that
condition~\ref{it:UO-ops-different-inputs} below is met.
The key insight towards \Cref{lem:cone-operad-pushout}
explains why condition~\ref{it:UO-factors-OO} holds.
All other conditions are immediate from the definition.

\begin{prop}
  \label{prop:operad-pushout}
  Let \(\UO\) be an operad with two full suboperads \(\NO\) and \(\HO\).
  Let us write \(\OO\coloneqq\HO\setminus(\NO\cap\HO)\),
  by which we mean the full suboperad of \(\HO\)
  spanned by those colors that do not lie in \(\NO\).
  We make the following assumptions:

  \begin{enumerate}[ref=(A\arabic*),label=(A\arabic*)]
  \item
    Every color of \(\UO\) lies in \(\HO\) or in \(\NO\),
    i.e., \(\UO=\HO\cup\NO\).
  \item
    \label{it:no-identity-composites}
    No color of \(\OO\) has any non-trivial automorphisms
    and the composite of two non-invertible \(1\)-ary operations is always non-invertible.
  \item
    The operad \(\UO\) is unital, i.e.,
    every color \(x\) admits a unique 0-ary operation \(\emptyset_x\colon ()\to x\).
  \item 
    \label{it:OO-ops-0-1}
    For every operation in \((x_1,\dots,x_n)\to x\) in \(\UO\),
    at most one input color \(x_i\) lies in \(\OO\).
    If an input color \(x_i\) lies in \(\OO\), then so does the output color \(x\).
  \item
    \label{it:UO-ops-different-inputs}
    For every operation \((x_1,\dots,x_n)\to x\) in \(\UO\),
    the input colors are pairwise non-isomorphic.
  \item
    \label{it:UO-factors-OO}
    Let \(\iota\colon \overline{o}\to o'\) be an \(n\)-ary operation in \(\OO\)
    (this implies \(n=0\) or \(n=1\) by Assumption~\ref{it:OO-ops-0-1}).
    We assume that we have a \emph{complement operation} for \(\iota\),
    by which we mean an (\(n+1\))-ary operation
    \(a_\iota\colon (\overline{o},o'')\to o'\) in \(\HO\)
    with the following properties: 
    \begin{itemize}
    \item
      The color \(o''\) lies in \(\NO\) (hence in \(\NO\cap \HO\)).
    \item
      We have \(\iota= a_\iota\circ(\overline{o}, \emptyset_{o''})\).
    \item
      Let \(f\colon (\overline{o},x_1\dots,x_n)\to o'\) be an operation in \(\UO\)
      with \(x_1,\dots,x_n\in \NO\) and
      with \(\iota=f \circ(\overline{o},\emptyset_{x_1},\dots,\emptyset_{x_n})\).
      Then there is a unique operation
      \(b_f\colon (x_1,\dots,x_n)\to o''\) in \(\NO\)
      with \(f = a_\iota\circ (\overline{o}, b_f)\).
    \end{itemize}
  \end{enumerate}
  Then the pullback square of operads
  \begin{equation}
    \label{eq:pushout-operads}
    \cdsquareNA[pb]
    {\NO\cap\HO}
    {\HO}
    {\NO}
    {\UO}
  \end{equation}
  is a pushout square of operads and even of \infy-operads.
\end{prop}

\begin{rem}
  In the setting of Assumption~\ref{it:UO-factors-OO}
  the complement operation \(a_\iota\) is uniquely determined up to unique isomorphism.
  Thus we may say that \(o''\) is the complement of \(\iota\),
  and we write \(o'-\overline{o}\coloneqq o''\).

  It might be tempting to believe that when \(\overline{o}=()\) is the empty tuple,
  the complement \(o'-()\) is just \(o'\);
  this is never the case because the latter lies in \(\OO\)
  while the former lies in \(\NO\).
  Explicitly, in the context of \Cref{lem:cone-operad-pushout}
  we have \(U'-()=U'\setminus \{0\}\neq U'\) whenever \(U'\in \open{\topcone{Z}}\)
  is any horizontal open that contains the cone point \(0\).
\end{rem}

From now on, we assume that we are in the setting of \Cref{prop:operad-pushout}.
We make some preliminary observations and definitions.

First, we can lift the square \eqref{eq:pushout-operads}
to a square in the \((1,1)\)-category \(\strictOpd\) of operads and functors,
where everything commutes strictly and there are no non-trivial natural isomorphisms.
Moreover, we may assume that
the operads \(\UO\), \(\HO\), \(\NO\), \(\NO\cap\HO\) and \(\OO\)
are skeletal,
and that all the inclusion functors are injective on objects\footnote{
  Note that an operad being skeletal or a functor being injective on objects
  is not a property that is invariant under equivalence,
  so it is meaningless for objects or morphisms in \(\Opd\).
  It only becomes meaningful after we lift our diagram
  from the theory of operads and equivalences
  to the theory of operads and isomorphisms.
};
in particular the complement operations of \ref{it:UO-factors-OO}
are unique on the nose\footnote{
  again, a statement that is meaningless before lifting to \(\strictOpd\)
}.
Indeed, to achieve all that we may simply choose a set of representatives
for the isomorphism classes of the colors of \(\UO\);
the other skeletal operads are then just given by appropriate subsets.%
\footnote{
  Note that in the setting of \Cref{sec:link-decomp},
  the operads in question
  are already presented as skeletal operads in \(\strictOpd\)
  since they arise from explicit inclusions of \posetsofopens{}.
  So in that case there is nothing that needs to be done.
}

By applying the dendroidal nerve \(\ND{}\colon\strictOpd\to\dSets\),
we may view all these operads as dendroidal sets;
going forward we will do so
without explicitly writing the dendroidal nerve \(\ND{}\) every time.

Explicitly, a dendrex of \(\UO\)
is a tree (in the sense of \cite{MW07}) where each edge is decorated with a color of \(\UO\)
and each vertex is decorated with an appropriate operation.
Such a dendrex is degenerate precisely if some \(1\)-ary operation
is decorated with an identity.
Since we assume that \(\UO\) is skeletal,
we may reformulate Assumption~\ref{it:no-identity-composites}
as saying that the composite of non-identities is always a non-identity,
which translates to saying that every face of a non-degenerate dendrex is non-degenerate.
\begin{defn}
  \label{def:degeneration}
  For every dendrex \(F\) we denote by \(m(F)\)
  the dendrex obtained from \(F\)
  by removing all \(1\)-ary operations labelled by identities.
  This means that we have \(m(F)=F\) if and only if \(F\) is non-degenerate;
  otherwise \(m(F)\) is the unique non-degenerate dendrex of which \(F\)
  is a degeneration.
\end{defn}
Assumption \ref{it:UO-ops-different-inputs}
implies that for every dendrex \(F\in\UO(T)\),
all of its input colors are pairwise different.
In particular, any non-identity automorphism of the tree \(T\)
will permute these input colors non-trivially;
hence it cannot fix the dendrex \(F\).
In other words, the action of \(\Aut(T)\) on \(\UO(T)\) is free,
which means that the dendroidal set \(\UO\) is cofibrant.
Therefore the same is true for the dendroidal subsets \(\NO\), \(\OO\) and \(\HO\);
moreover, all inclusions between them are cofibrations.
(Recall that the cofibrations of dendroidal sets are precisely the normal monomorphisms,
i.e., those monomorphisms \(A\to B\) where for all trees \(T\),
the automorphism group \(\Aut(T)\) acts freely on \(B(T)\setminus A(T)\).)

\begin{obs}
  Due to condition \ref{it:OO-ops-0-1},
  every dendrex in \(\UO\) is of one of the following types:

  \begin{itemize}
  \item
    A dendrex of \(\NO\), which means that no color in \(\OO\) appears.
  \item
    A tree of the form given in \Cref{eqref:dendrex-geq-1}
    \begin{figure}
      \includegraphics*[scale=0.8]{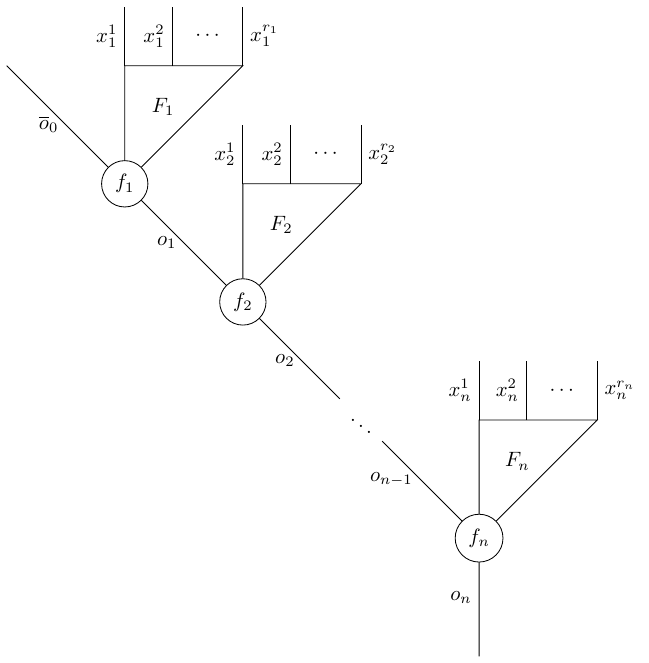}
      \caption{
        An illustration of a tree of shape \((n, d)\).
      }
      \label{eqref:dendrex-geq-1}
    \end{figure}
    with \(n\geq 0\) and the following conventions:
    \begin{itemize}
    \item 
      Each \(\overline{o}_{i}\) is always just a single color \(\overline{o}_{i}=o_i\in\OO_{1_+}\),
      except for \(0=i<n\)
      where we additionally allow the possibility \(\overline{o}_{0}=()\)
      (the empty tuple).
      The latter case happens precisely if the total dendrex has
      no input colors that lie in \(\OO\);
      pictorially we would then just not draw the upper-left-most edge.
    \item
      Each \((\overline{o}_{i-1},\overline{x}_{i})\to o_i\)
      is a tree obtained by grafting a forest \(F_i\) in \(\NO\)
      with roots \(y_1^1,\dots,y_i^{l_i}\)
      onto an operation
      \(f_i\colon (\overline{o}_{i-1},y_i^1,\dots,y_i^{l_i})\to o_i\)
      which definitely does not lie in \(\NO\)
      (since \(o_i\in\OO\)),
      but might also not lie in \(\HO\) if some \(y_i^j\notin \HO\).
    \end{itemize}
    We say that such a tree is of \emph{shape} \((n,d)\),
    where \(d=\sum_{i=1}^n|F_i|\) is the total number of vertices across all forests \(F_i\).
  \end{itemize}
\end{obs}
\begin{rem}
  By convention we usually denote an element \(\overline{y}\in \UO\) with an overline,
  but drop the overline and just write \(\overline{y}=y\) whenever we know for sure
  that it lies over \(1_+\).
  For example, in \Cref{eqref:dendrex-geq-1} we may write \(o_i\) for \(i\geq 1\),
  but must write \(\overline{o_i}\) when we do not know \(i\),
  thus allowing for the case \(\overline{o_0}=()\).
\end{rem}

Let \(F\) be a dendrex of \(\UO\) that does not lien in \(\NO\).
Consider one of the operations
\(f_i\colon (\overline{o}_{i-1},\overline{y}_{i})\to o_i\) described above.
By \ref{it:UO-factors-OO}, this operation factors uniquely as
\begin{equation}
  \label{eq:tree-special-factorization}
  (\overline{o}_{i-1},\overline{y}_{i})\xrightarrow{(\overline{o}_{i-1},b_{f_i})} (\overline{o}_{i-1},o_i-\overline{o}_{i-1})\xrightarrow{a_\iota} o_i
\end{equation}
where \(\iota=f_i\circ(\overline{o}_{i-1},\emptyset_{\overline{y}_{i}}): \overline{o}_{i-1}\to o_i\)
and \(a_\iota\) is the chosen complement operation.

\begin{defn}
  We say that the dendrex \(F\) is \emph{special at position} \(i\)
  if the resulting
  operation \(b_{f_i}\colon \overline{y}_{i}\to o_i-\overline{o}_{i-1}\) is an identity,
  i.e., if the forest \(F_i\) is actually a tree with root
  \(o_i-\overline{o}_{i-1}\) and
  \(f_i=a_\iota\) is the chosen complement operation.
  In this case we write \(a_i=a_\iota\).
  We say that \(F\) is \emph{special} if it is special at all \(1\leq i\leq n\).
\end{defn}

\begin{constr}
  Denote by \(R(F)\) the dendrex obtained from \(F\)
  by replacing the corolla
  \(f_i\colon(\overline{o}_{i-1},\overline{y}_{i})\to o_i\)
  with the 2-vertex tree \eqref{eq:tree-special-factorization}
  for each position \(i\) at which \(F\) is not special;
  by construction, \(R(F)\) is now special.
  Tautologically, we have \(R(F)=F\) if and only if \(F\) is already special.
\end{constr}

\begin{defn}
  We now define the precise \emph{spawn time} of any dendrex \(F\) of \(\UO\) as follows:
  \begin{itemize}
  \item
    If \(F\) fully lies in \(\NO\) or fully lies in \(\HO\)
    then we say that it spawns at time \(0\).
  \item
    If \(F\) is non-degenerate and special of shape \((n,d)\),
    we say that \(F\) spawns at time \((n,d)\).
  \item
    Otherwise, we say that \(F\) spawns at the same time as \(mR(F)=Rm(F)\).
  \end{itemize}
\end{defn}

The set of possible spawn times of a dendrex
is \(L\coloneqq\{0\}\amalg(\naturals_{\geq 1}\times\naturals_{\geq 1})\),
which we view well-ordered lexicographically:
we always have $0<(n,d)$;
and $(n,d)<(n',d')$ iff either $n<n'$ or both $n=n'$ and $d<d'$.
Note that every special tree of shape \((n,0)\) or \((0,d)\) is fully contained in \(\HO\);
which is why those spawn times do not appear.

We observe the following:
\begin{itemize}
\item
  A special dendrex of shape \((n,d)\) spawns at or before time \((n,d)\).
\item
  A dendrex \(F\) always spawns at the same time as \(R(F)\).
\item
  An arbitrary dendrex \(F\) of shape \((n,d)\) spawns at or before time \((n,d+n)\),
  since that is the largest possible shape of \(R(F)\).
\end{itemize}
Moreover, for each dendrex \(F\) of \(\UO\)
\begin{itemize}
\item
  applying a degeneracy does not change the spawn time,
\item
  applying a face map does not increase the spawn time.
\end{itemize}
Therefore for each \(\alpha\in L \), the collections
\begin{equation}
  \Qleq[\alpha]\coloneqq \{\text{dendrices of }\UO\text{ that spawn at time }\leq \alpha\}
\end{equation}
are dendroidal subsets of \(\UO\);
their ascending union is \(\UO\) because every dendrex spawns at some time.
Note that by definition we have \(\Qleq[0]=\HO\amalg_{\NO\cap\HO}\NO\)
(the strict pushout of dendroidal sets).
We also abbreviate \(\Qless[\alpha]\coloneqq \bigcup_{\beta<\alpha}\Qleq[\beta]\)
consisting of all dendrices that spawn before time \(\alpha\).

\begin{constr}
  For each non-degenerate special dendrex \(F\) that spawns at time \(\alpha\)
  we define the dendroidal set \(\spawndary{F}\) as the following pullback:
  \begin{equation}
    \label{eq:pullack-spawndary}
    \cdsquareOpt[pb]
    {\spawndary{F}}
    {\Qless[\alpha]}
    {\Omega[T]}
    {\Qleq[\alpha]}
    {}{hookrightarrow}{hookrightarrow}{"F"}
  \end{equation}
\end{constr}

\begin{lemma}
  \label{lem:char-spawndary}
  Let \(F\) be a non-degenerate special dendrex that spawns at time \(\alpha>0\).
  The inclusion \(\spawndary{F}\hookrightarrow \Omega[T]\)
  is the dendroidal subset spanned by the following faces:
  \begin{itemize}
  \item
    all outer faces
  \item
    all inner faces except for those which under \(F\)
    correspond to a root of one of the trees \(F_i\) in the picture
    \eqref{eqref:dendrex-geq-1}
    (recall that a priori each \(F_i\) is a forest, 
    but since \(F\) is special they are actually trees).
  \end{itemize}
\end{lemma}

\begin{figure}
  \centering
  \includegraphics*[scale=0.7]{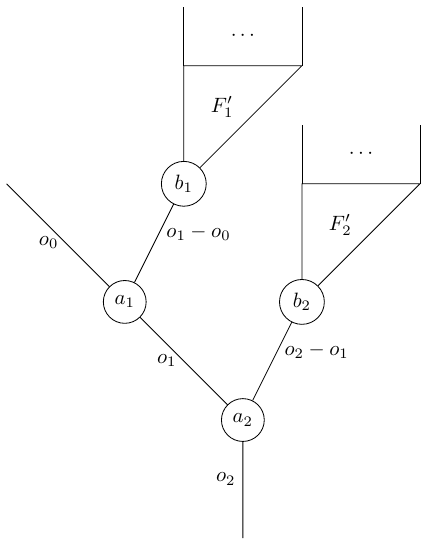}
  \quad$\leadsto$\quad
  \includegraphics*[scale=0.7]{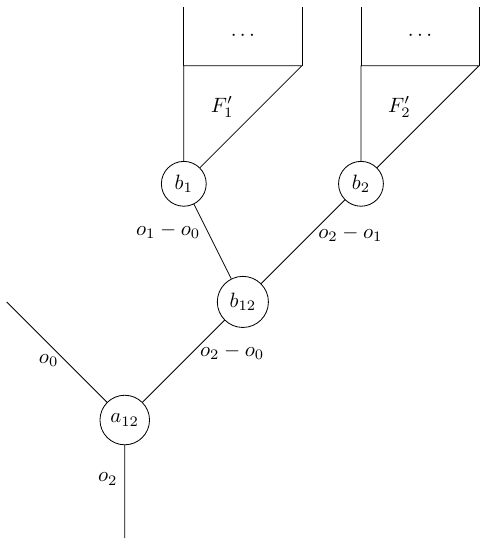}
  \caption{
    \label{fig:special-dendrex-collapse}
    On the left, a special dendrex \(F\) of shape \((2,2+d')\),
    where \(a_1\) and \(a_2\) are chosen complement operations
    and \(d'\coloneqq |F'_1|+|F'_2|\).
    After collapsing the edge decorated by \(o_1\)---%
    which means passing to the dendrex \(F\sigma\) for the appropriate inner face \(\sigma\)---%
    the tree is no longer special
    because the composition \(a_2\circ a_1\) of two complement operations
    is not a complement operation anymore.
    We make it special again by passing to the tree \(R(F\sigma)\),
    depicted on the right.
    The operation \(b_{12}\) is uniquely determined by the relation
    \(a_{12}\circ b_{12}=a_2\circ a_1\)
    due to the defining property of the complement operation \(a_{12}\).
    The tree \(R(F\sigma)\) has shape \((1,3+d')\) so it spawned before \(F\).
  }
\end{figure}

\begin{proof}
  The right vertical map in the square \eqref{eq:pullack-spawndary}
  is a levelwise inclusion, hence the same is true for the left one.
  Concretely, \(\spawndary{F}\)
  is the dendroidal subset of \(\Omega[T]\)
  spanned by those non-degenerate dendrices \(\sigma\colon \Omega[S]\to\Omega[T]\)
  where the dendrex \(F\sigma\) of \(\UO\) spawns before time \(\alpha=(n,d)\).

  First let \(\sigma\) be an outer face
  which means that the dendrex \(F\sigma\) is obtained from \(F\)
  by omitting an external vertex \(v\).
  There are the following possibilities:
  \begin{itemize}
  \item
    The vertex \(v\) is a leaf-vertex of one of the \(F_i\).
  \item
    The tree \(F_1\) consists of a single edge \(o_1-\overline{o}_{0}\)
    and \(v\) is the vertex \((\overline{o}_{0},o_1-\overline{o}_{0})\to o_1\).
  \item
    The tree \(F_n\) consists of a single edge \(o_n-\overline{o}_{n-1}\)
    and \(v\) is the vertex \((\overline{o}_{n-1},o_n-\overline{o}_{n-1})\to o_n\).
  \end{itemize}
  (Note that in the last two cases we must have \(n\geq 2\),
  because otherwise \(F\) would lie in \(\HO\),
  which is impossible since \(F\) spawns at time \(\alpha>0\).)
  Then the dendrex \(F\sigma\) is again special
  and is of shape \((n,d-1)\), \((n-1,d)\) or \((n-1,d)\), respectively.
  Hence it already spawned before time \((n,d)\).

  Next, let us consider inner faces,
  which means that \(F\sigma\) is obtained from \(F\)
  by contracting an inner edge \(e\).
  Let us consider the following possibilities:
  \begin{itemize}
  \item
    The edge \(e\) is a non-root edge of one of the \(F_i\).
    In this case the dendrex \(F\sigma\) is special of shape \((n,d-1)\)
    so that is the latest possible spawn time for it.
  \item
    The edge \(e\) is one of the edges labelled \(o_i\), for \(1\leq i < n\)
    (hence \(n\geq 2\)).
    In this case the dendrex \(F\sigma\) spawned at or before time
    \((n-1,d+1)\),
    because that is the shape of \(R(F\sigma)\).
    Indeed, contracting the edge reduces \(n\) by one;
    the resulting tree is not special at position \(i\) anymore,
    requiring \(R(F)\) to add one new vertex.
    See \Cref{fig:special-dendrex-collapse} for an illustration of this case.
  \end{itemize}
  In either of those cases we see that \(F\sigma\)
  has spawned before time \(\alpha=(n,d)\).

  Finally, let \(\sigma\colon \Omega[S]\to\Omega[T]\) be the intersection
  of all the inner faces of \(T\) that have not yet been listed above.
  Then the dendrex \(F\sigma\) is obtained from \(F\)
  by collapsing the root edge of each of the trees \(F_i\).
  Note that \(F\sigma\) cannot be special at any position \(i\)
  because otherwise the tree \(F_i\) would be of the form
  \begin{equation}
    ... \to o_i-\overline{o}_{i-1} \xrightarrow{v} o_i-\overline{o}_{i-1}
  \end{equation}
  with \(a_i\circ (\overline{o}_{i-1},v)=a_i\),
  which by the uniqueness in Assumption~\ref{it:UO-factors-OO}
  would imply that \(v=\id\) making \(F\) degenerate, which it is not.
  It follows that \(R(F\sigma)=F\),
  so that \(F\sigma\) spawns at the same time as \(F\) and not before.

  For every subdendrex \(\sigma\colon \Omega[S]\to\Omega[T]\) we have shown: 
  \begin{itemize}
  \item
    If \(\sigma\) is contained in one of the listed faces,
    then it lies in \(\spawndary{F}\).
  \item
    If \(\sigma\) contains the intersection of the unlisted faces
    then it does not lie in \(\spawndary{F}\).
  \end{itemize}
  This means that \(\spawndary{F}\)
  consists precisely of the union of the listed faces, as desired.
\end{proof}

\begin{lemma}
  \label{lem:unique-non-deg-special}
  Let \(G\colon \Omega[S]\to\UO\) be a dendrex that spawns at time \(\alpha>0\).
  Then there is a factorization
  \begin{equation}
    G\colon \Omega[S]\xrightarrow{\sigma}\Omega[T]\xrightarrow{F}\UO
  \end{equation}
  where \(F\colon \Omega[T]\to\UO\) is a special non-degenerate dendrex
  that spawns at time \(\alpha\).
  Moreover, the pair \((F, \sigma)\)
  is unique up to the tautological action of \(\Aut(T)\)
  given by \(\tau . (F, \sigma)\coloneqq (F\tau,\tau^{-1}\sigma)\).
\end{lemma}

\begin{proof}
  First, assume we had such a factorization \((F,\sigma)\) with \(F\sigma=G\).
  Factor \(\sigma=\tau\rho\), where \(\rho\) consists of degeneracies
  and \(\tau\) of face inclusions.
  Since \(F\) must be non-degenerate so must \(F\tau\);
  hence we must have \(F\tau=mG\) up to tree-automorphism, which also uniquely determines \(\rho\).
  Note that \(\tau\) cannot lie in \(\spawndary{F}\),
  since otherwise \(mG=F\tau\) would have already spawned before time \(\alpha\).
  Hence by \Cref{lem:char-spawndary},
  \(mG\) arises from \(F\) precisely by collapsing
  the root edges of the trees \(F_i\) in all positions \(i\) where \(G\) is not special.
  Therefore, up to tree-automorphism, \(F\) is uniquely determined to be \(F=Rm(G)=mR(G)\)
  and \(\tau\) corresponds to the aforementioned edge collapses.
\end{proof}

We also make use of the following result.

\begin{lemma}[\cite{MW09}, Lemma~5.1]
  \label{lem:CM-anodyne-outer-not-all-inner}
  Let \(\{\sigma_j\colon\Omega[T_j]\to\Omega[T]\}_j\)
  be a set of faces of \(\Omega[T]\) which
  contains all outer faces 
  and does \emph{not} contain all inner faces.
  Then the inclusion \(\bigcup_j\Omega[T_j]\to \Omega[T]\) is anodyne.
\end{lemma}

\begin{proof}[Proof of \Cref{prop:operad-pushout}]
  We claim that the square \eqref{eq:pushout-operads}
  is a homotopy pushout in the Cisinski-Moerdijk model structure~\cite{CM11},
  thus yielding the desired pushout of \infy-operads
  (and a fortiori, of operads).

  Fix a positive spawn time \(\alpha=(n,d)>0\).
  We make the following claims.

  \begin{enumerate}
  \item
    For each non-degenerate special dendrex \(F\) that spawns at time \(\alpha\),
    the inclusion \(\spawndary{F}\to\Omega[T]\)
    is an anodyne extension.
  \item
    Let \(\{F_j\colon \Omega[T_j]\to \UO\}_{j\in J}\)
    be a set of representatives up to tree-automorphism
    of all non-degenerate special dendrices of \(\UO\)
    that spawn at time \(\alpha\).
    Then we have a pushout square
    \begin{equation}
      \label{eq:attaching-square-Qleq}
      \cdsquare[po]
      {\coprod\limits_{j \in J}\spawndary{F_j}}
      {\Qless[\alpha]}
      {\coprod\limits_{j \in J}\Omega[T_j]}
      {\Qleq[\alpha]}
      {}{}{}{(F_j)}
    \end{equation}
  \end{enumerate}

  Jointly these claims imply that each inclusion
  \(\Qless[\alpha]\to\Qleq[\alpha]\)
  is anodyne, hence also their transfinite composition
  \(\HO\amalg_{\NO\cap\HO}\NO=\Qleq[0]\to \UO\).
  Since the vertical inclusions in the square 
  \eqref{eq:pushout-operads} are cofibrations,
  this implies that this square is a homotopy pushout, as desired.

  So it remains to show the claims:
  \begin{enumerate}
  \item
    The first claim follows directly from \Cref{lem:char-spawndary}
    and \Cref{lem:CM-anodyne-outer-not-all-inner}.
  \item
    To show that the square \eqref{eq:attaching-square-Qleq} is a pushout of dendroidal sets,
    we have to show that for each tree \(S\) the comparison map
    \begin{equation}
      \label{eq:gap-map-attaching-Qleq}
      \left(
        \Qless[\alpha]
        \amalg_{\coprod_{j}\spawndary{F_j}}
        \coprod_{j}\Omega[T_j]
      \right)
      (S)
      \to \Qleq[\alpha](S)
    \end{equation}
    is a bijection.
    For this, let \(G\colon \Omega[S]\to \Qleq[\alpha]\) be an arbitrary dendrex;
    we have to show that there is a unique dendrex in the pushout that hits it.
    We distinguish the following two cases:
    \begin{itemize}
    \item
      The dendrex \(G\) spawned before time \(\alpha\).
      In this case, \(G\) already lies in the component \(\Qless[\alpha]\).
      Moreover, if \(\sigma\colon \Omega[S]\to \Omega[T_j]\)
      is a dendrex in the other component of the pushout
      that hits \(G\), i.e., \(F_j\sigma=G\),
      then the pair \((G,\sigma)\) lies in \(\spawndary{F_j}\)
      and witnesses an equality \([\sigma]=[G]\) in the pushout.
    \item
      The dendrex \(G\) spawns at time \(\alpha\).
      In this case the only elements in the pushout that can possibly hit \(G\)
      come from the component \(\coprod_j\Omega[T_j]\).
      \Cref{lem:unique-non-deg-special} states precisely
      that this component has a unique dendrex
      \(\sigma \colon \Omega[S]\to \Omega[T_i]\) with \(G=F_i\sigma\).
      \qedhere
    \end{itemize}
  \end{enumerate}
\end{proof}

\subsection{Factorization algebras from factorization homology}

We show that fixing a smooth manifold $M$
(possibly stratified, with framings or other tangential structures),
taking factorization homology locally gives a factorization algebra.
This proof is a relatively straightforward consequence
of the results established in \Cref{sec:external-tools}
and does not rely on our gluing results.
In the locally constant case a proof relying on gluing of locally constant factorization algebras was given in \cite{GTZ}.

In the remainder of this section,
we adopt terminology similar to that of \cite{AFT-fh-stratified}, in particular:
\begin{itemize}
\item $\targetcatOT$ is a symmetric monoidal \(\otimes\)-presentable \infy-category;
\item $\cB$ is a notion of tangential structure on smooth conical manifolds,
  encoded as a presheaf on the \infy-category of basics;
\item $\mfldB$ and $\mmfldB$ are the category and \infy-category of $\cB$-manifolds
  (i.e.\ smooth conical manifolds equipped with a tangential structure of type \(\cB\)), respectively;
  the morphisms are open embeddings.
\item
  \(\diskB\subset \mfldB\) and \(\ddiskB\subset \mmfldB\)
  are the full sub-(\infy-)categories of multidisks with tangential
  structure of type \(\cB\). 
\item $\MfldB$ and $\MMfldB$ are the symmetric monoidal category and \infy-category of $\cB$-manifolds, respectively. 
\end{itemize}
If $\cB$ is the trivial presheaf (encoding no tangential structure)
then $\mfldB=\Snglr$ and $\mmfldB = \SSnglr$ (even as symmetric monoidal categories).
Two other cases to keep in mind are orientations and framings, defined compatibly with the stratifications.

We first recall the definition of factorization homology and a main computational tool
(\cite{AFT-fh-stratified}  Theorem 2.15, Lemmas 2.16 \& 2.17).
 The {\em  factorization homology functor $\int_{-} $} is a left adjoint to the restriction
 \begin{equation}
 \Fun^\otimes ( \DDiskB, \targetcatOT) \longleftarrow \Fun^\otimes ( \MMfldB, \targetcatOT) \, ,
  \end{equation}
 which, since $\targetcat$ is $\otimes$-sifted cocomplete, can be computed as the (ordinary) left Kan extension
 \begin{equation}
 \begin{tikzcd}
 \ddiskB \arrow[hook]{d} \arrow{rr}{\Aa} && \targetcat\ \,,\\
 \mmfldB \arrow[dashed, swap]{urr}{\int_{-} \Aa}
 \end{tikzcd}
\end{equation}
which is automatically an operadic left Kan extension;
the resulting functor \(\int_{-}\Aa\) is automatically symmetric monoidal.
 
The structure of a $\cB$-manifold is uniquely inherited along open embeddings.
Thus forgetting the tangential structure gives an equivalence
\begin{equation}
  \overcat{\mfldB}{M}\xrightarrow{\simeq}\overcat{\Snglr}{M}=\open{M}
  \quad
  \text{and}
  \quad
  \overcat{\diskB}{M}\xrightarrow{\simeq}\overcat{\mathrm{Disks}}{M}=\disk{M}
\end{equation}
of posets\footnote{Recall \Cref{not:snglr-cat}.} ;
see also \cite[Observation 2.7]{AFT-fh-stratified}.
Then we have an induced map of \infy-operads
\begin{equation}
\iota_M\colon \Open{M}  \to \MfldB \to \MMfldB \,.
\end{equation}
Hence, given a symmetric monoidal functor
$\Aa\colon \DDiskB\to \targetcatOT$, we obtain by precomposing,
a prefactorization algebra 
\begin{equation}
  \Fa_\Aa \colon \Open{M}   \to \MfldB \to
  \MMfldB \xrightarrow{\int_{-}\Aa} \targetcatOT\, .
\end{equation}

\begin{thm}\label{thm:fact_hom_is_fact_alg}
The prefactorization algebra $\Fa_\Aa$ is a constructible factorization algebra. Hence, for $M$ a $\cB$-manifold, the composite
\begin{equation}
  \Fa_{-}\colon
  \Fun^\otimes  ( \DDiskB, \targetcatOT)
  \xrightarrow{\int_{-}} \Fun^\otimes ( \MMfldB, \targetcatOT)
  \xrightarrow{-\circ \iota_M} \Alg[\targetcat]{\open{M}} 
\end{equation}
factors through $\FactCstr[\targetcat]{M}$.
\end{thm}

\begin{proof}
To see that $\Fa_\Aa$ is constructible, note that inclusions of disks in $M$ which are abstractly isomorphic become equivalences in $\mmfldB$, hence $\Fa_\Aa$ evaluates to an equivalence. Next, since factorization homology is symmetric monoidal, $\Fa_\Aa$ is multiplicative.

It remains to check that $\Fa_\Aa$ is a left Kan extension of its restriction to
$\disk{M}\simeq\overcat{\diskB}{M}$;
indeed, the Weiss condition will then be automatic by
\Cref{cor:constructible=>Weiss} applied to the tautological disk-basis 
\(\Bf\coloneqq \disk{M}\) of \(M\).
Consider the following diagram:
\begin{equation}
 \begin{tikzcd}[bezier bounding box]
   \disk{M} \arrow{r}\arrow{d}
   & \diskB \arrow{d}\arrow{r}{L}&
   \ddiskB\arrow{d} \arrow{r}{\Aa} & \targetcat\\
   \open{M} \arrow[controls={+(1.5,-0.5) and +(0,-4)}]{rrru}{\Fa_\Aa} \arrow{r}
   & \mfldB \arrow{r}{\ell}& \mmfldB
   \arrow[dashed, "\int_{-} \Aa" description]{ru}
 \end{tikzcd}
\end{equation}
We need to show that the induced map
\begin{equation}
  \colim_{D\in \overcat{\disk{M}}{U}} \Aa(L(D)) \longrightarrow \Fa_\Aa(U)
  = \colim_{D' \in \overcat{\ddiskB}{\ell(U)}} \Aa(D')
\end{equation}
is an equivalence. Here the right expression comes from the formula for factorization homology as a left Kan extension.

To see this, we observe that the composite map
\begin{equation}
  \disk{U} = \overcat{\disk{M}}{U} \longrightarrow \overcat{\ddiskB}{\ell(U)} \xrightarrow{\simeq} \overcat{\ddisk}{\ell(U)}
\end{equation}
is a localization by \Cref{thm:Disk-localized} for $X=U$ and hence cofinal.
\end{proof}

\newpage

%%% START APPENDIX

\appendix

\section{Category theory}
\label{app:categories}

\subsection{Quillen's Theorem B}

The following version of Quillen's theorem B is a well known fact of higher category theory;
we provide a proof for completeness and the convenience of the reader.

\begin{lemma}[Quillen B for cocartesian fibrations]
  \label{lem:Quillen-B-cocart}
  Let
  \begin{equation}
    \label{eq:pb-Quillen-B}
    \cdsquare[pb]
    {X'}{X}
    {W'} {W}
    {}{}{p}{f}
  \end{equation}
  be a cartesian square of \(\infty\)-categories
  where both vertical maps are cocartesian fibrations. 
  Assume that for each arrow \(d\to d'\) in \(W\),
  the induced map \(\classspace{X_d}\to \classspace{X_{d'}}\) is an equivalence
  of \(\infty\)-groupoids.
  Then the square \eqref{eq:pb-Quillen-B} remains cartesian
  after applying the classifying space functor \(\classspace{-}\).
\end{lemma}

\begin{proof}
  We write \(X_\bullet\colon W \to \Catinfty\)
  for the functor classifying the cocartesian fibration \(p\),
  so that we can write \(X\simeq \int_{d\in W} X_d\).
  Then the assumption of the lemma says
  that the composite functor
  \begin{equation}
    \classspace{X_\bullet}\coloneqq
    \left(
    W\xrightarrow{X_\bullet}\Catinfty\xrightarrow{\classspace{-}}\Gpdinfty
    \right)
  \end{equation}
  factors through the classifying space \(W \to \classspace{W}\)
  so that we can compute
  \begin{align}
    \classspace{X}
    &=
      \localize{X}{\mathrm{all}}
    \simeq
      \localize{
      \localize{
      \left(\int_{d\in W}X_d \right)
      }
      {\mathrm{fiberwise}}
      }{\mathrm{cocart}}
    \\
    &\overset{(a)}{\simeq}
      \localize{
      \left(
      \int_{d\in W} \classspace{X_d}\right)
      }{\mathrm{cocart}}
    \overset{(b)}{\simeq}
    \colim_{d\in W}\classspace{X_d}
    \\
    &\overset{(c)}{\simeq}
    \colim_{d\in \classspace{W}}\classspace{X_d}
    \overset{(b)}{\simeq}
    \int_{d\in \classspace{W}}\classspace{X_d}
  \end{align}
  (over \(\classspace{W}\)),
  where we used
  (a) that fiberwise localizations are compatible with unstraightening,
  (b) the explicit computation of colimits of \infy-groupoids via unstraightening
  and
  (c) that \(W\to \classspace{W}\) is cofinal.
  Together with a similar computation for \(X'_\bullet = X_{f(\bullet)}\),
  this shows that the square
  \begin{equation}
    \cdsquareNA
    {\classspace{X'}}
    {\classspace{X}}
    {\classspace{W'}}
    {\classspace{W}}
  \end{equation}
  can be identified with the square
  \begin{equation}
    \cdsquare[pb]
    {\int_{d\in \classspace{W'}}\classspace{X_{fd}}}
    {\int_{d\in \classspace{W}}\classspace{X_{d}}}
    {\classspace{W'}}
    {\classspace{W}}
    {}{}
    {}{\classspace{f}}
  \end{equation}
  which is manifestly cartesian.
\end{proof}

\subsection{A localization criterion}

The following lemma is a general criterion for detecting localizations
and underpins the computation in \Cref{sec:localizing-disk}.

We extracted this lemma as the abstract and reusable core of an argument
that was first sketched in \cite{AFT-fh-stratified} with a major gap
and then completed using a result of Mazel-Gee~\cite{Mazel-Gee}.
We learned about this argument from Berry's thesis~\cite{Eric}
and private communication with David Ayala;
it also appeared in \cite{Cepek}.

The same lemma with essentially the same proof
was also independently isolated by Arakawa
\cite[Proposition 2.27]{Arakawa}.

\begin{lemma}
  \label{lem:Berry-localization}
  Let \(F\colon D \to E\) be a functor between \(\infty\)-categories.
  Define \(W\coloneqq F^{-1}(E^\simeq)\)
  to be the wide subcategory of \(D\)
  containing those arrows which become invertible in \(E\).
  Assume that
  \begin{enumerate}
  \item
    \label{it:total-for-localization}
    the map of \(\infty\)-groupoids
    \begin{equation}
      \classspace{W} \xrightarrow{\simeq}E^\simeq
    \end{equation}
    induced by \(F\) is an equivalence and
  \item
    \label{it:slices-for-localization}
    for each \(d\in D\) the map
    \begin{equation}
      \classspace{W\times_{D}\overcat{D}{d}}
      \xrightarrow{\simeq}
      (\overcat{E}{Fd})^\simeq
    \end{equation}
    induced by \(F\) is an equivalence.
  \end{enumerate}
  Then \(F\) induces an equivalence
  \begin{equation}
    \localize{D}{W} \xrightarrow{\simeq} E.
  \end{equation}
\end{lemma}

\begin{proof}
  \newcommand\RezkDW{\mathrm{N}^{\mathrm{Rezk}}(D,W)}
  For this proof we will view \(\infty\)-categories as complete Segal simplicial spaces.
  We will consider the so-called Rezk nerve
  \(\RezkDW\) of the pair \((D,W)\)
  and show separately that the canonical structure map
  \(\RezkDW\to \localize{D}{W}\)
  and the composite
  \(\RezkDW \to\localize{D}{W} \to E \)
  are equivalences of simplicial spaces.

  We recall the definition of the Rezk-nerve as described in \cite[Section 3]{Mazel-Gee}: 
  For any category \(T\) write
  \(W(T,D)\subset \Fun(T,D)\)
  for the wide subcategory where we only allow transformations
  that are pointwise in \(W\).
  Letting \(T=[n]\in\Delta^\op\subseteq\Cat^\op\) vary and passing to classifying spaces,
  we obtain a simplicial space
  \begin{equation}
    \RezkDW\colon [n] \mapsto \classspace{W([n],D)}.
  \end{equation}
  It comes equipped with a canonical map of simplicial spaces
  \begin{equation}
    \label{eq:map-rezk-nerve}
    \RezkDW \to \localize{D}{W}\coloneqq \Fun([\bullet],\localize{D}{W})^\simeq.
  \end{equation}
  Theorem~3.8 of \cite{Mazel-Gee} states that up to Rezk-fibrant replacement,
  the Rezk nerve presents the
  \(\infty\)-categorical localization \(\localize{D}{W}\).
  In particular, the comparison map \eqref{eq:map-rezk-nerve}
  is an equivalence of spaces \emph{if} \(\RezkDW\) is already complete Segal.
  
  We claim that under the assumptions of \Cref{lem:Berry-localization},
  the simplicial space \(\RezkDW\) is indeed complete Segal
  so that \(\RezkDW\to \localize{D}{W}\) is an equivalence.
  \begin{itemize}
  \item
    For each \(n\geq 1\) we have a pullback square of categories
    \begin{equation}
      \label{eq:Segal-square-W}
      \cdsquareNA[pb]
      {W([n],D)}
      {W(\{0<1\},D)}
      {W(\{1<\dots<n\},D)}
      {W(\{1\},D)}
    \end{equation}
    in which the vertical maps are cocartesian fibrations
    classified by
    \begin{equation}
      W(\{1<\dots<n\},D) \to W(\{1\},D) = W
      \xrightarrow{ d\, \mapsto\, W\times_D\overcat{D}{d}} \Catinfty.
    \end{equation}
    For each map \(d\to d'\) in \(W\),
    we have a commutative square
    \begin{equation}
      \cdsquare
      {\classspace{W\times_D\overcat{D}{d}}}
      {\classspace{W\times_D\overcat{D}{d'}}}
      {(\overcat{E}{Fd})^\simeq}
      {(\overcat{E}{Fd'})^\simeq}
      {}{\simeq}{\simeq}{\simeq}
    \end{equation}
    in which the lower horizontal map is an equivalence
    (because \(Fd \to Fd'\) is an equivalence in \(E\))
    and the vertical maps are equivalences by the second assumption.
    It follows that the square \eqref{eq:Segal-square-W}
    satisfies the assumptions of \Cref{lem:Quillen-B-cocart}
    so that it induces the pullback square
    \begin{equation}
      \cdsquareNA[pb]
      {\RezkDW_n}
      {\RezkDW_1}
      {\RezkDW_{n-1}}
      {\RezkDW_0}
    \end{equation}
    of \(\infty\)-groupoids
    which exhibits \(\RezkDW\) as a Segal space.
  \item
    To show that \(\RezkDW\) is complete,
    let \(\alpha\colon d\to d'\) be any arrow in \(D\)
    such that \(\alpha\in \classspace{W([1],D)}\) is invertible
    (with respect to the composition in the Segal space \(\RezkDW\)).
    Since the two maps of Segal spaces
    \begin{equation}
      \RezkDW \to \localize{D}{W} \xrightarrow{F} E
    \end{equation}
    preserve invertible arrows,
    it follows that \(F\alpha\) is invertible in \(E\),
    which means that we have \(\alpha \in W\).
    But then the commutative square
    \begin{equation}
      \cdsquare d{d'}{d'}{d'}{\alpha}{\alpha}{=}{=}
    \end{equation}
    is a morphism \(\alpha \to \id_{d'}\)
    in \(W([1],D)\),
    which in \(\classspace{W([1],D)}\) becomes a homotopy from \(\alpha\)
    to a degenerate \(1\)-simplex.
    Thus we have shown that every invertible \(1\)-simplex
    of \(\RezkDW\) is equivalent to a degenerate one;
    this is precisely the completeness condition.
  \end{itemize}

  Since \(\RezkDW\) is a complete Segal space, it presents the localization
  \(\localize D W\), so we are reduced to showing that the map
  \begin{equation}
    \RezkDW \to E
  \end{equation}
  induced by \(F\) is an equivalence of complete Segal spaces.
  For this it suffices to check that it is an equivalence on \(0\) and \(1\)-simplices:
  \begin{itemize}
  \item
    On \(0\)-simplices, we have the map
    \begin{equation}
      \classspace{W}\cong \classspace{W([0], D)} = \RezkDW_0 \to E_0=E^\simeq,
    \end{equation}
    which is an equivalence of \(\infty\)-groupoids by the first assumption.
  \item
    We have a commutative square induced by \(F\)
    \begin{equation}
      \label{eq:pb-DW-E-0-1}
      \cdsquareNA
      {\classspace{W(\{0<1\},D)}}
      {\Fun(\{0<1\},E)^\simeq}
      {\classspace{W}}
      {E^\simeq}
    \end{equation}
    where by definition the right vertical map is classified by
    \begin{equation}
      e \mapsto (\overcat{E}e)^\simeq.
    \end{equation}
    Moreover, we established earlier that the left vertical map is classified by
    the functor
    \begin{equation}
      d\mapsto \classspace{W\times_D \overcat{D}{d}}.
    \end{equation}
    Since the comparison transformation
    \begin{equation}
      d\mapsto \classspace{W\times_D \overcat{D}{d}} \xrightarrow{\simeq}
      (\overcat{E}{Fd})^\simeq
    \end{equation}
    is an equivalence by the second assumption,
    it follows that the square \eqref{eq:pb-DW-E-0-1}
    is a pullback of \(\infty\)-groupoids.
    Since we have established the lower horizontal map to be an equivalence,
    so is the upper map, which is just the map
    \begin{equation}
      \RezkDW_1\to E_1
    \end{equation}
    on \(1\)-simplices.
    \qedhere
  \end{itemize}
\end{proof}

As an application of \Cref{lem:Berry-localization}
we give a simple proof of the following well-known fact.

\begin{prop}
  \label{lem:localize-to-Delta}
  Consider the poset \(\Ka\) of compact nonempty subsets of \(\reals\)
  with finitely many connected components.
  The connected-component functor
  \begin{equation}
    \pi_0\colon \Ka\to \Delta
  \end{equation}
  is an \(\infty\)-categorical localization at the \(\pi_0\)-equivalences.
\end{prop}

\begin{proof}
  We invoke \Cref{lem:Berry-localization}:
  \begin{enumerate}
  \item
    \label{it:K-Delta-localization-1}
    Denote by \(W\subset\Ka\) the wide subposet of \(\pi_0\)-equivalences
    and \(W_n\subseteq W\) the full subposets consisting of compacts
    with exactly \(n+1\) many connected components.
    We aim to show that the canonical map
    \begin{equation}
     \coprod_{n\in\naturals}  \classspace{W_n} \cong \classspace{W}\to\Delta^\simeq\cong\naturals
    \end{equation}
    is an equivalence,
    i.e., that each \(W_n\) is weakly contractible.
    For each \(K\in W_n\),
    write \(\smallconf{\reals\setminus K}\subset\Conf[n+2]{\reals}\)
    for the subspace of those configurations
    such that there is precisely one point in each of the \(n+2\)
    components of \(\reals\setminus K\).
    Observe that for each \(S=\{s_0<\cdots<s_{n+1}\}\in\Conf[n+2]{\reals}\) the poset
    \(\{K\in W_n\mid S\in\smallconf{\reals\setminus K}\}\)
    is the product over \(0\leq i\leq n\)
    of the posets of nonempty closed intervals in \((s_i,s_{i+1})\);
    each of these posets is filtered, hence weakly contractible.
    Therefore we may apply Lurie's Seifert--Van Kampen~\cite[Theorem~A.3.1]{LurHA}
    to obtain the desired equivalence
    \begin{equation}
      \classspace{W_n}
      =
      \colim_{K\in W_n^\op}*
      \simeq
      \colim_{K\in W_n^\op}\smallconf{\reals\setminus K}
      \simeq
      \Conf[n+2]{\reals}
      \simeq
      *
    \end{equation}
    because both
    \begin{equation}
      \Conf[n+2]{\reals}
      =\{(s_i)\in \reals^{n+2}\mid s_0< \dots < s_{n+1}\}
    \end{equation}
    and each \(\smallconf{\reals\setminus K}\cong\reals^{n+2}\)
    are convex, hence contractible.
  \item
    For each \(K\in \Ka\) we have to show that the map
    \begin{equation}
      \classspace{W\times_\Ka \overcat{\Ka}{K}}\to
      \left(\overcat{\Delta}{\pi_0(K)}\right)^\simeq
    \end{equation}
    is an equivalence.
    By working one connected component of \(K\) at a time,
    we may assume that \(K\) is connected;
    hence the target is \(\left(\overcat{\Delta}{[0]}\right)^\simeq\cong\Delta^\simeq\cong\naturals\),
    so it suffices to show that each
    \(W_n\times_\Ka\overcat{\Ka}{K}\) is weakly contractible.
    We do a very similar argument as in part~\ref{it:K-Delta-localization-1},
    but this time we only cover the subspace
    \(\Ra_K \subseteq \Conf[n+2]{\reals}\)
    consisting of those configurations which have at most one point
    in each of the two components of \(\reals\setminus K\);
    in formulas
    \begin{equation}
      \Ra_K = \{(s_i)\in (-\infty,\max K]\times K^n\times [\min K,\infty)
      \mid s_0<\dots <s_{n+1}\}.
    \end{equation}
    Then we observe that for each \(n\in\naturals\)
    and each \(S=\{s_0<\dots<s_{n+2}\}\in\Ra_K\),
    the poset
    \begin{equation}
      \{K\supseteq L\in W_n
      \mid S\in\smallconf{\reals\setminus L}\}
    \end{equation}
    is isomorphic to the product over \(0\leq i\leq n\) of the posets
    of nonempty closed intervals in \((s_i,s_{i+1})\cap K\).
    Since each \((s_i,s_{i+1})\cap K\) is a
    possibly closed or half-closed, but crucially \emph{nonempty} interval
    (by the assumption \(S\in\Ra_K\)),
    we again conclude that this poset is filtered.
    Hence the Seifert--Van Kampen theorem yields the desired equivalence
    \begin{equation}
      \classspace{W_n}
      =
      \colim_{K\supseteq L\in W_n^\op}
      *
      \simeq
      \colim_{K\supseteq L\in W_n^\op}
      \smallconf{\reals\setminus L}
      \simeq
      \Ra_K
      \simeq
      *
    \end{equation}
    because \(\Ra_K\) and each
    \(\smallconf{\reals\setminus L}\cong \reals^{n+2}\) are convex,
    hence contractible.
    \qedhere
  \end{enumerate}
\end{proof}

\subsection{Colimits of right fibrations}
The following is a basic fact of \infy-category theory
that relates colimits of right fibrations and cofinality.
We record it here for completeness.

\begin{lemma}
  \label{lem:cofinality-right-fibrations}
  Let \(\Da\) be an \(\infty\)-category and
  \begin{equation}
    F\colon P^\triangleright\to \RFib(\Da),
    \quad
    p\mapsto (F(p)\to \Da)
  \end{equation}
  a cone of right fibrations over \(\Da\).
  The following are equivalent
  \begin{enumerate}
  \item
    \label{it:colimit-right-fib}
    The map
    \begin{equation}
      \label{eq:colimit-right-fib}
      \colim_{p\in P}F(p)\to F(\infty)
    \end{equation}
    of \(\infty\)-categories is colimit cofinal.
  \item
    \label{it:colimit-core}
    The map
    \begin{equation}
      \label{eq:colimit-core}
      \colim_{p\in P}\core{F(p)}\to \core{F(\infty)}
    \end{equation}
    of \(\infty\)-groupoids is an equivalence,
    where \(\core{(-)}\) denotes the groupoid core.
  \item
    \label{it:colimit-objects}
    For every object \(d\in \Da\),
    the induced map
    \begin{equation}
      \label{eq:colimit-objects}
      \colim_{p\in P}F(p)_d\to F(\infty)_d
    \end{equation}
    of $\infty$-groupoids
    is an equivalence.
  \end{enumerate}   
\end{lemma}

\begin{proof}
  Since \(F(\infty)\to \Da\) is a right fibration
  it is a general fact that a functor \(X \to F(\infty)\)
  is colimit cofinal if and only if it is a contravariant equivalence over \(\Da\);
  see \cite[Proposition~4.1.11]{Cisinski}).
  Recall that the fully faithful inclusion has a left adjoint
  \begin{equation}
    Q\colon \overcat{\Catinfty}{\Da} \leftrightarrows \RFib(\Da),
  \end{equation}
  given by right fibrant replacement.
  It follows that the map \eqref{eq:colimit-right-fib}
  is a contravariant equivalence
  if and only if it induces an equivalence
  \begin{equation}
    Q(\colim_PF) = \colim_{p\in P}^{\RFib} F(p) \xrightarrow{\simeq} F(\infty),
  \end{equation}
  where \(\colim^{\RFib}\) denotes the colimit computed in right fibrations.
  Under the straightening/unstraightening equivalence
  \(\RFib(\Da)\simeq \Fun(\Da^\op,\Gpdinfty)\)
  this means precisely that we have an equivalence
  \begin{equation}
    \colim_{p\in P} F(p)_\bullet\xrightarrow{\simeq} F(\infty)_\bullet
  \end{equation}
  of presheaves \(\Da^\op \to \Gpdinfty\).
  Since colimits of presheaves are computed pointwise,
  we get the biimplication 
  \ref{it:colimit-right-fib}${\iff}$ \ref{it:colimit-objects}.

  To get the logical equivalence 
  \ref{it:colimit-core}${\iff}$ \ref{it:colimit-objects}
  we apply the biimplication we just showed to the diagram
  \begin{equation}
    \iota^*F\colon P^\triangleright\xrightarrow{F}\RFib(\Da)\xrightarrow{\iota^*}
    \RFib(\core{\Da})
  \end{equation}
  obtained by pulling back along the inclusion \(\iota\colon \core{\Da}\hookrightarrow \Da\).
  For this, observe the following:
  \begin{itemize}
  \item
    The objects of \(\core{\Da}\) are just those of \(\Da\),
    so that condition \ref{it:colimit-objects}
    is the same for \(F\) and for \(\iota^* F\).
  \item
    Every right fibration \(X\to \Da\) is conservative,
    which means that the square
    \begin{equation}
      \cdsquare[pb]
      {\core{X}}
      {X}
      {\core{\Da}}
      {\Da}
      {\iota_X}{}{}{\iota_\Da}
    \end{equation}
    is a pullback.
    Since between \infy-groupoids there is no distinction between cofinality and equivalence,
    this means that condition~\ref{it:colimit-right-fib} for \(\iota^*F\)
    is indeed just condition~\ref{it:colimit-core}.
    \qedhere
  \end{itemize}
\end{proof}

\subsection{Localizations of pre-cocartesian $\infty$-operads}
%local macros
\newcommand\exoo{\Oa^\otimes}
\newcommand\expp{\Pa^\otimes}
\newcommand\exqq{\Qa^\otimes}
\newcommand{\exWW}{W^\otimes}
\newcommand\composW{\Ka_W}

\newcommand\precoc{{!^\mathrm{p}}}%decoration for pre-cocartesian arrows
\newcommand{\overcatpCocart}[2]{{#1}_{/^{\precoc}_{#2}}}

Analogously to the setting of \infy-categories,
we call a map \(\exoo\to\expp\) of \infy-operads a \emph{localization}
at a class of 1-ary operations \(W\),
if for each operad \(\exqq\) the restriction functor
\begin{equation}
  \Opdinfty(\expp,\exqq)
  \to
  \Opdinfty(\exoo,\exqq)
\end{equation}
is fully faithful with essential image spanned by those
operad maps \(\exoo\to\exqq\) which invert the operations in \(W\).

\begin{rem}
  Given a class \(W\) of \(1\)-ary operations in \(\exoo\),
  i.e., a class of arrows in \(\Oa\),
  we can consider the class \(\exWW\) of those morphisms in the various fibers
  \(\exoo_{n_+}\simeq \Oa^n\) which are represented by tuples
  all of whose components lie in \(W\).
  Then an operad map \(\exoo\to\expp\)
  inverts the operations in \(W\) if and only if
  it sends the arrows in \(\exWW\) to equivalences.
  We use these two perspectives interchangeably.
\end{rem}

The goal of this section is to provide an operadic analog
of the localization criterion \Cref{lem:Berry-localization}.
We start with an easy lemma.

\begin{lemma}
  \label{lem:localization-of-operads-total}
  Let \(\exoo\to \expp\) be a map of operads.
  Assume that the map \(\exoo\to\expp\) of total \infy-categories
  (i.e. forgetting the structure map to \(\Fin\))
  is a localization of \infy-categories.
  Then \(\exoo\to\expp\) is a localization of \infy-operads.
\end{lemma}

\begin{proof}
  Denote by \(\exWW\) the class of arrows in \(\exoo\) that become equivalences in \(\expp\).
  Assume that the functor \(\exoo\to\expp\) is a localization of \infy-categories at \(\exWW\).
  Then for every functor \(\exqq\to \Fin\) we can consider the commutative cube
  \begin{equation}
    \cdcubeNA
    {\Fun_{\Fin}(\expp,\exqq)}
    {\Fun(\expp,\exqq)}
    {\Fun_{\Fin}(\exoo,\exqq)}
    {\Fun(\exoo,\exqq)}
    {*}
    {\Fun(\expp,\Fin)}
    {*}
    {\Fun(\exoo,\Fin)}
  \end{equation}
  where the front and the back are cartesian by definition.
  By the localization assumption,
  the two right diagonal restriction maps along \(\exoo\to\expp\)
  are fully faithful with essential image spanned by those functors that invert arros in \(\exWW\);
  hence by pulling back the same is true for the upper left diagonal map.

  Next, we claim that every cocartesian map \(p\xrightarrow{!} p'\) in \(\expp\)
  over some inert \(I_+\to I'_+\) admits
  a lift to an inert cocartesian map in \(\exoo\):
  indeed we can lift \(p\) to an object \(o\)
  (because we assumed that \(\exoo\to\expp\) is a localization,
  in particular essentially surjective)
  and then find an inert lift \(o\xrightarrow{!} o'\) of \(o\) along  \(I_+\to I'_+\);
  this lift \(o\xrightarrow{!} o'\) must then map to \(p\xrightarrow{!} p'\)
  since \(\exoo\to\expp\) preserves inert maps and the cocartesian lift
  of \(p\) along \(I_+\to I'_+\) is unique.
  Hence, when \(\exqq\) is an \infy-operad,
  to check that a map \(\expp\to\exqq\) over \(\Fin\) preserves inert maps,
  it suffices to check it for the ones that come from \(\exoo\);
  in other words, we have a cartesian square of inclusions
  \begin{equation}
    \cdsquareOpt[pb]
    {\Opdinfty(\expp,\exqq)}
    {\Opdinfty(\exoo,\exqq)}
    {\Fun_{\Fin}(\expp,\exqq)}
    {\Fun_{\Fin}(\exoo,\exqq)}
    {hookrightarrow}
    {hookrightarrow}
    {hookrightarrow}
    {hookrightarrow}
  \end{equation}

  Unraveling the above constructions it is immediate that the essential image
  of the resulting fully faithful restriction map
  \begin{equation}
    {\Opdinfty(\expp,\exqq)}
    \hookrightarrow
    {\Opdinfty(\exoo,\exqq)}
  \end{equation}
  consists precisely of those operad maps \(\exoo\to\exqq\)
  which invert all arrows in \(\exWW\),
  which is precisely what we needed to show.
\end{proof}

We will only concern ourselves with operads which
are very close to being symmetric monoidal \infy-categories;
most operads in this paper have this property.

\begin{defn}
  \label{def:pre-cocartesian}
  An \infy-operad \(\exoo\to \Fin\) is called \emph{pre-cocartesian}
  if it admits a class of active maps called \emph{pre-cocartesian}
  (which we decorate by the symbol ``\(\precoc\)'')
  such that each active map \(\alpha\colon\bar{o}\to \bar{u}\) in \(\exoo\)
  with codomain \(\bar{u}\) over \(I_+\in\Fin\)
  factors uniquely as the composite
  \begin{equation}
    \bar{o}\xrightarrow{\precoc}\bar{o'}\xrightarrow{\alpha'}\bar{u}
  \end{equation}
  of a pre-cocartesian arrow
  followed by an arrow in the fiber over \(I_+\).

  An operad map \(\exoo\to\expp\) between pre-cocartesian operads
  is called \emph{multiplicative}
  if it preserves pre-cocartesian arrows.
\end{defn}

\begin{rem}
  \label{rem:left-adjoint-partially-monoidal}
  If \(\exoo\) is pre-cocartesian then for every object
  \(\bar{u}\in\exoo_{I+}\)
  the inclusion
  \begin{equation}
    \overcat{\exoo_{I_+}}{\bar{u}}\hookrightarrow\overcatAct{\exoo}{\bar{u}}
  \end{equation}
  has a left adjoint \(\alpha\mapsto\alpha'\).
  Thus the class of pre-cocartesian arrows is uniquely determined
  rather than being additional data:
  it consists precisely of those \(\alpha\)
  such that \(\alpha'\) is an equivalence.
\end{rem}

\begin{ex}
  \label{ex:partially-monoidal-pre-cocartesian}
  If \(\exoo\) is a symmetric monoidal \infy-category
  (i.e., if it admits all cocartesian lifts),
  then it is pre-cocartesian with the pre-cocartesian arrows
  being precisely the cocartesian ones.

  More generally we can consider operads arising from
  \emph{partially} symmetric monoidal \infy-categories:
  for such operads, the operations \((a_1,\dots,a_n)\to a\)
  are still given as maps \(a_1\otimes\dots\otimes a_n\to a\)
  and only exist for those tuples \((a_i)\) for which the tensor product exists.
  We call such operads \emph{partially monoidal};
  they are precisely those which are pre-cocartesian
  with all pre-cocartesian arrows being cocartesian.

  The operad \(\operad{\PS}\) associated to a \posetofopens{} \(\PS\)
  is partially monodial if and only if \(\PS\) admits joins/suprema
  of finitely many pairwise disjoint opens.
  But note that this partial monoidal structure is in general not given by \(\disjun\);
  in fact, this happens
  if and only if \(\PS\) is closed under subordinate disjoint unions
  (i.e., \(\PS=\cdisj{\PS}\cap\PS\)).
\end{ex}

\begin{ex}
  \label{ex:slice-operad}
  Let \((\Ma,\otimes)\) be a symmetric monoidal \infy-category
  and assume that the monoidal unit is an initial object \(\emptyset\).
  Let \(\Da^\otimes\subseteq \Ma^\otimes\) be a full symmetric monoidal subcategory. 
  Then for each object \(x\in \Ma\) there is
  a canonical pre-cocartesian \infy-operad structure
  \(\operad{(\overcat{\Da}{x})}\)
  on the slice category \(\overcat{\Da}{x}\).
  We only provide an informal description of this \infy-operad
  and refer to \cite[Corollary~1.20]{AFT-fh-stratified}
  for a formal construction:

  An operation \((y_i\to x)_{i\in I}\to (z\to x)\) consists of
  \begin{enumerate}
  \item
    \label{it:slice-monoidal-tensor}
    an extension along the inclusions \(y_i\to y\coloneqq \bigotimes_{i\in I}y_i\)
    (which exist, using that the monoidal unit is initial)
    of the objects \(y_i\to x\)
    to an object \(y\to x\) in \(\overcat{\Da}{x}\), and
  \item
    a map \((y\to x)\to (z\to x)\) in \(\overcat{\Da}{x}\).
  \end{enumerate}
  The unique factorization of an operation
  \(\alpha\colon (y_i\to x)\to (z\to x)\)
  is precisely
  \((y_i\to x)\xrightarrow{\precoc} (y\to x) \xrightarrow{\alpha'} (z\to x)\)
  with the notation above.
  The \infy-operad
  \(\operad{(\overcat{\Da}{x})}\)
  comes equipped with the forgetful operad map
  \(\operad{(\overcat{\Da}{x})}\to \Da^\otimes\);
  the pre-cocartesian arrows in \(\operad{(\overcat{\Da}{x})}\)
  are precisely those which are sent to cocartesian arrows in \(\Da^\otimes\).
 
  If the tensor product \(\otimes\) is the coproduct,
  then the choice of auxiliary extension data
  \ref{it:slice-monoidal-tensor} is unique,
  and \(\operad{(\overcat{\Da}{x})}\) is again a symmetric monoidal \infy-category.
  However, in general this extra data is \emph{not} unique
  and the \infy-operad \(\operad{(\overcat{\Da}{x})}\)
  need not have any non-inert cocartesian arrows at all.
\end{ex}

We can now state the extension of \Cref{lem:Berry-localization}
to pre-cocartesian operads.

\begin{lemma}
  \label{lem:operadic-Berry-localization}
  Let \(f\colon\exoo\to\expp\) be a multiplicative map
  of pre-cocartesian \infy-operads.
  Denote by \(\exWW\) the class of those arrows in \(\exoo\)
  that are inverted by \(f\)
  and write \(W\coloneqq \exWW\cap \Oa\).
  Assume the following:
  \begin{enumerate}
  \item
    \label{it:op-localization-underlying}
    The functor \(f\restrict{1_+}\colon\Oa\to \Pa\) on the underlying \infy-categories
    satisfies the two assumptions of \Cref{lem:Berry-localization},
    i.e.,
    \begin{equation}
      \label{eq:op-localization-underlying-conditions}
      \classspace{W}\xrightarrow{\simeq}\Pa^\simeq
      \quad
      \text{and}
      \quad
      \classspace{W\times_\Oa \overcat{\Oa}{u}}
      \xrightarrow{\simeq}
      \left(\overcat{\Pa}{f(u)} \right)^\simeq
    \end{equation}
    for all \(u\in \Oa\).
  \item
    \label{it:op-localization-left-fib}
    Given a composition
    \(\bar{d}\xrightarrow{\precoc} o \xrightarrow{w} o'\)
    in \(\exoo\)
    with \(o,o' \in \Oa\)
    where the first arrow is pre-cocartesian and the second
    is in \(W\),
    there is a unique factorization of the form
    \begin{equation}
      \label{eq:op-localization-W-factorize}
      \cdsquare
      {\bar{d}}
      {\bar{d'}}
      {o}
      {o'}
      {\bar{w}}
      {\precoc}
      {\precoc}
      {w}
    \end{equation}
    where the top arrow is in \(\exWW\) and the right arrow is pre-cocartesian.
  \item
    \label{it:op-localization-fibers}
    For each object \(o\in \Oa\), the map \(f\) induces an equivalence
    \begin{equation}
      \label{eq:op-localization-fibers}
      \exWW\times_{\exoo}\overcatpCocart{\exoo}{o}
      \xrightarrow{\simeq}
      \left(
        \overcatpCocart{\expp}{f(o)}
      \right)^\simeq
    \end{equation}
    of \infy-groupoids\footnote{
      As explained in the proof, condition~\ref{it:op-localization-left-fib}
      implies that the left side of
      \eqref{eq:op-localization-fibers}
      is automatically an \infy-groupoid,
      so that we do not need to pass to its classifying space
      like we do in \eqref{eq:op-localization-underlying-conditions}.
    },
    where on both sides the symbol ``\(\precoc\)'' denotes the full subcategories
    of the slice spanned by the pre-cocartesian arrows.
  \end{enumerate}
  Then \(f\colon \exoo\to\expp\) is a localization of \infy-operads
  at \(W\).
\end{lemma}

\begin{proof} %[\Cref{lem:operadic-Berry-localization}]
  By \Cref{lem:localization-of-operads-total}
  it suffices that the functor \(f\colon \exoo\to\expp\)
  on total categories is a localization of \infy-categories at \(\exWW\).
  For this we aim to apply \Cref{lem:Berry-localization} to this total functor.
  Before checking its two assumptions, we perform a preliminary computation.

  First, note that we have an (a priori lax) square
  \begin{equation}
    \cdsquareNA
    {\overcatAct{\exoo}{u}}
    {\overcatAct{\expp}{f(u)}}
    {\overcat{\Oa}{u}}
    {\overcat{\Pa}{f(u)}}
  \end{equation}
  where the horizontal morphisms are induced by \(f\)
  and the vertical morphisms are induced by the
  left adjoints described in \Cref{rem:left-adjoint-partially-monoidal};
  the square actually commutes because \(f\) preserves pre-cocartesian arrows.

  After restricting only to those morphisms that become equivalences in \(\expp\),
  we obtain the commutative square
  \begin{equation}
    \label{eq:op-localization-comparison-square}
    \cdsquareNA
    {\exWW\times_{\exoo}\overcatAct{\exoo}{u}}
    {\left(\overcatAct{\expp}{f(u)}\right)^\simeq}
    {W\times_{\Oa}\overcat{\Oa}{u}}
    {\left(\overcat{\Pa}{f(u)}\right)^\simeq}
  \end{equation}

  We know the following facts about the square
  \eqref{eq:op-localization-comparison-square}:
  \begin{itemize}
  \item
    Condition \ref{it:op-localization-left-fib} states exactly
    that the left vertical map is a left fibration;
    in particular its fibers are \infy-groupoids.
  \item
    After factoring through the classifying space,
    the lower horizontal map becomes an equivalence
    by the second assumption of condition~\ref{it:op-localization-underlying}.
  \item
    For each object \((o\to u)\) in the lower left,
    the induced map on fibers is
    \begin{equation}
      \exWW\times_{\exoo}\overcatpCocart{\exoo}{o}
      \to
      \left(
        \overcatpCocart{\expp}{f(o)}
      \right)^\simeq,
    \end{equation}
    which is an equivalence of \infy-groupoids
    by condition~\ref{it:op-localization-fibers}.
  \end{itemize}
  It is a straightforward consequence of Quillen's theorem B
  (see \Cref{lem:Quillen-B-cocart})
  that for any such square the top horizontal map induces an equivalence
  \begin{equation}
    \label{eq:op-localization-key}
    \classspace{\exWW\times_{\exoo}\overcatAct{\exoo}{u}}
    \xrightarrow{\simeq}
    {\left(\overcatAct{\expp}{f(u)}\right)^\simeq}
  \end{equation}
  of \infy-groupoids.

  Now we can easily show the two conditions of \Cref{lem:Berry-localization}
  for the functor \(f\colon \exoo\to\expp\):
  \begin{enumerate}
  \item
    We have the equivalence
    \begin{equation}
      \classspace{\exWW} = \coprod_{n \in \naturals} \classspace{W}^n
      \xrightarrow{\simeq}
      \coprod_{n\in\naturals}(\Pa^\simeq)^n
      = \left(\expp\right)^\simeq
    \end{equation}
    using the first assumption of
    condition~\ref{it:op-localization-underlying}.
  \item
    Note that, fixing the codomain \(\bar{u}\),
    the data of a general map \(h\colon \bar{d}\to \bar{u}\)
    in an \infy-operad is uniquely described
    by a collection of active maps \(\bar{d}_i \to u_i\)
    (which jointly form the active part
    \(h^\act\colon \bigoplus_i\bar{d}_i\to u\) of \(h\))
    together with a finite tuple \((d_j)_{j\in J}\)
    (which correspond to components of the domain
    that are discarded by the inert part \(h^\inert\)).
    
    Hence for each \(\bar{u}=(u_i)\in \exoo_{I_+}\) we have the desired equivalence
    \begin{align}
      \classspace{\exWW\times_{\exoo}\overcat{\exoo}{\bar{u}}}
      &=
        \left(
        \prod_{i\in I}
        \classspace{\exWW\times_{\exoo}\overcatAct{\exoo}{u_i}}
        \right)
        \times
        \left(
        \coprod_{n\in\naturals}
        \classspace{W}^n
        \right)
      \\
      &\xrightarrow{\simeq}
        \left(
        \prod_{i\in I}
        \left(
        \overcatAct{\expp}{f(u_i)}
        \right)^\simeq
        \right)
        \times
        \left(
        \coprod_{n\in\naturals}
        (\Pa^\simeq)^n
        \right)
        =
        \left( \overcat{\expp}{f(\bar{u})} \right).
    \end{align}
    again using the equivalence \(\classspace{W}\to\Pa^\simeq\)
    and the key computation \eqref{eq:op-localization-key}.
    \qedhere
  \end{enumerate}
\end{proof}

\newpage

\section{Operadic colimits and left Kan extensions}
\label{app:OLKE}

We make liberal use of the theory of operadic colimits
and operadic left Kan extensions developed in Section~3.1 of \cite{LurHA}.
The results therein are written in much greater generality
than what is needed for our purposes.
For the convenience of the reader,
we summarize the relevant special cases used in this paper.

Throughout, let \((\targetcat,\otimes)\),
or more precisely, \(q\colon\targetcatOT\to \Fin\)
be a symmetric monoidal \(\infty\)-category;
the underlying \(\infty\)-category is
\(\targetcat\coloneqq\targetcatOT_{1+}\subset\targetcatOT\).
We sometimes write objects of \(\targetcatOT\) as tuples
\((I_+,(c_i)_i)\) or \((c_i)_{i\in I}\),
or abbreviate them as \(\bar{c}\).
We write \(c=(1_+,(c))\) for objects in the fiber
\(\targetcat=\targetcatOT_{1_+}\) over \(1_+\).

\begin{itemize}
\item
  We always assume that \((\targetcat,\otimes)\)
  is \emph{compatible with all (small) colimits}
  (\cite[Definition~3.1.1.18]{LurHA}),
  which means that the \(\infty\)-category \(\targetcat\)
  has all colimits and that the tensor product functor
  \begin{equation}
    \bigotimes_I\colon\targetcat^I\to \targetcat
  \end{equation}
  preserves them in each variable.
\item \cite[Definition~3.1.1.2]{LurHA}
  We write \emph{operadic colimit cone} for what Lurie calls
  ``operadic \(q\)-colimit diagram''
  since we never consider a case where \(q\) is anything other than
  the structure-defining cocartesian fibration to \(\Fin\).
  When we have an operadic colimit cone
  \(F^\triangleright\colon K^\triangleright \to \targetcatact\)
  whose value at the cone point \(\infty\in K^\triangleright\)
  lies in the fiber over \(I_+\in \Fin\),
  we write  \begin{equation}
    \ocolim_{k\in K}F(k)\xrightarrow{\simeq} F^\triangleright(\infty)\quad \in \targetcatOT_{I_+}.
  \end{equation}
  This is a mild abuse of notation, since to specify an operadic colimit
  one needs to specify not just the base diagram
  \(F\colon K\to \targetcatact\)
  but also a cone
  \(F_0^\triangleright\colon K^\triangleright \to \Finact\) of \(F_0\coloneqq qF\),
  which is left implicit in the notation.
\item\cite[Proposition~3.1.1.15 (2) and Proposition~3.1.1.20]{LurHA}
  Let \(F^\triangleright\colon K^\triangleright\to\targetcatact\)
  be a cone and let
  \(\bigotimes F^\triangleright\colon K^\triangleright \to \targetcat\)
 be the transferred cone along the terminal transformation
  \(qF^\triangleright(-)\Rightarrow 1_+\).
  Then \(F^\triangleright\) is an operadic colimit cone if and only if
  \(\bigotimes F^\triangleright\) is an (ordinary) colimit cone in \(\targetcat.\)
  In other words, we have
  \begin{equation}
    \ocolim_{k\in K}F(k)\xrightarrow{\simeq}F^\triangleright(\infty)
    \quad \in \targetcatOT_{qF^\triangleright(\infty)}
  \end{equation}
  if and only if we have
  \begin{equation}
    \colim_{k\in K}\bigotimes F(k) \xrightarrow{\simeq} \bigotimes F^\triangleright(\infty)
    \quad \in \targetcat.
  \end{equation}
  In particular, every diagram \(F\colon K\to \targetcatact\)
  admits an operadic colimit with cone point in \(\targetcat=\targetcatOT_{1_+}\)
  by the formula
  \begin{equation}
    \label{eq:opcolim-1-as-colim}
    \ocolim_{k\in K} F(k)
    \coloneqq \colim_{k\in K} \bigotimes F(k)
    \quad\in\targetcat
  \end{equation}
  (lying over the unique cone \(F_0^\triangleright\colon qF\Rightarrow 1_+\)
  which is left implicit).
\item \cite[Proposition~3.1.1.8]{LurHA}
  If \(\left(F_i^\triangleright\colon K_i^\triangleright\to \targetcatact\right)_{i\in I}\)
  is a finite tuple of operadic colimit cones,
  then the induced cone
  \begin{equation}
    \left( \prod_{i}K_i \right)^\triangleright
    \to
    \prod_iK_i^\triangleright
    \xrightarrow{\prod_iF_i^\triangleright}
    \prod_i\targetcatact
    \xrightarrow{\oplus_I}
    \targetcatact;
  \end{equation}
  is also an operadic colimit cone
  (where \(\oplus\) denotes concatenation of tuples).
  In particular, if all cone points \(F_i^\triangleright(\infty)\)
  take value in \(\targetcat=\targetcatOT_{1_+}\),
  we can write
  \begin{equation}
    \label{eq:tuple-oLKE}
    \ocolim_{k\in\prod_{i}K_i}(F_i(k_i))_{i\in I} = \left( \ocolim_{k_i\in K_i}F_i(k_i) \right)_{i\in I}
    \quad \in \targetcatOT_{I_+}.
  \end{equation}
\item
  \cite[Remark~3.1.1.4]{LurHA}
  Precomposition with colimit cofinal maps preserves operadic colimit cones.
  In other words, we have
  \begin{equation}
    \ocolim_{h\in H}F(\alpha(h))\simeq\ocolim_{k\in K}F(k)
    \quad \in \targetcatOT_{F_0^\triangleright(\infty)}
  \end{equation}
  whenever \(\alpha\colon H\to K\) is a colimit cofinal map.
\item
  \cite[Corollary~3.1.3.5]{LurHA}
  Given a map \(\iota\colon \Oa^\otimes\to \Pa^\otimes\) of small \(\infty\)-operads,
  we have an adjunction
  \begin{equation} \label{eq:OLKEAdjunction}
    \iota_!\colon \Alg[\targetcat]{\Oa}\leftrightarrows \Alg[\targetcat]{\Pa} : \iota^*
  \end{equation}
  where the right adjoint \(\iota^*\) is restriction
  and the left adjoint \(\iota_!\) is \emph{operadic left Kan extension}.
  The operadic left Kan extension is characterized uniquely
  by the pointwise operadic colimit formula
  \begin{equation}
    \label{eq:pointwise-oLKE}
    \ocolim \Aa\restrict{\overcatAct{\Oa^\otimes}{\bar{p}}}
    \xrightarrow{\simeq}
    (\iota_!\Aa)(\bar{p}).
  \end{equation}
\end{itemize}

\begin{cor}
  If \(\iota\colon \Oa^\otimes\to \Pa^\otimes\) is fully faithful
  then \(\iota_!\) is also fully faithful and induces an equivalence
  \begin{equation}
    \Alg[\targetcat]{\Oa}\xrightarrow{\simeq}
    \left\{
      \Aa \in \Alg[\targetcat]{\Pa}
      \,\,
      \middle|
      \,\,
      \forall p\in \Pa\colon \, \ocolim_{\bar{o}\xrightarrow{\act}p}\Aa\restrict{\Oa}(\bar{o})\xrightarrow{\simeq} \Aa(p)
    \right\}.
  \end{equation}
\end{cor}

\begin{proof}
  If \(\iota\) is fully faithful,
  then for all \(o\in\Oa\) the unit
  \begin{equation}
    \Aa(o) \xrightarrow{\simeq} \ocolim\Aa\restrict{\overcatAct{\Oa^\otimes}{o}}\simeq (\iota^*\iota_!\Aa)(o)
  \end{equation}
  is an equivalence because
  \(o\in \overcatAct{\Oa^\otimes}{o}\simeq\overcatAct{\Oa^\otimes}{i(o)}\)
  is a terminal object.
\end{proof}

\begin{cor}
  \label{cor:oLKE-pointwise-colim}
  An algebra \(\Aa\in \Alg[\targetcat]{\Pa}\) is the operadic left Kan extension
  of its restriction to \(\Oa\)
  if and only if it induces an equivalence
  \begin{equation}
    \colim_{\bar{o}\xrightarrow{act}p}\bigotimes\Aa(\bar{o}) \xrightarrow{\simeq} \Aa(p)
    \quad \in \targetcat
  \end{equation}
  for all \(p\in \Pa\).
\end{cor}

\begin{proof}
  It follows from \eqref{eq:tuple-oLKE}
  and the canonical identification
  \begin{equation}
    \overcatAct{\Oa^\otimes}{\bar{p}}\xrightarrow{\simeq}
    \prod_{i\in I}\overcatAct{\Oa^\otimes}{p_i}
  \end{equation}
  for each \(\bar{p}=(p_i)\in \Pa^\otimes_{I_+}\)
  that it suffices to check the operadic left Kan extension condition
  \eqref{eq:pointwise-oLKE}
  for only those \(\bar{p}=p\) that lie over \(1_+\in\Fin\).
  Then the claim follows from the identification of \eqref{eq:opcolim-1-as-colim}.
\end{proof}

In special situations, there is a very tight relationship between the operadic
and the ordinary left Kan extension.

\begin{lemma}
  \label{lem:oLKE-LKE}
  Let \(\iota\colon \Oa^\otimes\hookrightarrow \Pa^\otimes\)
  be a fully faithful inclusion of \infy-operads.
  Assume that for every object \(\bar{o}\in \Oa^\otimes_{I_+}\)
  and every active map
  \(\iota(\bar{o})\to p\) with codomain \(p\in \Pa\),
  there exists a cocartesian arrow \(\bar{o}\to o\)
  in \(\Oa^\otimes\) lifting the active map \(I_+\to 1_+\)
  which remains cocartesian in \(\Pa^\otimes\).
  Then the (a priory lax) square
  \begin{equation}
    \label{eq:oLKE-LKE-square}
    \cdsquareOpt[pb]
    {\Alg[\targetcat]{\Oa}}
    {\Alg[\targetcat]{\Pa}}
    {\Fun(\Oa,\targetcat)}
    {\Fun(\Pa,\targetcat)}
    {"\mathrm{oLKE}",hookrightarrow}
    {"(-)_{1_+}"'}
    {"(-)_{1_+}"}
    {"\mathrm{LKE}",hookrightarrow}
  \end{equation}
  commutes and is a pullback,
  where the upper and lower horizontal functors
  are operadic and ordinary left Kan extension,
  respectively.
\end{lemma}

\begin{rem}
  In this paper we are mostly concerned with operads
  arising from (subposets of) \(\open{X}\) together with its
  partial symmetric monoidal structure given by disjoint union.
  For an inclusion \(\Bf\hookrightarrow \Uf\) of such \posetofopens{},
  the associated operad inclusion \(\operad{\Bf}\hookrightarrow\operad\Uf\)
  satisfies the assumption of \Cref{lem:oLKE-LKE}
  if \(\Bf\) is closed under disjoint unions subordinate to \(\Uf\),
  i.e.,  \(\Bf=\cdisj{\Bf}\cap\cdown{\Uf}\).
\end{rem}

\begin{proof}[Proof of \Cref{lem:oLKE-LKE}]
  The assumption imply that every active arrow of the form
  \(\alpha\colon \iota(\bar{o})\to p\)
  admits a unique factorization
  \(\alpha= (\iota(\bar{o})\xrightarrow{!}\iota(o)\xrightarrow{\alpha'} p)\)
  where the first arrow is cocartesian
  and the second one lies in the \(\Pa\).
  Hence the canonical inclusion
  \begin{equation}
    \overcat{\Oa}{p} \hookrightarrow \overcatAct{\Oa^\otimes}{p}
  \end{equation}
  has a left adjoint given by \(\alpha \mapsto \alpha'\),
  in particular it is colimit cofinal.
  Thus for every algebra \(\Aa\in \Alg[\targetcat]{\Oa}\)
  and every \(p\in \Pa\) the canonical comparison map
  \begin{equation}
    \label{eq:opcolim-colim}
    \mathrm{LKE}\Aa_{1_+}(p)
    \simeq
    {\colim \Aa\restrict{\overcat{\Oa}{p}}}
    \xrightarrow{\simeq}
    {\colim \bigotimes\Aa\restrict{\overcatAct{\Oa^\otimes}{p}}}
    \simeq
    \mathrm{oLKE}\Aa(p)
  \end{equation}
  is an equivalence
  (using \Cref{cor:oLKE-pointwise-colim} for the rightmost equivalence).
  It follows directly that the diagram~\eqref{eq:oLKE-LKE-square}
  commutes and is a pullback.
\end{proof}

We also use the following natural extension of Lurie's terminology
even though it does not, strictly speaking, appear in this form.

\begin{defn}
  Given the solid outer commutative square
  \begin{equation}
    \begin{tikzcd}
      K\ar[r,"F"]
      \ar[d,hookrightarrow]
      &
      \targetcatact
      \ar[d,"q"]
      \\
      L\ar[r,"\overline{F_0}"']
      \ar[ur,dashed,"\overline{F}"]
      &
      \Finact
    \end{tikzcd}
  \end{equation}
  where the left vertical map is fully faithful\footnote{
    We stick to (operadic) left Kan extensions along fully faithful inclusions
    to avoid having to work with lax diagrams.
    This generality is sufficient for the purposes of this paper.
  }
  we say that a dashed lift \(\overline{F}\) is a/the
  operadic left Kan extension of \(F\)
  (again, with \(\overline{F_0}\) left implicit),
  if for each \(l\in L\) we have %the operadic colimit
  \begin{equation}
    \ocolim_{k\to l}F(k)\xrightarrow{\simeq} \overline{F}(l).
  \end{equation}
\end{defn}

\begin{rem}
  Equivalently, \(\overline{F}\) is a left Kan extension of \(F\)
  if the transferred diagram
  \(\bigotimes\overline{F}\colon L\to \targetcat\)
  is an (ordinary) left Kan extension of the transferred diagram
  \(\bigotimes F\colon K\to \targetcat\).
  In particular, operadic left Kan extensions are transitive
  and preserve the operadic colimits of a diagram.
\end{rem}

\begin{lemma} \label{lemma:CocartContractibleDiagramCatOColim}
  Let \(F\colon K\to \targetcatact\) be a diagram that sends all arrows in \(K\)
  to cocartesian morphisms. Assume that \(K\) is weakly contractible.
  Then for each extension
  \(F_0^\triangleright\colon K^\triangleright\to\Finact\) of \(qF\)
  there exists an operadic colimit cone \(F^\triangleright\)
  of \(F\) over \(F_0^\triangleright\) and all structure maps
  \begin{equation}
    F(k)\to \ocolim_{k\in K}F(k)
  \end{equation}
  are cocartesian.
\end{lemma}

\begin{proof}
  Let \(F(-)\Rightarrow F'(-)\colon K \to \targetcatact\)
  be a pointwise cocartesian transformation exhibiting a transfer along
  \(F_0(-)\Rightarrow F_0^\triangleright(\infty)\).
  Since \(F\) sends all arrows to cocartesian morphisms,
  it follows that \(F'\) sends all arrows to equivalences.
  Since \(K\) is weakly contractible, \(F'\) then admits a colimit cone
  \(F'(-)\Rightarrow F^\triangleright(\infty)\)
  all of whose arrows are equivalences.
  Such colimit cones are preserved by any functor,
  in particular by the transfer along the unique active map
  \(F_0^\triangleright(\infty)\to 1_+\).
  It follows that the composite transformation
  \begin{equation}
    F(-)\Rightarrow F'(-)\Rightarrow F^\triangleright(\infty)
  \end{equation}
  is an operadic colimit cone
  \(F^\triangleright\colon K^\triangleright \to \targetcatact\)
  over \(F_0^\triangleright\);
  all of its arrows are cocartesian by construction.
\end{proof}

\newpage

\section{Homotopy colimits of transfinite cubes}
\label{appendix:transfinite-cubes}

\newcommand{\MM}{\mathbf{M}}%local name for the model category

Let \(\MM\) be a model category presenting the \(\infty\)-category
\(\Ca\coloneqq \localize{\MM}{W}\) obtained by localizing at the weak equivalences.
We say that a conical diagram \(D^\triangleright\to \MM\)
is a homotopy colimit cone
if its image in \(\Ca\) is a colimit cone.

In \Cref{subsect:FirstStepGluing} we need that a sufficiently cofibrant colimit cube
in the model category of dendroidal sets
(presenting the \(\infty\)-category of \(\infty\)-operads)
is a homotopy colimit cube.
Since we could not find a reference in the literature for
the infinite version of this statement,
we provide an elementary proof by induction
which directly generalizes to the transfinite case.

\subsection{Transfinite sequences}

% local macros
\newcommand\Ord{\mathrm{Ord}}% the large ordinal of ordinals
\newcommand\ord[1]{\mathrm{ord}(#1)}%order type of a well-ordered set

\begin{itemize}
\item
  By default, all ordinals are assumed to be small.
  We denote by \(\Ord\) the first large ordinal.
  It is the large set of (small) ordinals.
\item
  By a \emph{transfinite sequence}
  \begin{equation}
    x_0\to x_1\to\cdots \to x_\alpha\to x_{\alpha+1}\to \cdots
  \end{equation}
  in \(\Ma\),
  we mean a diagram
  \(x\colon \Ord\to\Ma \)
  or a diagram
  \(x\colon \beta\to \Ma \)
  where \(\beta\) is a (small) ordinal.
\item
  A transfinite sequence \(x\) is called \emph{cocontinuous}
  if at each limit ordinal \(\lambda\),
  the comparison map
  \(\colim_{\alpha<\lambda}x_{\alpha}\xrightarrow{\simeq} x_{\lambda}\)
  is an equivalence.
\item
  We say that a cocontinuous sequence \(x\colon \beta+1\to \Ma\)
  exhibits the map \(x_0\to x_\beta\)
  as a \emph{transfinite composition} of the maps
  \(x_{\alpha}\to x_{\alpha+1}\)
  (for \(\alpha<\beta\)).
\end{itemize}

\begin{lemma}
  \label{lem:transfinite-extension}
  Let \(T\) be a partially ordered (small) set with all suprema;
  in particular, \(T\) has a unique maximal element \(\max T=\sup T\).
  Let \(F\colon T \to\Ca\) be an \(T\)-indexed diagram
  in some \(\infty\)-category \(\Ca\).
  Assume that
  \begin{enumerate}
  \item
    \label{it:can-extend-equivalence}
    for each \(t<\max T\) there exists a \(t<t'\) such that the map
    \(F(t)\to F(t')\) is an equivalence.
  \item
    \label{it:equivalence-to-supremum}
    for each subset \(S\subset T\) the induced map
    \(\colim_{s\in S}F(s)\xrightarrow{\simeq} F(\sup S)\) is an equivalence.
  \end{enumerate}
  Then for each \(t\in T\),
  the map \(F(t)\to F(\max T)\) is an equivalence.
\end{lemma}

\begin{proof}
  Fix \(t\in T\).
  By transfinite recursion, we find a successor ordinal \(\beta+1\)
  and a transfinite sequence \(p\colon \beta+1\to T\) as follows:
  \begin{itemize}
  \item
    We start with \(p_0\coloneqq t\).
  \item
    For each ordinal \(\alpha\),
    the successor value \(p_{\alpha+1}>p_{\alpha}\) is chosen
    using condition~\ref{it:can-extend-equivalence}
    such that
    \(F(p_\alpha)\to F(p_{\alpha+1})\) is an equivalence;
    unless \(p_\alpha=\max T\), in which case the recursion terminates
    and we set \(\beta\coloneqq \alpha\).
  \item
    For each limit ordinal \(\lambda\),
    we set \(p_{\lambda}\coloneqq\sup_{\alpha<\lambda}p_\alpha\).
  \end{itemize}
  Note that the recursion must eventually terminate
  because otherwise it would produce an injection \(\Ord\hookrightarrow T\)
  which is impossible due to size.
  By construction, the sequence \(p\) is cocontinuous,
  hence so is \(F\circ p\) by condition~\ref{it:equivalence-to-supremum}.
  Hence the map \(F(t)=F(p_0)\to F(p_\beta)=F(\max T)\)
  is the transfinite composition of the maps
  \(F(p_\alpha)\to F(p_{\alpha+1})\),
  which are all equivalences by construction.
  We are done, since the transfinite composition of equivalences is an equivalence.
\end{proof}

\subsection{Transfinite cubes} \label{sect:TransfiniteCubes}

% local macros for various constructions
\newcommand\bdr[3][]{\bd[#1]{#3}(#2)}%the restricted latching object at #2
\newcommand\bd[2][]{\clatch[#1]{#2}}%the restricted latching cube
\newcommand\cres[1]{\tau_{#1}}%restriction of cubes
\newcommand\clatch[1][]{\partial_{#1}}%
\newcommand\cubes[2][]{\mathrm{Cube}_{#1}{#2}}
\newcommand\cshift[1][]{\sigma_{#1}}%shift a #1-cube by a finite amount #2

Let \(\Ma\) be an \(\infty\)-category with all colimits.
Denote by \(\emptyset\) its initial object.

\begin{itemize}
\item
  For every ordinal \(\beta\)
  we denote by \(\Pfinop{\beta}\)
  the poset of finite subsets of \(\beta\) ordered by reverse inclusion.
\item
  A \emph{\(\beta\)-cube} in \(\Ma\) is a diagram
  \(\Pfinop{\beta}\to \Ma \).
\item
  We denote by
  \begin{equation}
    \cubes[\beta]{\Ma}\coloneqq \Fun(\Pfinop{\beta},\Ma)
  \end{equation}
  the \(\infty\)-category of \(\beta\)-cubes in \(\Ma\).
\item
  For each \(\alpha< \beta\),
  the \emph{restriction operator}
  \begin{equation}
    \cres{\alpha}\colon \cubes[\beta]{\Ma}\to\cubes[\alpha]{\Ma}
  \end{equation}
  has a fully faithful left adjoint given by left Kan extension.
  Via this inclusion, we identify \(\alpha\)-cubes with
  those \(\beta\)-cubes \(F\)
  where \(F(I)\simeq \emptyset\)
  for all \(I\not\subseteq \alpha\).
\item
  We write
  \begin{equation}
    \cubes{\Ma}=\colim\limits_{\beta\in\Ord}\cubes[\beta]{\Ma}
  \end{equation}
  for the \(\infty\)-category of all cubes;
  the colimit is taken along the aforementioned full embeddings.
  Note that the \(\infty\)-category \(\cubes{\Ma}\) is still locally small,
  since every pair of cubes lie in some common
  \(\cubes[\beta]{\Ma}\subset \cubes{\Ma}\).
\item 
  We can identify \(\cubes{\Ma}\) with the category of those (large) cubes
  \(F\colon \Pfinop{\Ord}\to \Ma\)
  which take the value \(F(I)=\emptyset\) for all except a small amount of \(I\)'s.
\item
  The counits of the adjunctions
  \begin{equation}
    \dots
    \rightleftarrows
    \cubes[\alpha]{\Ma}
    \rightleftarrows
    \cubes[\alpha+1]{\Ma}
    \rightleftarrows
    \dots
    \rightleftarrows
    \cubes{\Ma}
  \end{equation}
  induce a cocontinuous \(\Ord\)-indexed sequence
  \begin{equation}
    \label{eq:restriction-approximation}
    \cres{0}\to\cres{1}\dots\to\cres{\alpha}\to\cres{\alpha+1}\to\dots
    \to \id
  \end{equation}
  of restriction operators \(\cubes{\Ma}\to\cubes{\Ma}\).
\item
  For each (small) subset \(S\subset \Ord\)
  denote by \(\ord{S}\) the \emph{order type} of \(S\),
  i.e.\ the unique ordinal \(\beta\) isomorphic to \((S,<)\).
  Then we have the \emph{restriction operator}
  \begin{equation}\label{eq:SRestOperator}
    \cres{S}\colon \cubes{\Ma}\to\cubes[\ord\beta]{\Ma},
  \end{equation}
  given by precomposition with the embedding
  \begin{equation}
    \Pfinop{\ord{S}}\cong\Pfinop{S}\hookrightarrow\Pfinop{\Ord}.
  \end{equation}
  This notation agrees with the previous one in the case
  where \(S\) is an initial segment of \(\Ord\), i.e.\ \(S=\ord{S}\).
\item
  For each ordinal \(\alpha\),
  we have the cocontinuous \emph{back face operator}
  \begin{equation}
    \cshift[\alpha]\colon \cubes{\Ma}\to\cubes[\alpha]{\Ma}
  \end{equation}
  obtained by pulling back along the embedding
  \begin{equation}
    -\disjun\{\alpha\}\colon
    \Pfinop{\alpha}
    \hookrightarrow
    \Pfinop{\Ord}
  \end{equation}
  of posets.
\item
  The data of an \((\alpha+1)\)-cube \(F\),
  amounts precisely to
  \begin{itemize}
  \item
    the \(\alpha\)-cubes \(\cshift[\alpha]F\) and \(\cres{\alpha}F\)
    which are its back and front face
    ``in direction \(\alpha\)'' respectively, and
  \item
    the transformation \(\cshift[\alpha]F\to\cres{\alpha}F\)
    between those two faces.
  \end{itemize}
  More precisely, we have an equivalence of \(\infty\)-categories
  \begin{equation}
    (\cshift[\alpha]\to\cres{\alpha})\colon
    \cubes[\alpha+1]{\Ma}\xrightarrow{\simeq} \Fun(\Delta^1,\cubes[\alpha]{\Ma}).
  \end{equation}
\end{itemize}

\begin{ex}\label{ex:3Cube}
  Consider a \(3\)-cube \(F \in \cubes[3]{\Ma}\):
  \begin{equation}
    F=\quad\quad
    \cdcubeNA[small]
    {F(012)}{F(12)}
    {F(01)}{F(1)}
    {F(02)}{F(2)}
    {F(0)}{F(\emptyset)}
  \end{equation}  
  Then \(\cshift[2]F\to \cres{2}F\) are the \(2\)-cubes (=squares)
  \begin{equation}
    \cdsquareNA
    {F(012)}{F(12)}
    {F(02)}{F(2)}
    \longrightarrow
    \cdsquareNA
    {F(01)}{F(1)}
    {F(0)}{F(\emptyset)}
  \end{equation}
  comprising the back and front face of \(F\), respectively.
  
\end{ex}

\subsection{The latching operator}

\begin{itemize}
\item
  For each ordinal \(\beta\),
  we have the \emph{latching operator}
  \begin{equation}
    \clatch[\beta]\colon \cubes[\beta]{\Ma}\to\cubes[\beta]{\Ma}
  \end{equation}
  given by left Kan extension and then restriction along the cospan
  \begin{equation}
    \begin{tikzcd}
      \Pfinop{\beta}\times\{1\}
      \ar[r,hookrightarrow]
      &Q&
      \Pfinop{\beta}\times\{0\}
      \ar[from=l,hookleftarrow]
    \end{tikzcd}
  \end{equation}
  (read from left to right)
  where \(Q\) is the set \(\Pfinop{\beta}\times\{0,1\}\)
  equipped with the lexicographic ordering,
  i.e.\
  \((J,j)< (I,i)\) if and only if
  \(J\supsetneq I\)
  or
  both
  \(J=I\) and \(j<i\).
\item
  Explicitly, we have the pointwise formula
  \begin{equation}
    \label{eq:clatch-formula-small}
    \clatch[\beta]F(I)
    =\colim\limits_{(J,1)<(I,0)}F(J)
    =\colim\limits_{\beta\supseteq J\supsetneq I} F(J),
  \end{equation}
  where we keep implicit the fact that \(J\)
  only ranges over finite subsets\footnote{
    Even if \(F\) is defined on a larger class of subsets of \(\beta\)
    (not just the finite ones),
    this notation is still not ambiguous
    because the inclusion
    \(\{J \text{ finite}\}\hookrightarrow \{\beta\supseteq J\supsetneq I\}\)
    is colimit cofinal for every \(I\).
  }.
\item
  A straightforward computation shows
  that for any two ordinals \(\alpha<\beta\),
  the square
  \begin{equation}
    \begin{tikzcd}
      \cubes[\alpha]{\Ma}
      \ar[d,hookrightarrow]
      \ar[r,"{\clatch[\alpha]}"]
      &
      \cubes[\alpha]{\Ma}
      \ar[d,hookrightarrow]
      \\
      \cubes[\beta]{\Ma}
      \ar[r,"{\clatch[\beta]}"]
      &
      \cubes[\beta]{\Ma}
    \end{tikzcd}
  \end{equation}
  canonically commutes.
\item
  Therefore we may assemble all the \(\clatch[\beta]\)
  into a single cocontinuous latching operator
  \begin{equation}
    \clatch=\clatch[\Ord]\colon \cubes{\Ma}\to\cubes{\Ma},
  \end{equation}
  with the pointwise formula
  \begin{equation}
    \label{eq:clatch-formula}
    \clatch F(I)=\colim\limits_{\Ord\supset J\supsetneq I}F(J).
  \end{equation}
  This colimit is well defined,
  since it is equivalent to the
  small colimit \eqref{eq:clatch-formula-small}
  for any sufficiently large (but still small) ordinal \(\beta\).
\item
  The operator \(\clatch\) comes equipped
  with the tautological transformation \(\clatch{F}\to F\),
  which we call the \emph{latching map}.
  A cube \(F\) is \emph{cocartesian} (i.e.\ a colimit diagram in \(\Ma\))
  if and only if
  the latching map \(\clatch{F}(\emptyset)\to F(\emptyset)\)
  at \(\emptyset\) is an equivalence.
\item
  For each (small) subset \(S\subset \Ord\),
  we write\footnote{
    Note that for \(S=\ord{S}=\beta\), this notation is compatible
    with the previous one because for every \(\beta\)-cube \(F\) we have
    \(\cres{\beta}F=F\).
  }
  \(\clatch[S]\coloneqq  \clatch\circ \cres{S}\). 
  We call 
  \begin{equation} \label{eq:SLatchingMap}
    \clatch[S]F\to \cres{S}F\to F
  \end{equation}
  the \(S\)-\emph{restricted latching map}. The corresponding pointwise formula is simply that of \eqref{eq:clatch-formula} applied to the \(S\)-restricted cube, i.e.
  \begin{equation}\label{eq:SLatchingMapPointwise}
  	\clatch[S]F(I)
    =\colim\limits_{\Ord\supset J\supsetneq I}\cres{S}F(J)
    \simeq\colim\limits_{S\supseteq J\supsetneq I} F(J).
  \end{equation}
\item
  Composing the cocontinuous functor \(\clatch\)
  with the cocontinuous sequence
  \eqref{eq:restriction-approximation}
  yields the cocontinuous \(\Ord\)-indexed sequence
  \begin{equation}
    \label{eq:clatch-cocontinuous}
    \clatch[0]\to\clatch[1]\to\dots\to\clatch[\alpha]\to\clatch[\alpha+1]\to\dots \to \clatch\to \id
  \end{equation}
  of functors \(\cubes{\Ma}\to\cubes{\Ma}\).
\item
  The transformation \(\cshift[\alpha]\to\cres{\alpha}\)
  factors as
  \begin{equation}
    \cshift[\alpha]\to \cres{\alpha}\clatch[\alpha+1]\to\cres{\alpha}
  \end{equation}
  where the first map is pointwise just the structure map
  \begin{equation}
    \cshift[\alpha]F(I)=
    F(I\disjun\{\alpha\})\xrightarrow{J=I\disjun\{\alpha\}}
    \colim\limits_{\alpha+1\supseteq J\supsetneq I} F(J)
    =\clatch[\alpha+1]F(I)
  \end{equation}
  of the colimit (for \(I\subseteq \alpha\)).
\item
  For all ordinals \(\alpha\)
  and subsets \(S\subseteq\alpha\)
  we have the natural identification
  \begin{equation}
    \label{eq:clatchshift}
    \begin{tikzcd}
      \cshift[\alpha]\clatch[S\disjun\{\alpha\}]F(I)
      =
      \colim\limits_{S\disjun\{\alpha\}\supseteq J\supsetneq I\disjun\{\alpha\}}F(J)
      \ar[r,equals,"J=K\disjun\{\alpha\}"]
      &
      \colim\limits_{S\supsetneq K\supsetneq I}F(K\disjun\{\alpha\})
      =\clatch[S]\cshift[\alpha]F(I)
    \end{tikzcd}
  \end{equation}
  and in particular (for \(S=\alpha\))
  \begin{equation}
    \cshift[\alpha]\circ \clatch[\alpha+1]=\clatch[\alpha]\circ\cshift[\alpha].
  \end{equation}
\end{itemize}

\begin{ex}\label{ex:SRestrictionAndLatchingMap}
  Consider again the 3-cube from \Cref{ex:3Cube}.
  For the subset \(S=\{0,2\}\) the corresponding restriction operator from Equation \eqref{eq:SRestOperator} is
  \begin{equation}
    \cres{S} F=\quad\quad
    \scalebox{0.8}{
      \cdcubeNA[small]
      {\emptyset}{\emptyset}
      {\emptyset}{\emptyset}
      {F(02)}{F(2)}
      {F(0)}{F(\emptyset)}	}. 
  \end{equation}
  The corresponding \(S\)-restricted latching map is pointwise given by \eqref{eq:SLatchingMapPointwise}.
  For example for \(I=\emptyset\) we have
  \begin{equation}
    \clatch[S]F(\emptyset) = \colim_{\{0,1,2\} \supseteq J \supsetneq \emptyset} \cres{S}F(J) \simeq \colim \big(F(0) \leftarrow F(02) \rightarrow F(2)  \big).
  \end{equation}
  Note that a priori the colimit is over the entire punctured cube, i.e.\ with the 4 copies of \(\emptyset\), but these do not contribute to the colimit so we omit them for simplicity. 
\end{ex}

Now, let \(\La\colon \Ma\to \Ca\)
be a functor between cocomplete \(\infty\)-categories
that preserves the initial object.
We do \emph{not} assume that \(\La\) preserves colimits.
The typical situation is when \(\Ma\) is a model category
and \(\La\) is its \(\infty\)-categorical localization.

\begin{itemize}
\item
  The functor \(\La\) induces an operator \(\cubes{\Ma}\to\cubes{\Ca}\)
  by postcomposition, which we also denote \(\La\).
\item
  For each ordinal \(\alpha\)
  the operators \(\La\) and \(\cshift[\alpha]\) canonically commute,
  since one is given by postcomposition and the other by precomposition.
\item
  Since \(\La\) preserves the initial objects,
  it also commutes with the restriction operators \(\cres{\beta}\).
\item
  On the other hand,
  \(\La\) does not in general commute with the latching operator
  \(\clatch\).
  The universal property of the colimit \eqref{eq:clatch-formula}
  (or more precisely of the left Kan extension
  assembled from those pointwise colimits)
  yields a natural transformation
  \(\clatch \La \to \La\clatch\),
  which is not usually an equivalence.
\end{itemize}

\subsection{Very cofibrant cubes}

From now on let \(\Ma\) be equipped with a wide subcategory of morphisms
called cofibrations.
We call an object \(x\in \Ma\) \emph{cofibrant} if \(\emptyset\to x\) is a cofibration.

\begin{as}
  \label{as:nice-cofibrations}
  We require the following properties:
  \begin{enumerate}[ref= (\arabic*), label=(\arabic*)]
  \item
    \label{as:pushout-cofibration}
    Let
    \begin{equation}
      \label{eq:as:pushout-cofibration}
      \cdsquareNA[po]{x}{x'}{y}{y'}
    \end{equation}
    be a pushout square in \(\Ma\) where \(x\to y\) is a cofibration.
    Then \(x'\to y'\) is also a cofibration.
    Moreover, the square \eqref{eq:as:pushout-cofibration}
    remains a pushout in \(\Ca\) after applying \(\La\).
  \item
    \label{as:transfinite-cofibrations}
    Let \(x_0\to x_\beta\) be a transfinite composition in \(\Ma\)
    of cofibrations \(x_\alpha\to x_{\alpha+1}\) (for \(\alpha<\beta\))
    and assume that \(x_0\) is cofibrant.
    Then the comparison map
    \(\colim_{\alpha<\beta}\La x_{\alpha}\to \La x_\beta\)
    is an equivalence.
  \end{enumerate}
\end{as}

\begin{defn}\label{defn:VeryCofibrantCube}
  A cube \(F\) in \(\Ma\) is called \emph{very cofibrant}
  if for each subset \(S\subset \Ord\),
  the \(S\)-restricted latching map \(\clatch[S]F\to F\)
  is a levelwise cofibration.
\end{defn}

\begin{rem}
  When \(F\) is a \(\beta\)-cube,
  it suffices to check the levelwise cofibrancy of the maps
  \(\clatch[S]F\to F\)
  only for the subsets \(S\subseteq \beta\),
  since for every \(S\subset\Ord\),
  we have \(\clatch[S\cap\beta]F=\clatch[S]F\).
\end{rem}

\begin{lemma}
  If \(F\) is a very cofibrant cube,
  then so is the back face \(\cshift[\alpha]F\)
  for every ordinal \(\alpha\).
\end{lemma}
\begin{proof}
  Assume that \(F\) is very cofibrant.
  Since \(\cshift[\alpha]F\) is an \(\alpha\)-cube,
  it suffices to show that for each subset \(S\subseteq \alpha\),
  the restricted latching map \(\clatch[S]\cshift[\alpha]F\to \cshift[\alpha]F\)
  is a levelwise cofibration.
  But by \eqref{eq:clatchshift}, this map is identified with
  \begin{equation}
    \cshift[\alpha](\clatch[S\disjun\{\alpha\}]F\to F),
  \end{equation}
  which is a levelwise cofibration by assumption,
  and \(\cshift[\alpha]\) preserves levelwise cofibrations
  (since it is just given by a precomposition). 
\end{proof}

Our goal is to prove the following:

\begin{thm}
  \label{thm:cofibrant-cubes}
  Let \(\La\colon\Ma\to\Ca\) be a functor between cocomplete \(\infty\)-categories
  which preserves the initial object and
  satisfies \Cref{as:nice-cofibrations}.

  Then for every very cofibrant cube \(F\) in \(\Ma\)
  the comparison map
  \begin{equation}
    \clatch\La F (\emptyset)\to \La\clatch F(\emptyset)
  \end{equation}
  is an equivalence in \(\Ca\).
\end{thm}

More specifically, we care about the following consequence.
\begin{cor}\label{cor:VeryCofibrantCubeHomotopyColimit}
  Let \(\MM\) be a left proper model category and \(F\in \cubes{\MM}\)
  a very cofibrant cube in \(\MM\).
  Then the colimit
  \begin{equation}
    \clatch{F}(\emptyset)\coloneqq\colim\limits_{I\supsetneq \emptyset}F(I)
  \end{equation}
  is a homotopy colimit,
  i.e.\ remains a colimit in the \(\infty\)-category
  \(\Ca\coloneqq \localize{\MM}{W}\) presented by \(\MM\).
\end{cor}

\begin{proof}
  Follows immediately from \Cref{thm:cofibrant-cubes},
  because for every left proper model category \(\MM\),
  the localization functor
  \(\MM\to\localize{\MM}{W}\) satisfies \Cref{as:nice-cofibrations}:
  \begin{itemize}
  \item[\ref{as:pushout-cofibration}]
    Pushouts of cofibrations in a model category are always cofibrations
    because they are characterized by left lifting against trivial fibrations.
    In a left proper model category,
    a pushout square of the form \eqref{eq:as:pushout-cofibration}
    is automatically a homotopy pushout;
    for example see \cite[Proposition~A.2.4.4]{LurHTT}.
  \item[\ref{as:transfinite-cofibrations}]
    This follows from the existence of the (generalized) Reedy model structure
    on the category of transfinite sequences \(x\colon \beta+1\to \MM\)
    where the cofibrant objects are precisely those sequences
    where \(x_0\) is cofibrant and all latching maps
    \(\colim_{\alpha<\beta}x_\alpha \to x_\beta\)
    are cofibrations.
    For example, see \cite[Corollary~5.1.5]{Hovey} or \cite[Corollary~A.2.9.25]{LurHTT}.
    \qedhere
  \end{itemize}
\end{proof}

We will prove \Cref{thm:cofibrant-cubes}
by induction on the dimension of the cube.
The following is the key lemma used in the inductive step.

\begin{lemma}
  \label{lem:pushout-clatch}
  For every ordinal \(\alpha\) and every cube \(F\),
  we have a natural pushout square
  of \(\alpha\)-cubes
  \begin{equation}
    \label{eq:pushout-clatch}
    \begin{tikzcd}
      {\cshift[\alpha] \clatch[\alpha+1]}F
      \ar[r,equal]
      \ar[d]
      &
      {\clatch[\alpha] \cshift[\alpha]}F
      \ar[r]
      &
      {\clatch[\alpha]}F
      \ar[d]
      \\
      {\cshift[\alpha]}F
      \ar[rr]
      &
      &
      \cres{\alpha}{\clatch[\alpha+1]}F
    \end{tikzcd}
  \end{equation}
  Moreover, if \(F\) is very cofibrant,
  then the two vertical maps are levelwise cofibrations.
\end{lemma}

\begin{proof}
  For every \(I\subset \alpha\),
  the value in the bottom right is by definition the colimit
  \begin{equation}
    \clatch[\alpha+1]F(I)=\colim\limits_{\alpha+1\supseteq J \supsetneq I}F(J).
  \end{equation}
  We can compute this colimit
  by considering the cover
  \begin{equation}
    \{\alpha+1\supseteq J \supsetneq I\}
    =
    \{J\neq I\disjun\{\alpha\}\}
    \cup
    \{\alpha\in J\}
  \end{equation}
  of the indexing poset;
  with the intersection being the poset
  \begin{equation}
    \{J\neq I\disjun\{\alpha\}\}
    \cap
    \{\alpha\in J\}
    =
    \{\alpha+1\supseteq J\supsetneq I\disjun\{\alpha\}\}
    \xleftarrow[J=K\disjun\{\alpha\}]{\cong}
    \{\alpha\supseteq K\supsetneq I\}.
  \end{equation}
  The inclusion
  \begin{equation}
    \{\alpha\supset J\supsetneq I\}\hookrightarrow
    \{\alpha+1\supseteq J\supsetneq I\mid J\neq I\disjun\{\alpha\}\},
  \end{equation}
  has a left adjoint \(J\mapsto J\cap \alpha\),
  hence is colimit cofinal;
  the poset \(\{\alpha\in J\supsetneq I\}\)
  has a terminal object \(J=I\disjun\{\alpha\}\).
  All told we have a pushout square
  \begin{equation}
    \cdsquareNA[po]
    {
      \colim\limits_{\alpha\supseteq K\supsetneq I} F(K\disjun\{\alpha\})
    }
    {
      \colim\limits_{\alpha\supseteq J\supsetneq I}F(J)
    }
    {
      F(I\disjun\{\alpha\})
    }
    {
      \clatch[\alpha+1]F(I)
    }
  \end{equation}
  which is exactly the desired pushout square
  \eqref{eq:pushout-clatch}
  after substituting back in the definitions of
  \(\cshift[\alpha]\) and \(\clatch[\alpha]\).

  It remains to show that for a very cofibrant cube \(F\),
  the left vertical map is a levelwise cofibration;
  the same then follows for the right one by
  \Cref{as:nice-cofibrations}~\ref{as:pushout-cofibration}.
  Indeed, the comparison map \(\clatch[\alpha+1]F\to F\)
  is a levelwise cofibration by the definition of very cofibrant;
  and applying \(\cshift[\alpha]\) does not change that.
\end{proof}

\begin{proof}[Proof of \Cref{thm:cofibrant-cubes}.]
  By transfinite induction on \(\beta\),
  we prove that for every very cofibrant cube \(F\),
  the comparison map
  \(\clatch[\beta]\La F(\emptyset)\to \La\clatch[\beta]F(\emptyset)\)
  is an equivalence in \(\Ca\).
  This yields the claimed result,
  since for every fixed \(F\) we have
  \(\clatch[\beta]F=\clatch F\) for sufficiently large \(\beta\).

  \begin{itemize}
  \item
    The induction start \(\beta=0\) is trivial,
    since \(\clatch[0]F=\emptyset\) (independently of \(F\))
    and \(\La\) preserves initial objects.
  \item
    For the successor step \(\alpha\leadsto\alpha+1\)
    consider the commutative cube
    \begin{equation}
      \begin{tikzcd}
        \clatch[\alpha]\cshift[\alpha]{\La F}{(\emptyset)}
        \ar[rr]
        \ar[dd]
        \ar[dr]
        &&
        \clatch[\alpha]{\La F}{(\emptyset)}
        \ar[dd]
        \ar[dr]
        \\
        &
        \La\clatch[\alpha]\cshift[\alpha]{F}{(\emptyset)}
        \ar[rr,crossing over]
        &&
        \La\clatch[\alpha]{F}{(\emptyset)}
        \ar[dd]
        \\
        \cshift[\alpha]{\La F}{(\emptyset)}
        \ar[dr,equals]
        \ar[rr]
        &&
        \clatch[\alpha+1]{\La F}{(\emptyset)}
        \ar[dr]
        \\
        &
        \La\cshift[\alpha]{F}{(\emptyset)}
        \ar[rr]
        \ar[from=uu,crossing over]
        &&
        \La\clatch[\alpha+1]{F}{(\emptyset)}
      \end{tikzcd}
    \end{equation}
    The back face is the pushout square of \Cref{lem:pushout-clatch}
    for the cube \(\La F\) in \(\Ca\).
    The front face is obtained by applying \(\La\)
    to the pushout square of \Cref{lem:pushout-clatch}
    for the cube \(F\) itself;
    it remains a pushout in \(\Ca\)
    by \Cref{as:nice-cofibrations}~\ref{as:pushout-cofibration}.

    The two upper diagonal arrows
    are equivalences by applying the induction hypothesis
    to the cubes \(\cshift[\alpha]F\) and \(F\),
    respectively.
    From the pushout property,
    it follows that the lower right diagonal arrow
    is also an equivalence, as required.
  \item
    Finally, let \(\lambda\) be a limit ordinal.
    Then we can factor the comparison map
    \(\clatch[\lambda]\La F(\emptyset)\to \La \clatch[\lambda]F(\emptyset)\)
    as a composition of equivalences
    \begin{equation}
      \clatch[\lambda]\La F(\emptyset)
      =
      \colim\limits_{\alpha<\lambda}\clatch[\alpha]\La F(\emptyset)
      \xrightarrow{\simeq}
      \colim\limits_{\alpha<\lambda}\La\clatch[\alpha]F(\emptyset)
      \xrightarrow{\simeq}
      \La\clatch[\lambda] F(\emptyset),
    \end{equation}
    where
    \begin{itemize}
    \item
      the first \(\simeq\) is the induction hypothesis;
    \item
      the second \(\simeq\) follows from
      \Cref{as:nice-cofibrations}~\ref{as:transfinite-cofibrations} because
      \(\clatch[0]F(\emptyset)\to \clatch[\lambda]F(\emptyset)\)
      is a transfinite composition of cofibrations
      \(\clatch[\alpha]F(\emptyset)\to\clatch[\alpha+1]F(\emptyset)\)
      with cofibrant domain \(\clatch[0]F(\emptyset)=\emptyset\).
      \qedhere
    \end{itemize}
  \end{itemize}
\end{proof}

\newpage

\bibliographystyle{alpha}

\bibliography{mybib}

\end{document}